\documentclass[10pt, reqno]{amsart}
\usepackage{appendix}
\usepackage[a4paper, centering]{geometry}
\usepackage{xcolor}
\usepackage{comment}
\usepackage[T1]{fontenc}

\usepackage{xr-hyper}
\usepackage[colorlinks=true,hyperindex,pagebackref,linktocpage=true]{hyperref}
\usepackage{aliascnt}
\usepackage{cleveref}
\usepackage{amsfonts}
\usepackage{amssymb}
\usepackage{tikz-cd}
\usepackage{graphicx}
\usepackage{url}

\hypersetup{
	colorlinks,
	linkcolor={red!50!black},
	citecolor={blue!50!black},
	urlcolor={blue!80!black}
}

\usepackage{color}

\usepackage{stmaryrd}
\usepackage{mathrsfs}  
\usepackage{varwidth}
\theoremstyle{plain}
\newtheorem{theorem}{Theorem}[section]

\newaliascnt{proposition}{theorem}
\newtheorem{proposition}[proposition]{Proposition}
\aliascntresetthe{proposition}

\newaliascnt{lemma}{theorem}
\newtheorem{lemma}[lemma]{Lemma}
\aliascntresetthe{lemma}

\newaliascnt{lemma and definition}{theorem}
\newtheorem{lemma and definition}[lemma and definition]{Lemma and Definition}
\aliascntresetthe{lemma and definition}

\newaliascnt{corollary}{theorem}
\newtheorem{corollary}[corollary]{Corollary}
\aliascntresetthe{corollary}

\newaliascnt{conjecture}{theorem}
\newtheorem{conjecture}[conjecture]{Conjecture}
\aliascntresetthe{conjecture}

\theoremstyle{definition}
\newaliascnt{definition}{theorem}
\newtheorem{definition}[definition]{Definition}
\aliascntresetthe{definition}

\newaliascnt{example}{theorem}
\newtheorem{example}[example]{Example}
\aliascntresetthe{example}

\newaliascnt{remark}{theorem}
\newtheorem{remark}[remark]{Remark}
\aliascntresetthe{remark}

\newaliascnt{hypothesis}{theorem}

\aliascntresetthe{hypothesis}

\newaliascnt{condition}{theorem}
\newtheorem{condition}[condition]{Condition}
\aliascntresetthe{condition}

\theoremstyle{plain}

\newaliascnt{alphconjecture}{alphtheorem}

\aliascntresetthe{alphconjecture}

\crefname{theorem}{Theorem}{Theorems}
\Crefname{theorem}{Theorem}{Theorems}
\crefname{proposition}{Proposition}{Propositions}
\Crefname{proposition}{Proposition}{Propositions}
\crefname{lemma}{Lemma}{Lemmas}
\Crefname{lemma}{Lemma}{Lemmas}
\crefname{lemma and definition}{Lemma and Definition}{Lemmas and Definitions}
\Crefname{lemma and definition}{Lemma and Definition}{Lemmas and Definitions}
\crefname{condition}{Condition}{Conditions}
\Crefname{condition}{Condition}{Conditions}
\crefname{corollary}{Corollary}{Corollaries}
\Crefname{corollary}{Corollary}{Corollaries}
\crefname{conjecture}{Conjecture}{Conjectures}
\Crefname{conjecture}{Conjecture}{Conjectures}
\crefname{definition}{Definition}{Definitions}
\Crefname{definition}{Definition}{Definitions}
\crefname{example}{Example}{Examples}
\Crefname{example}{Example}{Examples}
\crefname{remark}{Remark}{Remarks}
\Crefname{remark}{Remark}{Remarks}
\crefname{hypothesis}{Hypothesis}{Hypotheses}
\Crefname{hypothesis}{Hypothesis}{Hypotheses}
\crefname{alphtheorem}{Theorem}{Theorems}
\Crefname{alphtheorem}{Theorem}{Theorems}
\crefname{alphconjecture}{Conjecture}{Conjectures}
\Crefname{alphconjecture}{Conjecture}{Conjectures}

\newcommand{\nc}{\newcommand}
\nc{\on}{\operatorname}

\nc{\Q}{\mathbb{Q}}
\nc{\Z}{\mathbb{Z}}
\nc{\cl}{\mathrm{cl}}

\nc{\fraka}{{\mathfrak a}} \nc{\bba}{{\mathbf a}}
\nc{\frakb}{{\mathfrak b}}
\nc{\frakc}{{\mathfrak c}}
\nc{\frakd}{{\mathfrak d}}
\nc{\frake}{{\mathfrak e}}
\nc{\frakf}{{\mathfrak f}}
\nc{\frakg}{{\mathfrak g}}
\nc{\frakh}{{\mathfrak h}}
\nc{\fraki}{{\mathfrak i}}
\nc{\frakj}{{\mathfrak j}}
\nc{\frakk}{{\mathfrak k}}
\nc{\frakl}{{\mathfrak l}}
\nc{\frakm}{{\mathfrak m}}
\nc{\frakn}{{\mathfrak n}}
\nc{\frako}{{\mathfrak o}}
\nc{\frakp}{{\mathfrak p}}
\nc{\frakq}{{\mathfrak q}}
\nc{\frakr}{{\mathfrak r}}
\nc{\fraks}{{\mathfrak s}}
\nc{\frakt}{{\mathfrak t}}
\nc{\fraku}{{\mathfrak u}}
\nc{\frakv}{{\mathfrak v}}
\nc{\frakw}{{\mathfrak w}}
\nc{\frakx}{{\mathfrak x}}
\nc{\fraky}{{\mathfrak y}}
\nc{\frakz}{{\mathfrak z}}
\nc{\frakA}{{\mathfrak A}}
\nc{\frakB}{{\mathfrak B}}
\nc{\frakC}{{\mathfrak C}}
\nc{\frakD}{{\mathfrak D}}
\nc{\frakE}{{\mathfrak E}}
\nc{\frakF}{{\mathfrak F}}
\nc{\frakG}{{\mathfrak G}}
\nc{\frakH}{{\mathfrak H}}
\nc{\frakI}{{\mathfrak I}}
\nc{\frakJ}{{\mathfrak J}}
\nc{\frakK}{{\mathfrak K}}
\nc{\frakL}{{\mathfrak L}}
\nc{\frakM}{{\mathfrak M}}
\nc{\frakN}{{\mathfrak N}}
\nc{\frakO}{{\mathfrak O}}
\nc{\frakP}{{\mathfrak P}}
\nc{\frakQ}{{\mathfrak Q}}
\nc{\frakR}{{\mathfrak R}}
\nc{\frakS}{{\mathfrak S}}
\nc{\frakT}{{\mathfrak T}}
\nc{\frakU}{{\mathfrak U}}
\nc{\frakV}{{\mathfrak V}}
\nc{\frakW}{{\mathfrak W}}
\nc{\frakX}{{\mathfrak X}}
\nc{\frakY}{{\mathfrak Y}}
\nc{\frakZ}{{\mathfrak Z}}
\nc{\bbA}{{\mathbb A}}
\nc{\bbB}{{\mathbb B}}
\nc{\bbC}{{\mathbb C}}
\nc{\bbD}{{\mathbb D}}
\nc{\bbE}{{\mathbb E}}
\nc{\bbF}{{\mathbb F}} \nc{\bbf}{{\mathbf f}}
\nc{\bbG}{{\mathbb G}}
\nc{\bbH}{{\mathbb H}}
\nc{\bbI}{{\mathbb I}}
\nc{\bbJ}{{\mathbb J}}
\nc{\bbK}{{\mathbb K}}
\nc{\bbL}{{\mathbb L}}
\nc{\bbM}{{\mathbb M}}
\nc{\bbN}{{\mathbb N}}
\nc{\bbO}{{\mathbb O}}
\nc{\bbP}{{\mathbb P}}
\nc{\bbQ}{{\mathbb Q}}
\nc{\bbR}{{\mathbb R}}
\nc{\bbS}{{\mathbb S}}
\nc{\bbT}{{\mathbb T}}
\nc{\bbU}{{\mathbb U}}
\nc{\bbV}{{\mathbb V}}
\nc{\bbW}{{\mathbb W}}
\nc{\bbX}{{\mathbb X}}
\nc{\bbY}{{\mathbb Y}}
\nc{\bbZ}{{\mathbb Z}}
\nc{\calA}{{\mathcal A}}
\nc{\calB}{{\mathcal B}}
\nc{\calC}{{\mathcal C}}
\nc{\calD}{{\mathcal D}}
\nc{\calE}{{\mathcal E}}
\nc{\calF}{{\mathcal F}}
\nc{\calG}{{\mathcal G}}
\nc{\calH}{{\mathcal H}}
\nc{\calI}{{\mathcal I}}
\nc{\calJ}{{\mathcal J}}
\nc{\calK}{{\mathcal K}}
\nc{\calL}{{\mathcal L}}
\nc{\calM}{{\mathcal M}}
\nc{\calN}{{\mathcal N}}
\nc{\calO}{{\mathcal O}}
\nc{\calP}{{\mathcal P}}
\nc{\calQ}{{\mathcal Q}}
\nc{\calR}{{\mathcal R}}
\nc{\calS}{{\mathcal S}}
\nc{\calT}{{\mathcal T}}
\nc{\calU}{{\mathcal U}}
\nc{\calV}{{\mathcal V}}
\nc{\calW}{{\mathcal W}}
\nc{\calX}{{\mathcal X}}
\nc{\calY}{{\mathcal Y}}
\nc{\calZ}{{\mathcal Z}}

\nc{\scrA}{{\mathscr A}}
\nc{\scrB}{{\mathscr B}}
\nc{\scrC}{{\mathscr C}}
\nc{\scrD}{{\mathscr D}}
\nc{\scrE}{{\mathscr E}}
\nc{\scrF}{{\mathscr F}}
\nc{\scrG}{{\mathscr G}}
\nc{\scrH}{{\mathscr H}}
\nc{\scrJ}{{\mathscr J}}
\nc{\scrI}{{\mathscr I}}
\nc{\scrK}{{\mathscr K}}
\nc{\scrL}{{\mathscr L}}
\nc{\scrM}{{\mathscr M}}
\nc{\scrN}{{\mathscr N}}
\nc{\scrO}{{\mathscr O}}
\nc{\scrP}{{\mathscr P}}
\nc{\scrQ}{{\mathscr Q}}
\nc{\scrR}{{\mathscr R}}
\nc{\scrS}{{\mathscr S}}
\nc{\scrT}{{\mathscr T}}
\nc{\scrU}{{\mathscr U}}
\nc{\scrV}{{\mathscr V}}
\nc{\scrW}{{\mathscr W}}
\nc{\scrX}{{\mathscr X}}
\nc{\scrY}{{\mathscr Y}}
\nc{\scrZ}{{\mathscr Z}}

\nc{\tG}{{\widetilde{G}}}
\nc{\tH}{{\widetilde{H}}}
\nc{\tM}{{\widetilde{M}}}
\nc{\tR}{{\widetilde{R}}}
\nc{\tS}{{\widetilde{S}}}
\nc{\tT}{{\widetilde{T}}}
\nc{\tX}{{\widetilde{X}}}
\nc{\tDelta}{{\widetilde{\Delta}}}

\nc{\tcalG}{{\widetilde{\calG}}}

\nc{\bnu}{{\bar{ \nu}}}
\nc{\olO}{\bar{\calO}}
\nc{\pf}{{\on{pf}}}

\nc{\al}{{\alpha}} 
\nc{\be}{{\beta}}
\nc{\ga}{{\gamma}} \nc{\Ga}{{\Gamma}}
\nc{\hGa}{\hat{\Gamma}}
\nc{\ve}{{\varepsilon}} 
\nc{\la}{{\lambda}} \nc{\La}{{\Lambda}}
\nc{\om}{\omega} \nc{\Om}{\Omega} 
\nc{\sig}{{\sigma}} \nc{\Sig}{{\Sigma}}
\nc{\dR}{{\mathrm{dR}}}
\nc{\Perf}{{\mathrm{Perf}}}
\nc{\Gm}{{\mathbb{G}_m}}

\nc{\Spa}{\on{{Spa}}}
\nc{\Spd}{\on{{Spd}}}
\nc{\tnb}{\psi_{\rm tame}}
\nc{\oM}{\overline{{M}}}
\nc{\op}{{\on{op}}}
\nc{\ad}{{\on{ad}}}
\nc{\alg}{{\on{alg}}}
\nc{\Ad}{{\on{Ad}}}
\nc{\Adm}{{\on{Adm}}}
\nc{\aff}{{\on{af}}}
\nc{\Aut}{{\on{Aut}}}
\nc{\Bun}{{\on{Bun}}}
\nc{\cha}{{\on{char}}}
\nc{\colim}{{\on{colim}}}
\nc{\der}{{\on{der}}}
\nc{\Der}{{\on{Der}}}
\nc{\diag}{{\on{diag}}}
\nc{\Dist}{{\on{Dist}}}
\nc{\End}{{\on{End}}}
\nc{\Fl}{{\calF\!\ell}}
\nc{\Tr}{{\on{Transp}}}
\nc{\TR}{{\calT\!\calR}}
\nc{\Gal}{{\on{Gal}}}
\nc{\Gr}{{\on{Gr}}}
\nc{\Hk}{{\on{Hk}}}
\nc{\rH}{{\on{H}}}
\nc{\Hom}{{\on{Hom}}}
\nc{\IC}{{\on{IC}}}
\nc{\id}{{\on{id}}}
\nc{\Id}{{\on{Id}}}
\nc{\ind}{{\on{ind}}}
\nc{\Ind}{{\on{Ind}}}
\nc{\Lie}{{\on{Lie}}}
\nc{\Pic}{{\on{Pic}}}
\nc{\pr}{{\on{pr}}}
\nc{\Res}{{\on{Res}}}
\nc{\res}{{\on{res}}}
\nc{\Sat}{{\on{Sat}}}
\nc{\spc}{{\on{sc}}}
\nc{\drv}{{\on{der}}}
\nc{\sgn}{{\on{sgn}}}
\nc{\Spec}{{\on{Spec}}}
\nc{\Spf}{\on{Spf}} 
\nc{\Sph}{\on{Sph}}
\nc{\St}{{\on{St}}}
\nc{\tr}{{\on{tr}}}
\nc{\Mod}{{\mathrm{-Mod}}}
\nc{\Hilb}{{\on{Hilb}}} 
\nc{\Ext}{{\on{Ext}}} 
\nc{\vs}{{\on{Vec}}}
\nc{\ev}{{\on{ev}}}
\nc{\nO}{{\breve{\calO}}}
\nc{\spe}{{\on{sp}}}
\nc{\loc}{{\on{loc}}}

\nc{\nscrR}{{\mathscr{R}^{\on{nr}}}}

\nc{\GL}{{\on{GL}}}

\nc{\Gl}{\on{Gl}} 
\nc{\GSp}{{\on{GSp}}}
\nc{\gl}{{\frakg\frakl}}
\nc{\SL}{{\on{SL}}} 
\nc{\SU}{{\on{SU}}} 
\nc{\SO}{{\on{SO}}}
\nc{\PGL}{{\on{PGL}}}

\nc{\Conv}{{\on{Conv}}}
\nc{\gr}{{\on{gr}}}
\nc{\Rep}{{\on{Rep}}}
\nc{\Dom}{{\on{Dom}}}
\nc{\red}{{\on{red}}}
\nc{\act}{{\on{act}}}
\nc{\nr}{{\on{nr}}}
\nc{\ctf}{{\on{ctf}}}

\nc{\Perm}{{\on{Perm}}} 
\nc{\str}{{\on{-}}} 
\nc{\os}{{\bar{s}}}
\nc{\oeta}{{\bar{\eta}}}

\nc{\hookto}{\hookrightarrow}
\nc{\longto}{\longrightarrow}
\nc{\leftto}{\leftarrow}
\nc{\onto}{\twoheadrightarrow}
\nc{\lonto}{\twoheadleftarrow}

\nc{\Bmu}{\mbox{$\raisebox{-0.59ex}{$l$}\hspace{-0.18em}\mu\hspace{-0.88em}\raisebox{-0.98ex}{\scalebox{2}{$\color{white}.$}}\hspace{-0.416em}\raisebox{+0.88ex}{$\color{white}.$}\hspace{0.46em}$}{}}

\nc{\pot}[1]{ [\![ {#1} ]\!] }
\nc{\rpot}[1]{ (\!( {#1} )\!) }

\numberwithin{equation}{section}

\begin{document}
	
	\title{Moduli descriptions of v-sheaf local models}
	
	\author[J.~Louren\c{c}o, T.~Richarz, E.~Viehmann, and T.~Wedhorn]{Jo\~ao Louren\c{c}o, Timo Richarz, Eva Viehmann, and Torsten Wedhorn}
	
 \address{L'Institut Galilée, Université Sorbonne Paris Nord, 99 avenue Jean Clément, Villetaneuse, France}
 \email{lourenco@math.univ-paris13.fr}
 	
	\address{Technische Universit\"at Darmstadt, Department of Mathematics, 64289 Darmstadt, Germany}
	\email{richarz@mathematik.tu-darmstadt.de}
	
	\address{Mathematisches Institut, Universit\"at M\"unster, Einsteinstrasse 62, M\"unster, Germany}
	\email{viehmann@uni-muenster.de}
	
	\address{Technische Universit\"at Darmstadt, Department of Mathematics, 64289 Darmstadt, Germany}
	\email{wedhorn@mathematik.tu-darmstadt.de}
	
	\begin{abstract}
We propose a uniform moduli-theoretic description of the v-sheaf local model $\mathbb{M}_{\mathcal G,\mu}$ for an arbitrary parahoric group scheme $\mathcal G$ and an arbitrary dominant cocharacter $\mu$. 
		Our approach is Tannakian and builds on a new theory of parahoric induced modules. 
		When $\mu$ is minuscule, this framework explains the linear-algebraic conditions arising in the theory of local models of Shimura varieties.

	\end{abstract}
	\date{\today}

	\maketitle
	\tableofcontents

\section{Introduction}

Local models are meant to describe the singularities of integral models of Shimura varieties with parahoric level structure at a prime $p$.  
In examples arising from moduli spaces of abelian varieties with additional structure, one is naturally led to linear-algebraic moduli problems over complete discrete valuation rings of mixed characteristic $(0,p)$.  
The evident ``naive'' moduli problem, however, is often not flat \cite{Pap00}.  
Recovering the correct local model then requires additional conditions, such as determinant, wedge, spin, or strengthened spin conditions, whose form can appear specific to the group and its chosen representation.

The interplay between intrinsic constructions and explicit moduli descriptions is a recurring theme throughout the theory.
The early examples \cite{DR73,Rap90,CN90,dJ93,DP94} were organized into a more general framework by Rapoport--Zink \cite{RZ96}, and their geometry was subsequently studied by Faltings, G\"ortz, Pappas, Rapoport, and others \cite{Fal97,Fal01,Pap00,Goertz:LinearCase,Goertz:SymplecticCase,Fal03,PR03,PR05,PR08,PR09}, see \cite{PRS13} for a survey.  
Pappas--Zhu \cite{PZ13,Zhu14} introduced a group-theoretic construction as a flat closure in an affine Grassmannian, later extended by Levin \cite{Lev16} and Lourenço \cite{Lou23}.
Scholze--Weinstein proposed an intrinsic construction in $p$-adic geometry \cite{SW20}, whose existence, uniqueness, and compatibility with the scheme-theoretic theory were established in \cite{AGLR22}.
Local models also enter the construction by Kisin, Pappas, and Zhou of integral models for Shimura varieties of abelian type with parahoric level structure \cite{kisin_pappas_integral_models_of_shimura_varieties_with_parahoric_level_structure,KisinPappasZhou_IntegralModels}.

What remained largely case-by-case is the passage from a canonical construction to concrete moduli conditions.  
Explicit descriptions are known for several classes of groups and level structures \cite{Goertz:LinearCase,Goertz:SymplecticCase,Goertz:RamifiedCase,PR03,PR05,PR09,SmithlingOddUnitary,Smithling:TopFlatnessEven,SmithlingEvenUnitary,Smithling15,Zachos,Pappas_Zachos:OrthogonalShimuraVarieties,Luo,Yang:SpinLocalModelsEven}, but their defining conditions have lacked a common conceptual explanation.
In the present manuscript, we provide such an explanation for the $v$-sheaf local models of Scholze--Weinstein at arbitrary parahoric level, using algebraic representations of parahoric group schemes and Tannakian formalism.
Examples are discussed in the companion paper \cite{LRVW:Minuscule}.

Our approach is inspired by the work of Finkelberg--Mirkovi\'c \cite[Section~10.3]{FM:SemiInfinite} on spherical Schubert varieties in affine Grassmannians; see also \Cref{remark reverse inclusion} for further references.
Earlier work of Hartl--Viehmann \cite[Proposition~5.5]{HartlViehmann1} gives a topological description of local models for equi-characteristic local shtukas at hyperspecial level, where the problem reduces to the spherical Schubert case.
Haines \cite{Haines:Pluecker} proposes a related construction via Drinfeld--Pl\"ucker relations at parahoric level, but the resulting moduli problems need not coincide with the local models.

Our results also imply topological flatness of the corresponding scheme-theoretic moduli problems in both equal and mixed characteristic.
We will address their scheme-theoretic structure in subsequent work.

\subsection*{Definition of the moduli problem}
Let $F$ be a non-archimedean local field with ring of integers $O$.
Let $G$ be a connected reductive $F$-group, let $\calG$ be a parahoric $O$-model of $G$, and let $\mu$ be a conjugacy class of geometric cocharacters of $G$, not necessarily minuscule.
We denote by $E$ the reflex field of $\mu$ and by $E^{\mathrm{ur}}$ its maximal unramified subextension over $F$. Let $\Rep(\calG_{O_E})$ be the category of algebraic representations $\rho\colon \calG_{O_E}\to \GL(\calV)$ on finite free $O_E$-modules.

The v-sheaf Beilinson--Drinfeld Grassmannian
$\Gr_\calG \to \Spd(O)$
parametrizes $\calG$-bundles $\calE$ on the de Rham period ring $B_{\mathrm{dR}}^+$ together with a trivialization after inverting $\xi$, see \cite[Section~19]{SW20}, or \Cref{recollection} for a quick reminder.
Every $(\rho,\calV) \in \Rep(\calG_{O_E})$ determines therefore an $O_E\otimes_{O_{E^{\mathrm{ur}}}}B_{\mathrm{dR}}^+$-lattice
$$
\tilde\rho_*\calE \;\subset\; \calV\otimes_{O_{E^{\mathrm{ur}}}} B^+_{\mathrm{dR}}[\xi^{-1}],
$$
where $\tilde \rho_*\colon \Gr_{\calG,O_E}\to \Gr_{\Res_{O_E/O_{E^{\mathrm{ur}}}}\GL(\calV), O_E}$ is the morphism induced by $\rho$ by restriction of scalars from $O_E$ to $O_{E^{\mathrm{ur}}}$.

\begin{definition}[{\Cref{definition local model}}]\label{definition local model intro}
	For any subset $J\subset \Rep(\calG_{O_E})$, let ${\mathbb M}_{\mathcal{G},\mu}^J\subset \Gr_{\mathcal{G},O_E}$ be the sub-functor parametrizing those bundles $\calE$ with a generic trivialization such that 
	\begin{equation}\label{Tannaka bound local model intro}
		\underline{\xi}^{{\underline{\mu}_{\max}}(\calV)}\cdot\big(\calV\otimes_{O_{E^{\mathrm{ur}}}}B^+_{\mathrm{dR}}\big)\subset \tilde\rho_*\calE, \quad \text{for every $(\rho,\calV)\in J$.}
	\end{equation}
	Here $\mu_{\mathrm{max}}(\calV)\in \mathbb Z$ denotes the maximal $\mu$-weight of $\calV_{\bar F}$ (\Cref{representation weight}), and the underlined notation records the different $E^{\mathrm{ur}}$-embeddings of $E$ into $\bar F$, see \Cref{preliminary lemma}.
\end{definition}

\begin{remark}
If $E$ is unramified over $F$, then \Cref{Tannaka bound local model intro} simplifies to the condition
\begin{equation}\label{defining lattice inclusion unramified}
\xi^{{{\mu}_{\max}}(\calV)}\cdot\big(\calV\otimes_{O_E}B^+_{\mathrm{dR}}\big)\subset \rho_*\calE, \quad \text{for every $(\rho,\calV)\in J$.}
\end{equation}
\end{remark}

Thus $\mathbb M_{\calG,\mu}^J$ is described by ordinary lattice containments in linear representations. 
We want to compare it with the v-sheaf local model $\mathbb M_{\calG,\mu}$ defined as the v-sheaf closure in $\Gr_{\calG,O_E}$ of the Schubert variety defined by $\mu$, see \cite{AGLR22}.
For every $J\subset\Rep(\calG_{O_E})$, one has an inclusion
\begin{equation}
\label{inclusion intro}
\mathbb M_{\calG,\mu}\subset\mathbb M_{\calG,\mu}^J,
\end{equation}
see \Cref{easy inclusion model}. 
In general, one cannot expect to have an equality, even for $J = \Rep(\calG_{O_E})$: 
first, algebraic groups whose derived group is not simply connected ``do not have enough linear representations''; 
second, the inertia-coinvariants of the algebraic fundamental group $\pi_1(G)_I$ need to be torsion-free. 
The first obstruction can already be seen for the quasi-minuscule Schubert variety for $G = \PGL_2$; see \Cref{pgl2 remark}. 
The second obstruction appears for ramified tori; see \Cref{torsion torus example}.
Accordingly, whenever needed, we will assume that $G_{\mathrm{der}}$ is simply connected and that $\pi_1(G)_I$ is torsion-free. 
Since our aim is to describe the v-sheaf local model $\mathbb M_{\calG,\mu}$, we may always reduce to this setting: after passing to a suitable central extension of $G$, both conditions hold without changing the associated v-sheaf local model; see \Cref{central extension lemma}.

\subsection*{Main result}
We aim to establish an integral equality $\mathbb M_{\calG,\mu}=\mathbb M_{\calG,\mu}^J$ for a finite, and ideally minimal, subset $J$ of $\calG_{O_E}$-representations. 
Our approach relies on the construction of suitable algebraic representations of parahoric group schemes.
More precisely, let $S$ be a maximal $\breve F$-split torus in $G$ defined over $F$ whose
apartment contains the facet defining $\calG$, let $T=Z_G(S)$ be its centralizer, and let $B_0$ be a Borel subgroup in $G_{\breve F}$ containing $S_{\breve F}$.
To a dominant character $\chi\in X^*(T)^+$ and a Borel subgroup $B$ containing $S_{\breve F}$, we associate the
\textit{parahoric induced module} (or \textit{parahoric dual Weyl module})
\begin{equation}
\label{parahoric weyl module intro}
\calV_{\chi,B}
 =
\Res_{\tilde O/O_E}\left(
\Ind_{\calB^-}^{\calG_{\tilde O}}({}^B\!\chi)\right),
\end{equation}
where $\tilde E$ is an auxiliary finite field extension of $E$ splitting $T$ with ring of integers $\tilde O$,
$\calB^-$ denotes the scheme-theoretic closure in $\calG_{\tilde O}$ of the opposite Borel subgroup $B^-$, 
and ${}^B\!\chi$ is the $B$-dominant representative in the Weyl orbit of $\chi$,
see \Cref{parahoric Weyl module sec}.
Even though the quotient $\calG_{\tilde O}/\calB^-$ need not be proper over $\tilde O$, the algebraic induction in \Cref{parahoric weyl module intro} is a finite free $\tilde O$-module and is a lattice in the
dual Weyl module of highest weight $\chi$, see \Cref{enough reps of parahorics}.

Choose a finite set $J_1\subset X^*(T)^+$ that generates $X^*(T)^{+}_{\mathbb Q}$ as a $\mathbb Q_{\geq0}$-monoid and that satisfies \Cref{J_1 condition}.
The latter condition is automatic if $E/F$ is unramified, and it can always be satisfied by adjoining inertia invariant characters to $J_1$.
Let $J_2$ be the finite set of $\calG$-relevant Borel subgroups, see \Cref{relevant Borel subgroups}.
We write
\begin{equation}
\label{definition of J intro}
J(\calG,\mu)
 :=
 J_1\times J_2,
\end{equation}
viewed as a subset of $\calG_{O_E}$-representations via $(\chi,B)\mapsto \calV_{\chi,B}$.

Recall that $G$ is \textit{essentially tamely ramified}, see \cite[Section~5]{PR21} and \cite[Definition~3.1.4]{KisinPappasZhou_IntegralModels}, if $G_{\breve F,\mathrm{ad}}\simeq\prod_j\Res_{F_j/\breve F}G_j$, where each $F_j$ is a finite separable field extension of $\breve F$ and each $G_j$ is an absolutely simple $F_j$-group that splits over a tamely ramified extension. 
In particular, if the residue characteristic of $F$ is at least $5$, then every reductive $F$-group is essentially tamely ramified.
Our main result is the following.

\begin{theorem}[{\Cref{cor:qbg-iwahori-moduli}}]\label{main theorem intro}
Assume that $G_\der$ is simply connected, that $\pi_1(G)_I$ is torsion-free, and that $G$ is essentially tamely ramified. 
Then, for every parahoric $O$-model $\calG$ and every conjugacy class $\mu$ of geometric cocharacters, not necessarily minuscule, one has $\mathbb M_{\calG,\mu}^{J}=\mathbb M_{\calG,\mu}$ for $J=J(\calG,\mu)$.
\end{theorem}

This is the first uniform group-theoretic result on moduli descriptions of local models.
In particular, the v-sheaf local model is cut out by finitely many explicit linear-algebraic conditions obtained directly from the root datum.
This gives a uniform explanation for several conditions that previously appeared separately in moduli descriptions of local models of Shimura varieties \cite{Goertz:LinearCase,Goertz:SymplecticCase,Goertz:RamifiedCase,PR03,PR09,SmithlingOddUnitary,Smithling:TopFlatnessEven,SmithlingEvenUnitary,Smithling15,Zachos,Pappas_Zachos:OrthogonalShimuraVarieties,Yang:SpinLocalModelsEven}.
Our results also imply topological flatness of the corresponding scheme-theoretic moduli problems in both equal and mixed characteristic.

The construction is effective and can be implemented by finite calculations in Weyl groups.
The proof relies on the quantum Bruhat graph and is outlined in the last part of this introduction, where we restrict for simplicity to groups that split over $\breve F$.

\subsection*{Reducing the number of representations}
The test set $J(\calG,\mu)$ is usually not minimal for the moduli description to hold.
Motivated by the known linear-algebraic conditions arising for local models of Shimura varieties, we seek a test set whose associated conditions exactly recover these linear-algebraic conditions in the corresponding cases.

First, we remove redundancies among the representations in $J(\calG,\mu)$.
For a vertex $a\in \mathrm{Vert}(\mathbf f)$ of the facet defining $\calG_{\breve O}$ and a weight $\chi$, two relevant Borel subgroups yield isomorphic representations if they have the same image in $W_a\backslash W_0/W_\chi$; see \Cref{redundancies between representations}.
Choosing a set $J_{a,\chi}$ of representatives for these double cosets, we define a subset of $J(\calG,\mu)$ by
\[
J(\calG,\mu)'
:= (J_1^{\mathrm{cen}}\times \{B_0\})\cup
\bigcup_{\substack{\chi\in J_1\setminus J_1^{\mathrm{cen}} \\
                   a\in\operatorname{Vert}(\mathbf f)}}
\{\chi\}\times J_{a,\chi},
\]
where $J_1^{\mathrm{cen}}=J_1\cap X^*(G)$ is the set of central characters in $J_1$.
Since every representation in $J(\calG,\mu)$ is isomorphic to one in $J(\calG,\mu)'$ by \Cref{rep reduction lemm}, the conclusion of \Cref{main theorem intro} remains valid for $J=J(\calG,\mu)'$.

Second, to compare our construction with known moduli descriptions of local models, we further reduce the test set.
For a fixed $\mu$, only certain faces of the dominant cone $X^*(T)^+_{\mathbb Q}$ are needed to recover the Schubert variety in the generic fiber; see \Cref{minimal test set prop}.
These faces determine a unique subset $J(\mu_{\mathrm{ad}})$ of the fundamental weights (\Cref{minimal test set}).
When $G_{\der}$ is simply connected, we choose arbitrary lifts of these weights to $X^*(T)^+$.
We also choose a minimal set of central characters in $X^*(G)$ that generates $X^*(G)_{\mathbb Q}$ as a $\mathbb Q_{\geq 0}$-monoid.
We denote the union of these central characters and the chosen lifts of $J(\mu_{\mathrm{ad}})$ by $J(\mu)\subset X^*(T)^+$.
In the definition of $J(\calG,\mu)'$, we then take $J_1\subset X^*(T)^+$ to be an inclusion-minimal subset satisfying the conditions in \Cref{J1 conditions} with respect to $J(\mu)$.
This yields a smaller test set ({\Cref{condition implies J1 conditions}}), denoted $J_{\min}(\calG,\mu)$.
If $G_{\breve F}$ is split, we may take $J_1=J(\mu)$.

The following theorem shows that $J_{\min}(\calG,\mu)$ gives the correct generic fiber and bounds the special fiber by the permissible locus of Kottwitz--Rapoport \cite{KR00}, without assuming that $G$ is essentially tamely ramified.
We expect $J_{\min}(\calG,\mu)$ to be minimal among the test sets for which both the generic-fiber equality and the permissible-locus bound hold.

\begin{theorem}[{\Cref{inclusion permissible locus theorem}}]
\label{permissible theorem intro}
Assume that $G_{\mathrm{der}}$ is simply connected and that $\pi_1(G)_I$ is torsion-free.
For $J=J_{\min}(\calG,\mu)$ as in \Cref{definition of J intro}, the inclusion $\mathbb M_{\calG,\mu}\subset\mathbb M_{\calG,\mu}^J$ is an equality on generic fibers, and the special fiber of $\mathbb M_{\calG,\mu}^J$ is contained in the permissible
locus.
Consequently, $\mathbb M_{\calG,\mu}^J=\mathbb M_{\calG,\mu}$ whenever the admissible and permissible loci agree.
\end{theorem}

Whenever the admissible and permissible loci agree, \Cref{permissible theorem intro} shows that $J_{\min}(\calG,\mu)$ suffices to define $\mathbb M_{\calG,\mu}$.
The required equality depends only on the échelonnage root system of $G_{\breve F}$, see \cite{Hai18}, the facet defining $\calG_{\breve O}$, and the image of the dominant representative of $\mu$ in the inertia coinvariants $X_*(T)_I^+$.
It is established in many cases \cite{Rap05,PRS13,HN02,Goertz:RamifiedCase,Smithling:PermissibleOn,SmithlingOddUnitary,SmithlingEvenUnitary,HH17,Schremmer24}, which we review in the companion paper \cite{LRVW:Minuscule}; see \Cref{permissible corollary}.
Under the remaining hypotheses of \Cref{permissible theorem intro}, this gives effective moduli descriptions even when $G$ is not essentially tamely ramified, and hence lies outside the scope of \Cref{main theorem intro}.

The conditions indexed by $J_{\min}(\calG,\mu)$ recover many of the linear-algebraic descriptions studied in the theory of local models of Shimura varieties.
These include linear and symplectic groups \cite{Goertz:LinearCase,Goertz:SymplecticCase}, split orthogonal groups \cite{Smithling:TopFlatnessEven,Pappas_Zachos:OrthogonalShimuraVarieties,Yang:SpinLocalModelsEven}, and ramified linear and symplectic groups \cite{PR03,Goertz:RamifiedCase}.
They also include ramified unitary groups in odd dimension and in even dimension with signature $(r,s)$ satisfying $\min(r,s)\leq 2$ \cite{PR09,SmithlingOddUnitary,SmithlingEvenUnitary}, as well as the ramified orthogonal cases of \cite{Pappas_Zachos:OrthogonalShimuraVarieties,YangZachosZhao:PELTypeD}.
Thus the test set $J_{\min}(\calG,\mu)$ gives a common representation-theoretic construction of these conditions.

The conclusion of \Cref{main theorem intro} need not hold for $J=J_{\min}(\calG,\mu)$: if $\mu$ is non-minuscule or $G_{\breve F}$ is not split, the special fiber of $\bbM_{\calG,\mu}^J$ may strictly contain the admissible locus.
When $\mu$ is minuscule and $G_{\breve F}$ is split, however, we show in the companion paper \cite{LRVW:Minuscule} that $\bbM_{\calG,\mu}^J=\bbM_{\calG,\mu}$ for $J=J_{\min}(\calG,\mu)$; see \Cref{refined conjecture minuscule}.
For example, we show that at Iwahori level, the exceptional groups $E_6$ and $E_7$ require $31$ and $64$ representation-theoretic conditions, respectively.

\begin{example}
Let $G=\GL_n$, $n\geq2$, and let
$\calG=\underline{\Aut}(\calV_\bullet)$ be the Iwahori group
scheme for the standard periodic lattice chain
\[
\calV_i
 =\mathrm{span}_O\{\pi^{-1}e_1,\ldots,\pi^{-1}e_i,e_{i+1},\ldots,e_n\}
 \quad(0\leq i<n),
\qquad
\calV_{i+n}=\pi^{-1}\calV_i,
\]
where $\pi$ is a chosen uniformizer of $O$.
Let $T$ be the standard diagonal torus and $B_0$ the upper
triangular Borel subgroup.
For the dominant characters $\omega_j=(1^{(j)},0^{(n-j)})$, $1\leq j<n$, and
$\delta=(1,\ldots,1)$, we put
$J_1:=\{\omega_1,\ldots,\omega_{n-1},\delta,-\delta\}
\subset X^*(T)^+$.
This choice of $J_1$ satisfies \Cref{J_1 condition} for every $\mu\in X_*(T)^+$.
The associated test set is
\[
J(\calG,\mu)=J_1\times\{B_0,\ldots,B_{n-1}\},
\]
where $B_i$ is the Borel subgroup determined by the cyclically
ordered basis $e_{i+1},\ldots,e_n,e_1,\ldots,e_i$.
It is the relevant Borel for the vertex corresponding to
$\calV_i$.
By \Cref{maximal parahoric rep cor}, the parahoric induced
module from \Cref{parahoric weyl module intro} satisfies
$\calV_{\omega_j,B_i}\simeq\bigwedge^j\calV_i$,
with $\calG$ acting through $\bigwedge^j\rho_i$, where
$\rho_i\colon\calG\to\GL(\calV_i)$ denotes the projection onto
the $i$-th lattice.
Thus the Borel determines the projection onto a lattice in
the chain, while the weight determines its exterior power.
For fixed $j$, the generic representations are all isomorphic
to $\bigwedge^jF^n$, but their integral forms
$\bigwedge^j\calV_i$ are non-isomorphic as
$\calG$-representations in general.

Points of $\Gr_\calG$ are periodic
$B_{\mathrm{dR}}^+$-lattice chains
$\calE_\bullet=\{\calE_i\}_{i\in\bbZ}$, and
$(\bigwedge^j\rho_i)_*\calE_\bullet=\bigwedge^j\calE_i$.
For any dominant cocharacter $\mu=(m_1,\ldots,m_n)$,
\Cref{Tannaka bound local model intro} applied to the
parahoric induced module $\calV_{\omega_j,B_i}$ for
$(\omega_j,B_i)\in J(\calG,\mu)$ gives the inclusion
\begin{equation}\label{general conditions}
\xi^{m_1+\cdots+m_j}\cdot
\bigwedge^j\bigl(\calV_i\otimes_OB_{\mathrm{dR}}^+\bigr)
 \subset\bigwedge^j\calE_i,
\qquad
0\leq i<n,\quad 1\leq j<n.
\end{equation}
The characters $\delta$ and $-\delta$ give the determinant
and its inverse, imposing equality for $j=n$.
Their integral representations are independent of $i$ up to
isomorphism, so this equality need only be imposed at $i=0$.
By \Cref{main theorem intro}, these conditions define
$\mathbb M_{\calG,\mu}$.

Assume that $\mu=(1^{(r)},0^{(n-r)})$ is minuscule, where
$1\leq r<n$, so $\mu_{\mathrm{ad}}=\omega_r^\vee$.
The classification of the minimal test set $J(\mu_\ad)$ in
\Cref{minuscule classification} implies that we may choose
\[
J_{\min}(\calG,\mu)
 \subset
\{(\delta,B_0),(-\delta,B_0)\}
\cup
\{(\omega_1,B_i),(\omega_{n-1}-\delta,B_i)\mid 0\leq i<n\}.
\]
The $\calG$-modules $\calV_{\omega_1,B_i}$ and
$\calV_{\omega_{n-1}-\delta,B_i}$ are $\calV_i$ and
its dual $\calV_i^\vee$, with actions through $\rho_i$ and
$\rho_i^\vee$, respectively.
The conditions in \Cref{general conditions}, together with
the determinant equality, simplify to
\begin{equation}
\label{inclusion GL intro}
\xi\cdot\bigl(\calV_i\otimes_OB_{\mathrm{dR}}^+\bigr)
 \subset\calE_i
 \subset\calV_i\otimes_OB_{\mathrm{dR}}^+,
\qquad
\bigwedge^n\calE_0
 =
\xi^r\cdot\bigwedge^n
 \bigl(\calV_0\otimes_OB_{\mathrm{dR}}^+\bigr),
\end{equation}
for $i=0,\ldots,n-1$.
By \Cref{permissible theorem intro} and the equality of the
admissible and permissible loci in this case
\cite[Theorem~1]{HN02}, these conditions still define
$\mathbb M_{\calG,\mu}$.
Reduction modulo $\xi$ translates \eqref{inclusion GL intro}
into the lattice-chain conditions of the Rapoport--Zink linear
local model \cite{RZ96} studied by G\"ortz
\cite{Goertz:LinearCase}.
\end{example}

\subsection*{Weight modules and the quantum Bruhat graph}
We now outline the proof of \Cref{main theorem intro}.
For simplicity, assume that $G_{\breve F}$ is split, and suppose that $G_{\der}$ is simply connected.
Then $\pi_1(G)_I=\pi_1(G)$ is torsion-free by \Cref{elementary group lemma}.
Let $\mu\in X_*(T)^+$ be the dominant representative.
Fix a hyperspecial vertex in the standard apartment over $\breve F$, and let $\mathbf a$ be the adjacent dominant alcove.
Assume that the facet $\mathbf f$ defining $\calG_{\breve O}$ lies the closure of $\mathbf a$.
For each vertex $a\in\operatorname{Vert}(\mathbf f)$, dominant weight $\chi\in X^*(T)^+$, and $v\in W_0$, consider the parahoric induced module
\[
\calV_{a,\chi,v}
=\Ind_{\calB^-_v}^{\calG_a}(v\chi)\in\Rep(\calG_a)
\]
over $\breve O$, associated with the Borel subgroup $B_v=vB_0v^{-1}$.
The extremal weight space of weight $v\chi$ defines a rank-one $\breve O$-submodule $\calV_{a,\chi,v}(v\chi)$.
For an element $x=\nu(\pi)u\in W$ of the Iwahori--Weyl group, define the elementary divisor $d_{a,\chi,v}(x,v\chi)\in\mathbb Z$ by
\[
x\cdot\calV_{a,\chi,v}(v\chi)
=\pi^{d_{a,\chi,v}(x,v\chi)}
 \cdot\calV_{a,\chi,v}(uv\chi),
\]
where $\pi$ is a fixed uniformizer of $F$.
The proof of \Cref{main theorem intro} rests on the following comparison with Schremmer's semi-affine quantum Bruhat graph weight \cite{Schremmer24}.

\begin{theorem}[{\Cref{prop:extremal-minor-reduction}}]\label{schremmer comparison intro}
Let $x=\nu(\pi)u\in W$ be an element in the Iwahori--Weyl group, and $\chi\in X^*(T)^+$. 
Then, for every vertex $a\in \mathrm{Vert}(\mathbf f)$ and every element $v\in W_0$ in the Weyl group, one has
\begin{equation}\label{schremmer comparison intro eq}
d_{a,\chi, v}(x,v\chi)
 =\left\langle\chi,\,
 v^{-1}u^{-1}\nu+
 \operatorname{wt}_a(v\Rightarrow uv)
 \right\rangle,
\end{equation}
where $\operatorname{wt}_a$ is Schremmer's semi-affine quantum Bruhat weight; see \Cref{lrvw-semi-affine-weight}.
\end{theorem}

The proof uses the graded character formula of Feigin--Makedonskyi \cite{FeiginMakedonskyi17,FeiginMakedonskyiTwisted17}. 
Let us explain how it implies \Cref{main theorem intro}.
We know the inclusion $\mathbb M_{\calG,\mu}\subset \mathbb M_{\calG,\mu}^J$ to be an isomorphism in the generic fiber. 
So, suppose we have an element $x\in W$ whose associated Schubert variety lies in the special fiber of $\mathbb M_{\calG,\mu}^J$. 
We want to show that the double $W_{\mathbf f}$-coset of $x$ is $\mu$-admissible. 
The defining lattice condition in \Cref{defining lattice inclusion unramified} applied to $\calV_{a,\chi,v}$ and the bundle represented by the element $x\in W$ gives after restriction to the weight module $\calV_{a,\chi,v}(v\chi)$:
$$
\pi^{\langle\chi,\mu\rangle} \cdot  \calV_{a,\chi,v}\left(uv\chi\right)\,\subset\, x\cdot  \calV_{a,\chi,v}(v\chi)=\pi^{d_{a,\chi, v}(x,v\chi)}\cdot  \calV_{a,\chi, v}\left(uv\chi\right),
$$
which is equivalent to $d_{a,\chi, v}(x,v\chi)\leq \langle\chi,\mu\rangle$.
Since the quantum Bruhat weight $\operatorname{wt}_{\mathbf f}$ at the facet $\mathbf f$ is the coefficient-wise supremum of the $\operatorname{wt}_a$, $a\in\mathrm{Vert}(\mathbf f)$, see \Cref{supremums formula}, this implies
\[
\left\langle\chi,
v^{-1}u^{-1}\nu+\operatorname{wt}_{\mathbf f}(v\Rightarrow uv)
\right\rangle=
\max_{a\in \mathrm{Vert}(\mathbf f)} d_{a,\chi,v}(x,v\chi)
\leq \langle\chi,\mu\rangle
\]
for $\chi$ central or a lift of a fundamental weight.
Since $G_\der$ is simply connected, these inequalities for a generating set $\chi\in X^*(T)^+$ imply
\[
v^{-1}u^{-1}\nu+
\operatorname{wt}_{\mathbf f}(v\Rightarrow uv)
\leq\mu, \quad \text{for all } v\in W_0.
\]
By Schremmer's admissibility criterion \cite[Proposition~4.35]{Schremmer24} at the vertices of $\mathbf f$ together with the vertexwise criterion of \cite{HH17}, see \Cref{schremmer-admissibility}, this is equivalent to the $\mu$-admissibility of the double $W_{\mathbf f}$-coset of $x$ 
and shows $\mathbb M_{\calG,\mu}=\mathbb M_{\calG,\mu}^J$.

\begin{remark}
The quantum Bruhat weights $\operatorname{wt}_a(v\Rightarrow uv)$ in \Cref{schremmer comparison intro} can be viewed as correction terms that account for the difference between admissibility and permissibility.
Replacing this formula by the weaker estimate of \Cref{estimate prop}, the defining lattice conditions still give
\[
\langle v\chi,\delta_a(x)\rangle
\leq d_{a,\chi,v}(x,v\chi)
\leq \langle\chi,\mu\rangle,
\qquad
\delta_a(x)=u^{-1}(a-x(a));
\]
see \Cref{key corollary permissible}.
These inequalities give the convex-hull bounds for permissibility required for \Cref{permissible theorem intro}, which we prove independently of \Cref{main theorem intro}, without assuming that $G$ is essentially tamely ramified.
\end{remark}

\subsection*{Acknowledgements}
This project grew out of the 2024 Oberwolfach Seminar \textit{Reduction of Arithmetic Varieties}, organized by three of the authors together with Tom Haines. 
During the seminar, we gave courses on the approaches to local models for Shimura varieties and moduli spaces of shtukas available at the time. 
We thank the MFO for its hospitality.

We thank Johannes Anschütz, Alexis Bouthier, Robert Cass, Ian Gleason, Tom Haines, Si Ying Lee, Michael Rapoport, and Ioannis Zachos for helpful discussions.
T.R.~thanks the research groups at ESAGA (University of Duisburg--Essen) and the organizers of the ReMoLD special semester in Padova for the opportunity to present preliminary results.

During the final stages of preparing this manuscript, we used generative AI to streamline our concepts, find literature, improve the exposition, and obtain feedback on earlier drafts.

The authors acknowledge funding by the Deutsche Forschungsgemeinschaft (DFG, German Research Foundation) TRR 326 \textit{Geometry and Arithmetic of Uniformized Structures}, project number 444845124.
T.R.~acknowledges funding by the European Research Council (ERC) under Horizon Europe (grant agreement 101040935), the LOEWE professorship in Algebra, project number LOEWE/4b/519/05/01.002(0004)/87, and the MSCA doctoral network ReMoLD (project number 101168795) under Horizon Europe.
E.V.~acknowledges funding by Germany's Excellence Strategy EXC 2044/2 -- 90685587, Mathematics Münster: Dynamics--Geometry--Structure and by a Leibniz prize.

\section{Notation and recollections}\label{section notation recollections}
We introduce our basic notation in \Cref{group theoretic notation} and recollect some facts on affine Grassmannians and local models in \Cref{recollection}. 
Our moduli description is based on the notion of maximal $\mu$-weights of representations defined in \Cref{maximal weights}.

Throughout, we fix a prime number $p$.  
Let $F$ be a non-archimedean local field, i.e., a finite field extension of $\Q_p$ or $\mathbb{F}_p\rpot{t}$, with ring of integers $O$, uniformizer $\pi$ and finite residue field $k=O/(\pi)$. 
Fix the completion of an algebraic closure $C/F$ with ring of integers $O_C$ and residue field $\bar k= k_C$. 
Denote by $\breve F\subset C$ the maximal unramified completed subextension with ring of integers $\breve O$ and the same residue field $\bar k= k_{\breve F}$.

The absolute Galois group of a subextension $K\subset C$ is denoted $\mathrm{Gal}_K$.
If $K=\breve F$, then we abbreviate $I=\mathrm{Gal}_{\breve F}$ for the inertia group of $F$.

\subsection{Group theoretic notation}\label{group theoretic notation}
Let $G$ be a connected reductive $F$-group, $\mu$ a conjugacy class of cocharacters of $G_C$ and let $\calG$ be a parahoric $O$-group scheme with generic fiber $G$.
The reflex field of $\mu$ is a finite, separable field extension $E=E_\mu$ of $F$ with ring of integers $O_E$ and residue field $k_E$.
Let $E^{\mathrm{ur}}$ be the maximal unramified intermediate extension of $F\subset E$ with ring of integers $O_{E^{\mathrm{ur}}}$ and the same residue field $k_{E^{\mathrm{ur}}}=k_E$.

\subsubsection{Auxiliary subgroups}
Also, we fix the following auxiliary subtori of $G$, defined over $F$:
\begin{equation}
\label{auxiliary.groups}
	A \subset S \subset T \subset G
\end{equation}
Here, $A$ is a maximal $F$-split torus in $G$, $S$ a maximal $\breve F$-split torus, $T$ the centralizer of $S$ in $G$. Such subtori exist \cite[Proposition~9.3.4]{KP23}, and $T$ is a maximal torus because $G_{\breve F} := G \otimes_F {\breve F}$ is quasi-split by Steinberg's theorem. Moreover, we can choose a Borel subgroup $B_0$ of $G_{\breve F}$ containing the maximal $\breve F$-split subtorus $S_{\breve F}$ of $G_{\breve F}$. Then we have
\begin{equation}
\label{auxiliary.groups2}
	T_{\breve F} \subset B_0 \subset G_{\breve F}
\end{equation}
by \cite[Proposition~2.4.9]{KP23}. 
In addition, $A\subset S$ and hence $T$ are chosen such that their connected Néron lft-models $\calA\subset \calS\subset \calT$ over $O$ are contained in $\calG$; see \cite[Proposition~5.1.10]{BT84}.

\subsubsection{Weyl groups}\label{Weyl group section}
The Iwahori-Weyl group for $(G,T)$ and its subgroup associated to $\calG$ are the groups 
$$
W = N_S(\breve F)/\calT(\breve O)\;\;\text{and}\;\; W_\calG=\big(N_S(\breve F)\cap \calG(\breve O)\big)/\calT(\breve O),
$$
where $N_S$ is the normalizer of $S$ in $G$. 
The Iwahori decomposition $G(\breve F)=\calG(\breve O)N_S(\breve F)\calG(\breve O)$ induces the bijection \cite[Proposition~8]{HR08}:
\begin{equation}\label{Iwahori decomposition}
	W_\calG\backslash W/W_\calG \cong \calG(\breve O)\backslash G(\breve F)/\calG(\breve O).
\end{equation}
The Iwahori-Weyl group admits two presentations \cite[Proposition~13, Lemma~14]{HR08} depending on the choice of an alcove, respectively a special vertex, in the apartment of the Bruhat--Tits building of $G_{\breve F}$ defined by $S_{\breve F}$:
\begin{align}
W=W_{\mathrm{aff}}\rtimes \pi_1(G)_I, \label{Bruhat presentation} \\
W=X_*(T)_I\rtimes W_{0} \label{translation presentation},
\end{align}
where $W_{\mathrm{aff}}$ is the affine Weyl group, $\pi_1(G)_I$ the inertia coinvariants of Borovoi's fundamental group, $X_*(T)_I$ the inertia coinvariants of the geometric cocharacter lattice which are identified with a subgroup of $W$ via Kottwitz' isomorphism \cite[Equation~(7.2.1)]{Kot97}
\begin{equation}\label{coinvariants identification}
	X_*(T)_I \overset{\cong}{\longto} T(\breve F)/\calT(\breve O), \;\; \nu_I\mapsto \nu_I(\pi),
\end{equation}
and $W_{0}=N_S(\breve F)/T(\breve F)$ the (relative) Weyl group of $(G_{\breve F},S_{\breve F})$. 
The Weyl group of $(G_C,T_C)$ is denoted $W_{0,C}=N_T(C)/T(C)$ where $N_T$ is the normalizer of $T$ in $G$. As $T$ is the centralizer of $S$ in $G$, one has an inclusion $N_S\subset N_T$ as subgroups in $G$.
In particular, $W_0$ is naturally a subgroup of $W_{0,C}$.

The following lemma is used in \Cref{root lemma twisted} below.

\begin{lemma}\label{elementary group lemma}
Assume that the derived group $G_{\mathrm{der}}$ is simply connected. 
Then, the following hold:
\begin{enumerate}
\item \label{elementary group lemma 1}
	The abelian group $\pi_1(G)$ is isomorphic to $X_*(G/G_\der)$, hence torsion-free. 
\item \label{elementary group lemma 2}
	If $\pi_1(G)_I$ is torsion-free, so is $X_*(T)_I$. 
	In this case, the following sequence is exact:
	\begin{equation}\label{inertia invariants sequence}
	0\to X^*(G)^I\to X^*(T)^I\to X^*(T_\der)^I\to 0,
	\end{equation}
	where $T_\der=T\cap G_\der$. 
\end{enumerate}
\end{lemma}
\begin{proof}
\Cref{elementary group lemma 1}:
Note that $T_{\mathrm{der}}$ is a maximal torus in $G_\der$; see \cite[XXII, Proposition~6.2.8]{SGA3}.
By simply connectedness of $G_{\mathrm{der}}$ and the definition of $\pi_1(G)$, the sequence of abelian groups
\begin{equation}\label{elementary group lemma ses}
0\to X_*(T_{\mathrm{der}})\to X_*(T)\to \pi_1(G)\to 0
\end{equation}
is exact. 
In particular, $\pi_1(G)=X_*(T/T_\der)=X_*(G/G_\der)$; see \cite[XXII, Section~6.2]{SGA3}.

\Cref{elementary group lemma 2}:
As the derived group is simply connected, $X_*(T_{\mathrm{der}})$ is an induced inertia module \cite[Proposition~4.4.16]{BT84} and hence the coinvariants $X_*(T_{\mathrm{der}})_I$ are torsion-free. 
Therefore, the exact sequence \eqref{elementary group lemma ses} induces a short exact sequence on $I$-coinvariants and if $\pi_1(G)_I$ is torsion-free, so is $X_*(T)_I$.
In addition, the sequence dual to \Cref{elementary group lemma ses} gives an exact sequence
\begin{equation}\label{exact sequence dual}
0\to X^*(G)\to X^*(T)\to X^*(T_{\mathrm{der}})\to 0,
\end{equation}
using $X^*(G)=X^*(G/G_\der)$.
Passing to $I$-invariants gives the exact sequence in \Cref{inertia invariants sequence} because $\mathrm{H}^1(I,M^\vee)=\Hom_{\mathbb Z}((M_I)_{\mathrm{tors}},\mathbb Q/\mathbb Z)$ vanishes when applied to $M=\pi_1(G)=X_*(G/G_\der)$.
\end{proof}

\subsubsection{Apartment} \label{apartment sec}
The enlarged apartment $\mathscr A=\mathscr A(G,S,\breve F)$ of the Bruhat--Tits building of $G_{\breve F}$ corresponding to $S_{\breve F}$ is a torsor under $X_*(S)_{\mathbb R} = X_*(T)_{I, \mathbb R}:=X_*(T)_I\otimes_{\mathbb Z}\mathbb R$.
We use the sign convention of Bruhat--Tits \cite[Section~4.2.3]{BT84} (cf.~also \cite[Axiom~4.1.4, A~1]{KP23}) for the affine action of $N_S(\breve F)$ on $\mathscr A$: under the
identification \eqref{coinvariants identification}, the element
$\nu_I(\pi)\in T(\breve F)/\calT(\breve O)$ acts by
\begin{equation}\label{BT translation convention}
 \nu_I(\pi)\cdot a=a-\nu_I,
\end{equation}
for all $a\in \mathscr A$, see \Cref{plus translation convention} for the plus sign convention. 
Thus, after choosing a special vertex in $\mathscr A$ and using it as the origin to identify $\mathscr A=X_*(T)_{I,\mathbb R}$ and $W=X_*(T)_I\ltimes W_{0}$, an element
$x=\nu_I(\pi)u\in X_*(T)_I\rtimes W_0$ acts by
\begin{equation}\label{BT affine Weyl action}
 x\cdot a=u(a)-\nu_I.
\end{equation}
This leads to the following correspondence between Weyl chambers and Borel subgroups:
Let $B\subset G_{\breve F}$ be a Borel subgroup containing $S_{\breve F}$, and let $\Phi_B^+$ be the subset of the relative roots $\Phi(G,S)$ occurring in the unipotent radical $U$ of $B$.  
For $a\in\mathscr A$, set
\begin{equation}\label{B dominant cone}
 \mathscr C_B(a)
 :=\{b\in \mathscr A\;\mid\; \alpha(b-a)\geq 0\text{ for all }\alpha\in\Phi_B^+\}.
\end{equation}
The interior $\mathscr C_B(a)^\circ$ is the Weyl chamber corresponding to $B$ with apex $a$. 
The assignment $B\mapsto\mathscr C_B(a)^\circ$ is a bijection from the set of Borel subgroups of $G_{\breve F}$ containing $S_{\breve F}$ onto the set of Weyl chambers in $\mathscr A$ with apex $a$. 

\begin{lemma}\label{unipotent inclusion lemma}
Let $a,b\in \mathscr A$. 
Denote by $\calU_a$ and $\calU_b$ the respective $\breve O$-models of $U$.
If $b\in \mathscr C_B(a)$, then $\calU_a(\breve O)\subset \calU_b(\breve O)$.
In particular, the identity on $U$ extends to a map $\calU_a\to \calU_b$ of $\breve O$-group schemes. 
\end{lemma}
\begin{proof}
Let $\Phi_{B,\mathrm{nd}}^+$ denote the set of non-divisible positive relative roots determined by $B$. 
For $c\in\mathscr A$, $\alpha\in\Phi_{B,\mathrm{nd}}^+$, and $r\in\mathbb R$, write $U_{\alpha,c,r}\subset U_\alpha(\breve F)$ for the corresponding Bruhat--Tits filtration subgroup \cite[Definition~6.1.2]{KP23}.
Put $v=b-a$. 
By \cite[Lemma~6.1.6]{KP23}, the dependence of the filtration on the point of the apartment is given by
\[
U_{\alpha,a,0}
 =
U_{\alpha,a+v,\alpha(v)}
 =
U_{\alpha,b,\alpha(b-a)}.
\]
Since $b\in\mathscr C_B(a)$, one has $\alpha(b-a)\geq 0$ for all $\alpha\in\Phi_B^+$.
The filtration $r\mapsto U_{\alpha,b,r}$ is decreasing by definition, and hence
\begin{equation}\label{inclusion root models}
U_{\alpha,a,0}
 =
U_{\alpha,b,\alpha(b-a)}
 \subset
U_{\alpha,b,0}
\qquad
\text{for every }\alpha\in\Phi_{B,\mathrm{nd}}^+.
\end{equation}
The groups $U_{\alpha,c,0}$ are schematic, with smooth affine $\breve O$-models $\calU_{\alpha,c,0}$, see \cite[Lemma~8.3.7]{KP23}. 
Moreover, \cite[Proposition~8.3.14(2)--(3)]{KP23} shows that, after fixing any ordering of $\Phi_{B,\mathrm{nd}}^+$, multiplication identifies $\calU_c$ with the ordered product of these root-group models:
\begin{equation}\label{root group isomorphism}
\prod_{\alpha\in\Phi_{B,\mathrm{nd}}^+}
\calU_{\alpha,c,0}
\;\xrightarrow{\ \sim\ }\;
\calU_c .
\end{equation}
Taking $\breve O$-points and using \Cref{inclusion root models} factor by factor therefore gives
\[
\calU_a(\breve O)
 =
\prod_{\alpha\in\Phi_{B,\mathrm{nd}}^+}
U_{\alpha,a,0}
\subset
\prod_{\alpha\in\Phi_{B,\mathrm{nd}}^+}
U_{\alpha,b,0}
 =
\calU_b(\breve O).
\]

Finally, $\calU_a$ and $\calU_b$ are smooth affine integral models of the same $\breve F$-group $U$. 
Since the identity of $U$ carries $\calU_a(\breve O)$ into $\calU_b(\breve O)$, it extends uniquely to $\calU_a\longrightarrow\calU_b$ by the étoffé property \cite[Corollary~2.10.10]{KP23}.
\end{proof}

\begin{remark}\label{split apartment conventions}
Suppose that $G$ is split, i.e., $A=S=T$ is a split maximal torus.  
Choose a hyperspecial vertex $o\in \mathscr A$ and identify $\mathscr A=o+X_*(T)_\mathbb R$. 
Then, the
affine roots are the functions $\alpha+n$, with $\alpha\in\Phi(G,T)$ and
$n\in\mathbb Z$.
For each $\al\in \Phi_B^+$, choose an isomorphism $x_\al\colon \mathbb G_{a,O}\cong \calU_{\alpha,o,0}$ with the $O$-model of $U_\al$ at the vertex $o$. 
Then, $x_\alpha(\pi^nO)$ fixes the half-apartment $\{\,a\in\mathscr A\mid \alpha(a)+n\geq0\,\}$ and we have
$$
U_{\alpha,a,0}(O)
=
U_{\alpha,o,-\alpha(a)}(O)
 =x_\alpha\!\left(\pi^{\lceil-\alpha(a)\rceil}O\right)
 =x_\alpha\!\left(\pi^{-\lfloor\alpha(a)\rfloor}O\right).
$$
\end{remark}

\begin{remark}\label{plus translation convention}
In parts of the literature (e.g., \cite{HN02, HH17}), the element $\nu_I(\pi)$ is
instead made to act by translation by $+\nu_I$.  To pass to that convention,
choose a special vertex $o$ and conjugate the above affine action by the
opposition involution
\[
 \iota_o\colon o+a\longmapsto o-a,
\]
i.e., define $n\star_o a:=\iota_o(n\cdot\iota_o(a))$ for all $n\in N_S(\breve F)$, $a\in \mathscr A$.  
Then $W_0$ has the same linear action and $ \nu_I(\pi)\star_o a=a+\nu_I$.
All affine-geometric statements transport along $\iota_o$: affine roots
are replaced by their composites with $\iota_o$, and
$ \iota_o\bigl(\mathscr C_B(a)\bigr)
 =\mathscr C_{B^-}\bigl(\iota_o(a)\bigr)$.
 Thus, if one keeps the same algebraic positive-root convention, every
chamber label $B$ must be replaced by $B^-$ when passing from the
$-\nu_I$ convention to the $+\nu_I$ convention.  Equivalently, one may keep
the labels and simultaneously replace $\Phi_B^+$ by $-\Phi_B^+$.
\end{remark}

\subsection{Recollections on $B_{\mathrm{dR}}^+$-affine Grassmannians and v-sheaf local models}\label{recollection}
	By the term ``v-sheaf'', we mean a set-valued functor on the site $\Perf$ of perfectoid spaces in characteristic $p$ equipped with the v-topology. 
	Since all functors considered below are v-sheaves, it is enough to define them on the basis of the topology given by affinoid perfectoid spaces $\Spa (R,R^+)$.
	We denote by $\Spd(O)$ the diamantine spectrum of $O$.

	\subsubsection{De Rham period rings}
	For any $\Spa(R,R^+)\in \Perf$ equipped with a map to $\Spd (O)$, we denote by 
	$$B_{\mathrm{dR}}^+(R^\sharp)\;\;\;\;\text{and}\;\;\;\; B_{\mathrm{dR}}(R^\sharp)=B_{\mathrm{dR}}^+(R^\sharp)[\xi^{-1}]$$ 
	the rings of de Rham periods formed using $O$-Witt vectors.
	Also, set $B_{\mathrm{dR}}^+=B_{\mathrm{dR}}^{+}(C)$ and $B_{\mathrm{dR}}=B_{\mathrm{dR}}(C)$.
	These rings are non-canonically isomorphic to $C\pot{\xi}$ and $C\rpot{\xi}$, respectively.  
	If $\Spa(R,R^+)\to \Spd(O)$ factors through $\Spd(k)$, then
$B_{\mathrm{dR}}^+(R^\sharp)=W_O(R)$ and
$B_{\mathrm{dR}}(R^\sharp)=W_O(R)[\pi^{-1}]$, where, for any
perfect $k$-algebra $R$, we set
$$
W_O(R)=
\begin{cases}
W(R)\otimes_{W(k)}O, & \operatorname{char}(F)=0,\\
R\pot{\pi}, & \operatorname{char}(F)=p.
\end{cases}
$$
In equal characteristic, we use the identification $O=k\pot{\pi}$
determined by the chosen uniformizer.

	\subsubsection{Beilinson--Drinfeld Grassmannian}
	The Beilinson--Drinfeld Grassmannian $\Gr_\calG$ is the set-valued functor that associates to any $\Spa(R,R^+)\to \Spd(O)$ in $\Perf$ isomorphism
classes of pairs $(\calE,\al)$ where $\calE $ is a $\calG$-bundle on $\Spec(B_{\mathrm{dR}}^+(R^\sharp))$ and $\al$ is a trivialization over $\Spec( B_{\mathrm{dR}}(R^\sharp))$. 
	This defines a v-sheaf, and the map $\Gr_\calG\to \Spd(O)$ is representable in ind-(proper spatial diamonds).
	Its generic and special fibers are given by the $B_{\mathrm{dR}}^+$-affine Grassmannian $\Gr_G$ and the v-sheaf associated with the Witt vector affine flag variety $\Fl_\calG$ respectively:
	\begin{equation}\label{BD_fibers}
	\Gr_\calG\times_{\Spd(O)} \Spd F \cong \Gr_G\;\;\;\;\text{and}\;\;\;\; \Gr_\calG\times_{\Spd(O)} \Spd(k) \cong \Fl_\calG^\diamondsuit.
	\end{equation}
	For a complete subextension $K\subset C$ with residue field $\kappa$, we denote by $\Gr_{\calG, O_K}$,  $\Gr_{\calG, K}=\Gr_{G,K}$ and $\Gr_{\calG,\kappa}=\Fl_{\calG,\kappa}^\diamondsuit$ the base changes to $\Spd(O_K)$, $\Spd(K)$ and $\Spd(\kappa)$, respectively.

	\subsubsection{Loop groups}
	The loop groups of $\calG$ are the group-valued functors over $\Spd(O)$ given by $L_{O}\calG(R,R^+)=\calG(B_{\dR}(R^\sharp))$ and $L^+_{O}\calG(R,R^+)=\calG(B_{\dR}^+(R^\sharp))$, respectively. 
	Then, $L_{O}^+\calG\subset L_{O}\calG$ is a subgroup functor and 
	$$
	\Gr_\calG\cong L_{O}\calG/L_{O}^+\calG,
	$$
	compatibly with the base changes in \Cref{BD_fibers}.
	The base change of the loop groups to $\Spd(F)$, respectively to $\Spd(k)$ is given by  
	$$
	L^{(+)}_FG=L_{O}^{(+)}\calG\times_{\Spd(O)}\Spd(F),\;\;\text{respectively}\;\; (L^{(+)}_{k}\calG)^\diamondsuit=L_{O}^{(+)}\calG\times_{\Spd(O)}\Spd(k),
	$$
	where $L^{(+)}_{k}\calG$ is the loop group defined on perfect
$k$-algebras $R$ by
$L^+_{k}\calG(R)=\calG(W_O(R))$ and
$L_{k}\calG(R)=\calG(W_O(R)[\pi^{-1}])$, respectively.

	\subsubsection{Schubert varieties in $\Gr_G$} \label{subsection schubert varieties in Gr_G}
	The $B_{\mathrm{dR}}^+$-affine Grassmannian $\Gr_G\to \Spd(F)$ is representable by an ind-(proper spatial diamond). 
	We denote by $L_CG$, $L_C^+G$ and $\Gr_{G,C}$ the base changes to $\Spd C$.
	The Cartan decomposition $G(B_{\mathrm{dR}})=G(B_{\mathrm{dR}}^+)T(B_{\mathrm{dR}})G(B_{\mathrm{dR}}^+)$ induces a bijection
	$$
	X_*(T)^+\overset{\la\mapsto \la(\xi)}{\cong} G(B^+_{\mathrm{dR}})\backslash G(B_{\mathrm{dR}})/G(B^+_{\mathrm{dR}}).	$$
	We obtain a stratification on the underlying topological spaces
	\begin{equation}
	\mathrm{Gr}_{G,C} = \bigcup_{\la \in X_*(T)^+} {\mathrm{Gr}}^\circ_{G,C,\la},
	\end{equation}
	where $\mathrm{Gr}_{G,\la,C}^\circ$ is the $L^+_CG$-orbit of $\la(\xi)$, i.e., the v-sheaf theoretic image of the orbit map.
	Its v-sheaf closure $\Gr_{G,\la, C}\subset \Gr_{G,C}$ is a union of orbits, which are locally closed. 
	For any $\la,\la'\in X_*(T)^+$, we have \cite[Corollary~4.6]{AGLR22}
	\begin{equation}\label{bruhat order on schubert varieties}
	\Gr_{G,\lambda',C}^\circ\subset \Gr_{G,\la, C} \iff \lambda' \leq \la
	\end{equation}
	in the Bruhat order on $X_*(T)^+$, i.e., $\la-\la'$ can be written as a sum of positive coroots with $\mathbb Z_{\geq 0}$-coefficients. 
	In particular, the (unique) representative of $\mu$ in $X_*(T)^+$ defines an orbit and its closure; these descend to the reflex field $\Spd(E)$ and are denoted $\Gr_{G,\mu}^\circ\subset \Gr_{G, \mu}$ inside $\Gr_{G,E}$.

	\subsubsection{Schubert varieties in $\Fl_\calG$}\label{witt vector affine flag section}
	The Witt vector affine flag variety $\Fl_\calG\to \Spec(k)$ is representable by an ind-(perfectly projective perfect variety) \cite{bhatt_scholze_projectivity_of_the_witt_vector_affine_grassmannian}. 
	We denote by $L_{\bar k}\calG$, $L_{\bar k}^+\calG$ and $\Fl_{\calG,\bar k}$ the base changes to $\Spec(\bar k)$.
	The Iwahori decomposition (\Cref{Iwahori decomposition}) induces a stratification on the underlying topological spaces
	$$
	\Fl_{\calG, \bar k}=\bigcup_{[x]\in W_\calG\backslash  W/W_\calG}\Fl_{\calG,\bar k, [x]}^\circ,
	$$
	where $\Fl_{\calG,\bar k, [x]}^\circ$ is the $L^+_{\bar k}\calG$-orbit of a representative of $[x]$ in $\Fl_\calG(\bar k)$. 
	The orbit closure $\Fl_{\calG,\bar k, [x]}$ is a perfectly projective perfect variety. 
	For any $[x],[y]\in W_\calG\backslash W/W_\calG$, we have
	\begin{equation}\label{orbit closure flag}
	\Fl_{\calG,\bar k, [y]}^\circ\subset \Fl_{\calG,\bar k, [x]} \iff [y] \leq [x]
	\end{equation} 
	in the Bruhat order on $W_\calG\backslash W/W_\calG$.
	The proof of \eqref{orbit closure flag} uses Demazure resolutions, as in the equal-characteristic case \cite[Proposition~2.8]{Ri13}.

	\subsubsection{Local models}\label{local models intro}
	The \textit{v-sheaf local model} $\mathbb M_{\calG,\mu}$ is the v-sheaf closure of $\Gr_{G,\mu}$ in $\Gr_{\calG,O_E}$, i.e., the smallest closed sub-v-sheaf containing $\Gr_{G,\mu}$.
	Then, $\mathbb M_{\calG,\mu}\to \Spd(O_E)$ is a flat, proper, $\pi$-adic kimberlite that carries an action of $L^+_{O_E}\calG$, see \cite[Propositions~4.13, 4.14]{AGLR22}. 
	Its generic and special fibers are given by 
	$$
	\mathbb M_{\calG,\mu, E} = \Gr_{G,\mu}\;\;\;\;\text{and}\;\;\;\; \mathbb M_{\calG,\mu, k_E}=A_{\calG,\mu}^\diamondsuit.
	$$
	Here, $A_{\calG,\mu}\subset \Fl_{\calG,k_E}$ is the admissible locus characterized by 
	$$
	A_{\calG,\mu,\bar k}=\bigcup_{\la} \Fl_{\calG,\bar k, [\lambda_I(\pi)]},
	$$
	where $\la$ runs through the elements in $X_*(T)$ in the conjugacy class $\mu$. 
	If $\mu$ is minuscule, then the v-sheaf local model $\mathbb M_{\calG,\mu}$ is representable by a flat, projective, and weakly normal $O_E$-scheme $M_{\calG,\mu}$, called simply the \textit{local model} \cite[Theorem~1.2]{AGLR22}.

	The following lemma allows us to assume that $G_\der$ is simply connected and that $\pi_1(G)_I$ is torsion-free without changing the v-sheaf local model.
	
	\begin{lemma} 	\label{central extension lemma}
	There exists a central extension
	\begin{equation}\label{coflasque z-extension}
	1\to Z\to \tilde G\to G\to 1
	\end{equation}
	of $F$-groups and a lift $\tilde \mu$ in $\tilde G$ of $\mu$ with the following properties: $Z$ is an induced torus, the derived group $\tilde G_\der$ is simply connected, $\pi_1(\tilde G)_I$ is torsion-free, and $\tilde \mu$ has the same reflex field as $\mu$. 
	In particular, the induced map $\tilde \calG\to\calG$ on matching parahoric $O$-models gives an isomorphism of v-sheaves:
	$$\mathbb M_{\tilde \calG,\tilde \mu}\overset{\sim}{\longto}\mathbb M_{\calG,\mu}.$$ 
	\end{lemma}
	\begin{proof}
	The extension in \Cref{coflasque z-extension} exists by \cite[Theorem~A.4.1(b)]{Cesnavicius2022Torsors} and the resulting map on v-sheaf local models is an isomorphism by \cite[Proposition~4.16]{AGLR22}.
	It remains to show that a lift $\tilde \mu$ with the same reflex field as $\mu$ exists:
	The set $\mathscr L_\mu$ of conjugacy classes of geometric cocharacters in $\tilde G$ lifting $\mu$ forms a $\Gal_E$-equivariant torsor under $X_*(Z)$. 
	This determines a class in Galois cohomology
	$$
	[\scrL_\mu] \in \mathrm{H}^1(\Gal_E,X_*(Z)),
	$$
	which vanishes because $Z$ is induced. 
	Hence, $\mathscr L_\mu$ contains a $\Gal_E$-fixed conjugacy class.
	Consequently, $\Gal_E\subset \Gal_{\tilde E}$, where $\tilde E$ is the reflex field of $\tilde \mu$.
	The reverse inclusion is automatic because $\tilde \mu$ lifts $\mu$, so $\Gal_E= \Gal_{\tilde E}$.
	This implies $E=\tilde E$ by Galois theory. 
	\end{proof}

\subsection{Maximal $\mu$-weights of representations}\label{maximal weights}
The moduli description of the v-sheaf local model (\Cref{definition local model}) is given in a Tannakian style and relies on the following notion. 
For any $O$-algebra $R$, we denote by $\Rep(\calG_R)$ the category of algebraic representations of $\calG_R:=\calG\otimes_{O}R$ on finite projective $R$-modules.  
		
\begin{definition}\label{representation weight}
Assume $\mathrm{Frac}(R)\subset C$, e.g., $R=K$ or $O_K$ for some intermediate field extension $F\subset K\subset C$.
For any nonzero representation $\calV\in \Rep(\calG_R)$, one defines the integer
\begin{equation}
	\mu_{\mathrm{max}}(\calV)=\max_{\chi\in X^*(T)}\{\langle\chi,\dot\mu\rangle \;|\; \text{$\chi$ is a $T_{C}$-weight of $\calV_{C}$} \}\in \mathbb Z,
\end{equation}
where $\dot{\mu}\in X_*(T)$ is a representative of $\mu$. 
This integer is well-defined independently of the choice of a representative $\dot{\mu}$ and a maximal torus in $G_{C}$.
\end{definition}

Here we collect some elementary properties for later use:

\begin{lemma}\label{elementary weight}
Let $R$ be as in \Cref{representation weight}, and $\calV, \calV', \calW\in \Rep(\calG_R)$ nonzero.
\begin{enumerate}
	\item \label{elementary weight 1}
		The integer $\mu_{\mathrm{max}}(\calV)$ only depends on the image of $\calV|_{T_C}$ in the Grothendieck group $K_0\Rep(T_{C})=\bbZ[X^*(T)]$.
	\item \label{elementary weight 2}
		There exists a highest $T_C$-weight $\chi$ of $\calV_{C}$ such that $$\mu_{\max}(\calV)=\langle\chi,\dot\mu\rangle,$$ where $\dot{\mu}\in X_*(T)^+$ is the dominant representative of $\mu$. 
	\item \label{elementary weight 4}
		If $R=O_K$ for a finite separable extension $K\subset C$ of the reflex field $E$, then 
		$$ \mu_{\max}(\calV)=\mu_{\max}(\Res_{O_K/O_E} \calV),$$
		where $\Res_{O_K/O_E} \calV$ is the $\calG_{O_E}$-representation on $\calV$ viewed as an $O_E$-module.
	\item \label{elementary weight 5}
		One has $$\mu_{\max}(\calV\otimes\calV')=\mu_{\max}(\calV)+\mu_{\max}(\calV').$$
	\item \label{elementary weight 6}
		If $0\to \calV\to \calW \to \calV'\to 0$ is a short exact sequence in $\Rep(\calG_R)$, then $$\mu_{\max}(\calW)=\max\{\mu_{\max}(\calV),\mu_{\max}(\calV')\}.$$ 
		In particular, this applies to $\calW=\calV\oplus \calV'$.
\end{enumerate}
\end{lemma}
\begin{proof}
\Cref{elementary weight 1}: This holds by definition. 

\Cref{elementary weight 2}: The weights of $\calV|_{T_C}$ and the pairing $\langle-,-\rangle$ are invariant under the action of the absolute Weyl group $W_{0,C}$. So, we may choose $\dot{\mu}\in X_*(T)^+$ to be the dominant representative of $\mu$. 
Then, \Cref{elementary weight 2} follows from $w\chi \leq \chi$ for any $\chi\in X^*(T)^+$ and $w\in W_{0,C}$.  

\Cref{elementary weight 4}:
Let $\Sigma$ be the set of highest weights of $\calV|_{T_C}$.
Then, the highest weights of $(\Res_{O_K/O_E}\calV)_{T_C}$ are the Galois orbits $\cup_{\gamma\in \mathrm{Gal}_E}\gamma*\Sigma$ for the Galois action on the based root datum.
Since the dominant representative $\dot\mu$ is fixed by $\mathrm{Gal}_E$, we obtain
\[
\mu_{\max}(\Res_{O_K/O_E}\calV)
=\max_{\chi\in\Sigma,\ \gamma\in\mathrm{Gal}_E}
 \{ \langle\gamma*\chi,\dot\mu\rangle\}
=\max_{\chi\in\Sigma}\{\langle\chi,\dot\mu\rangle\}
=\mu_{\max}(\calV).
\]

\Cref{elementary weight 5}: We use that the $T_{C}$-weights of $\calV_{C}\otimes\calV'_{C}$ are of the form $\chi +\chi'$ for $T_{C}$-weights $\chi$ and $\chi'$ of $\calV_{C}$ and $\calV'_{C}$, respectively. 

\Cref{elementary weight 6}: We restrict the short exact sequence to $T_{C}$, which then splits by semisimplicity of $\Rep(T_{C})$.
Then, each $T_{C}$-weight of $\calW_{C}$ appears in $\calV_{C}$ or $\calV'_{C}$, so \Cref{elementary weight 6} holds.
\end{proof}

\section{Moduli description of Schubert varieties}
In \Cref{moduli-description-generic-fiber}, we give a moduli description of the Schubert varieties $\Gr_{G, \mu}\to \Spd(E)$ in the $B_{\mathrm{dR}}^+$-affine Grassmannian $\Gr_G$ from \Cref{subsection schubert varieties in Gr_G}.
This will be used to identify the generic fiber of the moduli functor for the v-sheaf local models.

	\subsection{Moduli functor}
	Given a representation $\rho\colon G_{E}\to \GL(V)$ in $\Rep(G_{E})$, we denote the induced map of v-sheaves by 
	$$
	\rho_*\colon \Gr_{G,E}=\Gr_{G_E}\to \Gr_{\GL(V)},
	$$
	where we use that the formation of $B_{\mathrm{dR}}^+$-affine Grassmannians commutes with étale base change. 
	For any integer $n\in \mathbb Z$, we denote by $\Gr_{\GL(V)}^{\geq n}$ the closed sub-v-sheaf of $\Gr_{\GL(V)}$ parametrizing $B^+_{\mathrm{dR}}$-lattices inside $V\otimes_{E} B_{\mathrm{dR}}$ that contain $\xi^n\cdot \left( V\otimes_EB^+_{\mathrm{dR}}\right)$. 

\begin{definition}\label{Schubert functor definition}
	Let $J\subset \Rep(G_E)$ be a subset. 
	Define the sub-functor $\mathrm{Gr}_{G, \mu}^J\subset \Gr_{G,E}$ given on $\Spa(R,R^+)\to \Spd(E)$ by all $\underline{\calE}=[(\calE,\al)]\in \Gr_G(R,R^+)$ such that for all representations $\rho\colon G_{E}\to \GL(V)$ in $J$ one has 
	\begin{equation}\label{Tannaka.bound.gen}
		\xi^{\mu_{\max}(V)}\cdot \left(V\otimes_EB^+_{\mathrm{dR}}(R^\sharp)\right)\subset \rho_*\underline{\calE} 
	\end{equation}
	as $B^+_{\mathrm{dR}}(R^\sharp)$-lattices inside $V\otimes_{E} B_{\mathrm{dR}}(R^\sharp)$. 
	Equivalently, one has 
	\begin{equation}\label{tautological Schubert equality}
\mathrm{Gr}_{G, \mu}^J= \bigcap_{(\rho,V) \in J} \rho_*^{-1}\left(\mathrm{Gr}_{\mathrm{GL}(V)}^{\geq \mu_{\max}(V)}\right).
\end{equation}
\end{definition}

\begin{remark}\label{remark reverse inclusion}
	\Cref{Tannaka.bound.gen} applied to the dual representation $V^\vee$ gives the inclusion
	\begin{equation}\label{Tannaka.bound.dual}
		\rho_*\underline{\calE}  \subset {\xi}^{-{{\mu}_{\max}}(V^\vee)}\cdot\big(V\otimes_{E}B^+_{\mathrm{dR}}(R^\sharp)\big).
	\end{equation}
	Indeed, taking dual vector bundles reverses the inclusion and $\big((\rho^\vee)_*\underline{\calE}\big)^\vee=\rho_*\underline{\calE}$.
\end{remark}

\begin{remark}
In case $G$ is split reductive, analogues of the moduli functor appear in the equal characteristic setting in \cite[Section~10.3, Proposition]{FM:SemiInfinite}, \cite[Remark~2.1.4]{Zhu16}, \cite[Section~2.3]{KamnitzerMuthiahWeekes:ReducednessConjecture}, \cite[Definition~4.2.6]{KisinPappasZhou_IntegralModels} and \cite[Definition~4.1]{BessonHongYu:FMCounterexample}.
The aforementioned works require the condition in \eqref{Tannaka.bound.gen} only for the highest weight representations, cf.~\Cref{highest weight rep coro} below.
In fact, \eqref{Tannaka.bound.gen} specializes for these representations to the condition given in \cite[Section~2.3]{KamnitzerMuthiahWeekes:ReducednessConjecture}.
\end{remark}

\begin{lemma}\label{key lemma schubert}
Let $\nu\in X_*(T)$ and denote by $\nu(\xi)\in \Gr_G(C)$ the corresponding point. 
For any representation $\rho\colon G_{E}\to \GL(V)$, the following are equivalent:
\begin{enumerate}
	\item \label{key lemma schubert 1}
		The point $\rho_*(\nu(\xi))$ lies in $\mathrm{Gr}_{\mathrm{GL}(V)}^{\geq \mu_{\max}(V)}(C)$. 
	\item \label{key lemma schubert 2}
		One has $\nu_{\mathrm{max}}(V)\leq \mu_{\mathrm{max}}(V)$.
\end{enumerate}
\end{lemma}
\begin{proof}
By definition, \Cref{key lemma schubert 1} is equivalent to the following inclusion of $B^+_{\mathrm{dR}}$-lattices:
	$$
	\xi^{\mu_{\max}(V)}\cdot\left(V\otimes_EB^+_{\mathrm{dR}}\right) \subset \nu(\xi) \cdot\left(V\otimes_EB^+_{\mathrm{dR}}\right).
	$$
	Decomposing $V_C=\oplus_{\eta \in X^*(T)}V_C(\eta)$ into $T_C$-weight spaces, the inclusion is equivalent to the inclusions
	$$
	\xi^{\mu_{\max}(V)}\cdot\tilde V(\eta) \subset \nu(\xi) \cdot\tilde V(\eta)=\xi^{\langle\eta,\nu \rangle}\cdot \tilde V(\eta)\quad \text{for all $T_C$-weights $\eta$ of $V_C$},
	$$
	where $\tilde V(\eta)=V_C(\eta)\otimes_CB_{\mathrm{dR}}^+$.
	Equivalently, $\langle\eta,\nu \rangle \leq \mu_{\max}(V)$ for all $T_C$-weights $\eta$ of $V_C$, 
	which is equivalent to \Cref{key lemma schubert 2}.
\end{proof}

\begin{lemma}\label{easy inclusion}
For any subset $J\subset \Rep(G_E)$, the sub-functor $\mathrm{Gr}_{G, \mu}^J\subset \Gr_{G,E}$ defines a closed sub-v-sheaf that is $L^+_EG$-invariant and that contains the $\Spd(C)$-points $\dot{\mu}(\xi)$ for all representatives $\dot{\mu}$ of $\mu$.
In particular, one has $\mathrm{Gr}_{G, \mu}\subset \mathrm{Gr}_{G, \mu}^J$.
\end{lemma}
\begin{proof}
It follows from \Cref{tautological Schubert equality} that $\mathrm{Gr}_{G, \mu}^J$ is a closed sub-v-sheaf.
To check the invariance, let $g\in G(B_{\mathrm{dR}}(R^\sharp))$ and $\underline{\calE}=[(\calE,\al)]\in \Gr_G(R,R^+)$ for some affinoid perfectoid space $\Spa(R,R^+)\to \Spd(E)$.
Then, we have
$$
\rho_*(g\cdot \underline{\calE})=\rho_*[(\calE,g\cdot\alpha)]=[(\rho_*\calE,\rho(g)\cdot \rho_*\alpha)]=\rho(g)\cdot \rho_*\underline{\calE},
$$
for any representation $(\rho,V)\in \Rep(G_E)$.
Now, if $g\in G(B_{\mathrm{dR}}^+(R^\sharp))$ is integral, then $\rho(g)$ stabilizes $V\otimes_EB^+_{\mathrm{dR}}(R^\sharp)$, which implies the invariance, cf.~\Cref{Tannaka.bound.gen}.
That $\Gr_{G,\mu}^J(C)$ contains the points $\dot{\mu}(\xi)$ now follows from \Cref{key lemma schubert}. 
\end{proof}

For $\nu, \la\in X_*(T)^+$, we write $\nu\leq_{\mathbb Q}\la$ for the rational Bruhat order, i.e., $\la-\nu$ can be written as a sum of positive coroots with $\mathbb Q_{\geq 0}$-coefficients.

\begin{proposition}\label{moduli-description-generic-fiber}
	Let $J=\Rep(G_E)$.
	Then, one has
	\begin{equation}\label{union of schubert varieties}
	\mathrm{Gr}_{G, \mu,C}^J=\bigcup_{\substack{\nu\in X_*(T)^+\\ \nu\leq_{\mathbb Q}\dot{\mu}}}\mathrm{Gr}_{G, \nu,C}, 
	\end{equation}
	where $\dot{\mu}\in X_*(T)^+$ is the unique representative of the conjugacy class $\mu$.
	If $G_{\mathrm{der}}$ is simply connected, then $\mathrm{Gr}_{G, \mu}^J=\Gr_{G,\mu}$ .
\end{proposition}
\begin{proof}
Let $\tilde E$ be a finite Galois extension of $E$ splitting $T$.
Let $V_{\chi}$ be any $G_{\tilde E}$-represen\-tation with highest weight $\chi\in X^*(T)^+$ (e.g., the dual Weyl module $\Ind_{B^-}^{G_{\tilde E}}(\chi)$ or the associated irreducible representation) and obtain the $G_E$-representation $\Res_{\tilde E/E}V_\chi$.
        We compute
	$$
	\mu_{\max}(\Res_{\tilde E/E}V_\chi)\underset{\text{\Cref{elementary weight 4}}}{\overset{\text{\Cref{elementary weight}}}{=}}
	\mu_{\max}(V_\chi)\underset{\text{\Cref{elementary weight 2}}}{\overset{\text{\Cref{elementary weight}}}{=}}
	\langle\chi, \dot{\mu}\rangle.
	$$
	
	Now, $\mathrm{Gr}_{G, \mu, C}^J$ is $L^+_{C}G$-invariant (\Cref{easy inclusion}), hence a union of Schubert varieties. 
	So, let $\nu\in X_*(T)^+$. 
	Then, $\nu(\xi)\in \mathrm{Gr}_{G, \mu}^J(C)$ if and only if $\rho_*(\nu(\xi)) \in \mathrm{Gr}_{\mathrm{GL}(V)}^{\geq \mu_{\max}(V)}(C)$ for all $(\rho,V)\in J$ (\Cref{Schubert functor definition}) if and only if $\nu_{\max}(V)\leq \mu_{\max}(V)$ for all $V\in J$ (\Cref{key lemma schubert}).
	Applying this to the representations $V=\Res_{\tilde E/E}V_\chi$, the latter inequality implies
	$$
	\langle\chi,\nu\rangle \leq \max_{\gamma\in \mathrm{Gal}_E}\{\langle\gamma*\chi, \nu\rangle\}= \nu_{\max}(\Res_{\tilde E/E}V_\chi)\leq \mu_{\max}(\Res_{\tilde E/E}V_\chi)=\langle\chi,\dot{\mu}\rangle,
	$$
	i.e., $\langle \chi, \dot{\mu}-\nu\rangle\geq 0$ for all $\chi \in X^*(T)^+$.
	By \Cref{root lemma} below, this is equivalent to $\nu\leq_{\mathbb Q}\dot\mu$.
Conversely, if $\nu\leq_{\mathbb Q}\dot\mu$, then $\langle\chi,\nu\rangle\leq\langle\chi,\dot\mu\rangle$ for every $\chi\in X^*(T)^+$.
Taking the maximum over the highest weights of any $V\in\Rep(G_E)$ gives $\nu_{\max}(V)\leq\mu_{\max}(V)$ by \Cref{elementary weight 2}.
Thus \eqref{union of schubert varieties} holds. 
	If, in addition, $G_{\mathrm{der}}$ is simply connected, then $\nu\leq \dot{\mu}$ also by \Cref{root lemma}.
	So, the final claim follows from the closure relations of Schubert cells, see \Cref{bruhat order on schubert varieties}.
\end{proof}

\begin{lemma}\label{root lemma}
For all $\la, \nu\in X_*(T)^+$, the following are equivalent:
\begin{enumerate}
	\item \label{root lemma 1}
	$\nu \leq_{\mathbb Q} \la$, and
	\item \label{root lemma 2}
	$\langle \chi,\la-\nu\rangle \geq 0$ for all $\chi\in X^*(T)^+$ (equivalently, for all $\chi\in X^*(T)^+_{\mathbb Q}$). 
\end{enumerate}
If, in addition, $G_{\mathrm{der}}$ is simply connected, then the conditions are equivalent to $\nu\leq \la$.
\end{lemma}
\begin{proof}
	Assume \Cref{root lemma 2} holds (the converse is immediate).
	Let $Z_G$ be the center of $G$, and $\Delta$ the set of (absolute) simple roots in $B_0$.
	Using the isomorphism $X_*(T)_{\mathbb Q}\cong X_*(Z_G)_{\mathbb Q}\times \mathbb Q\Delta\!^\vee$ for the set of simple coroots $\Delta\!^\vee$, we can write 
	$$\la-\nu=\tau +\sum_{\alpha \in \Delta} c_\alpha\alpha^\vee$$ 
	for unique $\tau\in X_*(Z_G)_{\mathbb{Q}}$ and unique $c_\alpha \in \mathbb{Q}$ for all $\alpha \in \Delta$. 
	We have an isomorphism $X^*(T)^+_\mathbb Q\cong X^*(Z_G)_\mathbb Q\times (\mathbb Q\Delta)^+$.
	Pairing with $\chi \in X^*(Z_G)_{\mathbb{Q}}$, we deduce that $\langle \chi, \tau \rangle \geq 0$ from the non-negativity of the pairing, hence $\langle \chi, \tau \rangle =0$ by replacing $\chi$ with $-\chi$ (still dominant by centrality).
	The non-degeneracy of the pairing restricted to the center forces $\tau=0$. 
	Pairing with (positive multiples of) the fundamental weights, we deduce $c_\alpha\geq 0$ for all $\al\in \Delta$ and hence \Cref{root lemma 1} holds.
	
	If $G_{\mathrm{der}}$ is simply connected, then $\la-\nu \in X_*(T)\cap X_*(G_\der)_{\mathbb Q}=X_*(T\cap G_{\mathrm{der}})=\mathbb{Z}\Delta^\vee$, where the last equality uses the simply connectedness. 
	Thus, $c_\alpha\in \mathbb Z\cap \mathbb Q_{\geq 0}=\mathbb Z_{\geq 0}$ for all $\al \in \Delta$, hence $\nu \leq \la$.
\end{proof}

The following result shows that only finitely many representations are necessary to give moduli descriptions of Schubert varieties:

\begin{corollary}\label{highest weight rep coro}
	Assume $G_\der$ is simply connected. 
	Then, $\Gr_{G,\mu}=\Gr_{G,\mu}^J$
	where $J=\{\Res_{\tilde E/E}V_\chi\}_{\chi}$ with $\chi$ running through any subset of $X^*(T)^+$ whose $\mathrm{Gal}_E$-translates generate $X^*(T)_{\mathbb{Q}}^+$ as a $\mathbb{Q}_{\geq 0}$-monoid.
\end{corollary}
\begin{proof}
We follow the proof of \Cref{moduli-description-generic-fiber}.
The action of any $\gamma\in \mathrm{Gal}_E$ induces an isomorphism of $G_E$-representations $\Res_{\tilde E/E}V_\chi\cong\Res_{\tilde E/E}^\gamma V_{\chi}$. 
Here, $^\gamma V_{\chi}$ denotes the $\gamma$-twist of $V_\chi$ which is a $G_{\tilde E}$-representation of highest weight $\gamma*\chi$.
Hence, the choice of $J$ implies $\langle \chi,\dot{\mu}-\nu\rangle \geq 0$ for the elements $\chi$ in a generating set, hence for all $\chi\in X^*(T)^+_{\mathbb Q}$.
\Cref{root lemma} shows $\nu \leq \dot{\mu}$.
\end{proof}

\begin{remark}\label{pgl2 remark}
\Cref{union of schubert varieties} shows that the equality $\Gr_{G,\mu}^J=\Gr_{G,\mu}$ can fail when $G_\der$ is not simply connected, even when $J=\Rep(G_E)$.  
For example, let $G=\mathrm{PGL}_2$ and let $\mu$ be quasi-minuscule, so that $\dot{\mu}=2$ under the identification $X_*(T)=\mathbb Z$. 
Then, $\{\nu\in X_*(T)^+\mid \nu\leq_{\mathbb Q}\dot{\mu}\}=\{0,1,2\}$.
Hence, $\Gr_{G,\mu}^J$ is the disjoint union of $\Gr_{G,\mu}$ and the minuscule Schubert variety.
\end{remark}

\subsection{Minimal test sets for Schubert varieties}\label{minimal test sets}
We seek to find a minimal set $J=J(\mu)$ of representations that realizes the Schubert variety, i.e., $\Gr_{G,\mu}=\Gr_{G,\mu}^J$.

Let $\tilde E$ be a finite Galois extension of $E$ splitting $T$.
We have the $G_E$-representation $\Res_{\tilde E/E}V_\chi$ for any chosen $G_{\tilde E}$-representation $V_\chi$ with highest weight $\chi \in X^*(T)^+$.
This allows us to identify subsets of $X^*(T)^+$ with subsets of $\Rep(G_E)$.

Assume that $G_\der$ is simply connected.
Then $T_\der=G_\der\cap T$ is a maximal torus in $G_\der$ and $P=X^*(T_\der)$ is the weight lattice associated with the based root datum in $X^*(T_\der)_{\mathbb R}$ determined by $B_{0,C}\cap G_{\der, C}$. 
First assume that the absolute Dynkin diagram of $G$ is connected. In this case we use the labeling convention of \cite{Bou68}:
let $\omega_1,\ldots,\omega_n\in P$ be the fundamental weights, and $\omega_1^\vee,\ldots,\omega_n^\vee\in P^\vee$ the fundamental coweights.
We have $P^\vee=X_*(T_\der/Z_{G_\der})$ where $Z_{G_\der}$ is the center of $G_\der$.
Consider the exact sequence of abelian groups
$$
0\to X^*(G)\to X^*(T)\overset{\chi\mapsto \chi|_{T_\der}}{\to} P \to 0,
$$
see \Cref{exact sequence dual}. 
Let $\mu_{\mathrm{ad}} \in (P^\vee)^+$ be the image of the dominant representative $\dot\mu\in X_*(T)^+$ of $\mu$.
If $\mu_\ad\neq 0$, then we denote by $J(\mu_\ad)\subset \{\omega_1,\ldots,\omega_n\}$ the minimal test set as in \Cref{minimal test set definition}: 
that is, $J(\mu_\ad)$ is the smallest subset of $\{\omega_1,\ldots,\omega_n\}$ satisfying the following property: for every $\nu\in (P^\vee_{\mathbb Q})^+$, the inequalities
$$
\langle\chi,\mu_\ad-\nu\rangle\geq0
\qquad
\text{for all }\chi\in J(\mu_\ad)
$$
imply $\langle\omega_i,\mu_\ad-\nu\rangle\geq0$ for all $i=1,\ldots,n$.
By \Cref{minimal test set}, the set \(J(\mu_{\ad})\) is uniquely determined by the support of \(\mu_{\ad}\). 
If \(\mu_{\ad}=0\), we instead define \(J(\mu_{\ad})\) to consist of a single arbitrarily chosen fundamental weight. 
If the absolute Dynkin diagram of \(G\) is disconnected, we define \(J(\mu_{\ad})\) factorwise.

We consider 
\begin{equation}\label{minimal test set prop eq}
J(\mu):=\{\delta_1,\ldots,\delta_m\}\cup \tilde J(\mu_\ad) \subset X^*(T)^+,
\end{equation}
where $\delta_1,\ldots,\delta_m\in X^*(G)$ is any minimal generating set of the $\mathbb Q_{\geq 0}$-monoid $X^*(G)_{\mathbb Q}$ and $\tilde J(\mu_\ad)\subset X^*(T)^+$ are any lifts of the elements in $J(\mu_\ad)$.

For any subset $J\subset X^*(T)^+$, we denote by $\Gal_E*J$ its orbit for the Galois action on the based root datum. 

\begin{proposition}\label{minimal test set prop}
	Assume that $G_\der$ is simply connected. 
	Then, $\Gr_{G,\mu}=\Gr_{G,\mu}^{J}$ for any subset $J\subset X^*(T)^+$ such that $J(\mu)\subset \Gal_E*J$.
\end{proposition}
\begin{proof}
As in \Cref{highest weight rep coro}, we have to show that the inequality $\langle \chi,\dot{\mu}-\nu\rangle \geq 0$ for all elements $\chi$ in our test set $J$ already implies the same inequality for a generating set of the $\mathbb Q_{\geq 0}$-monoid $X^*(T)_{\mathbb Q}^+$.
Evidently, $\{\delta_1,\ldots,\delta_m\}$ together with preimages in $X^*(T)^+$ of all fundamental weights $\omega_1,\ldots,\omega_n$ is a generating set for $X^*(T)_{\mathbb Q}^+=X^*(G)_{\mathbb Q}\times P_{\mathbb Q}^+$.
Passing to the quotient of $X^*(T)$ by $X^*(G)$ and arguing factorwise, the description of $J$ follows from the definition of test sets when $\mu_\ad\neq 0$.
If $\mu_\ad=0$, let $J(\mu_\ad)=\{\omega_i\}$ for some $i$.
Since all entries of the inverse Cartan matrix of an irreducible root system are strictly positive, the inequality $\langle\omega_i,-\nu_\ad\rangle\geq 0$ forces the dominant coweight $\nu_\ad$ to vanish; see \Cref{inverse cartan multiplication} in \Cref{appendix minimal test sets}.
Thus, all remaining inequalities hold in this case as well.
This proves the proposition.
\end{proof}

\section{Moduli candidates for local models}
We formulate in \Cref{moduli conjecture} a moduli description of the v-sheaf local models $\mathbb M_{\calG,\mu}$ from \Cref{local models intro}.
Here, $\calG$ is a general parahoric $O$-group scheme and $\mu$ a conjugacy class of geometric cocharacters, not necessarily minuscule.

\subsection{Moduli functor}
If the reflex field $E$ is ramified over $F$, then the formulation of the moduli functor in \Cref{definition local model} is slightly more involved than the one for the Schubert varieties in \Cref{Schubert functor definition}, because the formation of Beilinson--Drinfeld Grassmannians only commutes with \textit{unramified} base change. 

\begin{lemma and definition}\label{preliminary lemma}
Let $V\in\Rep(G_{E})$, and write $V_\iota:=V\otimes_{E,\iota}C$ for any $E^{\mathrm{ur}}$-field embedding $\iota\colon E\hookrightarrow C$.
Then, the tuple of integers $\underline{\mu}_{\max}(V):=\bigl(\mu_{\max}(V_\iota)\bigr)_{\iota\in\Hom_{E^{\mathrm{ur}}}(E,C)}$ induces a unique element
\[
\underline{\xi}^{\underline{\mu}_{\max}(V)}
\in
\Gr_{\Res_{O_E/O_{E^{\mathrm{ur}}}}(\bbG_m)}\big(\Spd(O_E)\big),
\]
whose restriction to $\Spd(C)$ is the point corresponding to $\underline{\mu}_{\max}(V)$ under the identification of $\Gr_{\Res_{O_E/O_{E^{\mathrm{ur}}}}(\bbG_m)}(C)$ with $\mathbb Z^{\Hom_{E^{\mathrm{ur}}}(E,C)}$.
In addition, for $\calV\in \Rep(\calG_{O_E})$, one sets $\underline{\mu}_{\max}(\calV):=\underline{\mu}_{\max}(\calV\otimes_{O_E}E)$. 
\end{lemma and definition}
\begin{proof}
Put $\calH:=\Res_{O_E/O_{E^{\mathrm{ur}}}}\bbG_m$ and $H=\calH_{E^{\mathrm{ur}}}$ for its generic fiber.
The decomposition
\[
C\otimes_{E^{\mathrm{ur}}}E
\xrightarrow{\sim}
\prod_{\iota\in\Hom_{E^{\mathrm{ur}}}(E,C)}C,
\qquad
a\otimes b\longmapsto (a\iota(b))_\iota,
\]
identifies $\Gr_H(\Spd(C))$ with
$\prod_{\iota}\Gr_{\mathbb G_m}(C)=\mathbb Z^{\Hom_{E^{\mathrm{ur}}}(E,C)}$, see \cite[Proposition~21.3.1]{SW20}.
Let $\underline{\xi}^{\underline{\mu}_{\max}(V)}$
be the point corresponding to the displayed tuple.
Since ${}^\sigma\!\mu=\mu$ and ${}^\sigma V_\iota\simeq V_{\sigma\circ\iota}$, we have
\[
\mu_{\max}(V_{\sigma\circ\iota})
=
({}^\sigma\!\mu)_{\max}({}^\sigma V_\iota)
=
\mu_{\max}(V_\iota).
\]
Thus the point is invariant under $\Gal_E$, and descends to $\Spd(E)$.
By \cite[Proposition 4.14]{AGLR22}, the $\Spd(C)$-point extends uniquely to a $\Spd(O_C)$-point of $\Gr_{\calH}$.
By uniqueness of this extension, it is $\Gal_E$-equivariant, and since $\Spd(O_C)/\Gal_E=\Spd(O_E)$ by \cite[Lemma 18.1.2]{SW20}, it descends to a unique $\Spd(O_E)$-point.
\end{proof}

Given a representation $\rho\colon \mathcal{G}_{O_E}\to \GL(\calV)$, we denote by $\tilde\rho\colon \calG_{O_{E^{\mathrm{ur}}}}\to \Res_{O_E/O_{E^{\mathrm{ur}}}}\GL(\calV)$ the map corresponding to $\rho$ under the adjunction $(-\otimes_{O_{E^{\mathrm{ur}}}} O_E, \Res_{O_E/O_{E^{\mathrm{ur}}}}-)$.
Then, $\tilde \rho$ induces the following map on Beilinson--Drinfeld Grassmannians after base change $-\times_{\Spd(O_{E^{\mathrm{ur}}})}\Spd(O_E)$:
\begin{equation}\label{rho tilde}
\tilde\rho_*\colon \Gr_{\calG,O_E}\to \Gr_{\Res_{O_E/O_{E^{\mathrm{ur}}}}\GL(\calV), O_E},
\end{equation}
where we use the compatibility $\Gr_{\calG_{O_{E^{\mathrm{ur}}}}}=\Gr_{\calG,O_{E^{\mathrm{ur}}}}$ for the unramified extension $E^{\mathrm{ur}}/F$.
Further, we identify $\Res_{O_E/O_{E^{\mathrm{ur}}}}\mathbb G_m$ with the center of $\Res_{O_E/O_{E^{\mathrm{ur}}}}\GL(\calV)$, thus the former acts on the latter.

\begin{definition}\label{definition local model}
	Let $J\subset \Rep(\calG_{O_E})$ be a subset. 
	Define the sub-functor ${\mathbb M}_{\mathcal{G},\mu}^J\subset \Gr_{\mathcal{G},O_E}$ given on affinoid perfectoid spaces $\Spa(R,R^+)\to \Spd(O_E)$ by all $\underline{\calE}=[(\calE,\al)]\in \Gr_{\mathcal{G}}(R,R^+)$ such that for all representations $\rho\colon \mathcal{G}_{O_E}\to \GL(\calV)$ in $J$, one has
	\begin{equation}\label{Tannaka bound local model}
		\underline{\xi}^{{\underline{\mu}_{\max}}(\calV)}\cdot\big(\calV\otimes_{O_{E^{\mathrm{ur}}}}B^+_{\mathrm{dR}}(R^\sharp)\big)\subset \tilde\rho_*\underline{\calE} 
	\end{equation}
	as $O_E\otimes_{O_{E^{\mathrm{ur}}}}B^+_{\mathrm{dR}}(R^\sharp)$-lattices inside $\calV\otimes_{O_{E^{\mathrm{ur}}}} B^+_{\mathrm{dR}}(R^\sharp)[\xi^{-1}]$.
	Equivalently, one has 
	\begin{equation}\label{tautological local model equality}
{\mathbb M}_{\calG, \mu}^J= \bigcap_{(\rho,\calV) \in J} \tilde{\rho}_*^{-1}\left(\mathrm{Gr}_{\Res_{O_E/O_{E^{\mathrm{ur}}}}\mathrm{GL}(\calV), O_E}^{\geq \underline{\mu}_{\max}(\calV)}\right),
\end{equation}
where $\mathrm{Gr}_{\Res_{O_E/O_{E^{\mathrm{ur}}}}\mathrm{GL}(\calV), O_E}^{\geq \underline{\mu}_{\max}(\calV)}$ parametrizes lattices containing $\underline{\xi}^{{\underline{\mu}_{\max}}(\calV)}\cdot\big(\calV\otimes_{O_{E^{\mathrm{ur}}}}B^+_{\mathrm{dR}}(R^\sharp)\big)$.
\end{definition}

\begin{remark}\label{reverse inclusion model} 	
Equation \eqref{Tannaka bound local model} applied to the dual representation $\calV^\vee$ gives the inclusion
	\begin{equation}\label{Tannaka.bound.dual.model}
		\tilde\rho_*\underline{\calE}  \subset \underline{\xi}^{-{\underline{\mu}_{\max}}(\calV^\vee)}\cdot\big(\calV\otimes_{O_{E^{\mathrm{ur}}}}B^+_{\mathrm{dR}}(R^\sharp)\big),
	\end{equation}
	cf.~\Cref{remark reverse inclusion}.
\end{remark}

\begin{lemma}\label{easy inclusion model}
The functor $\mathbb M_{\calG,\mu}^J$ defines an
$L^+_{O_E}\calG$-invariant closed sub-v-sheaf of $\Gr_{\calG,O_E}$ such that on geometric generic fibers
$$\mathbb M_{\calG,\mu,C}^J=\Gr_{G,\mu,C}^{J'},$$ 
where $J'=\{\calV\otimes_{O_E,\iota}C\mid \calV\in J,\ \iota\in \Hom_{E^{\mathrm{ur}}}(E,C)\}$.
In particular, $\Gr_{G,\mu}\subset \mathbb M_{\calG,\mu,E}^J$ and so $\mathbb M_{\calG,\mu}\subset\mathbb M_{\calG,\mu}^J$.
\end{lemma}
\begin{proof}
Closedness and invariance follow as in \Cref{easy inclusion}.
For $(\rho,\calV)\in J$ and $\iota\in\Hom_{E^{\mathrm{ur}}}(E,C)$, write $V_{\iota}:=\calV\otimes_{O_E,\iota}C$ and $\rho_\iota\colon G_C\to\GL(V_{\iota})$.
Under the scalar-extension decomposition, the $\iota$-component of $\tilde\rho\otimes_{O_{E^{\mathrm{ur}}}}C$ is $\rho_\iota$.
Consequently,
\[
\mathbb M_{\calG,\mu,C}^J
=
\bigcap_{(\rho,\calV)\in J}
\ \bigcap_{\iota\in\Hom_{E^{\mathrm{ur}}}(E,C)}
(\rho_\iota)_*^{-1}
\left(
\Gr_{\GL(V_{\iota})}^{\geq\mu_{\max}(V_{\iota})}
\right)
=\Gr_{G,\mu,C}^{J'}.
\]
In particular, each factor contains $\Gr_{G,\mu,C}$ by \Cref{easy inclusion} and so $\Gr_{G,\mu}\subset \mathbb M_{\calG,\mu,E}^J$.
Since $\mathbb M_{\calG,\mu}^J$ is closed and contains $\Gr_{G,\mu}$, it contains its v-sheaf closure $\mathbb M_{\calG,\mu}$.
\end{proof}

\subsection{A conjectural moduli description}
This work revolves around the following conjecture, which we prove under a mild tameness assumption. 
A refined version, involving only finitely many representations, is studied in \Cref{refined conjecture section}.

\begin{conjecture}\label{moduli conjecture}
Assume that $G_\der$ is simply connected and that $\pi_1(G)_I$ is torsion-free. 
Then, $\mathbb M_{\calG,\mu}=\mathbb M^J_{\calG,\mu}$ for $J=\Rep(\calG_{O_E})$.
\end{conjecture}

Whether the inclusion $\mathbb M_{\calG,\mu} \subset \mathbb M^J_{\calG,\mu}$ is an isomorphism can be checked on points with values in algebraically closed perfectoid fields \cite[Lemma~12.5]{Sch17}. 
In particular it suffices that one has equality in the generic fiber and in the special fiber. 
The conjecture holds for the generic fiber:

\begin{proposition}\label{generic equality}
Assume that $G_\der$ is simply connected. 
Then, $\mathbb M_{\calG,\mu} \otimes E = \mathbb M^J_{\calG,\mu} \otimes E$ for $J=\Rep(\calG_{O_E})$.
\end{proposition}
\begin{proof}
Combine \Cref{easy inclusion model} and \Cref{moduli-description-generic-fiber}. 
Every finite-dimensional $G_E$-representation admits a $\calG_{O_E}$-stable lattice \cite[Section~I.10.3]{Jan03}. 
\end{proof}

Since the special fiber of $\mathbb M_{\calG,\mu}$ agrees with the v-sheaf associated with the admissible locus $A_{\calG,\mu}$ (\Cref{local models intro}), \Cref{moduli conjecture} is equivalent to the following conjecture:

\begin{conjecture}\label{moduli conjecture admissible}
Assume that $G_\der$ is simply connected and that $\pi_1(G)_I$ is torsion-free. 
Then, $A_{\calG,\mu}=A^J_{\calG,\mu}$ for $J=\Rep(\calG_{O_E})$.
\end{conjecture}

Here, $A_{\calG,\mu}^J\subset \Fl_{\calG, k_E}$ is the sub-functor given on perfect $k_E$-algebras $R$ by all $\underline{\calE}\in \Fl_\calG(R)$ such that for all representations $\rho\colon \calG_{O_E}\to \GL(\calV)$ in $J$ one has 
\begin{equation}\label{special fiber condition}
(\pi_E\otimes 1)^{\vert\underline{\mu}_{\max}(\calV)\vert}\cdot\big(\calV\otimes_{O_{E^{\mathrm{ur}}}} W_{O}(R)\big)\subset \tilde\rho_*\underline{\calE},
\end{equation}
where $\pi_E\in O_E$ is a uniformizer, $\vert\underline{\mu}_{\max}(\calV)\vert=\sum_{\iota\in \Hom_{E^{\mathrm{ur}}}(E,C)} {\mu}_{\max}(\calV\otimes_{O_E,\iota}C)$, and the map 
$$\tilde \rho_*\colon \Fl_{\calG,k_E}=\Fl_{\calG_{O_{E^{\mathrm{ur}}}}}\to \Fl_{\Res_{O_E/O_{E^{\mathrm{ur}}}}\GL(\calV)}= \Fl_{\GL(\calV)}$$ 
is induced by $\rho$ as in \Cref{rho tilde}.
(The final identification uses $W_{O_E}(R)=O_E\otimes_{O_{E^{\mathrm{ur}}}}W_{O}(R)$ for any perfect $k_E$-algebra $R$.)
Evidently, $(A_{\calG,\mu}^J)^\diamondsuit=\mathbb M_{\calG,\mu,k_E}^J$, so \Cref{moduli conjecture} is indeed equivalent to \Cref{moduli conjecture admissible}. 
Note that the inclusion $\mathbb M_{\calG,\mu}\subset \mathbb M^J_{\calG,\mu}$ (\Cref{easy inclusion model}) implies $A_{\calG,\mu}\subset A^J_{\calG,\mu}$.

\begin{remark}
For representations $\rho\colon \calG_{O_{E^{\mathrm{ur}}}}\to \GL(\calV)$ defined over $\calO_{E^{\mathrm{ur}}}$, the condition in \Cref{Tannaka bound local model} is equivalent to the condition
\begin{equation}\label{Tannaka bound local model simple} 
		\xi^{{{\mu}_{\max}}(\calV)}\cdot\big(\calV\otimes_{O_{E^{\mathrm{ur}}}}B^+_{\mathrm{dR}}(R^\sharp)\big)\subset \rho_*\underline{\calE}, 
	\end{equation}
where $\rho_*\colon \Gr_{\calG,O_E}\to \Gr_{\GL(\calV),O_E}$.
However, representations defined over $O_{E^{\mathrm{ur}}}$ are, in general, insufficient to detect the Schubert variety in the generic fiber when $E$ is ramified over $F$.
An example is given by $\calG=\Res_{O_K/O}\mathbb G_m$ with $K/F$ a ramified, quadratic extension, and $\mu$ corresponding to the geometric cocharacter of weight $(1,0)$, so $E=K$ and $E^{\mathrm{ur}}=F$. 
Representations defined over $O_{E^{\mathrm{ur}}}$ do not recover the generic fiber: their sets of weights are inertia-stable, so they cannot distinguish $(1,0)$ from $(0,1)$ whereas the generic Schubert variety contains only $(1,0)$.
\end{remark}

\begin{remark}\label{torsion torus example}
The torsion-freeness hypothesis in \Cref{moduli conjecture} and \Cref{moduli conjecture admissible} is already visible for tori.
Let $K/F$ be a ramified separable quadratic extension, and set $T$ to be the kernel of the norm $\mathrm{Norm}_{K/F}\colon\Res_{K/F}\mathbb G_m\to\mathbb G_m$.
Then inertia acts by $-1$ on $\pi_1(T)=X_*(T)\simeq\mathbb Z$, so that $X_*(T)_I\simeq\mathbb Z/2$, and thus $X_*(T)_{I,\mathbb Q}=0$.
Take $\mu=1$, so that $E=K$.
The admissible locus $A_{\calT,\mu}$ is the single point indexed by $1\in\mathbb Z/2$.
For any $\calV\in\Rep(\calT_{O_K})$, with generic weight set $S(\calV)\subset\mathbb Z$, one has $\vert\underline{\mu}_{\max}(\calV)\vert=\max S(\calV)-\min S(\calV)\geq 0$.
Hence the identity point also satisfies \Cref{special fiber condition}.
Together with $A_{\calT,\mu}\subset A_{\calT,\mu}^J$, this gives $A_{\calT,\mu}^J=\bbZ/2$ even for $J=\Rep(\calT_{O_K})$.
\end{remark}

\section{Parahoric induced modules}\label{parahoric Weyl module sec}
In the next \Cref{refined conjecture section}, we refine \Cref{moduli conjecture} to involve only finitely many representations of $\calG_{O_E}$ (\Cref{refined conjecture}) and single out a smaller set of representations whose moduli problem lies between the admissible and permissible loci (\Cref{definition min test set}). 
These conjectures are based on the construction of parahoric induced modules. 

\subsection{Induced representations}
In this subsection, let $R$ be a ring, $\calG$ any flat affine $R$-group scheme of finite presentation, and $\calH$ an $R$-flat, closed subgroup scheme in $\calG$. 
For a representation $\calH\to \GL(\calV)$, we denote by
\begin{equation}\label{induction of reps}
\Ind_\calH^{\calG}(\calV) =\{f\colon \calG\to \calV\;|\; f(gh)=h^{-1}f(g)\;\; \forall g\in \calG(S),h\in \calH(S)\text{, and }S\in \mathrm{Alg}_R\}
\end{equation}
the induced representation, viewed as a $\calG$-representation under left translation \cite[Section I.3.3]{Jan03}.  
It is the submodule of $\calH$-fixed elements in $\calV\otimes_R R[\calG]$ where we denote $R[\calG]:= \Gamma(\calG,\calO_{\calG})$.
Alternatively, let $\calV_{\mathrm{geo}}:=\calG\times^{\calH}\mathbb A(\calV)\to \calG/\calH$ be the geometric vector bundle defined by $\calV$.
Then \Cref{induction of reps} identifies with the global sections $\Gamma(\calG/\calH,\calV_{\mathrm{geo}})$.

The following result gives a criterion for finite generation of induced representations:

\begin{proposition}\label{finitely generated proposition}
Assume $R$ is a Dedekind ring and that the following properties hold:
\begin{enumerate}
\item \label{finitely generated proposition 1}
The generic fiber of the fppf quotient
$\calG/\calH\to\Spec(R)$ is proper.
\item \label{finitely generated proposition 2}
There exists an $R$-smooth, closed subscheme $\calU$ in $\calG$
with connected generic fiber and a section over $R$, such that
multiplication $\calU\times\calH\to\calG$ is a schematically dense open immersion.
\end{enumerate}
If $\calV$ is an $\calH$-representation on a finite projective
$R$-module, then $\Ind_{\calH}^{\calG}(\calV)$ is a finite
projective $R$-module.
\end{proposition}
\begin{proof}
The induced representation is an $R$-submodule of
$\calV\otimes_R R[\calG]$, which is $R$-flat if $\calV$ is so.
Over a Dedekind ring flatness is equivalent to torsion freeness, so every submodule of a
flat module is flat \cite[\href{https://stacks.math.columbia.edu/tag/0AUW}{Tag 0AUW}]{stacks-project}.
Hence, it suffices to show that $\Ind_{\calH}^\calG(\calV)$ is a finitely generated $R$-module \cite[\href{https://stacks.math.columbia.edu/tag/0AUW}{Tag 0AUW}]{stacks-project}. 
The induced map given by restriction of functions
\begin{equation}\label{restriction map reps}
\Ind_{\calH}^{\calG}(\calV)=\Gamma(\calG/\calH,\calV_{\mathrm{geo}})\overset{f\mapsto f|_{\calU}}{\longto} \Gamma(\calU,\calV_{\mathrm{geo}})= R[\calU] \otimes_R \calV,
\end{equation}
is injective by Assumption \eqref{finitely generated proposition 2}.
For the final equality, observe that the restricted bundle
$\calV_{\mathrm{geo}}|_\calU \cong \calU\times \mathbb A(\calV)$ is trivial. Indeed, the composite
$\calU\to \calG\to \calG/\calH$ is an open immersion and
the given $\calH$-torsor becomes trivial over the open locus $\calU$, by assumption.

For any $N\geq 1$, we denote by $\calU_N\subset \calU$ the $N$-th infinitesimal neighborhood of a chosen section of $\calU$ over $R$. 
The chosen section is a regular immersion by smoothness of $\calU$ over $R$ \cite[\href{https://stacks.math.columbia.edu/tag/067R}{Tag 067R}]{stacks-project}.
The successive quotients of its ideal powers are symmetric powers of the finite locally free conormal module, so $\calU_N$ is finite locally free over $\Spec R$.
We claim that the resulting map 
\begin{equation}\label{distribution_embedding}
\varphi_N\colon \Ind_{\calH}^{\calG}(\calV)\to R[\calU] \otimes_R \calV \to R[\calU_N] \otimes_R \calV
\end{equation}
is still injective for $N\gg 1$ which then shows that $\Ind_{\calH}^{\calG}(\calV)$ is finitely generated because $R$ is noetherian.
Again by $R$-flatness of $\Ind_{\calH}^{\calG}(\calV)$, it is enough to show injectivity of $\phi_N:=\varphi_N\otimes_R K$, where $K$ is the fraction field of $R$.
Flat base change identifies this $K$-vector space with $\Gamma((\calG/\calH)_K,\calV_{\mathrm{geo},K})$.
Assumption \eqref{finitely generated proposition 1} and proper pushforward of coherent sheaves \cite[\href{https://stacks.math.columbia.edu/tag/08AR}{Tag 08AR}]{stacks-project} therefore imply that it is finite-dimensional. 
Since $\calU_K$ is smooth and connected, it is integral.
Let $u\in\calU_K(K)$ be the point defined by the chosen section.
The maps $K[\calU_K]\to\calO_{\calU_K,u} \to\widehat{\calO}_{\calU_K,u}$ are injective: the first since $\calU_K$ is integral, and the second is faithfully flat.
Together with \Cref{restriction map reps}, this gives $\bigcap_{N\geq 1}\ker(\phi_N)=0$.
Since these kernels form a decreasing sequence of subspaces of the finite-dimensional $K$-vector space $\Ind_{\calH}^{\calG}(\calV)\otimes_R K$, we obtain $\ker(\phi_N)=0$ for $N\gg 1$.
This proves the proposition. 
\end{proof}

\subsection{Parahoric induced modules}\label{constructing reps section}
We return to the setting of \Cref{group theoretic notation}, i.e., $\calG$ is a parahoric $O$-model of the reductive group $G$.   
In addition, we assume $G$ to be residually split \cite[Definition~9.10.11]{KP23}.
So, $G$ is quasi-split and the subgroups $A=S\subset T\subset B_0$ are defined over $F$.  

Let $\mathrm{Bor}(G,S)$ be the set of Borel subgroups in $G$ defined over $F$ and containing $S$, hence containing $T$; see \cite[Proposition~2.4.9]{KP23}.
For $B\in \mathrm{Bor}(G,S)$, we denote by $\calB\subset \calG$ its scheme-theoretic closure, which is an $O$-smooth, closed subgroup scheme with connected special fiber admitting a semi-direct product decomposition $\calB=\calT\ltimes \calU$.
Indeed, as $S$ is $F$-split, a $B$-regular dominant cocharacter of $S$ extends to a cocharacter $\la\colon \bbG_m\to \calS\subset \calG$ over $O$, and $\calB$, $\calT$ and $\calU$ are the dynamic subgroups $\calP_\la$, $\calM_\la$ and $\calU_\la$ attached to $\la$.
Smoothness, connectedness of the fibers, and the asserted semi-direct product decomposition follow from \cite[Propositions~2.11.1--2.11.3]{KP23}.
Analogously, we have $\calB^-=\calT\ltimes \calU^-$ for the Borel subgroup $B^-$ opposite to $B$.

Any $\chi\in X^*(T)$ descends to a character $\chi\colon T_{F_\chi}\to \mathbb G_m$ over the field of definition $ F_\chi\subset C$ of $\chi$ over $F$. 
Let $O_\chi$ be its ring of integral elements.
The N\'eron mapping property of the full N\'eron lft-model of $\mathbb G_{m,F_\chi}$ induces a unique extension of $\chi$ to $\calT_{O_\chi}$.
Since $\calT_{O_\chi}$ has connected special fiber, this extension factors through the identity component $\mathbb G_{m,O_\chi}$ \cite[Section~10.1, Example~5]{bosch_luetkebohmert_raynaud_neron_models}, giving a character
$$
 \calB^-_{O_\chi}\longto \calT_{O_\chi} \overset{\chi}{\longto} \mathbb G_{m,O_\chi}.
$$

\begin{proposition}\label{enough reps of parahorics}
Let $\chi \in X^*(T)^+$ and $B\in \mathrm{Bor}(G,S)$. 
Then, the $O_\chi$-module 
\begin{equation}\label{lattice in induction}
\calV_{\chi, B}\: := \Ind_{{\calB}_{O_\chi}^{-}}^{\calG_{O_\chi}}\left({^B\!\chi}\right),
\end{equation}
is finite free, where ${^B\!\chi}\in X^*(T)$ denotes the ${B}$-dominant representative in $W_0\cdot \chi$.
In addition, the following properties hold:
\begin{enumerate}
	\item \label{enough reps of parahorics 1}
		One has $\calV_{\chi, B}\otimes F_\chi\cong \calV_{\chi, B'}\otimes F_\chi$ as $G_{F_\chi}$-representations for all $B'\in \mathrm{Bor}(G,S)$. 
		If $B'=nBn^{-1}$ for some $n\in N_S(F)\cap\calG(O_\chi)$, then $\calV_{\chi, B}\cong \calV_{\chi, B'}$ as $\calG_{O_\chi}$-representations. 
	\item \label{enough reps of parahorics 2}
		One has $\mu_{\max}(\calV_{\chi,B})=\langle{\chi},\dot{\mu}\rangle$, where $\dot{\mu}\in X_*(T)^+$ is the dominant representative of $\mu$.
\end{enumerate}
\end{proposition}
\begin{proof}
The $O_\chi$-module $\calV_{\chi, B}$ is finite free by \Cref{finitely generated proposition}: 
the generic fiber of $\calG/\calB^-\to \Spec(O)$ is the proper map $G/B^-\to \Spec(F)$;
the theory of the big cell \cite[Theorem~8.2.5]{KP23} gives a schematically dense, open immersion ${\calU}\times \calB^-\to \calG$.
The generic fiber $U$ of $\calU$ is the unipotent radical of $B$, hence geometrically connected, and the identity gives an $O$-section of $\calU$.
These properties remain valid after base change to $O_\chi$.
It remains to check the listed properties.

\Cref{enough reps of parahorics 1}: 
Choose $n\in N_T(F_\chi)$ such that $B'=nBn^{-1}$.
Right translation $f(g)\mapsto f(gn)$ commutes with the left $G_{F_\chi}$-action and induces an isomorphism
$$
\calV_{\chi, B}\otimes F_\chi=\Ind_{B_{F_\chi}^{-}}^{G_{F_\chi}}\left({^B\!\chi}\right) \cong \Ind_{({B'})_{F_\chi}^{-}}^{G_{F_\chi}}\big({^{B'}\!\chi}\big)=\calV_{\chi, B'}\otimes F_\chi.
$$
If, in addition, $n\in \calG(O_\chi)$, then the isomorphism holds integrally. 

\Cref{enough reps of parahorics 2}: 
This follows from \Cref{elementary weight}, \Cref{elementary weight 2} using that the extremal weights of $\calV_{\chi, B}\otimes C$ are given by $W_{0,C}\cdot \chi$ \cite[Lemma~22.15 and Proposition~22.27~(c)]{Milne:AlgebraicGroups}, where $W_{0,C}$ is the absolute Weyl group (\Cref{Weyl group section}).  
\end{proof}

The map in \Cref{restriction map reps} gives an embedding of $O_\chi$-modules
\begin{equation}\label{restriction of functions equation}
\calV_{\chi,B}\hookto O_{\chi}[\calU],\; f\mapsto f|_{\calU}.
\end{equation}

\begin{lemma}\label{embedding corollary}
Let $\chi \in X^*(T)^+$ and $B\in \mathrm{Bor}(G,S)$. 
\begin{enumerate}
\item\label{embedding corollary saturation}
The embedding in \Cref{restriction of functions equation} induces an equality of $O_\chi$-submodules in $F_\chi[U]$:
$$
\calV_{\chi, B}=V_{\chi}\cap O_\chi[\calU],
$$
where $V_\chi=\calV_{\chi,B}\otimes_{O_\chi}F_\chi$.
In particular, $\calV_{\chi, B}$ is saturated in $O_\chi[\calU]$, i.e., the quotient $O_\chi[\calU]/\calV_{\chi,B}$ is torsion-free.
\item\label{embedding corollary equivariance}
Let $\calT$ act on $\calU$ by conjugation, and
endow $O_\chi[\calU]$ with the $\calT_{O_\chi}$-action obtained from this one by twisting
with ${^B\!\chi}$, that is
$$
(t\star f)(u)={^B\!\chi}(t)\,f(t^{-1}ut),
\qquad u\in\calU(R),\ t\in \calT(R),\ R\in \mathrm{Alg}_{O_\chi}. 
$$
Then the embedding in \Cref{restriction of functions equation} is $\calT_{O_\chi}$-equivariant.
Equivalently, if $f\in \calV_{\chi,B}$ has weight $\eta$ for the left action of $\calT$, then
$f|_\calU$ has weight
$\eta-{^B\!\chi}$ for the untwisted action.
\item \label{embedding corollary invariants}
The embedding in \Cref{restriction of functions equation} is $\calU$-equivariant and induces an isomorphism 
$$
(\calV_{\chi, B})^\calU\overset{\cong}{\longto} O_\chi[\calU]^\calU=O_\chi, \quad f\mapsto f(1)
$$
on $\calU$-invariants.  Moreover $(\calV_{\chi, B})^\calU$ is the ${^B\!\chi}$-weight space of
$\calV_{\chi,B}$, and the group scheme $\calT_{O_\chi}$ acts on it via ${^B\!\chi}$.
\end{enumerate}
\end{lemma}
\begin{proof}
\Cref{embedding corollary saturation}:
Let $\calL_\chi:=\calG_{O_\chi}\times^{\calB_{O_\chi}^-}(\mathbb A^1_{O_\chi})_{{}^B\!\chi}$ be the geometric line bundle on $(\calG/\calB^-)_{O_\chi}$ defined by ${}^B\!\chi$.
The smooth $O_\chi$-scheme $X:=(\calG/\calB^-)_{O_\chi}$ is regular and integral.
Its special fiber is smooth and geometrically connected, hence geometrically irreducible, and $\calU_{O_\chi}$ meets it in a nonempty open subscheme.
Thus, the open subscheme $Y:=X_{F_\chi}\cup \calU_{O_\chi}$ has complement of codimension $\geq 2$.
Hartogs' theorem \cite[\href{https://stacks.math.columbia.edu/tag/0EBJ}{Tag 0EBJ}]{stacks-project} implies that the map
$$
\calV_{\chi, B}=\Gamma(X,\calL_\chi)\to \Gamma(Y,\calL_\chi)=V_\chi\cap O_{\chi}[\calU], \quad f\mapsto f|_Y
$$
is an isomorphism. 
In particular, $O_\chi[\calU]/\calV_{\chi,B}$ is torsion-free: given an element $f\in O_{\chi}[\calU]$ such that $\pi_\chi f\in \calV_{\chi, B}$, we have $f\in V_\chi$, and hence $f\in \calV_{\chi, B}$.

\Cref{embedding corollary equivariance}:
Let $R\in \mathrm{Alg}_{O_\chi}$, $f\in \calV_{\chi,B}$, $t\in \calT(R)$ and $u\in\calU(R)$.
Since $t^{-1}u=(t^{-1}ut)\cdot t^{-1}$ with $t^{-1}ut\in \calU(R)$, and since $t^{-1}\in \calT(R)\subset \calB^-(R)$, the defining equivariance in \Cref{induction of reps} gives
\begin{equation}\label{formula of item 2}
(t\cdot f)|_{\calU}(u)=f(t^{-1}u)={^B\!\chi}(t^{-1})^{-1}f\bigl(t^{-1}ut\bigr)={^B\!\chi}(t)\,(f|_\calU)(t^{-1}ut)
=\bigl(t\star (f|_\calU)\bigr)(u),
\end{equation}
which is the asserted $\calT_{O_\chi}$-equivariance.

\Cref{embedding corollary invariants}:
For $u_0\in \calU(R)$ we have $(u_0\cdot f)|_{\calU}(u)=f(u_0^{-1}u)=(f|_\calU)(u_0^{-1}u)$, so the embedding is $\calU$-equivariant for the left translation action of $\calU$ on itself.
In particular it carries $\calU$-invariants into constant functions $O_\chi[\calU]^\calU=O_\chi$.
Over $F_\chi$, the $U$-invariants of $V_\chi$ form its ${^B\!\chi}$-weight line and restrict to $F_\chi\cdot 1$, as may be checked after base change to $C$ using \cite[Propositions~22.22 and 22.27]{Milne:AlgebraicGroups}.
By \Cref{embedding corollary saturation}, the constant function $1$ therefore belongs to $\calV_{\chi,B}$, proving surjectivity on $\calU$-invariants.
The action on this line is ${^B\!\chi}$ by \Cref{formula of item 2}.
Finally, torsion freeness gives $(\calV_{\chi,B})^\calU=\calV_{\chi,B}\cap V_\chi^U$, so this is also the ${^B\!\chi}$-weight space of $\calV_{\chi,B}$.
\end{proof}

\begin{corollary}\label{parabolic induction corollary}
Let $\chi \in X^*(T)^+$ and $B\in \mathrm{Bor}(G,S)$. 
Let $P\subset G$ be a parabolic subgroup containing $B^-$ such that $\langle {^B\!\chi},\al^\vee\rangle =0$ for every root $\al$ of the Levi subgroup of $P_C$ containing $T_C$.  
Denote by $\calP \subset \calG$ the scheme-theoretic closure of $P$.
Then, ${^B\!\chi}$ extends to a character of $\calP_{O_\chi}$, and there is a natural isomorphism of $\calG_{O_\chi}$-representations
\begin{equation}
\calV_{\chi, B}\: \cong \Ind_{\calP_{O_\chi}}^{\calG_{O_\chi}}({^B\!\chi}).
\end{equation}
In particular, if $\chi$ is central, then $\calV_{\chi, B}={^B}\!\chi=\chi$.
\end{corollary}
\begin{proof}
Let $M$ be the Levi subgroup of $P$ containing $T$.
The assumption guarantees that ${^B\!\chi}$ extends to a character of $M_{F_\chi}$, hence by inflation to a character of $P_{F_\chi}$.
As above, $\calP$ is the dynamic subgroup attached to a suitable cocharacter of $\calS$, hence smooth with connected fibers.
The same N\'eron mapping argument therefore extends this character to $\calP_{O_\chi}$.
By transitivity of induction \cite[Section~I.3]{Jan03}, it suffices to show
$$
\Ind_{\calB^-_{O_\chi}}^{\calP_{O_\chi}}({^B\!\chi})\cong {^B\!\chi}
$$
as $\calP_{O_\chi}$-representations.
Twisting by the extended character reduces this to $\Gamma(Q,\calO)=O_\chi$, where $Q:=(\calP/\calB^-)_{O_\chi}$.
The generic fiber $Q_{F_\chi}$ is proper, geometrically connected, and geometrically reduced, so $\Gamma(Q_{F_\chi},\calO)=F_\chi$ by \cite[\href{https://stacks.math.columbia.edu/tag/0BUG}{Tag 0BUG}]{stacks-project}.
Flatness of $Q$ makes restriction to the generic fiber injective.
Every global function is therefore generically constant, and evaluation at the identity section forces this constant to lie in $O_\chi$.
\end{proof}

The question whether the representations $\calV_{\chi, B}$ and $\calV_{\chi, B'}$ are isomorphic depends on the relative position of the Borel subgroups $B$ and $B'$.
It is convenient to use the relative Weyl group $W_0$ introduced in \Cref{Weyl group section}.
Since $G$ is residually split, the set $\mathrm{Bor}(G,S)$ corresponds bijectively to the Borel subgroups over $\breve F$ containing $S_{\breve F}$.
Hence, it identifies with the relative Weyl group 
\begin{equation}\label{Weyl group identification}
W_0\overset{\cong}{\longto} \mathrm{Bor}(G,S), \quad w\mapsto wB_0w^{-1},
\end{equation}
using the choice of $B_0$.
Correspondingly, we denote the representation introduced in \Cref{lattice in induction} also by
\begin{equation}\label{alternative lattice notation}
\calV_{\chi, B} = \calV_{\chi, w},
\end{equation}
whenever $B=wB_0w^{-1}$ with $w\in W_0$, and so ${^B\!\chi}=w\chi$. 

The composition $W_\calG\subset W\to W_0$ is injective (cf.~with the proof of \Cref{bijection remark} below), and its image is denoted by $W_{\calG}$ as well. 
For $\chi\in X^*(T)^+$, we let $W_{\chi}$ be the stabilizer of $\chi$ in $W_0$, i.e., all elements $w\in W_0$ such that $w\chi=\chi$.

\begin{lemma}\label{redundancies between representations}
If $w,w'\in W_0$ have the same image in the double coset $W_{\calG}\backslash W_0/W_{\chi}$, then there exists an isomorphism
$
\calV_{\chi,w}\cong \calV_{\chi,w'}
$
as $\calG_{O_\chi}$-representations.
\end{lemma}

\begin{proof}
It is enough to treat separately the two operations $w\mapsto uw$ for $u\in W_{\calG}$ and $w\mapsto wv$ for $v\in W_\chi$.
\Cref{enough reps of parahorics}, \Cref{enough reps of parahorics 1} treats the first operation, using a representative of $u$ in $N_S(F)\cap\calG(O)$; such a representative exists as $G$ is residually split. 
For the second operation, let $v\in W_\chi$, so $wv\chi=w\chi$. 
For any $w'\in W_0$, let $\calB_{w'}\subset \calG$ be the scheme-theoretic closure of $B_{w'}=w'B_0w'^{-1}$.
Then, the opposite Borel subgroups $B_{wv}^-$ and $B_w^-$ lie in the same parabolic subgroup $P$; the Weyl group of its Levi quotient is $wW_\chi w^{-1}$.
\Cref{parabolic induction corollary} implies
\[
	\calV_{\chi,w}=\Ind_{(\calB^-_w)_{O_\chi}}^{\calG_{O_\chi}}(w\chi)
	\cong
	\Ind_{\calP_{O_\chi}}^{\calG_{O_\chi}}(w\chi)
	\cong
	\Ind_{(\calB^-_{wv})_{O_\chi}}^{\calG_{O_\chi}}(wv\chi)=\calV_{\chi, wv}.
\]
\end{proof}

\subsection{Relevant Borel subgroups}
The representations $\calV_{\chi,B}$ constructed in \Cref{lattice in induction} generally depend on the choice of Borel subgroup $B$.
For the moduli description, however, it suffices to consider certain Borel subgroups.
In this section, we specify these relevant Borel subgroups.

We continue to assume that $G$ is residually split.
So, the apartments $\mathscr A(G,S,F)=\mathscr A(G,S,\breve F)$ over $F$ and $\breve F$ agree, and we denote this common apartment by $\mathscr A$.
Let $\mathbf f \subset \mathscr A$ be the facet corresponding to the fixed parahoric subgroup $\mathcal G(O)$.
We denote by $\mathrm{Vert}(\mathbf f)$ the set of vertices contained in the closure of $\mathbf f$.
Fix an auxiliary alcove $\mathbf a$ containing $\mathbf f$ in its closure. 
For a vertex $a\in \mathrm{Vert}(\mathbf f)$, the tangent cone of $\mathbf a$ at $a$ is given by
\begin{equation}\label{tangent cone}
\operatorname{Tan}_a(\mathbf a)
:=
\bigcap_{\substack{\psi \in \Psi^+_{\mathbf a} \\ \psi(a)=0}}
\{\, b \in \mathscr A \mid \psi(b) \ge 0 \,\},
\end{equation}
where $\Psi=\Psi(G,S)$ is the set of affine roots and $\Psi^+_{\mathbf a}=\{\psi\in \Psi\mid \psi|_{\mathbf a}\geq 0\}$ its subset of $\mathbf a$-positive roots.
Equivalently, this is the intersection of all $\mathbf a$-positive root half spaces passing through the vertex $a$.

\begin{definition}\label{relevant Borel subgroups}
The set of \textit{$\mathbf f$-relevant Borel subgroups (with respect to $\mathbf a$, cf.~\Cref{independent of alcove})} is the subset
$$
J_{\mathbf f} = \bigcup_{a\in \mathrm{Vert}(\mathbf f)} J_a\subset \mathrm{Bor}(G,S),
$$
where $J_a=\{\, B \in \operatorname{Bor}(G,S)\;\mid\;\mathscr C_{B}(a) \subset \operatorname{Tan}_a(\mathbf a)\,\}$ and $\mathscr C_{B}(a)\subset \scrA$ is the $B$-dominant cone with apex $a$ in $\mathscr A$.
\end{definition}

\begin{lemma}\label{bijection remark}
The identification $\mathrm{Bor}(G,S)=W_0$ from \Cref{Weyl group identification} induces a bijection $J_a\overset{\sim }{\to} W_{\calG_a}\backslash W_0$.
In particular, $J_a$ is a singleton if $a$ is a special vertex.
\end{lemma}
\begin{proof}
Let $\Psi_a:=\{\,\psi\in\Psi\mid \psi(a)=0\,\}$ and $\Phi_a:=\nabla(\Psi_a)\subset\Phi(G,S)$, where $\nabla$ takes the linear parts.
By \cite[Proposition~1.3.35~(5)]{KP23}, $\nabla$ identifies the stabilizer $W(\Psi)_a$ with $W(\Phi_a)\subset W_0$.
Under $W_{\calG_a}\subset W\to W_0$, this identifies $W_{\calG_a}$ with $W(\Phi_a)$ by \cite[Axiom~4.1.9 and Lemma~1.3.17]{KP23}.
The alcove $\mathbf a$ determines the positive system $\Phi_a^+:=\nabla(\Psi_a\cap\Psi_{\mathbf a}^+)$ by \cite[Proposition~1.3.35~(6)]{KP23}.
Thus, \Cref{tangent cone} makes $\operatorname{Tan}_a(\mathbf a)-a$ the closed fundamental chamber of $W(\Phi_a)$.
Every $W(\Phi_a)$-orbit of Weyl chambers for $W_0$ has exactly one representative contained in this chamber, proving the bijection.
If $a$ is special, then $W_{\calG_a}=W_0$ by \cite[Lemma~1.3.42]{KP23}, so $J_a$ is a singleton.
\end{proof}

\begin{lemma}\label{isomorphism on open cell}
Let $a\in \on{Vert}(\mathbf f)$ and $B\in J_a$.
Then, the maps of parahoric group schemes $\calG_{\mathbf a}\to \calG=\calG_{\mathbf f}\to \calG_a$ induce isomorphisms 
$$
\calU_{\mathbf a}\overset{\cong}{\longto} \calU=\calU_{\mathbf f} \overset{\cong}{\longto} \calU_a,
$$
on the scheme-theoretic closures of the unipotent radical $U$ of $B$ in the respective parahoric group schemes. 
\end{lemma}

\begin{proof}
The decomposition $U\cong \prod_{\al\in \Phi_{\mathrm{nd}}^+} U_\al$ for any ordering on the set $\Phi_{\mathrm{nd}}^+$ of non-divisible $B$-positive roots extends to the corresponding smooth root-group models by \cite[Theorem~8.2.5~(2)]{KP23}.
It is therefore enough to prove that, for every $\al\in \Phi_{\mathrm{nd}}^+$, the induced maps
$
\calU_{\al,\mathbf a}\to \calU_{\al,\mathbf f}\to \calU_{\al,a}
$
are isomorphisms; equivalently \cite[Corollary~2.10.11]{KP23}, that the induced inclusions
\begin{equation}\label{inclusions relevance root}
\calU_{\al,\mathbf a}(\breve O)\subset \calU_{\al,\mathbf f}(\breve O)\subset \calU_{\al,a}(\breve O)
\end{equation}
are equalities.
These point groups are the root filtrations of \cite[Proposition~8.3.14~(2)]{KP23}.
Let $\Psi'$ be the auxiliary affine root system of filtration jumps from \cite[Definition~6.3.4]{KP23}.
If $\psi\in\Psi'$ has linear part $\al$ and $\psi(a)=0$, then $\theta:=\psi$ or $2\psi$ belongs to $\Psi$ by \cite[Proposition~6.3.8~(2)]{KP23}.
If $\psi$ were negative on $\mathbf a$, then $-\theta\in\Psi^+_{\mathbf a}$ would vanish at $a$, so relevance would force $\psi\leq 0$ on $\mathscr C_B(a)$.
This contradicts positivity of $\al$ on the interior of the $B$-dominant cone.
Thus, $\psi$ is positive on $\mathbf a$.
By \cite[Section~6.3.21]{KP23}, this gives $\calU_{\al,\mathbf a}(\breve O)=\calU_{\al,a}(\breve O)$; if no such $\psi$ exists, the same reference gives equality directly.
All inclusions in \Cref{inclusions relevance root} are therefore equalities.
\end{proof}

\begin{corollary}\label{maximal parahoric rep cor}
Let $a\in \on{Vert}(\mathbf f)$ and $B\in J_a$.
Then, for all $\chi\in X^*(T)^+$, the maps $\calG_{\mathbf a}\to \calG=\calG_{\mathbf f}\to \calG_a$ induce isomorphisms of $O_\chi$-modules
$$
\Ind_{({\calB}_{\mathbf a}^-)_{O_\chi}}^{(\calG_{\mathbf a})_{ O_\chi}}\left({^B\!\chi}\right)\overset{\cong}{\longleftarrow} \calV_{\chi, B}\overset{\cong}{\longleftarrow} \Ind_{({\calB}_{a}^-)_{O_\chi}}^{(\calG_{a})_{ O_\chi}}\left({^B\!\chi}\right)
$$
intertwining the actions of the respective parahoric group schemes. 
\end{corollary}
\begin{proof}
By \Cref{isomorphism on open cell}, the maps $\calG_{\mathbf a}\to\calG\to\calG_a$ identify $\calU_{\mathbf a}\cong\calU\cong\calU_a$.
Under restriction to $U$, \Cref{embedding corollary saturation} in \Cref{embedding corollary} identifies all three induced modules with $V_\chi\cap O_\chi[\calU]\subset F_\chi[U]$.
The maps are the identity in this common realization and intertwine the respective parahoric actions.
\end{proof}

\begin{remark}\label{independent of alcove}
Choosing another alcove $\mathbf a'$ containing $\mathbf f$ in its closure there exists a unique $x\in W_\calG$ with $x\cdot \mathbf a=\mathbf a'$.
In particular (\Cref{enough reps of parahorics}, \Cref{enough reps of parahorics 1}), for any $\chi \in X^*(T)^+$, the set of $\calG_{O_\chi}$-representations $\{\calV_{\chi, B}\;|\; B\in J_{\mathbf f}\}$ is (up to isomorphism) independent of the choice of the auxiliary alcove $\mathbf a$ used to define $J_{\mathbf f}$.
\end{remark}

\section{Refined moduli candidates}\label{refined conjecture section}

We return to the general setting of \Cref{group theoretic notation}. 
For a subset $J\subset \Rep(\calG_{O_E})$, recall from \Cref{definition local model} that $\mathbb M_{\calG,\mu}^J$ is the sub-functor of $\Gr_{\calG,O_E}$ given by 
	\begin{equation}\label{refined moduli candidate equation}
	{\mathbb M}_{\calG, \mu}^J= \bigcap_{(\rho,\calV)\in J} \tilde\rho_*^{-1}\left(\mathrm{Gr}_{\Res_{O_E/O_{E^{\mathrm{ur}}}}\mathrm{GL}(\calV), O_E}^{\geq \underline{\mu}_{\max}(\calV)}\right).
	\end{equation}
	
Let $\tilde F\subset C$ be a finite, unramified field extension of $F$ splitting $S$.
Then, $G_{\tilde F}$ is residually split \cite[Proposition~9.10.12]{KP23}, and the results from \Cref{constructing reps section} apply.
The set $\mathrm{Bor}(G,S)$ of Borel subgroups over $\tilde F$ containing $S$ corresponds bijectively to the ones over $\breve F$ containing $S_{\breve F}$.

Let $\tilde E\subset C$ be a finite field extension of the compositum $E \tilde F$ that splits $T$ and that is Galois over $F$.
In particular, any character $\chi \in X^*(T)$ is defined over $\tilde E$.
For any $B\in \mathrm{Bor}(G,S)$, \Cref{enough reps of parahorics} gives the $\calG_{O_E}$-representation 
\begin{equation}\label{parahoric test representations}
\calV_{\chi,B}:=\Res_{{\tilde O}/O_E}\left(\Ind_{{\calB}_{\tilde O}^{-}}^{\calG_{\tilde O}}({^{B}\!\chi})\right)\in \Rep(\calG_{O_E}),
\end{equation}
where $\tilde O$ denotes the ring of integral elements in $\tilde E$ and ${^{B}\!\chi}\in X^*(T)$ the ${B}$-dominant representative in $W_0\cdot \chi$. 
This allows us to identify subsets of $X^*(T)^+\times \textrm{Bor}(G,S)$ with subsets of $\Rep(\calG_{O_E})$.

Let $\mathbf f$ be the facet in the apartment $\mathscr A(G,S,\breve F)$ corresponding to $\mathcal G(\breve O)$.
We denote by $J_{\mathbf f}\subset \textrm{Bor}(G,S)$ the set of $\mathbf f$-relevant Borel subgroups with respect to a chosen auxiliary alcove (\Cref{relevant Borel subgroups} and \Cref{independent of alcove}). 

For any subset $J_1\subset X^*(T)^+$, we denote by $\Gal_E*J_1$ its orbit for the Galois action on the based root datum, and by $N_{\bar I}\colon X^*(T)^+\to X^*(T)^{I,+}$ the norm associated with any finite quotient $I\to \bar I$ through which its action factors.
Let $J_1^0\subset J_1$ be the subset of characters satisfying
\begin{equation}\label{J1 compatibility}
\langle \sigma\gamma*\chi,\dot\mu\rangle
=\langle \sigma*\chi,\dot\mu\rangle
\qquad\text{for every }\sigma\in I,\ \gamma\in\Gal_E.
\end{equation}

\begin{condition}\label{J_1 condition}
The set $J_1$ generates $X^*(T)^+_{\mathbb Q}$ as a $\mathbb Q_{\geq0}$-monoid, and the set $N_{\bar I}(J_1^0)$ generates $X^*(T)^{I,+}_{\mathbb Q}$ as a $\mathbb Q_{\geq0}$-monoid.
\end{condition}

If $E/F$ is unramified, then $J_1^0=J_1$ and so the additional condition on $J_1^0$ is automatic in this case. 
Further, it can always be achieved by adjoining to $J_1$ the norms $N_{\bar I}(\chi)$, $\chi\in J_1$.

\begin{definition}\label{definition test set}
Let $J_1\subset X^*(T)^+$ be a finite set satisfying \Cref{J_1 condition}.
Define the associated test set by
$$
J(\calG,\mu)
:=
J_1\times J_2,
$$
where $J_2=J_{\mathbf f}$ is the finite set of $\mathbf f$-relevant Borel subgroups.
\end{definition}

\begin{conjecture}\label{refined conjecture}
	Assume that $G_\der$ is simply connected and that $\pi_1(G)_I$ is torsion-free.
	Then, one has $\mathbb M_{\calG,\mu}=\mathbb M_{\calG,\mu}^J$ for $J=J(\calG,\mu)$.
\end{conjecture}

The conjecture is a strengthening of \Cref{moduli conjecture} involving only finitely many representations of $\mathcal G_{O_E}$.
\Cref{main theorem intro} proves the conjecture for essentially tamely ramified groups $G$.

\subsection{Reducing the number of test representations}\label{reduction number of reps section}
The test set $J=J(\calG,\mu)$ is not minimal in general.
This can be improved as follows. 

First, we remove redundancies among the representations in the set $J=J(\calG,\mu)=J_1\times J_2$ from \Cref{refined conjecture}, where $J_2=\cup_{a\in \mathrm{Vert}(\mathbf f)}J_a$ is the set of relevant Borel subgroups. 
We identify $\mathrm{Bor}(G,S)$ with the relative Weyl group $W_0$ via $w\mapsto B=wB_0w^{-1}$; see \Cref{Weyl group identification}.  
For every $\chi\in X^*(T)^+$ and $a\in \mathrm{Vert}(\mathbf f)$, we consider the map 
\begin{equation}\label{double class representatives eq}
J_a\subset \mathrm{Bor}(G,S)\overset{\sim}{\to} W_0 \to W_{a}\backslash W_0/W_\chi,
\end{equation}
where $W_a$ is the image of $W_{\calG_a}$ in $W_0$ and $W_\chi$ is the stabilizer of $\chi$ in $W_0$.
For every pair $(\chi,a)$, we choose representatives $J_{a,\chi}\subset J_a$ of the classes in $W_{a}\backslash W_0/W_\chi$.
We consider the set
\begin{equation}\label{finite test set redundant}
J'
:= (J_1^{\mathrm{cen}}\times \{B_0\})\cup
\bigcup_{\substack{\chi\in J_1\setminus J_1^{\mathrm{cen}} \\
                   a\in\operatorname{Vert}(\mathbf f)}}
\{\chi\}\times J_{a,\chi},
\end{equation}
where $J_1^{\mathrm{cen}}=J_1\cap X^*(G)$ are the central characters. 

\begin{lemma}\label{rep reduction lemm}
Every representation in $J$ is isomorphic to one of the representations in $J'$.
In particular, \Cref{refined conjecture} holds for $J$ if and only if it holds for the smaller test set $J'$.
\end{lemma}
\begin{proof}
If $\chi\in J_1^{\mathrm{cen}}$ is central, then $\calV_{\chi,B}=\Res_{\tilde O/O_E}{^B\!}\chi=\Res_{\tilde O/O_E}\chi$ by \Cref{parabolic induction corollary}.
So, every representation in $\{\chi\}\times J_2$ is isomorphic to $\Res_{\tilde O/O_E}\chi$. 
If $\chi\in J_1\setminus J_1^{\mathrm{cen}}$, denote $\calV_{\chi,w}=\calV_{\chi,B}$ for $B=wB_0w^{-1}$ and $w\in W_0$. 
Then, the $\calG_{O_E}$-representations $\calV_{\chi, w}$ and $\calV_{\chi, w'}$ are isomorphic if $w$ and $w'$ lie in the same double class in $W_{\calG_a}\backslash W_0/W_\chi$ by \Cref{maximal parahoric rep cor}.
This follows from \Cref{redundancies between representations,maximal parahoric rep cor}.
So, $J_{a,\chi}$ are representatives of these isomorphism classes. 
This proves the lemma. 
\end{proof}

Now assume that $G_{\der}$ is simply connected and that $\pi_1(G)_I$ is torsion-free.
We consider the exact sequence of abelian groups
\begin{equation}\label{exact sequence intro}
0\to X^*(G)\to X^*(T)\to P \to 0,
\end{equation}
and denote by $J(\mu_\ad)\subset \{\omega_1,\ldots, \omega_n\}$ the associated minimal test set from \Cref{minimal test sets}.
Here, $\omega_1,\ldots,\omega_n$ are the fundamental weights, which form a $\mathbb Z_{\geq 0}$-monoid basis for $P^+$ and are defined factorwise if the absolute Dynkin diagram of $G$ is not necessarily connected. 
We consider 
\begin{equation}
J(\mu):=\{\delta_1,\ldots,\delta_m\}\cup \tilde J(\mu_\ad) \subset X^*(T)^+,
\end{equation}
where $\delta_1,\ldots,\delta_m\in X^*(G)$ is a minimal set generating the $\mathbb Q_{\geq 0}$-monoid $X^*(G)_{\mathbb Q}$ and $\tilde J(\mu_\ad)\subset X^*(T)^+$ is a set of lifts of the elements in $J(\mu_\ad)$.
In addition, if $\mu_I\in X_*(T)_I^+$ is the image of $\dot \mu$ in the inertia coinvariants, then we can form the set $J(\mu_I)\subset X^*(T)^{I,+}$ in a similar way by using the échelonnage root system when $\pi_1(G)_I$ is torsion-free: 
the sequence in \Cref{exact sequence intro} induces a still exact sequence
\begin{equation}\label{exact sequence inertia}
0\to X^*(G)^I\to X^*(T)^I\to P^I \to 0,
\end{equation}
see \Cref{elementary group lemma}, \Cref{elementary group lemma 2}.
For every $I$-orbit $\calO$ on $\{\omega_1,\ldots,\omega_n\}$, we denote $\breve \omega_\calO=\sum_{\chi\in \calO}\chi\in P^I$, which are the fundamental weights for the échelonnage root system, as follows from \cite[Definition~3.3, Proposition~3.5, and Theorem~6.1]{Hai18}.
Let $\mu_{I,\ad}$ denote the image of $\mu_I$ in $X_*(T_\ad)_I$, which is a dominant coweight for $\Sigma$.
We denote by $J(\mu_{I,\ad})\subset \{\breve \omega_\calO\}_\calO$ the minimal test set for $\mu_{I,\ad}$ taken with respect to the échelonnage root system $\Sigma$, see \Cref{minimal test set definition} and \Cref{minimal test set}.
We then choose any lifts $J(\mu_I)\subset X^*(T)^{I,+}$ and note that it is not necessary to include in $J(\mu_I)$ monoid generators of $X^*(G)^I_{\mathbb Q}$.

\begin{lemma}\label{generating lemma J1}
Assume that $G_{\mathrm{der}}$ is simply connected and that $\pi_1(G)_I$ is torsion-free.
Let $J_1\subset X^*(T)^+$ be any subset satisfying
\begin{equation}\label{J1 conditions}
\begin{gathered}
J(\mu)\subset \operatorname{Cone}_{\mathbb Q_{\geq0}}(\Gal_E*J_1),
\qquad
J(\mu_I)\subset \operatorname{Cone}_{\mathbb Q_{\geq0}}\bigl(N_{\bar I}(J_1^0)\bigr),\\
\operatorname{Cone}_{\mathbb Q_{\geq0}}\bigl(N_{\bar I}(J_1^0\cap X^*(G))\bigr)
=X^*(G)^I_{\mathbb Q}.
\end{gathered}
\end{equation}
If $G_{\breve F}$ is split, then one may take $J_1=J(\mu)$.
Then, the following hold:
\begin{enumerate}
\item \label{generating lemma J1 1} Let $\nu\in X_*(T)_{\mathbb Q}^+$. If $\langle \gamma*\chi,\dot\mu-\nu\rangle\geq0$ for all $\chi\in J_1$ and $\gamma\in\Gal_E$, then
$$\langle\chi,\dot\mu-\nu\rangle\geq0\quad\text{for all }\chi\in X^*(T)^+_{\mathbb Q}.$$
\item \label{generating lemma J1 2} Let $\nu_I\in X_*(T)_{I,\mathbb Q}^+$. If $\langle N_{\bar I}\chi,\mu_I-\nu_I\rangle\geq0$ for all $\chi\in J_1^0$, then
$$\langle\chi,\mu_I-\nu_I\rangle\geq0\quad\text{for all }\chi\in X^*(T)^{I,+}_{\mathbb Q}.$$
\end{enumerate}
\end{lemma}
\begin{proof}
\Cref{generating lemma J1 1}:
The first containment in \Cref{J1 conditions} gives the inequalities for all elements in $J(\mu)$.
So, its central characters force equality on $X^*(G)_{\mathbb Q}$, and the remaining inequalities suffice on $P_{\mathbb Q}^+$ by the definition of $J(\mu_{\mathrm{ad}})$, as in \Cref{minimal test set prop}.

\Cref{generating lemma J1 2}:
The last condition in \Cref{J1 conditions} forces equality on $X^*(G)^I_{\mathbb Q}$.
The second containment gives the inequalities for $J(\mu_I)$, which suffice on $P_{\mathbb Q}^{I,+}$ by the definition of $J(\mu_{I,\mathrm{ad}})$.
The conclusion follows from \Cref{exact sequence inertia}.
\end{proof}

Finite sets satisfying \eqref{J1 conditions} exist: take $J(\mu)$, adjoin the $I$-invariant lifts $J(\mu_I)$, and adjoin central characters generating $X^*(G)^I_{\mathbb Q}$ as a $\mathbb Q_{\geq0}$-monoid.
The conditions in \eqref{J1 conditions} are weaker than \Cref{J_1 condition}:

\begin{lemma}\label{condition implies J1 conditions}
Assume that $G_{\mathrm{der}}$ is simply connected and that $\pi_1(G)_I$ is torsion-free.
Then every subset $J_1\subset X^*(T)^+$ satisfying \Cref{J_1 condition} satisfies \eqref{J1 conditions}.
\end{lemma}
\begin{proof}
As $J_1$ generates the $\mathbb Q_{\geq0}$-monoid $X^*(T)^+_{\mathbb Q}$, the first containment in \eqref{J1 conditions} holds.
As $N_{\bar I}(J_1^0)$ generates the $\mathbb Q_{\geq0}$-monoid $X^*(T)^{I,+}_{\mathbb Q}$, so does the second one.
For the last condition, let $\delta\in X^*(G)^I_{\mathbb Q}$.
As central characters are dominant, we write $\delta=\sum_{\chi\in J_1^0}c_\chi N_{\bar I}(\chi)$ for some $c_\chi\in\mathbb Q_{\geq 0}$.
We apply the projection $X^*(T)^I_{\mathbb Q}\to P^I_{\mathbb Q}$, $\chi\mapsto\bar\chi$, from \eqref{exact sequence inertia}, which commutes with the norm.
This gives $0=\sum_\chi c_\chi N_{\bar I}(\bar\chi)$, where each $N_{\bar I}(\bar\chi)=\sum_{\sigma\in\bar I}{}^\sigma\bar\chi$ is a sum of dominant elements of $P_{\mathbb Q}$, as $I$ preserves the dominant chamber of $B_0$.
Hence, whenever $c_\chi>0$, we get $N_{\bar I}\bar\chi=0$ and so $\bar\chi=0$, i.e., $\chi\in X^*(G)$.
Thus $\delta$ lies in the $\mathbb Q_{\geq0}$-monoid generated by $N_{\bar I}(J_1^0\cap X^*(G))$.
\end{proof}

\begin{definition}\label{definition min test set}
Let $J_1\subset X^*(T)^+$ be a finite inclusion-minimal set satisfying \Cref{J1 conditions}; when $G_{\breve F}$ is split, one takes $J_1=J(\mu)$.
Define the associated minimal test set by
$$
J_{\min}(\calG,\mu)
=
\bigl(J_1^{\mathrm{cen}}\times\{B_0\}\bigr)
\cup
\bigcup_{\substack{\chi\in J_1\setminus J_1^{\mathrm{cen}}\\
a\in\operatorname{Vert}(\mathbf f)}}
\{\chi\}\times J_{a,\chi},
$$
where $J_1^{\mathrm{cen}}=J_1\cap X^*(G)$.
\end{definition}

We prove in \Cref{inclusion permissible locus theorem} that the special fiber of $\bbM_{\calG,\mu}^J$ for $J=J_{\min}(\calG,\mu)$ is contained in the permissible locus for arbitrary $\calG$ and $\mu$.

\begin{remark}\label{refined conjecture minuscule} 
The companion paper \cite{LRVW:Minuscule} proves that, if
$G_{\mathrm{der}}$ is simply connected, $\pi_1(G)_I$ is torsion-free,
$\mu$ is minuscule, and $G_{\breve F}$ is split, then $\mathbb M_{\calG,\mu}^{J}=\mathbb M_{\calG,\mu}$ for $J=J_{\min}(\calG,\mu)$.
In the classical types, this follows from \Cref{inclusion permissible locus theorem} using the equality of the admissible
and permissible loci; see \Cref{permissible corollary}.
In types $E_6$ and $E_7$, these loci differ \cite{Sauermann}, and the proof relies on computer algebra calculations.
Neither the minuscule hypothesis nor the splitness hypothesis can be dropped.
The companion paper gives counterexamples for split
$\mathrm{GSpin}_7$ with $\mu_\ad=\omega_3^\vee$ and for a suitable central extension of ramified $\mathrm{SU}_6$ with minuscule $\mu$ of signature $(3,3)$.
\end{remark}

\section{Containment in the permissible locus}
The aim of this section is to show that the special fiber of the refined moduli candidate $\mathbb M_{\calG,\mu}^J$ for $J=J_{\mathrm min}(\calG,\mu)$ is contained in the permissible locus (\Cref{inclusion permissible locus theorem}).
Whenever the permissible locus agrees with the admissible one, this implies $\mathbb M_{\calG,\mu}^J=\mathbb M_{\calG,\mu}$, hence \Cref{refined conjecture}, \Cref{moduli conjecture}, and \Cref{moduli conjecture admissible} in these cases.
The known cases of this equality are listed in \Cref{permissible corollary}.

\subsection{Preliminary reduction}
For a subset $J\subset \Rep(\calG_{O_E})$, recall from \Cref{special fiber condition} the sub-functor $A_{G,\mu}^J$ of $\Fl_{\calG,k_E}$ given by
	$$
	{A}_{\calG, \mu}^J= \bigcap_{(\rho,\calV)\in J} \tilde\rho_*^{-1}\left(\Fl_{\mathrm{GL}(\calV)}^{\geq {\vert\underline{\mu}_{\max}(\calV)\vert}}\right).
	$$
With the notation as in \Cref{refined moduli candidate equation}, we have $({A}_{\calG, \mu}^J)^\diamondsuit=\mathbb M_{\calG,\mu, k_E}^J$.

\begin{lemma}\label{left to prove}
Assume that $G_\der$ is simply connected and that $\pi_1(G)_I$ is torsion-free.
Let $J=J_{\min}(\calG,\mu)$ be as in \Cref{definition min test set}.
Then, the inclusion $\mathbb M_{\calG,\mu}\subset \mathbb M_{\calG,\mu}^J$ is an equality on generic fibers. 
In particular, one has $\mathbb M_{\calG,\mu}^J=\mathbb M_{\calG,\mu}$ if and only if $A_{\calG,\mu}(\bar k)=A^J_{\calG, \mu}(\bar k)$. 
\end{lemma}
\begin{proof}
The inclusion is proven in \Cref{easy inclusion model}.
For the equality on generic fibers, it suffices to use the component of \Cref{easy inclusion model} corresponding to the fixed embedding $E\subset C$.
For every $\chi\in J_1$, the $B_{0,C}$-highest weights of the test representation form $\Gal_E*\chi$, and all have $\mu$-weight $\langle\chi,\dot\mu\rangle$ because $\Gal_E$ fixes $\dot\mu$.
Thus its lattice condition gives $\langle\gamma*\chi,\dot\mu-\nu\rangle\geq0$ for every $\gamma\in\Gal_E$ and every dominant Schubert parameter $\nu$ in the generic fiber.
\Cref{generating lemma J1}, \Cref{generating lemma J1 1}, and \Cref{root lemma} give the required equality, as in the proof of \Cref{minimal test set prop}. 
Finally, the inclusion $\mathbb M_{\calG,\mu}\subset \mathbb M_{\calG,\mu}^J$ is an equality if it is so fiberwise by \cite[Lemma~12.5]{Sch17}. 
If $A_{\calG,\mu}(\bar k)=A^J_{\calG, \mu}(\bar k)$, then $A_{\calG,\mu}=A^J_{\calG, \mu}$, and hence $A_{\calG,\mu}^\diamondsuit=(A^J_{\calG, \mu})^\diamondsuit$.
\end{proof}

\subsection{The permissible locus}\label{permissible sec}
Let $\mathbf f$ be the facet in the apartment $\mathscr A=\mathscr A(G,S,\breve F)$ given by $\calG(\breve O)$. 
Note that $\mathscr A$ is a torsor under $X_*(T)_{I, \mathbb R}$.
In particular, for any elements $a,b\in \mathscr A$, we denote by $a-b\in X_*(T)_{I, \mathbb R}$ the unique element such that $b + (a-b)=a$.
For $x\in W$, write $\bar x\in W_0$ for its image. The action of the Iwahori-Weyl group $W$ on such a difference through $W_0$ is denoted $\bar x(a-b)$.

Let $\dot\mu\in X_*(T)^+$ be the dominant representative of $\mu$, and $\mu_I$ its image in the coinvariants $X_*(T)_I$.  
Let $\Conv(W_0\cdot \mu_I)\subset X_*(T)_{I, \mathbb R}$ be the convex hull of the $W_0$-translates of $\mu_I$.
A double class $[x]\in W_{\mathbf f}\backslash W/W_{\mathbf f}$ is called \textit{$\mu$-permissible} if the following two conditions hold:
\begin{enumerate}
\item $[x]\equiv \mu_I$ in $\pi_1(G)_I$, and 
\item $\delta_a(x):=\bar x^{-1}\!\left(a-x(a)\right)\in \Conv(W_0\cdot \mu_I)$ for all $a\in \mathrm{Vert}(\mathbf f)$,
\end{enumerate}
where $x\in W$ is any representative of $[x]$.
The second condition is independent of the choice of representative. 

\begin{remark}\label{permissible remark}
The permissibility for split groups in \cite{HN02, HH17} is equivalent to our definition under the plus sign convention (\Cref{plus translation convention}) for the translation action on the apartment.
\begin{enumerate}
\item \label{permissible remark 1}
	We use the Bruhat--Tits sign convention from \Cref{BT translation convention}: $\nu_I(\pi)\in T(\breve F)/\calT(\breve O)$ acts as translation by $-\nu_I$ on $\mathscr A$.
	In particular, we see $\delta_a(\nu_I(\pi))=a-(a-\nu_I)=\nu_I$ for all $\nu_I\in X_*(T)_I$.
\item \label{permissible remark 2}
	Fix a special vertex $o\in \mathscr A$.
	Write $x=\nu_I(\pi)\cdot u$ with $\nu_I\in X_*(T)_I$ and $u\in W_0$ according to the decomposition $W=X_*(T)_I\rtimes W_0$ from \Cref{translation presentation}.   
	Then, we get
	$$
	\delta_a(x) = \bar x^{-1}\!(a - ua + \nu_I)=u^{-1}a-a+u^{-1}\nu_I.
	$$

\item \label{permissible remark 3}
	We have $\delta_a(x)\in \Conv(W_0\cdot \mu_I)$ if and only if $a-x(a)\in \Conv(W_0\cdot \mu_I)$ if and only if $a-x'(a)\in \Conv(W_0\cdot \mu_I)$ for any $x'\in xW_a$. 
	Twisting $a-x(a)$ by ${x^{-1}}$ in the definition of $\delta_a(x)$ is just a convenient normalization for our purposes.
\end{enumerate}
\end{remark}

\begin{definition}
The \textit{$\mu$-permissible locus $P_{\calG,\mu}$ in $\Fl_{\calG,k_E}$} is the $k_E$-descent of $\bigcup_{[x]}\Fl_{\calG,\bar k, [x]}$ where $[x]\in W_{\mathbf f}\backslash W/ W_{\mathbf f}$ runs through the $\mu$-permissible classes. 
\end{definition}

The permissible locus contains the admissible one but the converse does not hold in general \cite{HN02}.
The cases where the two loci are known to agree are listed in \Cref{permissible corollary}.

\begin{theorem}\label{inclusion permissible locus theorem}
Assume that $G_\der$ is simply connected and that $\pi_1(G)_I$ is torsion-free.
Let $J=J_{\min}(\calG,\mu)$ be as in \Cref{definition min test set}.
Then, one has inclusions 
\begin{equation}
A_{\calG,\mu}\subset A_{\calG,\mu}^J\subset P_{\calG,\mu}
\end{equation}
as sub-functors of $\Fl_{\calG,k_E}$.
In particular, if admissible and permissible loci agree, then $\mathbb M_{\calG,\mu}^J=\mathbb M_{\calG,\mu}$ by \Cref{left to prove}.
\end{theorem}

The proof of \Cref{inclusion permissible locus theorem} is given in \Cref{proof of permissible theorem sec}.
It implies the analogous statement for the test sets $J(\calG,\mu)$ by \Cref{condition implies J1 conditions}.

For the formulation of the next lemma, let us recall the (rational) Bruhat order in the ramified case.
Let $\Delta^\vee$ denote the (absolute) simple coroots for $B_0$.
The inertia coinvariants $(\mathbb Z\Delta\!^\vee)_I$ are torsion-free and admit a basis $\Delta^\vee_I$ in bijection with the $I$-orbits on $\Delta^\vee$. 
For $\nu_I, \la_I\in X_*(T)_I$, we have the Bruhat order
\begin{equation}\label{Bruhat order}
\nu_I\leq \la_I  \iff \la_I-\nu_I \in \mathbb Z_{\geq 0}\Delta_I^\vee,
\end{equation}
and for $\nu_I, \la_I\in X_*(T)_{I,\mathbb Q}$ its rationalized version
\begin{equation}\label{Bruhat order rational}
\nu_I\leq_{\mathbb Q} \la_I  \iff \la_I-\nu_I \in \mathbb Q_{\geq 0}\Delta_I^\vee,
\end{equation}
see \cite[Lemma~2.5, Remark~2.7~ff.]{ALRR24}.

\begin{lemma}\label{root lemma twisted}
For all $\nu_I,\la_I\in X_*(T)_{I, \mathbb Q}$, the following are equivalent:
\begin{itemize}
	\item $\nu_I \leq_{\mathbb Q} \la_I$, and
	\item $\langle \chi,\la_I-\nu_I\rangle \geq 0$ for all $\chi\in (X^*(T)^+)^I$ (equivalently, for all $\chi\in (X^*(T)^+_{\mathbb Q})^I$). 
\end{itemize}
Whenever these equivalent conditions hold, the following properties are satisfied:
\begin{enumerate}
	\item \label{root lemma twisted 1} If, in addition, $\nu_I,\lambda_I\in X_*(T)_{I, \mathbb Q}^+$ are dominant, then $\nu_I\in \mathrm{Conv}(W_0\cdot\la_I)$. 
	\item \label{root lemma twisted 2} 
		If $\nu_I,\la_I\in X_*(T)_{I}$ are integral, $G_\der$ is simply connected, and $\pi_1(G)_I$ is torsion-free, then $\nu_I\leq \la_I$. 
\end{enumerate}
\end{lemma}
\begin{proof}
The proof is analogous to that of \Cref{root lemma}:
Assume $\langle \chi,\la_I-\nu_I\rangle \geq 0$ for all $\chi\in (X^*(T)^+_{\mathbb Q})^I$.
The isomorphism $X_*(T)_{\mathbb Q}\cong X_*(Z_G)_{\mathbb Q}\times \mathbb Q\Delta\!^\vee$ is $I$-equivariant and induces an isomorphism $X_*(T)_{I, \mathbb Q}\cong X_*(Z_G)_{I,\mathbb Q}\times (\mathbb Q\Delta\!^\vee)_I$ on coinvariants.
So, we can write 
\begin{equation}\label{expand elements eq}
\la_I-\nu_I=\tau_I +\sum_{\alpha^\vee_I \in \Delta^\vee_I} c_{\alpha_I^\vee}\alpha^\vee_I
\end{equation}
for unique $\tau_I\in X_*(Z_G)_{I,\mathbb{Q}}$ and unique $c_{\alpha_I^\vee} \in \mathbb{Q}$ for all $\alpha_I^\vee \in \Delta^\vee_I$. 
Pairing with central elements and their inverses, we deduce $\langle \chi,\la_I-\nu_I\rangle = 0$ for all $\chi\in X^*(Z_G)_{\mathbb{Q}}^I$, hence $\tau_I=0$ by the non-degeneracy of the pairing on the center. 
Pairing with the $I$-averages  of the fundamental weights, we deduce $c_{\alpha_I^\vee}\geq 0$. 
Hence, $\nu_I\leq_{\mathbb Q} \la_I$ holds.
The converse direction is immediate from the definition of dominant elements. 
Now, we assume that the equivalent conditions are satisfied. 	
	
	\Cref{root lemma twisted 1}:
Applying \cite[Remark~4.6]{HH17} to the relative root system, one obtains
\[
\operatorname{Conv}(W_0\cdot\lambda_I)
\cap X_*(T)_{I,\mathbb Q}^{+}
=
\left\{
\la_I'\in X_*(T)_{I,\mathbb Q}^{+}
\;\middle|\;
\la'_I\leq_{\mathbb Q}\la_I
\right\},
\]
which implies \Cref{root lemma twisted 1}.
	
	\Cref{root lemma twisted 2}:
	Assume that we are in the situation of \Cref{root lemma twisted 2}.
	We have a short exact sequence of torsion-free abelian groups (\Cref{elementary group lemma} and its proof):
	$$
	0\to X_*(T_{\mathrm{der}})_I\to X_*(T)_I\to \pi_1(G)_I\to 0
	$$
	using the assumption that $G_{\mathrm{der}}$ is simply connected and $\pi_1(G)_I$ is torsion-free. 
	In particular, $X_*(T_{\mathrm{der}})_I$ is saturated in $X_*(T)_I$.
	Since $\nu_I$ and $\la_I$ are integral, this implies $\la_I-\nu_I \in X_*(T)_I\cap X_*(T_\der)_{I,\mathbb Q}=X_*(T_{\mathrm{der}})_I=(\mathbb{Z}\Delta^\vee)_I$. 
	Thus, all coefficients $c_\alpha$ are integral, i.e., $\nu_I \leq \la_I$.
\end{proof}

\begin{corollary}\label{central root lemma}
For all $\nu_I,\la_I\in X_*(T)_{I, \mathbb Q}$, the following are equivalent:
\begin{itemize}
	\item $\nu_I \equiv \la_I$ in $\pi_1(G)_{I,\mathbb Q}$, and
	\item $\langle \chi,\la_I-\nu_I\rangle \geq 0$ for all $\chi\in X^*(G)^I$ (equivalently, for all $\chi\in X^*(G)_{\mathbb Q}^I$). 
\end{itemize}
In particular, $\nu_I\equiv \la_I$ in $\pi_1(G)_I$ when $\nu_I,\la_I\in X_*(T)_I$ and $\pi_1(G)_I$ is torsion-free.
\end{corollary}
\begin{proof}
Writing $\la_I-\nu_I$ as in \Cref{expand elements eq}, the proof of \Cref{root lemma twisted} shows that both conditions are equivalent to $\tau_I=0$.
So, the corollary holds for $X^*(G)^I_{\mathbb Q}$ replaced by $X^*(Z_G)^I_{\mathbb Q}$. 
However, the central isogeny $Z_G\to G/G_\der$ has finite kernel, hence induces isomorphisms $X^*(Z_G)_{\mathbb Q}=X^*(G/G_\der)_{\mathbb Q}=X^*(G)_{\mathbb Q}$, which are Galois equivariant.
This implies the corollary. 
\end{proof}

\begin{remark}\label{permissible corollary}
We now list the cases in which the admissible and permissible loci coincide.
To this end, let $\Sigma=\coprod_j\Sigma_j$ be the decomposition of the
échelonnage root system of $G_{\breve F}$ into irreducible components,
and let $(\mu_{I,j})_j$ be the components of the image of $\mu_I$
in $X_*(T_\ad)_I$.
Then the admissible and permissible loci agree in each of the following cases:
\begin{enumerate}
\item\label{permissible corollary 1}
Every vertex of $\mathbf f$ is special, with no restriction on $\mu$.
\item\label{permissible corollary 2}
For each component $(\Sigma_j,\mu_{I,j})$, one of the following
conditions holds:
\begin{enumerate}
\item $\mu_{I,j}=0$, with no restriction on $\Sigma_j$;
\item\label{permissible corollary 3}
$\Sigma_j$ has type $A_n$, with no restriction on $\mu_{I,j}$;
\item $\Sigma_j$ has type $B_n$, $n\geq3$, and $\mu_{I,j}$ is
minuscule or quasi-minuscule;
\item\label{permissible corollary 4}
$\Sigma_j$ has type $C_n$, $n\geq2$, and
\[
\mu_{I,j}=m\omega_n^\vee\quad(m\in\mathbb Z_{\geq0}),
\qquad\text{or}\qquad
\mu_{I,j}=e_1+\cdots+e_s\quad(1\leq s<n);
\]
\item $\Sigma_j$ has type $D_n$, $n\geq4$, and $\mu_{I,j}$ is
minuscule or quasi-minuscule;
\item $\Sigma_j$ has type $B_2=C_2$, $C_3$, or $G_2$, with no
restriction on $\mu_{I,j}$.
\end{enumerate}
\end{enumerate}
For proofs and references, see the companion paper \cite{LRVW:Minuscule}.
In these cases, under the hypotheses of \Cref{inclusion permissible locus theorem} and for $J=J_{\min}(\calG,\mu)$, we obtain $\mathbb M_{\calG,\mu}^J=\mathbb M_{\calG,\mu}$.
\end{remark}

\subsection{Extremal weight modules in representations}\label{weight modules}
The proof of \Cref{inclusion permissible locus theorem} is based on a finer study of extremal weight modules in the representations $\calV_{\chi, B}$ defined in \Cref{parahoric test representations}. 
The key input for the theorem is \Cref{key corollary permissible} below.

Throughout this section, we fix a Borel subgroup $B\in\mathrm{Bor}(G_{\breve F},S_{\breve F})$, a $B$-dominant character $\chi\in X^*(T)$, and a finite field extension $\tilde E\subset C$ of $\breve E$ that splits $T$ and that is Galois over $\breve F$.
We denote by $\tilde O:=O_{\tilde E}$ its ring of integers. 

For any vertex in the apartment $a\in \mathscr A=\mathscr A(G,S,\breve F)$ with associated parahoric $\breve O$-model $\calG_a$, we consider the representation $\calG_a\to \Res_{O_{\breve E}/\breve O}\GL(\calV_a)$ given by
\begin{equation}\label{reps for proof}
\calV_{a} := \Res_{\tilde O/ O_{\breve E}}(\tilde \calV_{a})\quad\text{with}\quad \tilde \calV_{a}:= \Ind_{{\calB}_{a,\tilde O}^{-}}^{\calG_{a, \tilde O}}({\chi}),
\end{equation}
where $\calB_a^-$ is the closure of the Borel subgroup opposite to $B$ inside $\calG_a$.

The generic fiber $\tilde V_{} := \tilde\calV_{a}\otimes_{\tilde O}\tilde E$ is the dual Weyl module of the split group $G_{\tilde E}$ of $B$-highest weight $\chi$, and $V := \calV_{a}\otimes_{O_{\breve E}}\breve E$ its restriction of scalars to $\breve E=E\otimes_{E^{\mathrm{ur}}}\breve F$. 
For an extremal weight $\eta\in W_{0,C} \cdot{\chi}$, the resulting weight space $\tilde V(\eta)$ is one-dimensional \cite[Section~II.2.4]{Jan03}.
Hence, the intersection
$
\tilde \calV_{a}(\eta) := \tilde \calV_{a}\cap \tilde V(\eta)
$
is a free $\tilde O$-submodule in $\tilde \calV_a$ of rank $1$ that is stable under the action of $\calT_{\tilde O}$.
We consider the induced $O_{\breve E}$-submodule in $\calV_a$ and its highest exterior power
\begin{equation}\label{weight for proof}
\calV_a(\eta):= \Res_{\tilde O/ O_{\breve E}}\big(\tilde \calV_a(\eta)\big) \quad\text{and}\quad \wedge^{\mathrm{top}}\calV_a(\eta)
\end{equation}
which are stable under $\calT_{O_{\breve E}}$, and free over $O_{\breve E}$ of rank $[\tilde E:\breve E]$ and of rank $1$, respectively.

The absolute finite Weyl group $W_{0,C}$ acts on $X^*(T)$ by the formula $(w\eta)(t)=\eta(w^{-1}tw)$. 
So, we get an action of the Iwahori-Weyl group by using the composition $W\to W_0\subset W_{0,C}$.
For this reason, we now restrict our attention to $\eta \in W_0\cdot \chi$, i.e., 
the orbit of $\chi$ under the action of the relative finite Weyl group $W_0\subset W_{0,C}$ sitting inside the absolute one.
The Iwahori--Weyl action gives an inclusion of $O_{\breve E}$-submodules
\begin{equation}
x \cdot\calV_a(\eta)\subset V(\bar x\eta)
\end{equation}
stable under the action of $\calT_{O_{\breve E}}$.

\begin{lemma}\label{weight module properties}
Let $x\in W$ and $\eta\in W_{0}\cdot \chi$.
Assume either that $x$ is a translation element or that $a\in \scrA$ is a special vertex (or both).
Then, one has
$$
	\wedge^{\mathrm{top}}(x\cdot \calV_a(\eta))=\pi_E^{\langle N{\eta},\delta_a(x)\rangle}\cdot\wedge^{\mathrm{top}}\calV_a(\bar x\eta)
$$
where $\pi_E$ is a uniformizer of $E$ and $N\eta= \sum_{\gamma\in \Gal(\tilde E/\breve F)}{^\gamma\eta}\in X^*(T)^I$ is the norm of $\tilde E$ over $\breve F$.   
\end{lemma}
\begin{proof}
We apply the Kottwitz isomorphism from \Cref{coinvariants identification} for the fields $\tilde E$ and $\breve F$.
By \cite[Equation~(7.3.2)]{Kot97}, we have $\nu_I(\pi)=(N\nu_I)(\pi_{\tilde E})$ where $N\colon X_*(T)_I\to X_*(T)$ is the norm of $\tilde E$ over $\breve F$ and $\pi_{\tilde E}$ is a uniformizer of $\tilde E$. 

First, assume $x=\nu_I(\pi)\in T(\breve F)/\calT(\breve O)$ is a translation element. 
Then, we get
$$
(N\nu_I)(\pi_{\tilde E})\cdot \tilde \calV_a(\eta)=\pi_{\tilde E}^{\langle \eta, N\nu_I\rangle}\cdot \tilde \calV_a(\eta)
$$ 
by the definition of weight spaces. 
Note that $\langle \eta, N\nu_I\rangle=\langle N\eta, \nu_I\rangle$ by Galois invariance of the pairing.
So, taking highest exterior powers implies
$$
\wedge^{\mathrm{top}}\left(\nu_I(\pi)\cdot \calV_a(\eta)\right)=\pi_E^{\langle N\eta, \nu_I\rangle}\cdot \wedge^{\mathrm{top}}\calV_a(\eta),
$$ 
using that $\mathrm{Norm}_{\tilde E/\breve E}(\pi_{\tilde E})$ and $\pi_E$ generate the same ideals in $O_{\breve E}$. 
The lemma follows from $\delta_a(\nu_I(\pi))=\nu_I$ in this case (\Cref{permissible remark}, \Cref{permissible remark 1}).

Next, assume that $a=o$ is a special vertex.
Write $x=\nu_I(\pi)\cdot u$ with $\nu_I\in X_*(T)_I$ and $u\in W_0$ accordingly.   
Since the elements in $W_0$ admit representatives in $\calG_a(\breve O)$ and $\bar x\eta=u\eta$, we get $u\cdot \tilde\calV_a(\eta)=\tilde\calV_a(\bar x\eta)$.
As in the case of translation elements above, we compute
\begin{equation}\label{weight module properties computation}
x\cdot \tilde\calV_a(\eta)= \nu_I(\pi)\cdot \tilde \calV_a(\bar x\eta)=\pi_{\tilde E}^{\langle N(\bar x\eta), \nu_I\rangle}\cdot \tilde \calV_a(\bar x\eta).
\end{equation}
Since $W_0$ fixes the special vertex $a=o$, we have $\delta_a(x)=u^{-1}\nu_I=\bar x^{-1}\nu_I$ by \Cref{permissible remark}, \Cref{permissible remark 2}.
By invariance of the pairing under the Weyl group $W_0$, we have $\langle N(\bar x\eta), \nu_I\rangle=\langle N\eta, \delta_a(x)\rangle$.
Hence, the lemma follows from \Cref{weight module properties computation} by taking highest exterior powers. 
\end{proof}

 \begin{definition}\label{elementary divisor definition}
For any $x\in W$ and $\eta\in W_0\cdot \chi$, define the integer
$$
d=d_a(x,\eta)\in \mathbb Z
$$
as the elementary divisor of the rank one $O_{\breve E}$-lattices $\wedge^{\mathrm{top}}\calV_{a}(\bar x\eta)$ and $\wedge^{\mathrm{top}}(x\cdot \calV_{a}(\eta))$ inside $\wedge^{\mathrm{top}}V(\bar x\eta)$, i.e., 
$$
\wedge^{\mathrm{top}}(x\cdot \calV_{a}(\eta))\;=\;\pi_E^d\cdot \wedge^{\mathrm{top}}\calV_{a}(\bar x\eta).
$$
 \end{definition}

 \begin{proposition}\label{estimate prop}
 Assume that $\eta$ is $B$-dominant, i.e., ${\eta}=\chi$. 
 Then
 \begin{equation}\label{key inequality}
 \langle N\eta,\delta_a(x)\rangle \leq d_a(x,\eta).
 \end{equation}
 In addition, equality holds if $x$ is a translation or if $a\in \mathscr A$ is a special vertex. 
 \end{proposition}
 \begin{proof}
 If $x$ is a translation or $a\in \mathscr A$ is special, then $d_a(x,\eta)= \langle N\eta,\delta_a(x)\rangle$ by \Cref{weight module properties}.
 The inequality in \Cref{key inequality} requires more elaboration and follows from \Cref{key lemma} below. 
 \end{proof}
 
 The following result is the key input for the proof of \Cref{inclusion permissible locus theorem}.
 
 \begin{corollary}\label{key corollary permissible}
 Let $d\in\mathbb Z$. 
 If one has an inclusion $\pi_E^d\calV_a\subset x\calV_a$ of $O_{\breve E}$-lattices in $V$, then
 \[
 \langle N\chi,\delta_a(x)\rangle\leq[\tilde E:\breve E]d.
 \]
 \end{corollary}
 \begin{proof}
Intersecting with $V(\bar x\chi)$ and taking highest exterior powers gives
\[
\pi_E^{[\tilde E:\breve E]d}\wedge^{\mathrm{top}}\calV_a(\bar x\chi)
\subset\wedge^{\mathrm{top}}(x\calV_a(\chi))
\subset\pi_E^{\langle N\chi,\delta_a(x)\rangle}\wedge^{\mathrm{top}}\calV_a(\bar x\chi),
\]
where the second inclusion is \Cref{estimate prop}.
Comparing exponents proves the assertion.
 \end{proof}
 
In order to show the inequality in \Cref{key inequality}, it is convenient to fix a special vertex $o\in \mathscr A$ whose image under the scalar extension map
\begin{equation}\label{scalar extension map}
\iota\colon \mathscr A=\mathscr A(G,S,\breve F)\hookto \mathscr A(G,T,\tilde E)=\tilde{\mathscr A}
\end{equation}
is hyperspecial, see \cite[Construction~6.1.21]{KP23}.
The choice of $o$ induces an identification $\mathscr A=X_*(T)_{I, \mathbb R}$ and a special parahoric $\breve O$-model $\calG_o$ of $G_{\breve F}$.
We consider the representation $\calG_o\to \Res_{O_{\breve E}/\breve O}(\GL(\calV_o))$ associated with the vertex $o$ as above. 
Let $h=h_{a}(\eta)\in \mathbb Z$ be the elementary divisor of $\wedge^{\mathrm{top}} \calV_{o}(\eta)$ and $ \wedge^{\mathrm{top}}\calV_{a}(\eta)$ inside $\wedge^{\mathrm{top}}V(\eta)$, i.e., $\wedge^{\mathrm{top}}\calV_a(\eta)=\pi_E^h\cdot\wedge^{\mathrm{top}}\calV_o(\eta)$.

\begin{lemma}\label{h computation}
The following properties hold:
\begin{enumerate}
	\item \label{h computation 1} $h_a(\chi)=0$
		\item \label{h computation 2} If $a\in{\mathscr C_B(o)}$ (resp.~$a\in \mathscr C_{B^-}(o)$), then $h_a(\eta)\geq 0$ (resp.~$h_a(\eta)\leq 0$) for all $\eta\in W_{0,C}\cdot \chi$. 
	\item \label{h computation 3} 
		$N(\chi-\eta)(a)\geq h_a(\eta)$ for all $\eta\in W_0\cdot \chi$ with equality if $a-o\in X_*(T)_{I,\mathbb R}$ lifts to an element in $X_*(T)_I$.
\end{enumerate}
\end{lemma}
\begin{proof}
\Cref{h computation 1}:
\Cref{embedding corollary}, \Cref{embedding corollary invariants} gives a canonical isomorphism $\tilde \calV_b(\chi)=\tilde O$ for any vertex $b\in \mathscr A$. 
Hence, $\tilde\calV_a(\chi)=\tilde O=\tilde\calV_o(\chi)$, and restriction of scalars gives $h_a(\chi)=0$.

\Cref{h computation 2}:
For every vertex $b\in \mathscr A$, the big cell $\calU_b\times \calB_b^-\to \calG_b$ induces an embedding $\tilde\calV_b\subset\tilde O[\calU_b]$ such that $\tilde \calV_b(\eta)=\tilde V(\eta)\cap \tilde O[\calU_b]$, see \Cref{embedding corollary}, \Cref{embedding corollary saturation}.
The generic fiber $U:=\calU_b\otimes_{\breve O}\breve F$ is the unipotent radical of $B$, thus independent of the choice of $b$.
If $a\in {\mathscr C_B(o)}$, then we have $\calU_o(\breve O)\subset \calU_a(\breve O)$ (\Cref{unipotent inclusion lemma}) and hence $\tilde O[\calU_a]\subset \tilde O[\calU_o]$ as $\tilde O$-submodules in $\tilde E[U]$.
Intersecting with $\tilde V(\eta)$, we get $\tilde \calV_a(\eta)\subset \tilde \calV_o(\eta)$ as $\tilde O$-submodules in $\tilde V(\eta)$. 
This implies $\tilde \calV_a(\eta)=\pi_{\tilde E}^h\cdot \tilde \calV_o(\eta)$ for some $h\in \mathbb Z_{\geq 0}$.
Taking restriction of scalars and highest exterior powers replaces $\pi_{\tilde E}$ by its norm to $O_{\breve E}$. 
As $\tilde E$ is totally ramified over $\breve E$, we get
\begin{equation}\label{taking highest exterior powers}
\wedge^{\mathrm{top}} \calV_a(\eta)=\pi_{E}^h\cdot \wedge^{\mathrm{top}} \calV_o(\eta),
\end{equation}
and hence $h_a(\eta)=h\geq 0$.
Conversely, if $a\in {\mathscr C_{B^-}(o)}$, then $\calU_a(\breve O)\subset \calU_o(\breve O)$ and so $\tilde \calV_o(\eta)\subset \tilde \calV_a(\eta)$, which implies $h_a(\eta)\leq 0$.

\Cref{h computation 3}:
Let $\eta\in W_0\cdot\chi$ and write $\tilde\calV_a(\eta)=\pi_{\tilde E}^h\cdot \tilde\calV_o(\eta)$ for some $h\in\bbZ$. 
By \Cref{inclusion lemma key},
\[
 h\leq\lfloor N(\chi-\eta)(a)\rfloor\leq N(\chi-\eta)(a).
\]
As in \Cref{taking highest exterior powers}, we have $h=h_a(\eta)$. 
If $a-o$ lifts to an element of $X_*(T)_I$, then \Cref{inclusion lemma key} gives $h=N(\chi-\eta)(a)$, proving the equality assertion.
\end{proof}

The following lemma is used in the proof of \Cref{h computation 3} in \Cref{h computation}.

\begin{lemma}\label{inclusion lemma key}
For every $\eta \in W_0\cdot \chi$, one has an inclusion
\begin{equation}\label{inclusion equation integral lattices extremal weights}
 \pi_{\tilde E}^{\lfloor N(\chi-\eta)(a)\rfloor}\cdot \tilde \calV_o(\eta)
 \subset \tilde\calV_a(\eta).
\end{equation}
If, in addition, $a-o\in X_*(T)_{I,\mathbb R}$ lifts to an element in $X_*(T)_I$, then $N(\chi-\eta)(a)\in \mathbb Z$ and equality holds in \Cref{inclusion equation integral lattices extremal weights}.
\end{lemma}

\begin{proof} 
Choose a generator $P_\eta$ of the free rank-one $\tilde O$-module $\tilde\calV_o(\eta)$.
Note that under the embedding $\tilde\calV_o(\eta)\subset \tilde E[U]$, the element $P_\eta$ is homogeneous of weight $-\zeta$, where $\zeta:=\chi-\eta$, for the conjugation action of $T$ on $ \tilde E[U]$;
see \Cref{embedding corollary} \Cref{embedding corollary equivariance}.
Here, the conjugation action is normalized by $(t\cdot g)(u)=g(t^{-1}ut)$ throughout, so that the
coordinate $X_\al$ below has weight $-\al$.
Equivalently, $\sum_\al n_\al \al=\zeta$ for every monomial occurring in $P_\eta$ below.
Then, \Cref{inclusion equation integral lattices extremal weights} is equivalent to 
\begin{equation}\label{containment of Peta}
\pi_{\tilde E}^{\lfloor (N\zeta)(a)\rfloor}\cdot P_\eta\in \tilde \calV_a(\eta).
\end{equation}
We prove this in two steps: first assuming that \(T_{\breve F}\) is split, and then reducing the general case to this one.

\emph{Proof in the split case.}
Assume \(T_{\breve F}\) is split and $\tilde E=\breve F$. 
As in \Cref{split apartment conventions}, choose isomorphisms $\calU_{\al,o}\cong \mathbb G_{a,\breve O}$ for every $B$-positive root $\al\in \Phi_B^+$, and denote by $X_\al$ the induced coordinate.
After fixing an order on $\Phi_B^+$, we can write $\breve O[\calU_o]=\breve O[X_\al\mid \al\in \Phi_B^+]$ (cf.~\Cref{root group isomorphism}), and hence 
\[
        P_{\eta}
        =
        \sum_{\mathbf n=(n_\al)_\al} c_{\mathbf n}
        \prod_{\alpha\in\Phi_B^+}X_\alpha^{n_\alpha},
        \qquad c_{\mathbf n}\in \breve O .
\]
Then, one has
\begin{equation}\label{h formula}
h_a(\eta) = \max_{c_{\mathbf n}\neq 0} \big\{\textstyle{\sum_{\al}}n_\al\lfloor \al(a)\rfloor-\mathrm{val}_F(c_{\mathbf n})\big\}.
\end{equation}
Indeed, on root group $\breve O$-models, we have $\calU_{\al, a}(\breve O)=\pi^{-\lfloor\al(a)\rfloor}\calU_{\al,o}(\breve O)$ by \Cref{split apartment conventions}, so on coordinate rings $\breve O[\calU_{\al, a}]=\breve O[Y_\al]$ with $Y_\al=\pi^{\lfloor \al(a)\rfloor}X_\al$. 
This gives 
$$
P_\eta
=
    \sum_{\mathbf n} c_{\mathbf n}\pi^{-\sum_\al n_\al\lfloor \al(a)\rfloor}
        \prod_{\alpha\in\Phi_B^+}Y_\alpha^{n_\alpha}.
$$
Using the saturation property $\calV_a=V\cap \breve O[\calU_a]$ from \Cref{embedding corollary saturation} in \Cref{embedding corollary}, we see $\pi^hP_\eta\in \calV_a(\eta)$ if and only if $\pi^hP_\eta\in \breve O[\calU_a]$ if and only if 
$$
h-\sum_\al n_\al\lfloor \al(a)\rfloor+\mathrm{val}_F(c_{\mathbf n})\geq 0,
$$
for all $c_{\mathbf n}\neq 0$, where $\mathrm{val}_F(\pi)=1$.
Passing to the maximum proves \Cref{h formula}. 
In particular, $h_a(\eta)\leq \zeta(a)$ follows from the tautological estimate $\lfloor \al(a)\rfloor \leq \al(a)$ using that $\zeta=\chi-\eta=\sum_{\al\in \Phi_B^+} n_\al \al$ for any $c_{\mathbf n}\neq 0$.
This proves \Cref{containment of Peta} in the split case. 

\emph{Proof in the general case.}
We consider the scalar extension map $\iota\colon \mathscr A\to \tilde {\mathscr A}$ from \Cref{scalar extension map}. 
Let $\Phi_{\mathrm{abs}}$ be the (absolute) roots of $G_{\tilde E}$ with respect to $T_{\tilde E}$.
Since $\tilde E/\breve F$ is totally ramified, the comparison of relative and absolute apartments in \cite[Construction~6.1.21]{KP23} gives
\begin{equation}\label{split apartment norm comparison}
 \al\bigl(\iota(b)-\iota(o)\bigr)=(N\al)(b),
\end{equation}
for every $\al$ in the absolute root lattice $\bbZ\Phi_{\mathrm{abs}}$ and every $b\in\mathscr A$.
Further, the identity map on $G_{\tilde E}$ extends to a map of $\tilde O$-group schemes 
\begin{equation}\label{map to split model}
(\calG_b)_{\tilde O}\to \calG_{\iota(b)},
\end{equation}
which is compatible with the big open cell decomposition. 
Indeed, \cite[Construction~6.1.21]{KP23} gives an inclusion $\calG_b(\breve O)\subset \calG_{\iota(b)}(\tilde O)$ compatible with the big open cell decomposition. 
The étoffé property \cite[\S1.7.6]{BT84} shows that the diagonal $G\to \Res_{\tilde E/\breve F}(G_{\tilde E})$ extends to a morphism of parahoric $\breve O$-group schemes $\calG_b\to \Res_{\tilde O/\breve O}(\calG_{\iota(b)})$, hence the desired result by adjunction. 
In particular, \Cref{map to split model} induces a map $(\calU_b)_{\tilde O}\to \calU_{\iota(b)}$.
Since it is the identity on generic fibers and the associated rings of functions are torsion-free, we get an inclusion
\begin{equation}\label{map to split model unipotent}
\tilde O[\calU_{\iota(b)}] \subset \tilde O[\calU_{b}].
\end{equation}
Intersecting with the weight space $\tilde V(\eta)$, we get $\tilde \calV_{\iota(b)}(\eta)\subset \tilde \calV_{b}(\eta)$ by the saturation property from \Cref{embedding corollary saturation} in \Cref{embedding corollary}. 
We apply this to the special vertex $b=o$ and claim that $\tilde \calV_{\iota(o)}(\eta)= \tilde \calV_{o}(\eta)$ for every $\eta\in W_0\cdot \chi$. 
Indeed, equality holds by \Cref{embedding corollary invariants} in \Cref{embedding corollary} for $\eta=\chi$. 
For a general $\eta=w\chi$ with $w\in W_0$, choose a representative $n\in N_S(\breve F)$
of $w$.  As $o$ is a special vertex, $n$ may be chosen in $\calG_o(\breve O)$, and the
inclusion $\calG_o(\breve O)\subset \calG_{\iota(o)}(\tilde O)$ shows that $n$ remains a representative of $w$
in the larger group.  Hence $n$ carries
$\tilde \calV_{o}(\chi)$ onto $\tilde \calV_{o}(\eta)$ and $\tilde \calV_{\iota(o)}(\chi)$
onto $\tilde \calV_{\iota(o)}(\eta)$ simultaneously, and the equality for $\eta=\chi$
propagates to every $\eta\in W_0\cdot\chi$. 
In particular, $P_\eta$ is a generator of the $\tilde O$-module $\tilde \calV_{\iota(o)}(\eta)$. 
Since $\iota(o)$ is special, the computation in the split case above applies to show
$$
\pi_{\tilde E}^{\lfloor \zeta(\iota(a)-\iota(o)) \rfloor}\cdot P_\eta \in \tilde \calV_{\iota(a)}(\eta) \subset \tilde \calV_a(\eta).
$$
Since $\zeta(\iota(a)-\iota(o))=(N\zeta)(a)$ by \Cref{split apartment norm comparison}, this completes the proof of \Cref{inclusion equation integral lattices extremal weights}. 
 
\emph{Equality when $a-o$ lifts to $X_*(T)_I$.} 
Finally, suppose that $a-o$ lifts to $\lambda_a\in X_*(T)_I$. 
Then, $N(\chi-\eta)(a)\in\bbZ$, and conjugation by $\lambda_a(\pi)$ gives an isomorphism $\calU_a\xrightarrow{\sim}\calU_o$ by \cite[Remark~8.3.8]{KP23}. Over $\tilde E$, the element $\lambda_a(\pi)$ equals $(N\lambda_a)(\pi_{\tilde E})$ up to an element of $\calT(\tilde O)$, see \cite[Equation~(7.3.2)]{Kot97}. 
Hence, on the subspace in $\tilde E[U]$ of weight $-\zeta=\eta-\chi$ its pullback is multiplication, up to a unit, by
\[
 \pi_{\tilde E}^{\langle\chi-\eta,N\lambda_a\rangle}
 =
 \pi_{\tilde E}^{N(\chi-\eta)(a)}.
\]
Intersecting with $\tilde V(\eta)$ and using the saturation property from \Cref{embedding corollary saturation} in \Cref{embedding corollary} finishes the proof of \Cref{inclusion lemma key}.
\end{proof}

For our purposes, it is convenient to use the sum
$$
\epsilon_a(\eta) := (N\eta)(a)+h_{a}(\eta) \in \mathbb Q
$$
which is not integral in general.
Here, $(N\eta)(a)=\langle N\eta,a-o\rangle$ under our identification $\mathscr A=X_*(T)_{I,\mathbb R}$.
The following lemma finishes the proof of \Cref{estimate prop}. 

\begin{lemma}\label{key lemma}
For any $x\in W$ and $\eta\in W_{0}\cdot \chi$, the following properties hold:
\begin{enumerate}
\item \label{key lemma 1} $d_a(x,\eta) = \langle N\eta,\delta_a(x)\rangle + \epsilon_a(\eta)-\epsilon_a(\bar x\eta)$, and
\item \label{key lemma 2} $\epsilon_a({\chi})\geq \epsilon_a(\eta)$.
\end{enumerate}
\end{lemma}
\begin{proof}
\Cref{key lemma 1}: 
Write $x=\nu_I(\pi)\cdot u$ according to the decomposition $W=X_*(T)_I\rtimes W_0$ using the choice of the special vertex $o$.
As in the proof of \Cref{weight module properties}, we compute
$$
x\cdot \tilde\calV_a(\eta) 
=
x\cdot \pi_{\tilde E}^{h_a(\eta)}\cdot \tilde \calV_o(\eta)
\overset{\eqref{weight module properties computation}}{=}
\pi_{\tilde E}^{h_a(\eta) + \langle N(\bar x\eta), \nu_I\rangle}\cdot\tilde \calV_o(\bar x\eta)
=
\pi_{\tilde E}^{h_a(\eta) + \langle N(\bar x\eta), \nu_I\rangle - h_a(\bar x\eta)}\cdot\tilde \calV_a(\bar x\eta)
$$
using that multiplication by $\pi_{\tilde E}$ is a central endomorphism of $\tilde V$. 
Taking highest exterior powers, it follows that
\begin{equation}\label{formula for elementary div special}
d_a(x,\eta) = \langle N(\bar x\eta),\nu_I\rangle - h_a(\bar x\eta)+h_a(\eta).
\end{equation}
It remains to prove that
$$
\langle N(\bar x\eta),\nu_I\rangle - h_a(\bar x\eta)+h_a(\eta) = \langle N\eta,\delta_a(x)\rangle + \epsilon_a(\eta)-\epsilon_a(\bar x\eta).
$$
Since $\delta_o(x)=\bar x^{-1}\nu_I$ by \Cref{permissible remark}, \Cref{permissible remark 2}, this follows from the $W_0$-invariance of the pairing using $(N\eta)(a)=\langle N\eta,a-o\rangle$. 
This proves \Cref{key lemma 1}.

\Cref{key lemma 2}:
By \Cref{h computation 1,h computation 3} in \Cref{h computation},
$$
\epsilon_a(\eta)=(N\eta)(a)+h_a(\eta)
\leq (N\eta)(a)+N(\chi-\eta)(a)
=(N\chi)(a)=\epsilon_a(\chi).
$$
This proves the second assertion.
\end{proof}

\begin{corollary}\label{key corollary permissible cc}
Assume $\chi\in X^*(G)$ is central and write $x=\nu_I(\pi)\cdot u$ using the choice of $o$.
If $\pi_E^d\calV_a\subset x\calV_a$ for $d\in\mathbb Z$, then
\[
\langle N\chi,\bar x^{-1}\nu_I\rangle\leq[\tilde E:\breve E]d.
\]
\end{corollary}
\begin{proof}
Taking highest exterior powers as in \Cref{key corollary permissible} gives $d_a(x,\chi)\leq[\tilde E:\breve E]d$.
Since $\bar x\chi=\chi$, \eqref{formula for elementary div special} gives
\[
d_a(x,\chi)=\langle N(\bar x\chi),\nu_I\rangle-h_a(\bar x\chi)+h_a(\chi)
=\langle N\chi,\bar x^{-1}\nu_I\rangle.
\]
This proves the corollary.
\end{proof}

\subsection{Proof of \Cref{inclusion permissible locus theorem}}\label{proof of permissible theorem sec}
The inclusion $A_{\calG,\mu}\subset A_{\calG,\mu}^J$ and the final assertion follow from \Cref{left to prove}. 
In order to prove the inclusion $A_{\calG,\mu}^J\subset P_{\calG,\mu}$, let $[x]\in W_{\mathbf f}\backslash W/W_{\mathbf f}$ with $\Fl_{\calG,\bar k, [x]}\subset A_{\calG,\mu,\bar k}^J$, and choose a representative $x\in W$.
The defining lattice condition in \Cref{special fiber condition} for $\calV\in \Rep(\calG_{O_E})$ is given by the inclusion
\begin{equation}\label{lattice condition permissible}
(\pi_E\otimes 1)^{\vert\underline{\mu}_{\max}(\calV)\vert} \cdot (\calV\otimes_{O_{E^{\mathrm{ur}}}}\breve O)\,\subset\, x\cdot  (\calV\otimes_{O_{E^{\mathrm{ur}}}}\breve O)
\end{equation}
of $O_E\otimes_{O_{E^{\mathrm{ur}}}}\breve O$-lattices in $\calV\otimes_{O_{E^{\mathrm{ur}}}}\breve F$. 
We will show in \Cref{proof of permissible theorem sec ch} and \Cref{proof of permissible theorem sec cc} that the lattice inclusion in \Cref{lattice condition permissible} tested for the representations $\calV=\calV_{\chi, B}$ with $(\chi,B)$ in the set $J=J_{\min}(\calG,\mu)$ from \Cref{definition min test set} implies permissibility of $x$. 
Recall that $J_1\subset X^*(T)^+$ satisfies the conditions in \eqref{J1 conditions} and that $J_{a,\chi}\subset \mathrm{Bor}(G,S)$ are representatives of the double classes $W_{\calG_a}\backslash W_0/W_\chi$.  

\subsubsection{Proof of $\delta_a(x)\in \Conv(W_0\mu_I)$ for all $a\in \mathrm{Vert}(\mathbf f)$.} \label{proof of permissible theorem sec ch}
Let $a\in \mathrm{Vert}(\mathbf f)$ be a vertex, and choose an alcove $\mathbf a$ containing $\mathbf f$ in its closure.
For any $\chi\in J_1^0$ and any $B\in J_{a,\chi}$, the $\calG_{O_E}$-representation 
$$
\calV:=\calV_{\chi,B}=\Res_{{\tilde O}/O_E}\left(\Ind_{{\calB}_{\tilde O}^{-}}^{\calG_{\tilde O}}({^{B}\!\chi})\right)
$$
from \Cref{parahoric test representations} belongs to the test set when $\chi$ is noncentral.
For central $\chi$, its lattice condition follows from the retained $B_0$-test, since $\calV_{\chi,B}=\calV_{\chi,B_0}$ by \Cref{parabolic induction corollary}; this representation extends to every $\calG_a$.
Since $O_{\breve E}=O_{E}\otimes_{O_{E^{\mathrm{ur}}}}\breve O$, the base change of $\calG_{O_{E^{\mathrm{ur}}}}\to \Res_{O_E/O_{E^{\mathrm{ur}}}}\GL(\calV)$ to $\breve O$ has the form
$$
\calG_{\breve O} \to \Res_{O_{\breve E}/\breve O}\GL(\breve\calV),
$$
where $\breve \calV:=\calV\otimes_{O_E}{O_{\breve E}}$.
This representation naturally extends to $\calG_a\to \Res_{O_{\breve E}/\breve O}\GL(\breve \calV)$ by \Cref{maximal parahoric rep cor}, where $\calG_a$ denotes the parahoric $\breve O$-model of $G_{\breve F}$ associated with $a\in \mathrm{Vert}(\mathbf f)$.
If $B,B'\in J_a$ lie in the same double coset in $W_{\calG_a}\backslash W_0/W_\chi$, then one has a $\mathcal G_a$-equivariant isomorphism $V_{\chi,B}\simeq V_{\chi,B'}$. 
Hence, the lattice conditions imposed for the representatives $B'\in J_{a,\chi}$ hold for every $B\in J_a$.
In order to show $\delta_a(x)\in \Conv(W_0\cdot \mu_I)$, we may replace $x$ by any element in the right coset $x W_a$ (\Cref{permissible remark}, \Cref{permissible remark 3}), and assume
$$
{\delta_a(x)}\in \mathrm{Tan}_a(\mathbf a)-a=\bigcup_{B\in J_a}\left({\mathscr C_{B}(a)}-a\right).
$$
Indeed, $a-x(a)$ stays invariant by this replacement and $\delta_a(x)=\bar x^{-1}\!\left(a-x(a)\right)$ can be conjugated by $W_a$ into the tangent cone. 
So, choose $B\in J_a$ such that ${\delta_a(x)}\in {\mathscr C_{B}(a)}-a$, i.e., $\delta_a(x)$ is $B$-dominant.
We consider $\calV=\calV_{\chi,B}$ and $\breve \calV=\calV\otimes_{O_E}{O_{\breve E}}$ for varying $\chi\in J_1^0$ and our fixed choice of $B\in J_a$. 

We now rewrite \Cref{lattice condition permissible} in the form
\begin{equation}\label{lattice condition permissible rewritten}
\pi_E^{\vert\underline{\mu}_{\max}(\breve\calV)\vert} \cdot \breve \calV\,\subset\, x\cdot  \breve \calV,
\end{equation}
where we note that $\vert\underline{\mu}_{\max}(\breve\calV)\vert=\vert\underline{\mu}_{\max}(\calV)\vert$.
Let $\tilde E_1=\tilde E\breve F\subset C$ be the distinguished factor in the decomposition of $\tilde E\otimes_E\breve E$, and let $\tilde O_1$ be its ring of integers.
Projection to this factor gives
\begin{equation}\label{lattice condition permissible rewritten i}
\pi_E^d\breve\calV_1\subset x\breve\calV_1,
\qquad
\breve\calV_1=\Res_{\tilde O_1/O_{\breve E}}
\Ind_{{\calB}_{\tilde O_1}^{-}}^{\calG_{\tilde O_1}}({}^B\!\chi),
\end{equation}
where $d=|\underline\mu_{\max}(\calV_{\chi,B})|$ is the original exponent.
We compute it before applying \Cref{key corollary permissible}.
Put $H=\Gal_E$ and $I_E=I\cap H$.
For $\sigma$ representing a class in $I/I_E=\Hom_{E^{\mathrm{ur}}}(E,C)$, the $B_{0,C}$-highest weights of $\calV_{\chi,B}\otimes_{O_E,\sigma}C$ are $\sigma H*\chi$.
Since $\chi\in J_1^0$, their maximal $\mu$-weight is $\langle\sigma*\chi,\dot\mu\rangle$ by \eqref{J1 compatibility}.
Taking $N$ to be the norm for $\tilde E_1/\breve F$, we obtain
\begin{equation}\label{sharp maximal weight sum}
[\tilde E_1:\breve E]d
=[\tilde E_1:\breve E]\sum_{\sigma\in I/I_E}
\langle\sigma*\chi,\dot\mu\rangle
=\langle N\chi,\mu_I\rangle
=\langle N({}^B\!\chi),{}^B\!\mu_I\rangle.
\end{equation}
Here \eqref{J1 compatibility} makes each pairing constant on the coset $\sigma I_E$, and $W_0$ commutes with inertia.
Consequently, \Cref{lattice condition permissible rewritten i,sharp maximal weight sum,key corollary permissible} imply
$\langle N({}^B\!\chi),{}^B\!\mu_I-\delta_a(x)\rangle\geq0$.
Varying $\chi$ through the set $J_1^0$ therefore gives
\[
\langle N({^B\!\chi}),{^B\!\mu}_I-\delta_a(x)\rangle \geq 0\quad \text{for all } \chi\in J_1^0.
\]
Our choice of $J_1$ then implies
\[
\langle \tilde\chi,{^B\!\mu_I}-\delta_a(x)\rangle\geq 0 \quad\text{for all ${B}$-dominant }\tilde\chi\in X^*(T)^I_{\mathbb Q},
\]
see \Cref{generating lemma J1}, \Cref{generating lemma J1 2}.
Since $\delta_a(x)$ is ${B}$-dominant, \Cref{root lemma twisted} shows $\delta_a(x)\in \operatorname{Conv}(W_0\cdot{\mu_I})$.

\subsubsection{Proof of $[x]\equiv \mu_I$ in $\pi_1(G)_I$.}\label{proof of permissible theorem sec cc}
In addition, choose a special vertex $o\in \mathscr A$.
Write $x=\nu_I(\pi) u$ according to the decomposition $W=X_*(T)_I\rtimes W_0$ induced by $o$.
According to \Cref{key corollary permissible cc} and \eqref{sharp maximal weight sum}, \Cref{lattice condition permissible rewritten i} then implies $\langle N({^B\!\chi}),{^B\!\mu}_I-\bar x^{-1}\nu_I\rangle\geq 0$ when $\chi$ is central. 
Varying $\chi$ through the set $J_1^0\cap X^*(G)$ and noting that their norms generate $X^*(G)^I_{\mathbb Q}$ as a $\mathbb Q_{\geq 0}$-monoid, we see as above that
\[
\langle \tilde\chi,{^B\!\mu_I}-\bar x^{-1}\nu_I\rangle\geq 0 \quad\text{for all }\tilde\chi\in X^*(G)^I_{\mathbb Q}.
\]
Since $\pi_1(G)_I$ is torsion-free, \Cref{central root lemma} implies $\mu_I\equiv {^B\!\mu_I}\equiv \bar x^{-1}\nu_I\equiv  \nu_I\equiv [x]$ in $\pi_1(G)_I$. 
This finishes the proof of \Cref{inclusion permissible locus theorem}.

\section{Admissibility via quantum Bruhat weights}
\label{sec:qbg-elementary-divisors}
In this section, we compute the elementary divisors from \Cref{weight modules} in terms of Schremmer's semi-affine quantum Bruhat weights \cite{Schremmer24}.
\Cref{prop:extremal-minor-reduction} provides a vertexwise formula for every dominant character; its proof is given in \Cref{sec:proof-extremal-minor-reduction}.
In \Cref{moduli proof sec}, we explain how to deduce \Cref{cor:qbg-iwahori-moduli} stated in the introduction.

Throughout this section, we assume that $G_\der$ is simply connected and that $\pi_1(G)_I$ is torsion-free.
Choose a finite Galois extension $\tilde E/\breve F$ splitting $G$, as in \Cref{weight modules}, and set $\tilde O:=O_{\tilde E}$.
For a dominant character $\chi\in X^*(T)^+$, we use the notation $N\chi:=\sum_{\gamma\in\Gal(\tilde E/\breve F)}{}^\gamma\chi$
for its norm.
Let $o\in \scrA=\scrA(G,S,\breve F)$ be a special vertex whose image under the scalar extension map \eqref{scalar extension map} is hyperspecial; see \Cref{weight modules}. 
We identify
\begin{equation}\label{special vertex in 8}
\scrA=X_*(T)_{I,\mathbb R},
\qquad
W=X_*(T)_I\rtimes W_0,
\end{equation}
using throughout the Bruhat--Tits sign convention from \Cref{BT translation convention}. Let moreover $\mathbf a\subset \scrA$ be the unique alcove contained in the dominant cone and adjacent to $o$.
If $G_{\breve F}$ is split, we may take $\tilde E=\breve F$, in which case $N\chi=\chi$ and the notation simplifies considerably. The reader may find it helpful to keep this $\breve F$-split case in mind on a first reading.

Let $\Phi$ denote the relative root system of $G_{\breve F}$ with respect to $S_{\breve F}$, with positive roots $\Phi^+$ determined by the fixed Borel subgroup $B_0\subset G_{\breve F}$; see \Cref{auxiliary.groups2}.
Our constructions are based on the associated \'echelonnage root system $\Sigma$; see \cite[\S1.4]{BT72} and \cite[\S3]{Hai18}. 
We write $\Sigma^+$ for its set of positive roots.
Recall that the roots of $\Sigma$ lie in $X^*(T)^I_{\mathbb R}$, while its coroots lie in $X_*(T)_I$, and that the affine root hyperplanes in $\scrA$ are the zero loci of the affine functions $\beta+n$, with $\beta\in\Sigma$ and $n\in\bbZ$.
Thus $\Sigma$ is reduced in all cases, in contrast to $\Phi$, which need not be. 
The Weyl group $W_0$ of $\Phi$ identifies canonically with the one of $\Sigma$. 
If $G_{\breve F}$ is split, then $\Sigma=\Phi$.

\subsection{Quantum Bruhat weights}
First, assume that the \'echelonnage root system $\Sigma$ is irreducible.
For $v,w\in W_0$, the \textit{quantum Bruhat weight}
\begin{equation}\label{quantum weight}
\operatorname{wt}(v\Rightarrow w)\in\mathbb Z\Sigma^\vee
\end{equation}
is the weight of a shortest directed path from $v$ to $w$ in the quantum Bruhat graph of $\Sigma$ \cite[Definition~3.2, Proposition~3.4]{Schremmer24}.
For a vertex $a\in \mathrm{Vert}(\mathbf a)$, we denote by $K_a=\mathrm{Vert}(\mathbf a)\backslash \{a\}$ the complement.
If $w_0\in W_0$ is the longest element, then the action of $-w_0$ on $\scrA$ stabilizes $\mathbf a$, and we denote by $(-)^\star$ the resulting involution on the vertices $\mathrm{Vert}(\mathbf a)$.
We define the \textit{vertexwise quantum Bruhat weight at $a$} by
\begin{equation}
\label{lrvw-semi-affine-weight}
\operatorname{wt}_{a}(v\Rightarrow w)
:=
\operatorname{wt}_{K_{a^\star}}(w_0w\Rightarrow w_0v)\in\mathbb Z\Sigma^\vee,
\end{equation}
where the right-hand side is Schremmer's semi-affine weight for the subset $K_{a^\star}$, see \cite[Definition~3.26]{Schremmer24}.
Let $\mathbf f\subset \scrA$ be a facet adjacent to $\mathbf a$.
We define the \textit{quantum Bruhat weight at $\mathbf f$} by
\begin{equation}\label{supremums formula}
 \operatorname{wt}_{\mathbf f}(v\Rightarrow w):=\sup_{a\in \mathrm{Vert}(\mathbf f)}\mathrm{wt}_a(v\Rightarrow w)\in\mathbb Z\Sigma^\vee,
\end{equation}
where the supremum is taken in the coroot order, i.e., the coefficient-wise maximum in the expansion in simple coroots of $\Sigma$.
When $\mathbf f=\mathbf a$, we usually drop the subscript, as the formula recovers the usual quantum Bruhat weight by \cite[Lemma~3.8, Lemma~4.34]{Schremmer24}.

Next, if $\Sigma=\coprod_i\Sigma_i$ is reducible, then the standard alcove is the product of the factor alcoves.
Similarly, any adjacent facet $\mathbf f$ decomposes as a product of factor facets $(\mathbf{f}_i)_i$.
For a vertex $a=(a_i)_i$ and $v=(v_i)_i$, $w=(w_i)_i$, we set
\begin{equation}\label{componentwise weight function}
\operatorname{wt}_a(v\Rightarrow w)
:= \sum_i\operatorname{wt}_{a_i}(v_i\Rightarrow w_i),
\end{equation}
and define $\mathrm{wt}_{\mathbf f}(v\Rightarrow w)$ by \Cref{supremums formula} as before; one has $\mathrm{wt}_{\mathbf f}(v\Rightarrow w)=\sum_i\operatorname{wt}_{\mathbf{f}_i}(v_i\Rightarrow w_i)$, as the vertices of $\mathbf f$ are the tuples of vertices of the $\mathbf f_i$ and the coroot order is the product order.
The following result, essentially due to Schremmer, gives an admissibility criterion in terms of the quantum Bruhat weight at $\mathbf f$:

\begin{proposition}\label{schremmer-admissibility}
Let $\la_I\in X_*(T)_I^+$ and $x=\nu_I(\pi)\cdot u\in W$.
Then, the image of $x$ in $W_{\mathbf f}\backslash W/W_{\mathbf f}$ lies in the $\lambda_I$-admissible set if and only if
$$
v^{-1}u^{-1}\nu_I+\operatorname{wt}_{\mathbf f}(v\Rightarrow uv)\leq\la_I,
$$
for every $v\in W_0$.
\end{proposition}
\begin{proof}
First, let $\mathbf f=\{a\}$ be a vertex.
Then the proposition is \cite[Proposition~4.35]{Schremmer24} for the subset $K_{a^\star}$, applied to the root system $\Sigma$ and the lattice $X_*(T)_I$, see also \Cref{rem:schremmer-comparison}.

For a general facet $\mathbf f$, the image of $x$ in $W_{\mathbf f}\backslash W/W_{\mathbf f}$ is $\la_I$-admissible if and only if, for every $a\in\mathrm{Vert}(\mathbf f)$, its image in $W_{\calG_a}\backslash W/W_{\calG_a}$ is $\la_I$-admissible, by the vertexwise criterion \cite[Theorem~1.5]{HH17}.
By the vertex case, this holds if and only if $\eta_{a,v}:=v^{-1}u^{-1}\nu_I+\operatorname{wt}_{a}(v\Rightarrow uv)\leq\la_I$ for all $a\in\mathrm{Vert}(\mathbf f)$ and $v\in W_0$.
For fixed $v$, the elements $\eta_{a,v}$, $a\in\mathrm{Vert}(\mathbf f)$, differ only by elements of $\mathbb Z\Sigma^\vee$, so the inequalities $\eta_{a,v}\leq\la_I$ for all $a$ are equivalent to the coefficient-wise inequality $\sup_a\eta_{a,v}\leq\la_I$, that is, to $v^{-1}u^{-1}\nu_I+\operatorname{wt}_{\mathbf f}(v\Rightarrow uv)\leq\la_I$ by \Cref{supremums formula}.
\end{proof}

\begin{remark}\label{rem:schremmer-comparison}
The normalization in \Cref{lrvw-semi-affine-weight} is due to different conventions for the action of $T(\breve F)$ on $\scrA$ and different choices of base alcove $\mathbf a$:
we use the minus sign convention and the dominant base alcove, whereas \cite{Schremmer24} uses the plus sign convention and the dominant base alcove, see \Cref{plus translation convention}.
Conjugating the affine action by $a\mapsto-w_0a$ sends $x=\nu_I(\pi)u$ to $(w_0uw_0)\varepsilon^{w_0u^{-1}\nu_I}$ in Schremmer's notation and sends $\mathbf f$ to $\mathbf f^\star$; reindexing his criterion by $w_0v$ gives \Cref{schremmer-admissibility}.
\end{remark}

\subsection{Relation with extremal weight modules}
For a dominant character $\chi\in X^*(T)^+$, we denote the dual Weyl module over $\tilde E$ of highest weight $\chi$ by
$$\tilde V_\chi=\Ind_{B_{0,\tilde E}^-}^{G_{\tilde E}}(\chi).$$ 
Its weight subspace $\tilde V_\chi(\eta)$, for any extremal weight $\eta\in W_0\cdot \chi$, is $1$-dimensional.
For $v\in W_0$, let $B_v:=vB_0v^{-1}\in \mathrm{Bor}(G,S)$ be the associated Borel subgroup.
Recall from \Cref{reps for proof} the $\calG_{a,\tilde O}$-stable lattice inside $\tilde V_\chi$ attached to a vertex $a\in \mathscr A$:
\begin{equation}\label{def-integral-lattice-quantum-section}
\tilde\calV_{a,\chi,v} := \Ind_{\calB_{v,\tilde O}^-}^{\calG_{a,\tilde O}}(v\chi),
\end{equation}
where $\calB_v^-$ is the scheme-theoretic closure in $\calG_a$ of the Borel subgroup in $G_{\breve F}$ opposite to $B_v$.
Any extremal weight $\eta\in W_0\cdot \chi$ defines a rank one $\tilde O$-lattice $\tilde\calV_{a,\chi,v}(\eta)$ in $\tilde V_\chi(\eta)$, and we consider
\begin{equation}\label{def-extremal-weight-lattice-quantum-section}
\tilde\calV_{a,\chi,v}(\eta)=\pi_{\tilde E}^{h_{a,\chi, v}(\eta)} \cdot\tilde\calV_{o,\chi,v}(\eta)
\end{equation}
comparing the weight lattice at $a$ with the one at $o$.
By \Cref{taking highest exterior powers}, $h_{a,\chi,v}(\eta)$ is the integer $h_a(\eta)$ of \Cref{h computation}, formed for the Borel subgroup $B_v$.

\begin{theorem}
\label{prop:extremal-minor-reduction}
Assume that $G$ is essentially tamely ramified, and let $\chi\in X^*(T)^+$.
Then, for every $a\in\operatorname{Vert}({\mathbf a})$ and all $v,w\in W_0$, one has
\begin{equation}
\label{vertexwise-extremal-minor-formula}
-h_{a,\chi, v}(w\chi)
=
\left\langle N\chi,
\operatorname{wt}_a(v\Rightarrow w)
\right\rangle.
\end{equation}
\end{theorem}

The proof is given in \Cref{sec:proof-extremal-minor-reduction}.

\subsection{Proof of moduli descriptions}\label{moduli proof sec}
We now explain how \Cref{prop:extremal-minor-reduction} implies \Cref{main theorem intro} from the introduction; see \Cref{cor:qbg-iwahori-moduli}.
For $a\in \scrA$, $\chi\in X^*(T)^+$, and $v\in W_0$, we have the associated $\calG_a$-representation $\calV_{a,\chi, v}=\Res_{\tilde O/O_{\breve E}}\tilde\calV_{a,\chi,v}$ over $O_{\breve E}$.
For every $x\in W$ and $\eta\in W_0\cdot \chi$, we consider the elementary divisor $d_{a,\chi,v}(x,\eta)$ from \Cref{elementary divisor definition}, formed for the Borel subgroup $B_v$.
The formula from \Cref{formula for elementary div special}, together with $h_{a,\chi,v}(v\chi)=0$ (\Cref{h computation}), takes the following form:
\begin{equation}
\label{elementary-divisor-h-formula}
d_{a,\chi, v}(x,v\chi)
=
\left\langle N\chi,v^{-1}u^{-1}\nu_I\right\rangle
-h_{a,\chi, v}(uv\chi),
\end{equation}
for $x=\nu_I(\pi) u\in W$, $v\in W_0$, and $a\in \mathrm{Vert}(\mathbf a)$.

\begin{corollary}\label{cor:main-conjecture-admissibility}
Assume that $G$ is essentially tamely ramified.
Let $J_1\subset X^*(T)^+$ be any subset such that $N_{\bar I}(J_1)$ generates the $\mathbb Q_{\geq 0}$-monoid $X^*(T)^{I,+}_{\mathbb Q}$; cf.~\Cref{J_1 condition}.
Let $\mu\in X_*(T)^+$ and $x=\nu_I(\pi)u\in W$.
Then, the image of $x$ in $W_{\mathbf f}\backslash W/W_{\mathbf f}$ is $\mu_I$-admissible if and only if
\begin{equation}\label{cor:main-conjecture-admissibility eq}
d_{a,\chi,v}(x,v\chi)
\leq
\langle N\chi,\mu_I\rangle
\end{equation}
for all $a\in\operatorname{Vert}({\mathbf f})$, all $v\in W_0$, and all $\chi\in J_1$.
\end{corollary}
\begin{proof}
By \Cref{elementary-divisor-h-formula,prop:extremal-minor-reduction}, for fixed $a$ and $v$ one has
\[
d_{a,\chi,v}(x,v\chi)
=\left\langle N\chi,v^{-1}u^{-1}\nu_I+\mathrm{wt}_a(v\Rightarrow uv)\right\rangle.
\]
Thus each inequality in \eqref{cor:main-conjecture-admissibility eq} depends linearly only on $N\chi$.
Since $N(J_1)$ generates the invariant dominant cone, we may replace the tests by central weights and lifts of all fundamental weights before taking the maximum over vertices. 
The norms of the latter are positive multiples of the fundamental weights $\breve\omega_\calO$ of $\Sigma$, see \eqref{exact sequence inertia} and the discussion following it.
Under this assumption, \Cref{elementary-divisor-h-formula} implies
\[
\max_{a\in\operatorname{Vert}({\mathbf f})}
d_{a,\chi,v}(x,v\chi)
=
\left\langle
N\chi,v^{-1}u^{-1}\nu_I
\right\rangle + \max_{a\in\operatorname{Vert}({\mathbf f})}\left(-h_{a,\chi,v}(uv\chi)\right).
\]
\Cref{prop:extremal-minor-reduction} gives $-h_{a,\chi,v}(uv\chi)
=\langle N\chi,\mathrm{wt}_a(v\Rightarrow uv)\rangle$.
Now, if $\chi\in X^*(G)$ is central, then $\langle N\chi, \mathrm{wt}_a(v\Rightarrow uv)\rangle=0$ and the maximum vanishes as well.
If $\chi$ is a lift of a fundamental weight, then $\langle N\chi,-\rangle$ reads off a positive multiple of the coefficient of a single simple coroot of $\Sigma$, so that the maximum equals $\langle N\chi, \mathrm{wt}_{\mathbf f}(v\Rightarrow uv)\rangle$ by \eqref{supremums formula}.
In both cases, we get the formula
$$
\max_{a\in\operatorname{Vert}({\mathbf f})}
d_{a,\chi,v}(x,v\chi)
=
\left\langle
N\chi,v^{-1}u^{-1}\nu_I+\mathrm{wt}_{\mathbf f}(v\Rightarrow uv)\right\rangle.
$$
Since $N(J_1)$ generates $X^*(T)^{I,+}_{\mathbb Q}$, the inequalities in \Cref{cor:main-conjecture-admissibility eq}, for all $a\in \mathrm{Vert}(\mathbf f)$ and $\chi\in J_1$, are equivalent to
\[
v^{-1}u^{-1}\nu_I+
\operatorname{wt}_{\mathbf f}(v\Rightarrow uv)
\leq\mu_I, \quad \text{for all } v\in W_0,
\]
see \Cref{root lemma twisted} and \Cref{central root lemma}, which use that $G_\der$ is simply connected and that $\pi_1(G)_I$ is torsion-free.
So, the corollary follows from \Cref{schremmer-admissibility}.
\end{proof}

The following result proves \Cref{main theorem intro} from the introduction. 

\begin{theorem}\label{cor:qbg-iwahori-moduli}
Assume that $G$ is essentially tamely ramified.
Let $\calG$ be the parahoric $O$-model attached to $\mathbf f$, and let $\mu\in X_*(T)^+$ be arbitrary.
Then, one has $\mathbb M_{\calG,\mu}^{J}=\mathbb M_{\calG,\mu}$ for $J=J_1\times J_2$,
where $J_1\subset X^*(T)^+$ is any finite subset satisfying \Cref{J_1 condition}, and $J_2=J_{\mathbf f}$ the set of $\mathbf f$-relevant Borel subgroups.
\end{theorem}
\begin{proof}
In view of \Cref{cor:main-conjecture-admissibility}, the proof is analogous to that of \Cref{inclusion permissible locus theorem} in \Cref{proof of permissible theorem sec}:
By the generic-fiber argument in \Cref{left to prove}, it is enough to show that the inclusion on special fibers $A_{\calG,\mu}\subset A_{\calG,\mu}^J$ is an equality.
For the reverse inclusion, let $x\in W$ with $\Fl_{\calG,\bar k,x}\subset A_{\calG,\mu,\bar k}^J$.
The defining lattice condition in \Cref{special fiber condition} for $\calV\in \Rep(\calG_{O_E})$ is then given by the inclusion
\begin{equation}\label{lattice condition permissible split}
(\pi_E\otimes 1)^{\vert\underline{\mu}_{\max}(\calV)\vert} \cdot (\calV\otimes_{O_{E^{\mathrm{ur}}}}\breve O)\,\subset\, x\cdot  (\calV\otimes_{O_{E^{\mathrm{ur}}}}\breve O).
\end{equation}
We evaluate the condition on the representations
$$
\calV_{\chi, B}=\Res_{O_{\tilde E}/O_E}\left(\Ind_{\calB^-_{O_{\tilde E}}}^{\calG_{O_{\tilde E}}}({}^B\chi)\right)
$$
for $(\chi, B)\in J$.
As in \Cref{lattice condition permissible rewritten i}, we project to the distinguished factor $\tilde E_1=\tilde E\breve F\subset C$.
For $\chi\in J_1^0$, this gives
\begin{equation}\label{lattice condition permissible split 2}
\pi_E^{d_\chi}\breve\calV_1\subset x\breve\calV_1,
\qquad
\breve\calV_1=\Res_{\tilde O_1/O_{\breve E}}
\Ind_{{\calB}_{\tilde O_1}^{-}}^{\calG_{\tilde O_1}}({}^B\!\chi),
\end{equation}
where $d_\chi=|\underline\mu_{\max}(\calV_{\chi,B})|$ and
$[\tilde E_1:\breve E]d_\chi=\langle N\chi,\mu_I\rangle$ by \Cref{sharp maximal weight sum}.
We now apply the preceding results.

Let $a\in\operatorname{Vert}({\mathbf f})$, $\chi\in J_1^0$, and $v\in W_0$.
Choose $B\in J_a\subset J_{2}$, with $B=wB_0w^{-1}$, such that $[w]=[v]$ in $W_{\calG_a}\backslash W_0/W_{\chi}$.
By \Cref{maximal parahoric rep cor}, the $\calG_{\breve O}$-representation $\breve\calV_1$ extends to $\calG_a$, and we denote the resulting representation by $\breve \calV_{a,\chi,w}$.
\Cref{redundancies between representations} identifies this $\calG_a$-representation with $\breve \calV_{a,\chi,v}=\calV_{a,\chi,v}$, so \Cref{lattice condition permissible split 2} can be rewritten as
$$
\pi_E^{d_\chi} \cdot \calV_{a,\chi,v}\,\subset\, x\cdot  \calV_{a,\chi,v}.
$$
Passing to extremal weight modules and taking highest exterior powers, exactly as in the proof of \Cref{key corollary permissible}, gives
\begin{equation}\label{inclusion extremal weight modules}
\pi_E^{\langle N\chi,\mu_I\rangle} \cdot \wedge^{\mathrm{top}}\calV_{a,\chi,v}\left(\overline{xv}\chi\right)\,\subset\, \wedge^{\mathrm{top}}\left(x\cdot  \calV_{a,\chi,v}(v\chi)\right)=\pi_E^{d_{a,\chi, v}(x,v\chi)}\cdot  \wedge^{\mathrm{top}}\calV_{a,\chi, v}\left(\overline{xv}\chi\right),
\end{equation}
where we used $[\tilde E_1:\breve E]d_\chi=\langle N\chi,\mu_I\rangle$ from \Cref{sharp maximal weight sum}.
This is equivalent to $d_{a,\chi, v}(x,v\chi)\leq \langle N\chi,\mu_I\rangle$.
Since $N(J_1^0)$ generates $X^*(T)^{I,+}_{\mathbb Q}$, \Cref{cor:main-conjecture-admissibility} applied to $J_1^0$ implies that the double $W_{\mathbf f}$-coset of $x$ is $\mu_I$-admissible.
Since $\mu$ is dominant for the $\breve F$-rational Borel subgroup $B_0$, the $\mu_I$-admissible set agrees with the $\mu$-admissible set \cite[Theorem~4.2]{Hai18}.
This finishes the proof.
\end{proof}

\section{The elementary divisor formula}
\label{sec:proof-extremal-minor-reduction}
In this section, we prove \Cref{prop:extremal-minor-reduction},
which gives a formula for the elementary divisor.
The proof uses the graded character formula of Feigin--Makedonskyi
\cite{FeiginMakedonskyi17} and its analogue for twisted current algebras
\cite{FeiginMakedonskyiTwisted17}.

We continue with the notation of \Cref{sec:qbg-elementary-divisors}.
In particular, we assume that $G_\der$ is simply connected, that $\pi_1(G)_I$ is torsion-free, and we fix a special vertex $o\in \scrA=\scrA(G,S,\breve F)$ whose image under the scalar extension map \eqref{scalar extension map} is hyperspecial; see \Cref{special vertex in 8}.

In \Cref{weight linearity sec,reduction to simple simply connected},
we reduce to the case where $G_{\breve F}$ is simply connected,
absolutely almost simple, and tamely ramified, $\chi$ is a fundamental
weight, and $\tilde E$ is the splitting field of $G_{\breve F}$.
The essential tameness assumption is used in
\Cref{corollary reduction to simple simply connected} to ensure that
the resulting absolutely almost simple group is tamely ramified.
From \Cref{coordinates sec} onward, we work in this setting.

\subsection{Weight linearity of elementary divisors}\label{weight linearity sec}
The right-hand side of \Cref{vertexwise-extremal-minor-formula} is $\bbZ_{\geq 0}$-linear in $\chi$.
In this subsection, we check that the left-hand side is so as well; see \Cref{weight module monoidal}.

First, we fix explicit generators of the extremal weight lattices $\tilde\calV_{o,\chi,v}(w\chi)\subset \tilde V_\chi(w\chi)$ at the chosen origin $o\in \scrA$.
Restricting along the open embedding $U_v\to G/B_v^-$ yields an inclusion:
\begin{equation}
\tilde V_\chi \to  \tilde E[U_v],\ f\mapsto f|_{U_v}.
\end{equation}
This integrally extends to an inclusion
\begin{equation}
 \tilde\calV_{o,\chi,v} \to \tilde O[\calU_{o,v}], \ f\mapsto f|_{\calU_{o,v}}
 \end{equation}
with torsion-free cokernel given by restricting functions along the open immersion $\calU_{o,v} \to \calG_{o}/\calB_{o,v}^-$, see \Cref{embedding corollary}.
The highest $T_{\tilde E}$-weight of $\tilde V_\chi$ relative to $B_{v,\tilde E}$ is $v\chi$, and we can identify the image of $\tilde\calV_{o,\chi,v}(v\chi)$ under the above map with the constant functions $\tilde O$ by \Cref{embedding corollary}, \Cref{embedding corollary invariants}.

\begin{definition}\label{P generators}
For every $\chi\in X^*(T)^+$ and $v,w\in W_0$, choose a generator $f_{v,w}^\chi$ of the rank $1$ free $\tilde O$-module $\tilde\calV_{o,\chi,v}(w\chi)$. 
Define $P_{v,w}^\chi:=f_{v,w}^\chi|_{\calU_{o,v}}\in \tilde O[\calU_{o,v}]$.
\end{definition}

The generators $f_{v,w}^\chi$ are unique and conjugate under the $W_0$-action up to $\tilde O$-units, so $P_{v,w}^\chi$ are unique up to $\tilde O$-units as well. 
We suppress the superscript when $\chi$ is fixed.

\begin{proposition}\label{weight module monoidal}
The integer $h_{a,\chi,v}(w\chi)$ is $\bbZ_{\geq 0}$-linear in $\chi\in X^*(T)^+$.
\end{proposition}
\begin{proof}
By \Cref{embedding corollary saturation} in \Cref{embedding corollary}, one has
$\tilde\calV_{b,\chi,v}(w\chi)=\tilde V_\chi(w\chi)\cap\tilde O[\calU_{b,v}]$
for every $b\in\scrA$.
Hence $-h_{a,\chi,v}(w\chi)$ is the Gauss valuation of $P_{v,w}^\chi$ in
$\tilde O[\calU_{a,v}]$, i.e., the largest $d\in\bbZ$ with
$P^\chi_{v,w}\in\pi_{\tilde E}^d\tilde O[\calU_{a,v}]$.
The special fiber of $\calU_{a,v}$ is integral, as $\calU_{a,v}$ is smooth
with connected fibers \cite[\S4.3]{BT84}.
Thus Gauss valuations are multiplicative, and it suffices to show that
$P_{v,w}^{\chi+\chi'}=P_{v,w}^\chi P_{v,w}^{\chi'}$ up to $\tilde O$-units for every $\chi,\chi'\in X^*(T)^+$, $v,w\in W_0$.
For this, consider the $\calG_{o,\tilde O}$-equivariant map
\[
\tilde\calV_{o,\chi,v}\otimes_{\tilde O}\tilde\calV_{o,\chi',v}
\longto \tilde\calV_{o,\chi+\chi',v}, \, f\otimes f'\mapsto f f'.
\]
The map is an isomorphism on highest weight lattices, so the claim holds for $v=w$.
In general, since $o$ is special, we may choose a representative of $wv^{-1}$ in $\calG_o(\breve O)$.
By equivariance, we get isomorphisms on all extremal weight lattices, and the proposition follows.
\end{proof}

\subsection{Reduction to simply connected almost simple groups and fundamental weights}
\label{reduction to simple simply connected}
In this subsection, we reduce \Cref{prop:extremal-minor-reduction} to the case where $G$ is simply connected and absolutely almost simple, $\tilde E$ is the splitting field of $G$, and $\chi$ is a fundamental weight, which we assume from \Cref{coordinates sec} on.

Since $G_\der$ is simply connected, we have an exact sequence
$$
0\longto X^*(G)\longto X^*(T)\overset{\chi\mapsto \bar\chi}{\longto} P\longto 0,
$$
where $P=X^*(T\cap G_\der)$ is the weight lattice of the absolute root system.
The fundamental weights in $P$ form the basis dual to the absolute simple coroots.
We say that a dominant character $\chi\in X^*(T)^+$ is a lift of a fundamental weight if its image $\bar\chi\in P$ is a fundamental weight.

\begin{lemma}\label{lifts of fundamental weights suffice}
If \Cref{vertexwise-extremal-minor-formula} holds for all lifts of fundamental weights, then it holds for all $\chi\in X^*(T)^+$.
\end{lemma}

\begin{proof}
If $\chi\in X^*(G)$ is central, then both sides of \eqref{vertexwise-extremal-minor-formula} vanish: $\tilde V_\chi$ is a line, so $P_{v,w}^\chi\in\tilde O^\times$, and $N\chi$ is orthogonal to $\bbZ\Sigma^\vee$.
In general, write $\bar\chi=\sum_in_i\omega_i$ with $n_i\geq0$, and choose (necessarily dominant) lifts $\chi_i$ of the $\omega_i$.
So $\delta:=\chi-\sum_in_i\chi_i\in X^*(G)$, and both $\pm\delta$ are dominant.
Both sides of \eqref{vertexwise-extremal-minor-formula} are $\bbZ_{\geq0}$-linear on $X^*(T)^+$ by \Cref{weight module monoidal}.
Applying this to $\chi+(-\delta)=\sum_in_i\chi_i$ and using the central case gives the claim.
\end{proof}

Let $G_\der=\prod_iG_i$ be the decomposition into $\breve F$-almost simple factors, so that $\Sigma=\coprod_i\Sigma_i$, and decompose $a=(a_i)_i$, $v=(v_i)_i$ and $w=(w_i)_i$ accordingly.

\begin{lemma}\label{reduction almost simple factors}
Let $\chi\in X^*(T)^+$ lift a fundamental weight $\bar\chi$ of the $\breve F$-almost simple factor $G_{i_0}$, and let $\bar\chi_{i_0}$ be the corresponding character of the maximal torus of $G_{i_0}$.
Then, one has
\[
h^G_{a,\chi,v}(w\chi)=h^{G_{i_0}}_{a_{i_0},\bar\chi_{i_0},v_{i_0}}(w_{i_0}\bar\chi_{i_0}),
\qquad
\left\langle N\chi,\operatorname{wt}_a(v\Rightarrow w)\right\rangle=\left\langle N\bar\chi_{i_0},\operatorname{wt}_{a_{i_0}}(v_{i_0}\Rightarrow w_{i_0})\right\rangle.
\]
\end{lemma}
\begin{proof}
The second identity follows from \eqref{componentwise weight function}, as $N(\chi-\bar\chi)$ is orthogonal to $\bbZ\Sigma^\vee$ and $N\bar\chi$ vanishes on $\Sigma_i^\vee$ for $i\neq i_0$.
For the first one, consider the central isogeny
\begin{equation}\label{central isogeny}
\varphi\colon \tilde G:=Z_G^\circ\times G_{\der}\longrightarrow G,
\end{equation}
where $Z_G^\circ$ is the connected center of $G$.
It identifies the root subgroups of source and target together with their Bruhat--Tits models $\tilde\calU_{b,v}\simeq\calU_{b,v}$ at every $b\in \scrA$ \cite[Lemma~8.3.7]{KP23}.
Moreover $\varphi^*\tilde V_\chi=\delta\otimes V_{\bar\chi}$ for some $\delta\in X^*(Z_G^\circ)$, so \Cref{embedding corollary}, \Cref{embedding corollary saturation} yields weight-compatible isomorphisms
$$
\delta\otimes \tilde\calV^{\tilde G}_{b,\bar\chi,v}=\varphi^*\tilde V_\chi\cap\tilde O[\tilde \calU_{b,v}]\simeq\varphi^*\big(\tilde V_\chi\cap\tilde O[\calU_{b,v}]\big)=\varphi^*\tilde\calV_{b,\chi,v},
\qquad b\in\{a,o\}.
$$
Hence $h$ is unchanged by passing to $\tilde G$, and then to $G_\der=\prod_iG_i$.
There, the parahoric models, the base alcove and the root group models are products, and $V_{\bar\chi}$ is the external tensor product of $V_{\bar\chi_{i_0}}$ with trivial representations.
\end{proof}

\begin{lemma}\label{reduction restriction of scalars}
Let $G=\Res_{K/\breve F}H$ for a finite separable extension $K/\breve F$ and an absolutely almost simple, simply connected $K$-group $H$.
Fix an embedding $\iota_0\colon K\to\tilde E$, and let $\chi$ be a fundamental weight of $H\otimes_{K,\iota_0}\tilde E$.
View $\chi$ as a weight of $G_{\tilde E}$ via the projection
\[
G_{\tilde E}=\prod_{\iota\in \Hom_{\breve F}(K,\tilde E)}H\otimes_{K,\iota}\tilde E
\longrightarrow H\otimes_{K,\iota_0}\tilde E.
\]
Then the respective left-hand sides and the respective right-hand sides of \Cref{vertexwise-extremal-minor-formula} agree for both $(G,\breve F,\chi)$ and $(H,K,\chi)$, computed with the same $\tilde E$.
\end{lemma}
\begin{proof}
For the left-hand side, the parahoric and root group models of $G$ are $\calG_b=\Res_{O_K/\breve O}\calH_b$ and $\calU_{b,v}=\Res_{O_K/\breve O}\calU_{H,b,v}$, for $b$ in the common apartment.
The multiplication map $O_K\otimes_{\breve O}\tilde O\to\tilde O$, $x\otimes y\mapsto\iota_0(x)y$, induces a projection
\[
p_{\iota_0}\colon\calU_{b,v}\otimes_{\breve O}\tilde O
\longrightarrow\calU_{H,b,v}\otimes_{O_K,\iota_0}\tilde O,
\]
which is a morphism of smooth $\tilde O$-group schemes.
It is surjective on $\bar k$-points of the special fibers because $\calU_{H,b,v}$ is smooth and $O_K\otimes_{\breve O}\bar k\to\bar k$ has nilpotent kernel.
Since the special fiber of the target is reduced, $p_{\iota_0}^*$ preserves Gauss valuations.
Indeed, if $p_{\iota_0}^*f\in\pi_{\tilde E}\tilde O[\calU_{b,v}]$, then $f$ vanishes on the special fiber of the target, so $f\in\pi_{\tilde E}\tilde O[\calU_{H,b,v}]$.
Since $\tilde V_\chi$ is inflated along $p_{\iota_0}$ on the generic fibers, we have $P^{G,\chi}_{v,w}=p_{\iota_0}^*P^{H,\chi}_{v,w}$ up to units.
The proof of \Cref{weight module monoidal} therefore gives $h^G_{a,\chi,v}(w\chi)=h^H_{a,\chi,v}(w\chi)$.

For the right-hand side, put $\Gamma_K:=\Gal(\tilde E/\iota_0K)$.
Since $T=\Res_{K/\breve F}T_H$, we have $X_*(T)=\Ind_{\Gamma_K}^{\Gamma}X_*(T_H)$, where $\Gamma$ acts through the image of $I$.
Summing the $H$-coordinates defines a map $p\colon X_*(T)\to X_*(T_H)$ inducing the Shapiro isomorphism
\[
p_*\colon X_*(T)_I\xrightarrow{\ \sim\ }X_*(T_H)_{I_K}.
\]
Similarly, the projection $\mathrm{pr}_0$ onto the $\iota_0$-component induces an isomorphism $X_*(T)^I\simeq X_*(T_H)^{I_K}$.
These isomorphisms are compatible with the norm maps: $\mathrm{pr}_0(N^G_{\tilde E/\breve F}\nu)=N^H_{\tilde E/K}(p\nu)$ for $\nu\in X_*(T)$.
Moreover, $p_*$ identifies the coroots of the two \'echelonnage root systems, viewed as images of absolute coroots, as well as their Weyl groups and the vertices of their base alcoves.
Hence $p_*\operatorname{wt}^G_a(v\Rightarrow w)=\operatorname{wt}^H_a(v\Rightarrow w)$.
Since $\chi$ is inflated from the $\iota_0$-component, we have $\langle\chi,x\rangle=\langle\chi,\mathrm{pr}_0x\rangle$ for $x\in X_*(T)$.
The compatibility with the norm maps therefore gives
\[
\begin{aligned}
\langle N^G_{\tilde E/\breve F}\chi,
        \operatorname{wt}^G_a(v\Rightarrow w)\rangle
&=\langle\chi,
         N^G_{\tilde E/\breve F}\operatorname{wt}^G_a(v\Rightarrow w)\rangle\\
&=\langle\chi,
         N^H_{\tilde E/K}\operatorname{wt}^H_a(v\Rightarrow w)\rangle\\
&=\langle N^H_{\tilde E/K}\chi,
         \operatorname{wt}^H_a(v\Rightarrow w)\rangle.
\end{aligned}
\]
\end{proof}

\begin{lemma}\label{reduction splitting field}
Let $\tilde E'/\tilde E$ be a finite extension such that $\tilde E'/\breve F$ is Galois.
Then replacing $\tilde E$ by $\tilde E'$ multiplies both sides of \eqref{vertexwise-extremal-minor-formula} by $[\tilde E':\tilde E]$.
\end{lemma}
\begin{proof}
Put $\tilde O'=O_{\tilde E'}$, and use primes for the lattices, elementary divisors and norm computed with $\tilde E'$.
The parahoric induced modules commute with flat base change:
$\tilde\calV'_{b,\chi,v}=\tilde\calV_{b,\chi,v}\otimes_{\tilde O}\tilde O'$.
Since $\tilde E'/\tilde E$ is totally ramified, we have $\pi_{\tilde E}\in\pi_{\tilde E'}^{[\tilde E':\tilde E]}\tilde O'^\times$, and hence
$$h'_{a,\chi,v}(w\chi)=[\tilde E':\tilde E]\,h_{a,\chi,v}(w\chi).$$
Since $\Gal(\tilde E'/\tilde E)$ acts trivially on $X^*(T)$, we also have $N'\chi=[\tilde E':\tilde E]\,N\chi$.
\end{proof}

\begin{corollary}\label{corollary reduction to simple simply connected}
To prove \Cref{prop:extremal-minor-reduction}, it suffices to treat the case where $G_{\breve F}$ is simply connected, absolutely almost simple and tamely ramified, $\chi$ is a fundamental weight, and $\tilde E$ is the splitting field of $G_{\breve F}$.
In this case, $\Gal(\tilde E/\breve F)$ is cyclic of order $1,2$, or $3$.
\end{corollary}
\begin{proof}
By \Cref{lifts of fundamental weights suffice,reduction almost simple factors,reduction restriction of scalars,reduction splitting field}, it suffices to treat the case where $G_{\breve F}$ is simply connected and absolutely almost simple, $\chi$ is a fundamental weight, and $\tilde E$ is its splitting field.
The assumption that $G$ is essentially tamely ramified ensures that the resulting absolutely almost simple group is tamely ramified.
Moreover, $\Gamma$ identifies with the image of $I$ in the automorphism group of the absolute Dynkin diagram.
Since tame inertia is procyclic, this image is cyclic.
By the classification of connected Dynkin diagrams \cite[Chap.~VI, \S4.2]{Bou68}, their cyclic automorphism groups have order $1$, $2$ or $3$.
\end{proof}

\subsection{Coordinates on integral root groups}\label{coordinates sec}
In view of \Cref{corollary reduction to simple simply connected}, we assume from now on that $G_{\breve F}$ is simply connected, absolutely almost simple and tamely ramified, $\chi$ is a fundamental weight, and $\tilde E$ is the splitting field of $G_{\breve F}$.
So, $\Gamma=\Gal(\tilde E/\breve F)$ is cyclic of order $e\in\{1,2,3\}$, and we choose $\pi_{\tilde E}\in O_{\tilde E}$ with $\pi_{\tilde E}^e=\pi$.

We further fix a Chevalley--Steinberg system $(y_\gamma)$ for $G_{\breve F}$ \cite[4.1.3]{BT84}, compatible with the chosen special vertex $o$.
Recall that $\Phi$ denotes the relative root system of $G_{\breve F}$ with positive roots $\Phi^+$, and that $\Sigma$ is the associated \'echelonnage root system with positive roots $\Sigma^+$.
If $G_{\breve F}$ is split, then $\Sigma=\Phi$.

For $\alpha\in\Phi$, let $m_\alpha$ denote the number of absolute roots restricting to $\alpha$; these form a single $I$-orbit.
If $\alpha$ is neither multipliable nor divisible, then
\[
m_\alpha\alpha\in\Sigma,
\qquad
(m_\alpha\alpha)^\vee=\frac{\alpha^\vee}{m_\alpha},
\]
where the latter is the common image in $X_*(T)_I$ of the corresponding absolute coroots; see \cite[Definition~3.3, Proposition~3.5 and Theorem~6.1]{Hai18}.
If $\alpha$ is multipliable, i.e., if $2\alpha\in\Phi$, then instead
\[
2m_\alpha\alpha\in\Sigma,
\qquad
(2m_\alpha\alpha)^\vee=\frac{\alpha^\vee}{2m_\alpha}.
\]
Again, this coroot is the common image in $X_*(T)_I$ of the absolute coroots of the roots lying above $\alpha$, while the divisible roots yield no additional roots of $\Sigma$; see \cite[Definitions~3.1(3), 3.3 and Theorem~6.1]{Hai18}.

Next, we choose coordinates on the root groups.
For every $\alpha\in\Phi$, we fix a lift $\tilde\alpha\in X^*(T)$, so that $\calO_\alpha:=\Gamma\cdot\tilde\alpha$ is the set of the absolute roots restricting to $\alpha$, of cardinality $m_\alpha$.
Let $L_\alpha:=\tilde E^{\Gamma_{\tilde\alpha}}$ be the fixed field of the stabilizer $\Gamma_{\tilde\alpha}$, a totally ramified extension of $\breve F$ of degree $m_\alpha$.
This degree is equal to either $1$ or $e$, so $L_\alpha$ equals $\breve F$ or $\tilde E$, respectively. 
%It depends on the lift $\tilde\alpha$ up to $\Gamma$-conjugation, and so does the parametrization below.
We also choose the uniformizer $\varpi_\alpha:=\pi_{\tilde E}^{e_\alpha}$ of $L_\alpha$, where $e_\alpha:=e/m_\alpha$ is the ramification index of $\tilde E/L_\alpha$, so that $\varpi_\alpha^{m_\alpha}=\pi$.
To index the coordinates uniformly, put
\[
\bar\Psi:=\{(\alpha,j)\in \Phi\times \Z\mid \alpha\text{ nondivisible},\ 0\leq j<m_\alpha\}
\cup\{(\alpha,1)\mid\alpha\in\Phi\text{ divisible}\},
\]
and let $v\bar\Psi^+\subset\bar\Psi$ be the subset of pairs with $\alpha\in v\Phi^+$.
For $(\alpha,j)\in\bar\Psi$, put $\psi_{\alpha,j}:=\alpha+j/m_\alpha$ if $\alpha$ is nondivisible, and $\psi_{\alpha,1}:=\alpha+\frac12$ and $e_\alpha:=e_{\alpha/2}$ if $\alpha$ is divisible.

\begin{proposition}\label{proposition uniform rescaled coordinates}
For every $v\in W_0$, there are $\breve F$-coordinates $(X_{\alpha,j})_{(\alpha,j)\in v\bar\Psi^+}$ on $U_v$ such that, for every $a\in\scrA$,
\begin{equation}\label{uniform rescaled coordinates}
\tilde O[\calU_{a,v}]
=\tilde O\big[Y_{\alpha,j}\ \big|\ (\alpha,j)\in v\bar\Psi^+\big],
\qquad
Y_{\alpha,j}
=\pi_{\tilde E}^{e\lfloor\psi_{\alpha,j}(a)\rfloor}X_{\alpha,j}.
\end{equation}
\end{proposition}

\begin{proof}
We construct the coordinates separately for the two types of nondivisible roots.

\subsubsection*{Case: $\al\in\Phi$ is neither multipliable nor divisible.}
Then we have an isomorphism
\[
x_\alpha\colon\Res_{L_\alpha/\breve F}\mathbb G_a\xrightarrow{\ \sim\ }U_\alpha,
\qquad
x_\alpha(r)=\prod_{\sigma\in\Gamma/\Gamma_{\tilde\alpha}}
y_{\sigma\tilde\alpha}\big(\sigma(r)\big),
\]
whose factors commute \cite[4.1.5]{BT84}.
Writing $r=\sum_{j=0}^{m_\alpha-1}r_j\varpi_\alpha^j$ with $r_j\in\breve F$, we get
\[
x_\alpha(r)=\prod_{j=0}^{m_\alpha-1}x_{\alpha,j}(r_j),
\qquad
x_{\alpha,j}(z):=\prod_{\sigma\in\Gamma/\Gamma_{\tilde\alpha}}
y_{\sigma\tilde\alpha}\big(\sigma(\varpi_\alpha)^jz\big).
\]
We denote by $X_{\alpha,j}\in\breve F[U_\alpha]$ the coordinate corresponding to $x_{\alpha,j}$.
These form a basis of the $\breve F$-dual of $\Lie(U_\alpha)$, and can be expressed over $\tilde E$ in terms of the coordinates corresponding to the $y_{\sigma\tilde\alpha}$.
Since $o$ is the Chevalley--Steinberg point, one has $\calU_{\alpha,a}(\breve O)=\{x_\alpha(r)\mid\mathrm{val}(r)\geq-\alpha(a)\}$ \cite[4.2.2]{BT84}, where $\mathrm{val}(\pi)=1$.
As $\mathrm{val}(r)=\min_j\{\mathrm{val}(r_j)+j/m_\alpha\}$, comparing with the valuations of the $r_j$, we obtain
\begin{equation}\label{tame root coordinates}
\tilde O[\calU_{\alpha,a}]
=\tilde O[Y_{\alpha,j}\mid 0\leq j<m_\alpha],
\qquad
Y_{\alpha,j}
=\pi_{\tilde E}^{e\left\lfloor\alpha(a)+\frac{j}{m_\alpha}\right\rfloor}X_{\alpha,j}.
\end{equation}
In the split case, $m_\alpha=e=1$ and this reads $Y_\alpha=\pi^{\lfloor\alpha(a)\rfloor}X_\alpha$, as in the proof of \Cref{key lemma}.

\subsubsection*{Case: $\al\in\Phi$ is multipliable.}
Then $G$ is of ramified type ${}^2A_{2n}$ and $e=m_\alpha=2$; see \cite[\S2.6(b), p.~97]{KP23}.
Thus $e_\alpha=1$, $L_\alpha=\tilde E$ and $\varpi_\alpha=\pi_{\tilde E}$.
Since $\pi_{\tilde E}^2=\pi$, the nontrivial element $\sigma_0$ of $\Gamma$ satisfies $\sigma_0(\pi_{\tilde E})=-\pi_{\tilde E}$.
Tameness means $p\neq2$.
Moreover, $\tilde\alpha':=\sigma_0\tilde\alpha\neq\tilde\alpha$, and $\tilde\alpha+\tilde\alpha'$ is an absolute root, which restricts to $2\alpha$.
As in \cite[\S2]{Lou23}, we choose the lift $\tilde\alpha$ and the Chevalley--Steinberg system so that $[E_{\tilde\alpha},E_{\tilde\alpha'}]=E_{\tilde\alpha+\tilde\alpha'}$ for the associated root vectors, as for the standard pinning of $\mathrm{SL}_3$.
Then the root group $U_\alpha\supset U_{2\alpha}$ is a Heisenberg group \cite[4.1.9]{BT84}, \cite[D\'efinition~2.7]{Lou23}:
\begin{equation}\label{Heisenberg parametrization}
U_\alpha(\breve F)
=\Big\{x_\alpha(z,w):=
y_{\tilde\alpha'}(\sigma_0(z))\,y_{\tilde\alpha+\tilde\alpha'}(w)\,y_{\tilde\alpha}(z)
\ \Big|\ z,w\in\tilde E,\ w+\sigma_0(w)=z\sigma_0(z)\Big\},
\end{equation}
and $U_{2\alpha}(\breve F)=\{x_\alpha(0,w)\}$.
Moreover, $\calU_{\alpha,a}(\breve O)=\{x_\alpha(z,w)\mid\mathrm{val}(w)\geq-2\alpha(a)\}$ \cite[4.2.2]{BT84}, which forces $\mathrm{val}(z)\geq-\alpha(a)$.
Every pair $(z,w)$ in \Cref{Heisenberg parametrization} can be written uniquely as
\[
z=r_0+\pi_{\tilde E}r_1,
\qquad
w=\tfrac12z\sigma_0(z)+\pi_{\tilde E}r_2,
\qquad r_0,r_1,r_2\in\breve F,
\]
as $w-\frac12z\sigma_0(z)$ has trace zero.
Let $X_{\alpha,0},X_{\alpha,1},X_{2\alpha,1}\in\breve F[U_\alpha]$ be the coordinates $r_0,r_1,r_2$.
As $z\sigma_0(z)=r_0^2-\pi r_1^2$, the valuations of the three summands $r_0^2$, $\pi r_1^2$ and $\pi_{\tilde E}r_2$ lie in $2\bbZ$, $1+2\bbZ$ and $\frac12+\bbZ$, so that $\mathrm{val}(w)=\min\{2\,\mathrm{val}(r_0),1+2\,\mathrm{val}(r_1),\frac12+\mathrm{val}(r_2)\}$.
Hence the condition $\mathrm{val}(w)\geq-2\alpha(a)$ yields
\begin{equation}\label{Heisenberg rescaled coordinates}
\tilde O[\calU_{\alpha,a}]
=\tilde O[Y_{\alpha,0},Y_{\alpha,1},Y_{2\alpha,1}],
\quad
Y_{\alpha,j}
=\pi_{\tilde E}^{e\left\lfloor\alpha(a)+\frac{j}{2}\right\rfloor}X_{\alpha,j},
\quad
Y_{2\alpha,1}
=\pi_{\tilde E}^{e\left\lfloor2\alpha(a)+\frac12\right\rfloor}X_{2\alpha,1}.
\end{equation}

Taking products of the root group models over the nondivisible roots in $v\Phi^+$, Equations~\eqref{tame root coordinates} and \eqref{Heisenberg rescaled coordinates} give \Cref{uniform rescaled coordinates}.
\end{proof}

\begin{remark}\label{periodic coordinates remark}
It is convenient to extend the coordinates of \Cref{proposition uniform rescaled coordinates} periodically, with notation as in its proof.

If $\alpha$ is neither multipliable nor divisible, then for $n\in\bbZ$, let $X_{\alpha,n}\in\breve F[U_\alpha]$ be the coordinate corresponding to
\[
x_{\alpha,n}(z):=x_\alpha(\varpi_\alpha^nz),
\qquad z\in\breve F,
\]
so that $X_{\alpha,n+m_\alpha}=\pi^{-1}X_{\alpha,n}$ for all $n$.
As $X_{\alpha,j-m_\alpha N}=\pi^NX_{\alpha,j}$, \Cref{tame root coordinates} can be written as
\begin{equation}\label{periodic root coordinates}
\tilde O[\calU_{\alpha,a}]
=\tilde O\Big[X_{\alpha,n}\ \Big|\ 0\leq\alpha(a)+\frac{n}{m_\alpha}<1\Big].
\end{equation}
Thus $X_{\alpha,n}$ corresponds to the affine function $\alpha+n/m_\alpha$, or equivalently to the affine root $m_\alpha\alpha+n$ associated with $\Sigma$.
The generators at $a$ are those whose affine function takes values in $[0,1)$ at $a$.

If $\alpha$ is multipliable, we extend the coordinates of \eqref{Heisenberg rescaled coordinates} by $X_{\alpha,n+2}=\pi^{-1}X_{\alpha,n}$ for $n\in\bbZ$ and $X_{2\alpha,n+2}=\pi^{-1}X_{2\alpha,n}$ for odd $n$, so that \eqref{Heisenberg rescaled coordinates} reads
\begin{equation}\label{Heisenberg coordinates}
\tilde O[\calU_{\alpha,a}]
=\tilde O\Big[X_{\alpha,n},X_{2\alpha,n'}\ \Big|\
0\leq\alpha(a)+\frac n2<1,\ n'\text{ odd},\
0\leq2\alpha(a)+\frac{n'}2<1\Big].
\end{equation}
Thus $X_{\alpha,n}$ and $X_{2\alpha,n'}$ correspond to the affine functions $\alpha+n/2$ and $2\alpha+n'/2$, or equivalently to the affine roots $4\alpha+2n$ and $4\alpha+n'$ associated with $\Sigma$.
The generators at $a$ are again those whose affine function takes values in $[0,1)$ at $a$.

By Equations~\eqref{periodic root coordinates} and \eqref{Heisenberg coordinates}, $(\alpha,j)\mapsto\psi_{\alpha,j}+\bbZ$ identifies $\bar\Psi$ with the set of $\bbZ$-orbits of the set $\Psi$ of affine roots, see \Cref{tangent cone}, and $v\bar\Psi^+$ with those whose gradient lies in $v\Phi^+$.
\end{remark}

\begin{corollary}\label{Iwahori-minor-expansion recall coro}
Let $b$ be an interior point of $\mathbf a$.
For $(\alpha,j)\in\bar\Psi$, put
\[
n_b(\alpha,j):=
j-\frac{e}{e_\alpha}\lfloor\psi_{\alpha,j}(b)\rfloor.
\]
Then $n_b(\alpha,j)=j$ unless $j=0$ and $\alpha\in\Phi^-$, in which case $n_b(\alpha,0)=m_\alpha$, and
\begin{equation}\label{Iwahori-minor-expansion recall}
\tilde O[\calU_{b,v}]
=
\tilde O\big[X_{\alpha,n_b(\alpha,j)}
\ \big|\ (\alpha,j)\in v\bar\Psi^+\big].
\end{equation}
\end{corollary}

\begin{proof}
For every nondivisible $\alpha\in\Phi^+$, one has
$0<\alpha(b)<1/m_\alpha$, and even
$0<\alpha(b)<1/(2m_\alpha)$ if $\alpha$ is multipliable,
since $m_\alpha\alpha$, respectively $2m_\alpha\alpha$,
lies in $\Sigma$.
Put
$\epsilon(\alpha)=0$ if $\alpha\in\Phi^+$ and
$\epsilon(\alpha)=1$ if $\alpha\in\Phi^-$.
Then, for $(\alpha,j)\in\bar\Psi$, one has
$\lfloor\psi_{\alpha,j}(b)\rfloor=-\epsilon(\alpha)$ if $j=0$,
and $\lfloor\psi_{\alpha,j}(b)\rfloor=0$ otherwise
(for divisible $\alpha$, this uses $|\alpha(b)|<\frac12$).
This gives the asserted formula for $n_b(\alpha,j)$.
Since
\[
je_\alpha-e\lfloor\psi_{\alpha,j}(b)\rfloor
=
e_\alpha n_b(\alpha,j),
\]
\Cref{Iwahori-minor-expansion recall} follows from \eqref{uniform rescaled coordinates} and the periodicity of the coordinates; see \Cref{periodic coordinates remark}.
\end{proof}

\subsection{Coordinates on the split form}\label{split coordinates sec}
We continue with the notation and assumptions of \Cref{coordinates sec}.
In particular, $G_{\breve F}$ is simply connected, absolutely almost simple and tamely ramified with splitting field $\tilde E$.
The proof of \Cref{prop:extremal-minor-reduction} proceeds by reduction to equicharacteristic zero, achieved by comparing the integral root group coordinates of $G$ with those of its split form.

Choose a generator $\sigma$ of $\Gamma$ and put
\[
\eta_e:=\frac{\sigma(\pi_{\tilde E})}{\pi_{\tilde E}}
\in\tilde O^\times.
\]
Since $\pi_{\tilde E}^e=\pi$ and $e$ is prime to the residue characteristic, $\eta_e$ is an $e$-th root of unity, and all such roots lie in $\breve O^\times$.
Thus $\sigma^r(\pi_{\tilde E})=\eta_e^r\pi_{\tilde E}$.
As $\tilde E=\breve F(\pi_{\tilde E})$ and $\sigma$ has order $e$, the root $\eta_e$ has exact order $e$.

Put $\Lambda:=\bbZ[\zeta_e,1/e]$, where $\zeta_e$ is an abstract primitive $e$-th root of unity.
We regard $\breve O$, and hence $\tilde O$ and $\tilde E$, as $\Lambda$-algebras via $\zeta_e\mapsto\eta_e$; all base changes from $\Lambda$ to these rings below are formed in this way.
For $u\in\Lambda$, we write $\bar u\in\breve O$ for its image.
The map $u\mapsto\bar u$ is injective in characteristic zero; in equal characteristic $p$, its kernel contains $p$.
For the characteristic-zero comparison, we also fix an embedding $\Lambda\subset\bbC$, and regard $\bbC$ as a $\Lambda$-algebra through it.

Let $H$ be the split, simply connected Chevalley group over $\Lambda$ with the same absolute root datum as $G$, see \cite{SGA3}, with maximal split torus $T_H$ and a Chevalley system $(y^H_\gamma\colon\mathbb G_{a,\Lambda}\xrightarrow{\sim}U_{H,\gamma})_{\gamma\in\Phi_{\mathrm{abs}}}$.
For every $\gamma\in\Phi_{\mathrm{abs}}$, we write $Z_\gamma\in\Lambda[U_{H,\gamma}]$ for the coordinate corresponding to $ y^H_\gamma$, and $E_\gamma:=\mathrm{d}y^H_\gamma(1)\in\Lie(U_{H,\gamma})$ for the dual root vector.
As the Chevalley--Steinberg system $(y_\gamma)$ is a Chevalley system of $G_{\tilde E}$, see \cite[4.1.3]{BT84}, there is an isomorphism $H_{\tilde E}\simeq G_{\tilde E}$ carrying $T_{H,\tilde E}$ to $T_{\tilde E}$ and $y^H_\gamma$ to $y_\gamma$ for all $\gamma$, after possibly changing the signs of the $y^H_\gamma$.
We fix such an isomorphism, and view the $Z_\gamma$ and $E_\gamma$ also as coordinates and root vectors of $G_{\tilde E}$.
For $v\in W_0$, let $U_{H,v}\subset H$ be the subgroup corresponding to $U_{v,\tilde E}$, that is, the product of the $U_{H,\gamma}$ over the absolute roots $\gamma$ restricting to elements of $v\Phi^+$.
We order $v\Phi^+$ and the absolute roots restricting to it so that $\tilde\alpha',\tilde\alpha+\tilde\alpha',\tilde\alpha$ are consecutive, in this order, for multipliable $\alpha$; then $\prod_\gamma y^H_\gamma$ is an isomorphism onto $U_{H,v}$, and we view the $Z_\gamma$ as coordinates on $U_{H,v}$.

Now let $\alpha\in \Phi$ be nondivisible.
Assume first that $\alpha$ is not multipliable.
For $0\leq r,j<m_\alpha$, put
\begin{equation}\label{Fourier vectors}
c_{\sigma^r\tilde\alpha,j}:=\zeta_e^{r e_\alpha j}\in\Lambda^\times,
\qquad
E_{\alpha,j}:=\sum_{r=0}^{m_\alpha-1}
c_{\sigma^r\tilde\alpha,j}E_{\sigma^r\tilde\alpha}.
\end{equation}
These coefficients specialize to the arithmetic ones:
\[
\bar c_{\sigma^r\tilde\alpha,j}
=\eta_e^{re_\alpha j}
=\frac{\sigma^r(\varpi_\alpha)^j}{\varpi_\alpha^j},
\qquad
\varpi_\alpha=\pi_{\tilde E}^{e_\alpha}.
\]
They are well defined on the orbit, since changing $r$ by $m_\alpha$ does not change $\zeta_e^{r e_\alpha j}$.
The Fourier matrix $(c_{\sigma^r\tilde\alpha,j})_{r,j}$ has inverse
$\bigl(m_\alpha^{-1}\zeta_e^{-r e_\alpha j}\bigr)_{j,r}$, since $m_\alpha\mid e$ and the sums of nontrivial characters of $\bbZ/m_\alpha\bbZ$ vanish.
It is therefore invertible over $\Lambda$ and after every base change.
So the $E_{\alpha,j}$, $0\leq j<m_\alpha$, form a $\Lambda$-basis of $\bigoplus_{\gamma\in\calO_\alpha}\Lie(U_{H,\gamma})$, and we let $\xi_{\alpha,j}$ be the dual linear coordinates on the commutative group $\prod_{\gamma\in\calO_\alpha}U_{H,\gamma}$, so that
\[
Z_{\gamma}=\sum_{j=0}^{m_\alpha-1}c_{\gamma,j}\,\xi_{\alpha,j},\qquad\gamma\in\calO_\alpha.
\]

Now assume that $\alpha$ is multipliable, so that $e=m_\alpha=2$ and $\calO_\alpha=\{\tilde\alpha,\tilde\alpha'\}$.
Then \Cref{Fourier vectors} reads $E_{\alpha,j}=E_{\tilde\alpha}+(-1)^jE_{\tilde\alpha'}$ for $j=0,1$, and we put $E_{2\alpha,1}:=E_{\tilde\alpha+\tilde\alpha'}$.
We define coordinates $\xi_{\alpha,0},\xi_{\alpha,1},\xi_{2\alpha,1}$ on $U_{H,\tilde\alpha'}U_{H,\tilde\alpha+\tilde\alpha'}U_{H,\tilde\alpha}$ over $\Lambda$ by
\begin{equation}\label{Heisenberg Fourier}
Z_{\tilde\alpha}=\xi_{\alpha,0}+\xi_{\alpha,1},\qquad
Z_{\tilde\alpha'}=\xi_{\alpha,0}-\xi_{\alpha,1},\qquad
Z_{\tilde\alpha+\tilde\alpha'}=\xi_{2\alpha,1}+\tfrac12\big(\xi_{\alpha,0}^2-\xi_{\alpha,1}^2\big).
\end{equation}
Their differentials at the identity are dual to $E_{\alpha,0},E_{\alpha,1},E_{2\alpha,1}$.

In both cases, we obtain coordinates $\xi_{\alpha,j}$, $(\alpha,j)\in v\bar\Psi^+$, on $U_{H,v}$ over $\Lambda$, whose differentials at the identity are dual to the basis $(E_{\alpha,j})_{(\alpha,j)\in v\bar\Psi^+}$ of $\Lie(U_{H,v})$.
The following lemma compares them with the coordinates of \Cref{proposition uniform rescaled coordinates}.

\begin{lemma}\label{comparison of coordinates}
Under the identification $H_{\tilde E}\simeq G_{\tilde E}$, one has
\[
\xi_{\alpha,j}=\pi_{\tilde E}^{je_\alpha}X_{\alpha,j}
\qquad\text{on }U_{v,\tilde E},\text{ for all }(\alpha,j)\in v\bar\Psi^+.
\]
Dually, the tangent vector dual to the periodic coordinate $X_{\alpha,n}$ of \Cref{periodic coordinates remark} is $\pi_{\tilde E}^{e_\alpha n}E_{\alpha,j}$ for $n\equiv j \mod m_\alpha$, and the one dual to $X_{2\alpha,n'}$ is $\pi_{\tilde E}^{n'}E_{2\alpha,1}$ for odd $n'$.
\end{lemma}
\begin{proof}
It suffices to treat a single nondivisible $\alpha\in v\Phi^+$.
If $\alpha$ is not multipliable, then by the definition of the $x_{\alpha,j}$ in the proof of \Cref{proposition uniform rescaled coordinates}, one has
\[
Z_{\sigma^r\tilde\alpha}=\sum_j\sigma^r(\varpi_\alpha)^jX_{\alpha,j}=\sum_j\bar c_{\sigma^r\tilde\alpha,j}\,\varpi_\alpha^jX_{\alpha,j}
\]
on $U_{\alpha,\tilde E}$, so the invertibility of the specialized Fourier matrix gives $\xi_{\alpha,j}=\varpi_\alpha^jX_{\alpha,j}=\pi_{\tilde E}^{je_\alpha}X_{\alpha,j}$.
If $\alpha$ is multipliable, then $Z_{\tilde\alpha}=z$, $Z_{\tilde\alpha'}=\sigma_0(z)$ and $Z_{\tilde\alpha+\tilde\alpha'}=w$ at the point $x_\alpha(z,w)$ of \Cref{Heisenberg parametrization}.
Substituting $z=X_{\alpha,0}+\pi_{\tilde E}X_{\alpha,1}$ and $w=\frac12z\sigma_0(z)+\pi_{\tilde E}X_{2\alpha,1}$, where $z\sigma_0(z)=X_{\alpha,0}^2-\pi X_{\alpha,1}^2$, into \Cref{Heisenberg Fourier} gives $\xi_{\alpha,j}=\pi_{\tilde E}^jX_{\alpha,j}$ for $j=0,1$ and $\xi_{2\alpha,1}=\pi_{\tilde E}X_{2\alpha,1}$; here $e_\alpha=e_{2\alpha}=1$.
The dual statement follows by passing to differentials at the identity, using $X_{\alpha,n+m_\alpha}=\pi^{-1}X_{\alpha,n}$, $X_{2\alpha,n'+2}=\pi^{-1}X_{2\alpha,n'}$ and $\pi=\varpi_\alpha^{m_\alpha}$.
\end{proof}

\subsection{Lower bound at Iwahori level}\label{iwahori calc sec}
For every $\chi\in X^*(T)^+$ and $v,w\in W_0$, we expand the function $P^\chi_{v,w}$ (\Cref{P generators}) in the coordinates on $U_v$ constructed in \Cref{proposition uniform rescaled coordinates}:
\begin{equation}
\label{Iwahori-minor-expansion}
P_{v,w}^{\chi}
=
\sum_{\mathbf n}c_{\mathbf n}
\prod_{(\alpha,j)\in v\bar\Psi^+}X_{\alpha,j}^{n_{\alpha,j}}\in \tilde E[U_v],
\end{equation}
where $\mathbf n=(n_{\alpha,j})$ runs through $\bbZ_{\geq0}^{v\bar\Psi^+}$.

\begin{lemma}
For every $\chi\in X^*(T)^+$ and $v,w\in W_0$, one has
\begin{equation}\label{tame Gauss formula}
-h_{a,\chi,v}(w\chi)
=\min_{c_{\mathbf n}\ne0}
\left\{\operatorname{val}_{\tilde E}(c_{\mathbf n})-e
\sum_{(\alpha,j)\in v\bar\Psi^+}\left\lfloor\psi_{\alpha,j}(a)\right\rfloor n_{\alpha,j}\right\},
\end{equation}
where $\mathrm{val}_{\tilde E}(\pi_{\tilde E})=1$.
\end{lemma}
\begin{proof}
By the proof of \Cref{weight module monoidal}, the elementary divisor $-h_{a,\chi,v}(w\chi)$ is the Gauss valuation of $P^\chi_{v,w}$ in $\tilde O[\calU_{a,v}]$. 
So, \Cref{uniform rescaled coordinates} proves the lemma. 
\end{proof}

We describe the coefficients $c_{\mathbf n}$ of \Cref{Iwahori-minor-expansion}, with the aim of providing a universal formula: they are images in $\tilde O$ of universal elements of $\Lambda[t]$ under the map $u\mapsto\bar u$ extended by $t\mapsto\pi_{\tilde E}$, up to a common unit of $\tilde O$, see \Cref{rem:integrality-of-coefficients} below.

\begin{lemma}
\label{rem:integrality-of-coefficients}
With the split root datum, Chevalley system, generator of $\Gamma$, and root lifts fixed as above, there are universal elements $u_{\mathbf n}\in\Lambda$, depending on $\chi,v,w$, and a unit $\kappa\in\tilde O^\times$, such that the coefficients $c_{\mathbf n}$ of \Cref{Iwahori-minor-expansion} are given by
\[
c_{\mathbf n}=\kappa\,\bar u_{\mathbf n}\,\pi_{\tilde E}^{\sum_{(\alpha,j)\in v\bar\Psi^+}je_\alpha n_{\alpha,j}}.
\]
More precisely, let $B^-_{H,v}\subset H$ be the Borel subgroup containing $T_H$ that corresponds to $B^-_{v,\tilde E}$, let $V^H_{\chi,v}:=\Ind_{B^-_{H,v}}^{H}(v\chi)$ be the induced $H$-module over $\Lambda$, and let $f$ be a generator of its weight space $V^H_{\chi,v}(w\chi)$.
Then $f|_{U_{H,v}}=\sum_{\mathbf n}u_{\mathbf n}\prod_{(\alpha,j)\in v\bar\Psi^+}\xi_{\alpha,j}^{n_{\alpha,j}}$.
\end{lemma}

\begin{proof}
The $\Lambda$-module $V^H_{\chi,v}$ is free of finite rank, its formation commutes with base change, and its extremal weight spaces are direct summands of rank one.
Choose such a generator $f$, which is well defined up to $\Lambda^\times$. Its restriction to $U_{H,v}$ is a polynomial over $\Lambda$ in the coordinates $Z_\gamma$, hence also in the coordinates $\xi_{\alpha,j}$ of \Cref{comparison of coordinates}. 
This defines the coefficients $u_{\mathbf n}\in\Lambda$.
Let $\iota(o)$ be the hyperspecial point associated with the chosen isomorphism $H_{\tilde E}\simeq G_{\tilde E}$.
By the proof of \Cref{inclusion lemma key}, the extremal lines of $\tilde\calV_{o,\chi,v}$ and of $\tilde\calV_{\iota(o),\chi,v}=V^H_{\chi,v}\otimes_\Lambda\tilde O$ agree, so $P_{v,w}=\kappa\cdot(f\otimes_\Lambda1)|_{U_{v,\tilde E}}$ for some $\kappa\in\tilde O^\times$.
Substituting $\xi_{\alpha,j}=\pi_{\tilde E}^{je_\alpha}X_{\alpha,j}$ from \Cref{comparison of coordinates} proves the claim.
\end{proof}

\begin{definition}
For $\chi\in X^*(T)^+$ and $v,w\in W_0$, put
\[
\deg_{v,w}^\chi(\mathbf n)
:=\sum_{(\alpha,j)\in v\bar\Psi^+}
e_\alpha\,n_b(\alpha,j)\,n_{\alpha,j},
\]
and define
\[
d_\chi(v,w)
:=\min_{\bar u_{\mathbf n}\ne0}
\left\{\operatorname{val}_{\tilde E}(\bar u_{\mathbf n})
+\deg_{v,w}^\chi(\mathbf n)\right\}.
\]
Here $u_{\mathbf n}\in\Lambda$ are the universal coefficients of
\Cref{rem:integrality-of-coefficients}, and the minimum is taken over their nonzero images $\bar u_{\mathbf n}\in\breve O$.
\end{definition}

If $G_{\breve F}$ splits, then $m_\alpha=e_\alpha=1$ and $n_b(\alpha,0)=\epsilon(\alpha)$, so $\deg_{v,w}^\chi(\mathbf n)=\sum_{\alpha\in v\Phi^+\cap\Phi^-}n_\alpha$.

\begin{lemma}\label{interior point lemma}
At an interior point $b$ of the alcove $\mathbf a$, one has $h_{b,\chi,v}(w\chi)=-d_\chi(v,w)$.
\end{lemma}
\begin{proof}
By \Cref{tame Gauss formula}, the integer $-h_{b,\chi,v}(w\chi)$ is the
minimum valuation of the coefficients of $P^\chi_{v,w}$ in the
integral coordinates
$Y_{\alpha,j}
=\pi_{\tilde E}^{e\lfloor\psi_{\alpha,j}(b)\rfloor}X_{\alpha,j}$
of \Cref{proposition uniform rescaled coordinates}.
By \Cref{rem:integrality-of-coefficients}, the coefficient of
$\prod_{(\alpha,j)}Y_{\alpha,j}^{n_{\alpha,j}}$ has valuation
\[
\operatorname{val}_{\tilde E}(\bar u_{\mathbf n})
+\sum_{(\alpha,j)\in v\bar\Psi^+}
 \bigl(je_\alpha-e\lfloor\psi_{\alpha,j}(b)\rfloor\bigr)
 n_{\alpha,j}.
\]
Since $b$ is interior to $\mathbf a$,
\Cref{Iwahori-minor-expansion recall coro} gives
$je_\alpha-e\lfloor\psi_{\alpha,j}(b)\rfloor
=e_\alpha n_b(\alpha,j)$.
Taking the minimum therefore yields
$-h_{b,\chi,v}(w\chi)=d_\chi(v,w)$.
\end{proof}

The quantities $\deg_{v,w}^\chi(\mathbf n)$ are related to the quantum Bruhat graph by the following result.

\begin{proposition}[{\Cref{prop:characteristic-zero-extremal-degree}}]
For a fundamental weight $\chi$ and $v,w\in W_0$, one has
\[
\min_{u_{\mathbf n}\neq0}\{\deg_{v,w}^\chi(\mathbf n)\}
=\left\langle N\chi,\operatorname{wt}(v\Rightarrow w)\right\rangle.
\]
In particular, $\left\langle N\chi,\operatorname{wt}(v\Rightarrow w)\right\rangle\leq d_\chi(v,w)$.
\end{proposition}
\begin{proof}
On the left-hand side, the condition $u_{\mathbf n}\neq0$ in the minimum is tested in $\Lambda$.
As $\Lambda\subset\bbC$, this is the same as testing it in $\bbC$.
So the equality is \Cref{prop:characteristic-zero-extremal-degree}, which is proved in \Cref{appendix feigin makedonskyi} using the graded characters of generalized Weyl modules of Feigin--Makedonskyi \cite{FeiginMakedonskyi17,FeiginMakedonskyiTwisted17}.
For the inequality, let $\mathbf n$ be such that $\bar u_{\mathbf n}\neq0$.
Then $u_{\mathbf n}\neq0$, so $\deg_{v,w}^\chi(\mathbf n)\geq\langle N\chi,\operatorname{wt}(v\Rightarrow w)\rangle$ by the equality.
Moreover, $\operatorname{val}_{\tilde E}(\bar u_{\mathbf n})\geq0$, as $\bar u_{\mathbf n}\in\breve O$.
Hence every term in the minimum defining $d_\chi(v,w)$ is at least $\langle N\chi,\operatorname{wt}(v\Rightarrow w)\rangle$.
\end{proof}

\subsection{Upper bound at Iwahori level}
The previous section provided a lower bound for the elementary divisor $-h_{b,\chi,v}(w\chi)=d_\chi(v,w)$ at Iwahori level.
In this subsection, we continue our calculation of $d_{\chi}(v,w)$ by providing an upper bound.

First, we recall some facts on distribution algebras \cite[\S\S I.7, II.1.12]{Jan03}.
For a flat affine group scheme $\calH$ over a ring $R$ with augmentation ideal $\frakm_e\subset R[\calH]$, the distribution algebra is
\[
\Dist(\calH)=\bigcup_{n\geq0}\Hom_R\big(R[\calH]/\frakm_e^{n+1},R\big),
\]
with multiplication dual to the comultiplication of $R[\calH]$.
Every representation of $\calH$ on an $R$-module is a $\Dist(\calH)$-module, by $D\cdot m=(\id\otimes D)(\Delta(m))$ for the comodule map $\Delta$ \cite[I.7.11]{Jan03}.
If $R$ is a field of characteristic $0$, then $\Dist(\calH)=U(\Lie\calH)$ \cite[I.7.10]{Jan03}.
For $\calH=\mathbb G_a$ with coordinate $Y$, $\Dist(\calH)$ is free with basis the divided powers $E^{(n)}$ dual to $Y^n$, and $E^{(n)}E^{(m)}=\binom{n+m}{n}E^{(n+m)}$.
If $\calH=\calH_1\times\calH_2$ as schemes, then $\Dist(\calH)=\Dist(\calH_1)\otimes_R\Dist(\calH_2)$ \cite[I.7.4]{Jan03}.
For the split group $H_\bbZ$, this gives the Kostant $\bbZ$-form $\Dist(H_\bbZ)\subset U(\frakg_{\bbQ})$, generated by the $E_\gamma^{(n)}$ and $\Dist(T_{H,\bbZ})$, and a decomposition
\[
\Dist(U^-_{H,v})\otimes\Dist(T_H)\otimes\Dist(U_{H,v})\xrightarrow{\sim}\Dist(H)
\]
induced by the big cell \cite[II.1.12]{Jan03}.
The same decomposition holds for the parahoric group schemes $\calG_b$ over $\breve O$, using the integral big cell \cite[Theorem~8.2.5]{KP23}.

We will show that a coefficient of minimal universal degree remains a unit after specialization.
The first ingredient is the following finite-sum formula for root-group actions.

\begin{lemma}
\label{lem:finite-sum-action}
Let $R$ be a ring, let $\calU=\mathbb G_{a,R}$ with coordinate $Y$, and let $V$ be a representation of $\calU$
on an $R$-module, with comodule structure map $\Delta_V\colon V\to V\otimes_RR[Y]$.
For $m\in V$ write $\Delta_V(m)=\sum_{n\geq0}m_n\otimes Y^{n}$, a finite sum.  Then
$m_n=E^{(n)}m$, and for every $c\in R$ the element $u\in\calU(R)$ with $Y(u)=c$
acts by
\[
u\cdot m=\sum_{n\geq0}c^{n}E^{(n)}m.
\]
In particular every $\Dist(\calU)$-stable $R$-submodule of $V$
is stable under $\calU(R)$.
\end{lemma}

\begin{proof}
By definition $E^{(n)}\in\Dist(\calU)$ is the linear form on $R[\calU]=R[Y]$ dual
to $Y^{n}$, and $\Dist(\calU)$ acts on $V$ by
$D\cdot m=(\mathrm{id}\otimes D)\Delta_V(m)$, whence $E^{(n)}m=m_n$.  The point $u$
is the $R$-algebra map $\mathrm{ev}_u\colon R[Y]\to R$ with $Y\mapsto c$, and it acts
on the comodule $V$ by $u\cdot m=(\mathrm{id}\otimes\mathrm{ev}_u)\Delta_V(m)$,
which is $\sum_{n}c^{n}m_n$.
\end{proof}

\begin{lemma}
\label{lem:root-string}
Let $\eta\in W_0\chi$ be an extremal weight,
let $(\beta,j)\in\bar\Psi$, and assume that $r:=-\langle\eta,\beta^\vee\rangle\geq0$.
Then, the divided power $E_{\beta,j}^{(r)}$ carries $\tilde\calV_{o,\chi,v}(\eta)$ isomorphically to $\tilde\calV_{o,\chi,v}(s_\beta\eta)$.
\end{lemma}

\begin{proof}
We first treat an absolute root $\gamma$, in the split Chevalley lattice $\tilde\calV_{\iota(o),\chi,v}$ attached to the hyperspecial point $\iota(o)$, and put $m_\gamma:=-\langle\eta,\gamma^\vee\rangle\geq0$.
It is well known that the weights of
$\tilde\calV_{\iota(o),\chi,v}$ lying in
  $\eta+\bbZ\gamma$ are exactly
  $\eta,\eta+\gamma,\ldots,\eta+m_\gamma\gamma=s_\gamma\eta$, as this is the intersection of $\chi + \bbZ\Phi_{\mathrm{abs}}$ with the convex hull of $W_{0,C}\chi$.
  These weight considerations imply that
$E_{-\gamma}^{(n)}f=0$ for $n\geq1$ and $E_{\gamma}^{(n)}f=0$ for $n>m_\gamma$, where $f$ generates the $\eta$-line.
Let $\calH_\gamma$ be the reductive $\tilde O$-group scheme of semisimple rank one generated by the root groups of $\pm\gamma$ in the split model at $\iota(o)$.
So the $\Dist(\calH_\gamma)$-submodule generated by $f$ is
\[
M:=\Dist(\calH_\gamma)\cdot f=\sum_{n=0}^{m_\gamma}\tilde O\,E_\gamma^{(n)}f,
\qquad
M(s_\gamma\eta)=\tilde O\,E_\gamma^{(m_\gamma)}f .
\]
By \Cref{lem:finite-sum-action} the elements $y_{\pm\gamma}(c)$, $c\in\tilde O$, preserve $M$.
The same therefore holds for the lift $\dot s_\gamma$ of $s_\gamma$ in the subgroup
generated by them.
Consequently $\dot s_\gamma f\in M(s_\gamma\eta)=\tilde O\,E_\gamma^{(m_\gamma)}f$.
Since $\dot s_\gamma$ is invertible, it carries the $\eta$-line of $\tilde\calV_{\iota(o),\chi,v}$ isomorphically onto the $s_\gamma\eta$-line, which is therefore generated by $E_\gamma^{(m_\gamma)}f$.

Finally, we handle the relative case.
Write $\beta\in\{\alpha,2\alpha\}$ with $\alpha$ nondivisible, and put $m:=-\langle\eta,\alpha^\vee\rangle\geq0$, so that $r=m$ if $\beta=\alpha$ and $r=m/2$ if $\beta=2\alpha$; note $s_\beta=s_\alpha$.
Now let $\alpha\in\Phi$ be neither divisible nor multipliable.
Since $\eta=u\chi$ for some $u\in W_0$ and $u^{-1}\calO_\alpha=\calO_{u^{-1}\alpha}$ consists of absolute roots of the same sign, the integers $m_\gamma$, $\gamma\in\calO_\alpha$, are nonnegative, with sum $m$.
Since the roots in $\calO_\alpha$ are pairwise orthogonal, their root vectors commute.
Hence
\[
E_{\alpha,j}^{(m)}f=\sum_{\sum n_\gamma=m}\prod_{\gamma\in\calO_\alpha}\bar c_{\gamma,j}^{\,n_\gamma}E_\gamma^{(n_\gamma)}f
=\Big(\prod_{\gamma\in\calO_\alpha}\bar c_{\gamma,j}^{\,m_\gamma}\Big)\prod_{\gamma\in\calO_\alpha}E_\gamma^{(m_\gamma)}f,
\]
because every other term involves some $n_\gamma>m_\gamma$.
By the first part, applied successively, the right side generates the $s_\alpha\eta$-line of $\tilde\calV_{\iota(o),\chi,v}$, as $s_\alpha=\prod_{\gamma\in\calO_\alpha}s_\gamma$.

Now let $\alpha$ be multipliable, so that $\calO_\alpha=\{\gamma,\gamma'\}$ with $\gamma+\gamma'$ a root and $\langle\gamma',\gamma^\vee\rangle=-1$, and $\alpha^\vee=2(\gamma^\vee+\gamma'^\vee)$.
Then $m=2m_{\gamma+\gamma'}$ is even, and $s_\alpha=s_{\gamma+\gamma'}$, so the argument in the paragraph above gives the claim for $E_{2\alpha,1}^{(m/2)}=E_{\gamma+\gamma'}^{(m_{\gamma+\gamma'})}$.
For the claim about $E_{\alpha,0}$, consider the Levi subgroup $L\subset H_{\tilde O}$ of type $A_2$ generated by the root groups of $\pm\gamma,\pm\gamma'$, whose derived group is $\mathrm{SL}_3$ as $H$ is simply connected.
Put $E:=E_{\alpha,0}=E_\gamma+E_{\gamma'}$ and $F:=2(E_{-\gamma}+E_{-\gamma'})$.
In the standard matrix realization with $E_\gamma=E_{12}$ and $E_{\gamma'}=E_{23}$, these are $E_{12}+E_{23}$ and $2E_{21}+2E_{32}$, and $\alpha^\vee=2(\gamma^\vee+\gamma'^\vee)=\mathrm{diag}(2,0,-2)$.
So $(E,\alpha^\vee,F)$ is the image of the standard basis of $\mathfrak{sl}_2$ under the homomorphism $\varphi\colon\mathrm{SL}_{2,\tilde O}\to L$ given by the action of $\mathrm{SL}_2$ on $\mathrm{Sym}^2(\tilde O x\oplus\tilde O y)$ in the basis $(x^2,xy,\frac12y^2)$, which is defined as $2\in\tilde O^\times$.
Let $\calH_\alpha$ be its image, and note that $\varphi$ carries the element $\left(\begin{smallmatrix}0&1\\-1&0\end{smallmatrix}\right)$ to an antidiagonal matrix, which is a lift $\dot s_\alpha\in\varphi(\mathrm{SL}_2(\tilde O))$ of $s_{\gamma+\gamma'}=s_\alpha$.
As $m_\gamma,m_{\gamma'}\geq0$ and $\eta$ is extremal, $\eta-\gamma$ and $\eta-\gamma'$ are not weights, so $Ff=0$ and $\alpha^\vee f=-mf$.
Hence, by the $\mathfrak{sl}_2$-relations in divided powers,
\[
M:=\Dist(\calH_\alpha)\cdot f=\sum_{n=0}^{m}\tilde O\,E^{(n)}f,
\qquad
M(s_\alpha\eta)=\tilde O\,E^{(m)}f,
\]
where $M(s_\alpha\eta)$ is the part of $\alpha^\vee$-weight $m$.
As in the first paragraph, \Cref{lem:finite-sum-action}, applied to the one-parameter subgroups $\varphi$ of the upper and lower unipotent matrices, shows that $M$ is stable under $\varphi(\mathrm{SL}_2(\tilde O))$, which is generated by these subgroups as $\tilde O$ is local, so $\dot s_\alpha f\in M(s_\alpha\eta)$ generates the $s_\alpha\eta$-line, which proves the claim for $j=0$.
For $j=1$, take $t=\gamma^\vee(-1)\in T_H(\bbZ)$.
Then $\gamma(t)=1$ and $\gamma'(t)=-1$, so conjugation by $t$ carries $E_{\alpha,0}$ to $E_{\alpha,1}=E_\gamma-E_{\gamma'}$ and preserves all weight lines.
Finally, the extremal lines of $\tilde\calV_{\iota(o),\chi,v}$ and of $\tilde\calV_{o,\chi,v}$ agree, as shown in the proof of \Cref{inclusion lemma key}.
\end{proof}

\begin{lemma}
\label{lem:Iwahori-initial-form-primitive}
Let $\chi$ be a fundamental weight, let $v,w\in W_0$ and put
$d:=\langle N\chi,\operatorname{wt}(v\Rightarrow w)\rangle$.  Then some
$\bar u_{\mathbf n}$ with $\deg_{v,w}^\chi(\mathbf n)=d$ is a unit in $\tilde O$.
\end{lemma}

\begin{proof}
For a nondivisible $\beta\in\Phi$, write $\beta_\Sigma\in\Sigma$ for its positive multiple in $\Sigma$, i.e., $m_\beta\beta$ if $\beta$ is not multipliable, respectively $2m_\beta\beta$ if $\beta$ is multipliable.
Choose a shortest directed path in the quantum Bruhat graph of $\Sigma$
\[
\tau_0=v\xrightarrow{(\beta_1)_\Sigma}\tau_1\xrightarrow{(\beta_2)_\Sigma}\cdots
\xrightarrow{(\beta_s)_\Sigma}\tau_s=w,
\qquad
\tau_i=\tau_{i-1}s_{\beta_i},
\]
where $\beta_i\in\Phi^+$ is nondivisible.
Put
\[
\alpha_i=\tau_{i-1}\beta_i,
\qquad
m_i=\langle\chi,\beta_i^\vee\rangle,
\]
where $\beta_i^\vee$ is the relative coroot.
By the definition of the quantum Bruhat graph, $\epsilon(\alpha_i)=0$ on a Bruhat edge and
$\epsilon(\alpha_i)=1$ on a quantum edge.
Let $(\gamma_i,j_i)\in\bar\Psi$ be the pair with $\gamma_i\in\{\alpha_i,2\alpha_i\}$ such that $\psi_{\gamma_i,j_i}+\bbZ$ contains the affine root vanishing on the affine hyperplane $\{(\alpha_i)_\Sigma+\epsilon(\alpha_i)=0\}$, and put
\[
r_i:=-\langle\tau_i\chi,\gamma_i^\vee\rangle,
\qquad
q_i:=e_{\gamma_i}n_b(\gamma_i,j_i).
\]
Explicitly, if $\beta_i$ is not multipliable, then $\gamma_i=\alpha_i$, $j_i\equiv\epsilon(\alpha_i)$ modulo $m_{\beta_i}$, $r_i=m_i$ and $q_i=e_{\beta_i}\epsilon(\alpha_i)$.
If $\beta_i$ is multipliable, then $(\gamma_i,j_i,r_i,q_i)=(\alpha_i,0,m_i,0)$ on a Bruhat edge and $(\gamma_i,j_i,r_i,q_i)=(2\alpha_i,1,m_i/2,1)$ on a quantum edge, using the identification of \eqref{Heisenberg coordinates} with the affine roots $4\alpha+2n$ and $4\alpha+n'$ of $\Sigma$; here $m_i$ is even by \Cref{lem:root-string}.
As $\beta_i^\vee$ is $\Gal(\tilde E/\breve F)$-invariant, one has $\langle N\chi,\beta_i^\vee\rangle=e\,m_i$, and hence $\langle N\chi,(\beta_i)_\Sigma^\vee\rangle=e_{\beta_i}m_i$, respectively $m_i/2$ if $\beta_i$ is multipliable, as then $e=m_{\beta_i}=2$.
In both cases, this is $q_ir_i$ on a quantum edge, while $q_i=0$ on a Bruhat edge.
Therefore
\begin{equation}
\label{path-Iwahori-degree}
\sum_iq_ir_i
=\left\langle N\chi,\sum_{i:\,\mathrm{quantum}}(\beta_i)_\Sigma^\vee\right\rangle
=\langle N\chi,\operatorname{wt}(v\Rightarrow w)\rangle=d.
\end{equation}

Choose generators $f_i$ of the extremal lines $\tilde\calV_{o,\chi,v}(\tau_i\chi)$
so that
$f_s|_{U_v}=P_{v,w}$ up to $\tilde O$-units.  Since
$\langle\tau_i\chi,\alpha_i^\vee\rangle=-\langle\chi,\beta_i^\vee\rangle=-m_i$ and
$s_{\alpha_i}(\tau_i\chi)=\tau_{i-1}\chi$, \Cref{lem:root-string} applied to $\eta=\tau_i\chi$ and $(\gamma_i,j_i)$ gives
\begin{equation}
\label{path-divided-power}
E_{\gamma_i,j_i}^{(r_i)}f_i=u_i f_{i-1},
\end{equation}
for some unit $u_i\in \tilde O^\times$.
Let $b$ be an interior point of $\mathbf a$
and let $\calI:=\calG_b$ be the corresponding Iwahori $\breve O$-model, so that
$\tilde O[\calU_{b,v}]$ is generated by the $X_{\alpha,n_b(\alpha,j)}$ for $(\alpha,j)\in v\bar\Psi^+$, by \eqref{Iwahori-minor-expansion recall}.
Let $p_i\colon\mathbb G_{a,\tilde O}\to\calU_{\alpha_i,b}$ be the homomorphism with tangent vector $\pi_{\tilde E}^{q_i}E_{\gamma_i,j_i}$, that is, $p_i(c)=\exp(c\,\pi_{\tilde E}^{q_i}E_{\gamma_i,j_i})$.
Explicitly, in the notation of \eqref{Heisenberg parametrization} for multipliable $\beta_i$,
\[
p_i(c)=
\begin{cases}
x_{\alpha_i,\epsilon(\alpha_i)}(c), & \beta_i\text{ not multipliable},\\
x_{\alpha_i}(c,\tfrac12c^2), & \beta_i\text{ multipliable, Bruhat edge},\\
x_{\alpha_i}(0,\pi_{\tilde E}c), & \beta_i\text{ multipliable, quantum edge},
\end{cases}
\]
by \eqref{Fourier vectors}, respectively by the computation of $\exp$ in the Heisenberg group in the proof of \Cref{lem:minor-degree-one-factor-filtration}.
In all cases, the coordinate $X_{\gamma_i,n_b(\gamma_i,j_i)}$ takes the value $c$ at $p_i(c)$, and the other integral coordinates $X_{\alpha,n_b(\alpha,j)}$ of this point vanish: in the multipliable case, $x_{\alpha_i}(c,\frac12c^2)$ has $r_0=c$ and $r_1=r_2=0$, and $x_{\alpha_i}(0,\pi_{\tilde E}c)$ has $r_0=r_1=0$ and $r_2=c$ in the notation preceding \eqref{Heisenberg coordinates}.
For independent variables
$t_1,\ldots,t_s$ set
\[
h(\mathbf t)=\prod_{i=1}^s p_i(t_i),
\qquad
g(\mathbf t)=h(\mathbf t)^{-1}.
\]
Each factor lies in $\calI(\tilde O\pot{\mathbf t})$ and reduces to $1$ modulo
$\mathbf t$, hence so do $h(\mathbf t)$ and $g(\mathbf t)$.  For a formal series
$F\in\tilde O\pot{\mathbf t}$ we write $[\mathbf t^{\mathbf n}]F$ for the coefficient
of the monomial $\mathbf t^{\mathbf n}$. 
By the explicit formulas for $p_i$ and \Cref{lem:finite-sum-action}, applied to $\tilde\calV_{o,\chi,v}$ after base change to $\tilde O[\mathbf t]$, one has
\[
h(\mathbf t)\cdot f_s
=\left(\prod_{i=1}^s\sum_{n\geq0}
\pi_{\tilde E}^{nq_i}t_i^{\,n}E_{\gamma_i,j_i}^{(n)}\right)f_s,
\]
where the sums act finitely.
Since the variables $t_i$ are independent and
$f_s(g(\mathbf t))=(h(\mathbf t)\cdot f_s)(1)$ for the left action,
\eqref{path-Iwahori-degree} and \eqref{path-divided-power} give
\begin{equation}
\label{primitive-path-coefficient}
[t_1^{r_1}\cdots t_s^{r_s}]\,f_s(g(\mathbf t))
=\pi_{\tilde E}^{d}\bigl(E_{\gamma_1,j_1}^{(r_1)}\cdots E_{\gamma_s,j_s}^{(r_s)}f_s\bigr)(1)
=u\pi_{\tilde E}^{d},
\end{equation}
for some unit $u\in \tilde O^\times$: iterating \eqref{path-divided-power} carries $f_s$ to a unit multiple of $f_0$, whose restriction to $U_v$ is the unit $P_{v,v}\in\tilde O^\times$.
The integral big cell \cite[Theorem~8.2.5]{KP23} is open
at the identity, so we get a factorization
\begin{equation}
\label{formal-Iwahori-Gauss}
g(\mathbf t)=u(\mathbf t)c(\mathbf t),
\quad
u(\mathbf t)\in\calU_{b,v}(\tilde O\pot{\mathbf t}),
\quad
c(\mathbf t)\in\calB^-_{b,v}(\tilde O\pot{\mathbf t}).
\end{equation}
By \eqref{Iwahori-minor-expansion recall} and $X_{\alpha,n+m_\alpha}=\pi^{-1}X_{\alpha,n}$, we have
\begin{equation}
\label{Iwahori-Gauss-coordinate-divisibility}
X_{\alpha,j}(u(\mathbf t))
=\pi_{\tilde E}^{e_\alpha(n_b(\alpha,j)-j)}W_{\alpha,j}(\mathbf t),
\qquad W_{\alpha,j}(\mathbf t):=X_{\alpha,n_b(\alpha,j)}(u(\mathbf t))\in\tilde O\pot{\mathbf t},
\end{equation}
for all $(\alpha,j)\in v\bar\Psi^+$.  Right
equivariance and $f_s|_{U_v}=P_{v,w}$ yield
\begin{equation}
\label{path-pullback-big-cell}
f_s(g(\mathbf t))
=(v\chi)(c(\mathbf t))^{-1}P_{v,w}(u(\mathbf t)),
\end{equation}
where the first factor on the right is a unit with constant term $1$.

Substituting \eqref{Iwahori-Gauss-coordinate-divisibility} into
\eqref{Iwahori-minor-expansion} and using \Cref{rem:integrality-of-coefficients} gives
\[
P_{v,w}(u(\mathbf t))
=\kappa\sum_{\mathbf n}\bar u_{\mathbf n}\,\pi_{\tilde E}^{\deg_{v,w}^\chi(\mathbf n)}
\prod_{(\alpha,j)\in v\bar\Psi^+}W_{\alpha,j}(\mathbf t)^{n_{\alpha,j}}.
\]
If $\bar u_{\mathbf n}\neq0$, then
\Cref{prop:characteristic-zero-extremal-degree} gives $\deg_{v,w}^\chi(\mathbf n)\geq d$.
Let
\[
\fraka_d:=
\bigl(\bar u_{\mathbf n}\mid\deg_{v,w}^\chi(\mathbf n)=d\bigr)
\subseteq\tilde O.
\]
Since the specialized coefficients have nonnegative valuation, we conclude
\[
P_{v,w}(u(\mathbf t))\in
\pi_{\tilde E}^{d}\fraka_d\,\tilde O\pot{\mathbf t}+\pi_{\tilde E}^{d+1}\tilde O\pot{\mathbf t},
\]
and multiplying by the unit $(v\chi)(c(\mathbf t))^{-1}$ does not change this.
Comparing with \eqref{primitive-path-coefficient},
we deduce $\fraka_d=\tilde O$. Since $\tilde O$ is local, some $\bar u_{\mathbf n}$ of degree $d$ is a unit.
\end{proof}

We can now link the quantity $d_\chi(v,w)$ to the quantum Bruhat graph.

\begin{proposition}
\label{prop:Iwahori-extremal-minor-formula}
For every fundamental weight $\chi$ and every $v,w\in W_0$,
one has
$$
d_\chi(v,w)
=
\left\langle N\chi,
\operatorname{wt}(v\Rightarrow w)\right\rangle.
$$\end{proposition}

\begin{proof}
Put $d:=\langle N\chi,\operatorname{wt}(v\Rightarrow w)\rangle$.
If $\bar u_{\mathbf n}\neq0$, then by
\Cref{prop:characteristic-zero-extremal-degree} it satisfies $\deg_{v,w}^\chi(\mathbf n)\geq d$,
whence $d_\chi(v,w)\geq d$.
By \Cref{lem:Iwahori-initial-form-primitive} some $\bar u_{\mathbf n}$ with
$\deg_{v,w}^\chi(\mathbf n)=d$ is a unit in $\tilde O$, and therefore $d_\chi(v,w)\leq d$.
This yields the claim.
\end{proof}

\subsection{Passage from Iwahori level to maximal parahoric level}\label{proof for simple simply connected}
We first make the notation for the vertex calculation explicit.  Write
$\Delta=\{\alpha_1,\ldots,\alpha_r\}$ and
$\theta=\sum_{j=1}^r m_j\alpha_j$ for the simple roots and the highest
root of $\Sigma$.  The vertices of $\mathbf a$ are
\[
a_0=0,
\qquad
a_i=\frac{\omega_i^\vee}{m_i}\quad(1\leq i\leq r),
\]
where $\omega_i^\vee$ are the fundamental coweights of $\Sigma$.
We index the affine nodes by $\widetilde\Delta=\{0,1,\ldots,r\}$ and put
$K_i:=\widetilde\Delta\setminus\{i\}=K_{a_i}$, so that $K_i$ is the type of
$\calG_{a_i}$. We denote by $W_{K}\subset W$ the subgroup generated by the
reflections in the affine walls indexed by a subset $K\subset\widetilde\Delta$.
The involution $i\mapsto i^\star$ on $\widetilde\Delta$ is the one determined by
$a_{i^\star}=-w_0a_i$, as in \eqref{lrvw-semi-affine-weight}.  Throughout we write
elements of $W$ as $\nu_I(\pi)u$ with $\nu_I\in X_\ast(T)_I$ and $u\in W_0$: the translation to the conventions of
\cite{Schremmer24} is recorded in \Cref{rem:schremmer-comparison}.  Finally, for
$v\in W_0$ and an extremal weight $\eta$, we abbreviate
\[
\tilde\calV_{a_i}^{B_v}(\eta):=\tilde\calV_{a_i,\chi,v}(\eta),
\qquad
h_{a_i}^{B_v}(\eta):=h_{a_i,\chi,v}(\eta).
\]

Fix an affine node $i$ and $v,w\in W_0$, and put $K=K_{i^\star}$.  Schremmer's
semi-affine projection \cite[Lemma~3.24]{Schremmer24} attaches to every $u\in W_0$ a
pair ${}^K\!\pi(u)=(u',\mu)$ consisting of a representative
$u'\in{}^K\!W_0$ and a coweight $\mu\in\mathbb Z\Sigma^\vee$.  We write
\begin{equation}\label{schremmer projection identities}
{}^K\!\pi(w_0w)=(w_i,\mu_w),
\qquad
{}^K\!\pi(w_0v)=(v_i,\mu_v),
\end{equation}
and abbreviate $r_v:=w_0v_i$ and $r_w:=w_0w_i$.
We now compare the corresponding extremal lattices via these quantities.

\begin{proposition}
\label{prop:semi-affine-projection-extremal-lattices}
With the preceding notation, one has
\begin{equation}
\label{parahoric-Iwahori-extremal-reduction}
-h_{a_i}^{B_v}(w\chi)
=d_\chi(w_0v_i,w_0w_i)
-\langle N\chi,\mu_w\rangle
+\langle N\chi,\mu_v\rangle.
\end{equation}
\end{proposition}

\begin{proof}
For $s\in\{v,w\}$ we write $s_i$ and $\mu_s$ for the corresponding data in
\eqref{schremmer projection identities}.  By the proof of
\cite[Lemma~3.24]{Schremmer24}, transported to our conventions as in \Cref{rem:schremmer-comparison}, the element
\begin{equation}
\label{semi-affine-normalizer-element}
z_s:=(s\mu_s)(\pi)\cdot s(w_0s_i)^{-1}
\end{equation}
lies in $W_{K_i}\subset\operatorname{Stab}_W(a_i)$.
Its finite Weyl group-part is
$
u_s=s(w_0s_i)^{-1}=sr_s^{-1},
$
and hence
$
u_s(r_s\chi)=s\chi$. Similarly, we get
$u_sB_{r_s}u_s^{-1}=B_s.
$

We write $R_zf(g)=f(gz)$ for the right translation of induced functions and
$(z\cdot f)(g)=f(z^{-1}g)$ for the left action, and we fix a representative
$\dot z_s\in\calG_{a_i}(\breve O)$ of $z_s$ normalizing $\calT$.
Since $\dot z_vB_{r_v}^-\dot z_v^{-1}=B_v^-$, it induces an
isomorphism
\[
R_{\dot z_v}\colon
\Ind_{\calB_{r_v,a_i,\tilde O}^-}^{\calG_{a_i,\tilde O}}(r_v\chi)
\xrightarrow{\ \sim\ }
\Ind_{\calB_{v,a_i,\tilde O}^-}^{\calG_{a_i,\tilde O}}(v\chi).
\]
It commutes with the left $\calG_{a_i,\tilde O}$-action and therefore preserves the
$T$-weight spaces.  Then, left action by
$\dot z_w\in\calG_{a_i}(\breve O)$
carries the $r_w\chi$-weight lattice to the $w\chi$-one.  Thus
$
T_{v,w}:=\dot z_w\circ R_{\dot z_v}
$
identifies the integral lines
\begin{equation}
\label{semi-affine-integral-line-isomorphism}
\tilde\calV_{a_i,\chi, r_v}(r_w\chi)
\xrightarrow{\ \sim\ }
\tilde\calV_{a_i,\chi,v}(w\chi).
\end{equation}

It remains to compare this isomorphism with the lattices
$\tilde\calV_{o,\chi, r_v}(r_w\chi)$ and
$\tilde\calV_{o,\chi,v}(w\chi)$.  By \eqref{semi-affine-normalizer-element} and
$u_v=vr_v^{-1}$ we have
\[
z_v=(v\mu_v)(\pi)\,u_v=u_v\cdot(r_v\mu_v)(\pi),
\]
so we may take $\dot z_v=\dot u_v\cdot (r_v\mu_v)(\pi)$.  For an induced function,
right equivariance $f(gt)=(r_v\chi)(t)^{-1}f(g)$ for $t\in T(\breve F)$
and the Kottwitz isomorphism, as in \eqref{weight module properties computation},
show that $R_{\dot z_v}$ contributes the scalar
$
\pi_{\tilde E}^{-\langle N\chi,\mu_v\rangle}
$
up to a unit.
Choose a representative
$\dot u_v\in\calG_o(\breve O)$.  Its
right translation identifies the two lattices $\tilde\calV_{o,\chi, r_v}(r_w\chi)$ and
$\tilde\calV_{o,\chi,v}(r_w\chi)$.
Similarly, the left action of $\dot z_w$, lifting
$z_w=(w\mu_w)(\pi)\,u_w$, contributes
$
\pi_{\tilde E}^{\langle N\chi,\mu_w\rangle}
$
up to a unit. Thus,
\begin{equation}
\label{semi-affine-reference-scalar}
T_{v,w}(P_{r_v,r_w})
=\omega\pi_{\tilde E}^{\langle N\chi,\mu_w\rangle-
           \langle N\chi,\mu_v\rangle}P_{v,w},
\qquad \omega \in \tilde O^\times.
\end{equation}

Write
\[
h'=h_{a_i}^{B_{r_v}}(r_w\chi),
\qquad
h=h_{a_i}^{B_v}(w\chi).
\]
Comparing \eqref{semi-affine-integral-line-isomorphism} with
\eqref{semi-affine-reference-scalar} yields
\begin{equation}
\label{semi-affine-h-shift}
h=h'+\langle N\chi,\mu_w\rangle
       -\langle N\chi,\mu_v\rangle.
\end{equation}
It remains to identify $h'$ with the Iwahori degree. Let $b$ be an interior point of $\mathbf a$.
For every gradient $\gamma$ of an affine root $\psi\in\Psi$ that vanishes at $a_i$ and is positive on $\mathbf a$, the defining sign condition for $v_i,w_i\in{}^K\!W_0$ gives $r_v^{-1}\gamma,r_w^{-1}\gamma\in\Phi^-$.
Thus the root groups $U_{-\gamma}$ lie in $U_{r_v}$ and fix the extremal vector $f_{r_v,r_w}^\chi$.
These are precisely the root directions in which the models can enlarge from $b$ to $a_i$.
Order the corresponding nondivisible root-group factors first in the root-group decomposition of $U_{r_v}$.
Left invariance makes $P_{r_v,r_w}$ independent of their coordinates, and the remaining root-group models agree at $b$ and $a_i$.
By the saturation statement in \Cref{embedding corollary}, the $r_w\chi$-extremal lattices at $b$ and $a_i$ therefore agree.
Therefore $h'=h_{b,\chi,r_v}(r_w\chi)=-d_\chi(r_v,r_w)$.
Substitution into \eqref{semi-affine-h-shift} proves
\eqref{parahoric-Iwahori-extremal-reduction}.
\end{proof}

\begin{proof}[Proof of \Cref{prop:extremal-minor-reduction}]
By \Cref{weight module monoidal} and \Cref{reduction to simple simply connected}, we may assume that $G_{\breve F}$ is absolutely almost simple
and simply connected, that $\tilde E$ is its splitting field, and that $\chi$ is a fundamental weight.
Write $a=a_i$.
Substituting \Cref{prop:Iwahori-extremal-minor-formula} and the identity
$\operatorname{wt}(w_0v_i\Rightarrow w_0w_i)=\operatorname{wt}(w_i\Rightarrow v_i)$
of \cite[Lemma~3.8]{Schremmer24} into
\eqref{parahoric-Iwahori-extremal-reduction} gives
\begin{equation}\label{h formula weight plus translation}
-h_{a_i,\chi, v}(w\chi)
=
\left\langle N\chi,
\operatorname{wt}(w_i\Rightarrow v_i)-\mu_w+\mu_v
\right\rangle.
\end{equation}
By \cite[Definition~3.26]{Schremmer24} and the projection
identities \eqref{schremmer projection identities}, the expression in the right side of \eqref{h formula weight plus translation} is
\[
\langle N\chi, \operatorname{wt}_{K_{i^\star}}(w_0w\Rightarrow w_0v)\rangle
=
\langle N\chi,\operatorname{wt}_{a_i}(v\Rightarrow w)\rangle,
\]
by \eqref{lrvw-semi-affine-weight}.  This proves the formula for fundamental weights of absolutely almost simple, simply connected groups, and hence in general.
\end{proof}

\begin{appendix}

\section{Minimal test sets}\label{appendix minimal test sets}
We consider an irreducible, based root system in an $n$-dimensional $\mathbb R$-vector space and use the labeling convention of \cite{Bou68}.
In particular, we denote by $\al_1^{(\vee)},\ldots,\al_n^{(\vee)}$ the simple (co)roots, and by $\omega_1^{(\vee)},\ldots,\omega_n^{(\vee)}$ the fundamental (co)weights. 
Let $A=(\langle \alpha_i,\alpha_j^\vee\rangle)_{ij}\in \GL_n(\mathbb R)$ be the Cartan matrix.
Its inverse matrix is given by $A^{-1}=(\langle \omega_i,\omega_j^\vee\rangle)_{ij}$ because $\langle \al_i,\omega_j^\vee\rangle=\delta_{ij}$.

Let $P^\vee$ be the coweight lattice and denote by $P^\vee_{\mathbb Q}$ its rationalization. 
For dominant rational coweights $\nu,\mu\in (P^\vee_{\mathbb Q})^+$, we write $\nu\leq_{\mathbb Q}\mu$ if and only if $\mu-\nu\in \mathbb Q_{\geq 0} \{\al_1^\vee,\ldots,\al_n^\vee\}$ if and only if $\langle\omega_i,\mu-\nu\rangle\geq 0$ for all $i=1,\ldots,n$.

\begin{definition}\label{minimal test set definition}
Let $\mu\in (P^\vee_{\mathbb Q})^+$ be nonzero.
A subset $J\subset\{\omega_1,\ldots,\omega_n\}$ is called a \emph{test set for $\mu$}, if for every dominant rational coweight $\nu\in (P^\vee_{\mathbb Q})^+$, the inequalities
$$
\langle\chi,\mu-\nu\rangle\geq0
\qquad
\text{for all }\chi\in J
$$
imply $\langle\omega_i,\mu-\nu\rangle\geq0$ for all $i=1,\ldots,n$.
A test set for $\mu$ is called \emph{minimal} if none of its proper subsets is a test set for $\mu$.
\end{definition}

\begin{lemma}\label{minimal test set}
Every nonzero dominant coweight $\mu$ admits a unique minimal test set, denoted by $J(\mu)$.
Moreover, $J(\mu)$ depends only on the support
$$
\operatorname{supp}(\mu)
:=
\{\,i\in\{1,\ldots,n\}\mid
\langle\al_i,\mu\rangle>0\,\}.
$$
\end{lemma}
\begin{proof}
Identify $(P^\vee_{\mathbb Q})^+=\mathbb Q_{\geq 0}^n$ using the monoid basis $\omega_1^\vee,\ldots, \omega_n^\vee$, i.e., the dominant rational coweight $\nu$ corresponds to the vector $ \nu=(\nu_1,\ldots,\nu_n)$ with $\nu_i=\langle\al_i,\nu\rangle\geq 0$.
Then, the matrix vector multiplication gives 
\begin{equation}\label{inverse cartan multiplication}
A^{-1} \cdot \nu= \left(\langle\omega_1,\nu\rangle,\ldots,\langle\omega_n,\nu\rangle\right)^t. 
\end{equation}
Since the root system is irreducible, all entries of $A^{-1}$ are strictly positive \cite[Theorem~4.3 (Fin)]{Kac90}, and so are all entries of $A^{-1} \cdot \nu$ when $\nu\neq 0$.
We apply this as follows. 

For every $i=1,\ldots,n$, the inequality $\langle\omega_i,\mu-\nu\rangle\geq 0$ is equivalent to the inequality $r_i\cdot \nu \leq b_i$, where $r_i:=(A^{-1})_{i\bullet}$ is the $i$-th row and $b_i:=\langle \omega_i,\mu\rangle >0$.  
For a subset $J\subset\{\omega_1,\ldots,\omega_n\}$, we define the real polyhedron
$$
\mathscr P_{\mu,J}
:=
\left\{
x\in\mathbb R_{\geq0}^n
\ \middle|\
r_j\cdot  x\leq b_j
\text{ for every }\omega_j\in J
\right\},
$$
and abbreviate $\mathscr P_\mu:=\mathscr P_{\mu,\{\omega_1,\ldots,\omega_n\}}$. 
The latter is full dimensional because every $b_i$ is strictly positive.
Then, $J$ is a test set for $\mu$ if and only if $\mathscr P_{\mu,J}=\mathscr P_\mu$: indeed, if these real polyhedra were different, their difference would contain a rational point.

We claim that the minimal test set consists precisely of the fundamental weights whose inequalities define facets of $\mathscr P_\mu$. 
In particular, the minimal test set is unique.

Now, let $J\subset \{\omega_1,\ldots,\omega_n\}$ be any minimal test set. 
We first show that every inequality indexed by $J$ defines a facet of
$\mathscr P_\mu$. Fix $\omega_i\in J$. 
By the minimality of $J$, there exists $y\in
\mathscr P_{\mu,J\backslash\{\omega_i\}}
\backslash\mathscr P_\mu$.
Since $J$ is a test set, $y$ must violate the omitted
inequality:
$
r_i\cdot y>b_i
$.
Choose an interior point $x_0\in\mathscr P_\mu$, and let $z$ be the unique point of the segment joining
$ x_0$ to $y$ for which
$
r_i\cdot z=b_i$.
Because $x_0$ satisfies all inequalities strictly and
$y$ satisfies those indexed by
$J\backslash\{\omega_i\}$, one has
$
r_j\cdot z<b_j
$
for every $\omega_j\in J\backslash \{\omega_i\}$ and every coordinate of $z$ is positive. 
Hence, a relative neighborhood of $ z$ in the hyperplane $r_i\cdot x=b_i$
is contained in $\mathscr P_{\mu,J}=\mathscr P_\mu$. Therefore
$r_i\cdot x\leq b_i$ defines a facet of $\mathscr P_\mu$.

Conversely, suppose that $r_i\cdot x\leq b_i$ defines a facet of
$\mathscr P_\mu$, and let $J'$ be any test set for $\mu$.
We claim that $\omega_i\in J'$. Otherwise, choose a point $z$
in the relative interior of this facet that lies on none of the other
hyperplanes
$
r_j \cdot x=b_j
$
for $j\neq i$
and on none of the coordinate hyperplanes. Such a point exists because
the rows of $A^{-1}$ are pairwise non-proportional, and $b_i>0$ ensures that this supporting hyperplane is not a coordinate hyperplane. Thus
$
 z\in\mathbb R_{>0}^n$ and 
$
r_j\cdot z<b_j
$
for every $\omega_j\in J'$.
A sufficiently small displacement of $ z$ across the hyperplane
$r_i\cdot x=b_i$ would then remain in $\mathscr P_{\mu,J'}$ but leave
$\mathscr P_\mu$, contradicting
$
\mathscr P_{\mu,J'}=\mathscr P_\mu
$.
Thus every rational test set contains every fundamental weight whose
inequality defines a facet of $\mathscr P_\mu$.
This proves the claim. 

It remains to prove the dependence on the support. Let
$J\subset\{\omega_1,\ldots,\omega_n\}$ be arbitrary and fix
$i\in\{1,\ldots,n\}$. By the affine form of Farkas' lemma
\cite[Corollary~7.1h]{Schrijver:LinearIntegerProgramming}, the
inequalities
$$
r_j\cdot \nu\leq r_j\cdot \mu
\quad
(\omega_j\in J),
\qquad
\nu\geq0,
$$
imply
$
r_i\cdot \nu\leq r_i\cdot \mu
$
if and only if there exist numbers
$\lambda_j\in\mathbb R_{\geq0}$ such that
$
\sum_{\omega_j\in J}\lambda_jr_j\geq r_i
$
componentwise and
$
\sum_{\omega_j\in J}\lambda_jr_j\cdot\mu
\leq r_i\cdot\mu
$.
Set
$$
d:=\sum_{\omega_j\in J}\lambda_jr_j-r_i.
$$
The first condition says that $d\geq0$ componentwise, while the second
is equivalent to
$
d\cdot\mu\leq0$.
Since both $d$ and $\mu$ have nonnegative entries, this holds if
and only if
$
d_l=0
$
for every $l\in\operatorname{supp}(\mu)$.
Thus, the implication of the $i$-th inequality by those indexed by $J$
depends on $\mu$ only through $\operatorname{supp}(\mu)$. 
This shows the lemma. 
\end{proof}

\begin{remark}\label{minuscule classification}
	In \cite{LRVW:Minuscule}, we classify the minimal test sets for minuscule coweights.
	More precisely, if $\mu$ is a nonzero minuscule dominant coweight, then $\mu=\omega_i^\vee$ is a minuscule fundamental coweight, and we show that $J(\mu)$ belongs to the following list:
			\begin{itemize}
				\item Type $A_n$, $n\geq 1$: $J(\omega_1^\vee)=\{\omega_n\}$, $J(\omega_i^\vee)=\{\omega_1,\omega_n\}$ for $i=2,\ldots,n-1$, $J(\omega_n^\vee)=\{\omega_1\}$
				\item Type $B_n$, $n\geq 3$:
					$
						J(\omega_1^\vee)=\{\omega_n\}
					$
				\item Type $C_n$, $n\geq 2$:
					$
						J(\omega_n^\vee)=\{\omega_1\}
					$
				\item Type $D_n$, $n\geq 4$: $J(\omega_1^\vee)=\{\omega_{n-1},\omega_n\}$, $J(\omega_{n-1}^\vee)=\{\omega_{1},\omega_n\}$,  $J(\omega_n^\vee)=\{\omega_{1},\omega_{n-1}\}$
					
				\item Type $E_6$: $J(\omega_1^\vee)=\{\omega_2,\omega_6\}$, $J(\omega_{6}^\vee)=\{\omega_{1},\omega_2\}$
				\item Type $E_7$: $J(\omega_7^\vee)=\{\omega_1,\omega_2\}$
			\end{itemize}
	A posteriori, when the absolute Dynkin diagram has at least two nodes, the set \(J(\omega_i^\vee)\) can be recovered directly from the diagram as follows: 
start at the $i$-th vertex of the Dynkin diagram and include the fundamental weights corresponding to all other terminal vertices.  

In the classical types, all fundamental weights in $J(\omega_i^\vee)$ are minuscule but not so in the exceptional types.
The appearance of non-minuscule fundamental weights is closely related to the difference between $\mu$-admissibility and $\mu$-permissibility in exceptional types \cite{Sauermann}.
\end{remark}

\section{Root vector filtrations}
\label{appendix feigin makedonskyi}

The purpose of this appendix is to prove
\Cref{prop:characteristic-zero-extremal-degree}, which identifies the minimal
degree occurring in the expansion \eqref{Iwahori-minor-expansion} with a quantum
Bruhat weight.  The proof relies on the graded character formula for generalized
Weyl modules of Feigin--Makedonskyi \cite{FeiginMakedonskyi17} and its analogue for twisted current algebras \cite{FeiginMakedonskyiTwisted17}, together with the
shellability of the quantum Bruhat graph in the form given by
Lenart--Naito--Sagaki--Schilling--Shimozono \cite{LenartEtAl17}.  Since the
conventions of \cite{FeiginMakedonskyi17,FeiginMakedonskyiTwisted17} differ from ours in ways that are easy
to conflate, we begin by fixing a dictionary.

\subsection{Conventions and dictionary}
\label{FM dictionary}

We work over $\bbC$ via the embedding $\Lambda\subset\bbC$ fixed in \Cref{split coordinates sec}, and otherwise retain the hypotheses and notation of \Cref{iwahori calc sec}.
Thus $\frakg=\Lie(H\otimes_\Lambda\bbC)$ is simple.
Let $\tau_\Gamma$ be the pinned automorphism corresponding to the chosen Galois generator $\sigma$ on the based root datum, using the compatible Chevalley system of \Cref{split coordinates sec}, and put $\tau:=\tau_\Gamma^{-1}$.

For $k\in\bbZ/e\bbZ$, let $\frakg_k$ be the eigenspace on which $\tau$ acts by $\zeta_e^k$.
Then $\frakg=\bigoplus_k\frakg_k$.
We check that the vectors $E_{\alpha,k}$ of \eqref{Fourier vectors} and \eqref{Heisenberg Fourier} are eigenvectors of $\tau$, namely
\begin{equation}\label{tau eigenvectors}
\tau(E_{\alpha,k})=\zeta_e^{ke_\alpha}E_{\alpha,k},
\quad\text{that is,}\quad
E_{\alpha,k}\in\frakg_{ke_\alpha},
\qquad
(\alpha,k)\in v\bar\Psi^+,
\end{equation}
for every $v\in W_0$.
By the choice of the Chevalley--Steinberg system \cite[4.1.3]{BT84} and of the isomorphism $H_{\tilde E}\simeq G_{\tilde E}$ in \Cref{split coordinates sec}, the action of $\sigma$ on $G(\tilde E)$ is $\tau_\Gamma$ composed with the action of $\sigma$ on coefficients, and $\tau_\Gamma(E_\gamma)=E_{\sigma\gamma}$ for every absolute root $\gamma$ whose restriction is not divisible.
Let first $\alpha$ be nondivisible and not multipliable.
Then $\tau(E_{\sigma^r\tilde\alpha})=E_{\sigma^{r-1}\tilde\alpha}$, and substituting $r+1$ for $r$ in the sum \eqref{Fourier vectors} gives
\[
\tau(E_{\alpha,k})
=\sum_{r=0}^{m_\alpha-1}\zeta_e^{re_\alpha k}E_{\sigma^{r-1}\tilde\alpha}
=\zeta_e^{e_\alpha k}E_{\alpha,k}.
\]
Let now $\alpha$ be multipliable, so that $e=2$, $\zeta_e=-1$, $e_\alpha=1$, and $\tau=\tau_\Gamma$ exchanges $E_{\tilde\alpha}$ and $E_{\tilde\alpha'}$.
As $\tau$ is a Lie algebra automorphism and $[E_{\tilde\alpha},E_{\tilde\alpha'}]=E_{\tilde\alpha+\tilde\alpha'}$, we get
\[
\tau(E_{\alpha,0})=E_{\tilde\alpha'}+E_{\tilde\alpha}=E_{\alpha,0},
\qquad
\tau(E_{\alpha,1})=E_{\tilde\alpha'}-E_{\tilde\alpha}=-E_{\alpha,1},
\qquad
\tau(E_{2\alpha,1})=[E_{\tilde\alpha'},E_{\tilde\alpha}]=-E_{2\alpha,1}.
\]
This proves \eqref{tau eigenvectors}, matching $ke_\alpha$ with the degree of $E_{\alpha,k}$ in the twisted Iwahori current algebra, see \Cref{lem:minor-degree-one-factor-filtration}.

The fixed algebra $\frakg_0$ is simple with Cartan subalgebra $\frakh_0=\frakh^\tau$, its root system is the set $\Phi_0:=\Phi_{\mathrm{nd}}$ of nondivisible relative roots, with simple roots $\Delta$, and its Weyl group is $W_0$.
Here $\Phi_0=\Phi$, unless $\frakg$ is of type $A_{2n}$ and $e=2$, i.e., of twisted affine type $A_{2n}^{(2)}$, in which case $\Phi$ is of type $BC_n$ and $\frakg_0$ is of type $B_n$ \cite[\S2.2]{FeiginMakedonskyiTwisted17}.
The inclusion $W_0\subset W_H$ of the relative Weyl group into the absolute one matches longest elements $w_0$.
We write $P_0$ for the weight lattice of $\Phi_0$, $(-)^\star=-w_0$ for
the induced involution on $P_0$ and on $\Delta$, and
$\varpi_1,\ldots,\varpi_r$ for the fundamental weights in $P_0$.
For an absolute weight $\chi$, we write $\chi_0=\chi|_{\frakh_0}\in P_0$.
The coroot of $\alpha\in\Phi_0$ in $\frakh_0$ is the sum of the absolute coroots restricting to it, so restriction along $\frakh_0\to \frakg$ preserves fundamental weights, with one exception: in type $A_{2n}^{(2)}$, the coroot of a short root $\alpha$ of $B_n$ is $2(\tilde\alpha^\vee+\tilde\alpha'^\vee)$, as $\langle\tilde\alpha',\tilde\alpha^\vee\rangle=-1$, so that the absolute fundamental weights at the nodes $n,n+1$ of $A_{2n}$ restrict to $2\varpi_n$, compare \cite[Lemma~2.6]{FeiginMakedonskyiTwisted17}.
We write $\chi_0=c_\chi\varpi_i$ for an absolute fundamental weight $\chi$, where $c_\chi\in\{1,2\}$ and $c_\chi=2$ only in this exceptional case.
If $e=1$, then $\frakg_0=\frakg$, $\chi_0=\chi$, and $\Phi_0=\Sigma$.

For $\alpha\in\Phi_0$ we write $\operatorname{len}(\alpha)$ for its normalized length as in \cite[\S1.2, Figure~1]{FeiginMakedonskyiTwisted17}, that is, $1$ for short roots and $e$ for long roots, respectively $2$ for long roots in type $A_{2n}^{(2)}$.
Outside type $A_{2n}^{(2)}$, one has $\operatorname{len}(\alpha)=e_\alpha$.
We put $c:=1/2$ in type $A_{2n}^{(2)}$ and $c:=1$ otherwise, and write $\alpha_\Sigma\in\Sigma$ for the positive multiple of $\alpha$ in $\Sigma$, that is, $m_\alpha\alpha$, respectively $2m_\alpha\alpha$ if $\alpha$ is multipliable.
Then
\begin{equation}\label{Sigma coroot}
\alpha_\Sigma^\vee=\frac{c\operatorname{len}(\alpha)}{e}\,\alpha^\vee,
\end{equation}
as $(m_\alpha\alpha)^\vee=\alpha^\vee/m_\alpha$, and in type $A_{2n}^{(2)}$ the coroot of $4\alpha$ for short $\alpha$ is the image of $\tilde\alpha^\vee$, i.e., $\alpha^\vee/4$ in $\frakh_0$, see \Cref{sec:qbg-elementary-divisors}.
Moreover, the Galois group acts on $X^*(T)$ through the cyclic group generated by $\tau$, and the coroot $\alpha^\vee$ of $\alpha\in\Phi_0$ is $\tau$-invariant, so that $\langle N\chi,\alpha^\vee\rangle=e\langle\chi_0,\alpha^\vee\rangle$.
By \eqref{Sigma coroot}, we get
\begin{equation}\label{twisted pairing}
\langle N\chi,\alpha_\Sigma^\vee\rangle=c\operatorname{len}(\alpha)\,\langle\chi_0,\alpha^\vee\rangle,
\qquad \alpha\in\Phi_0.
\end{equation}
In particular, if $\chi_0=c_\chi\varpi_i$, then $N\chi$ is a positive multiple of the fundamental weight of $\Sigma$ at the node $i$.

The results of \cite{FeiginMakedonskyi17} also concern the Lie algebra
$\frakg$, and if $e=1$, then the quantum Bruhat graph of \cite[\S1.1]{FeiginMakedonskyi17} has vertex set
  $W_0$, edge labels in $\Phi_0^\vee$, and weights in $\mathbb Z\Phi_0^\vee$, so it agrees with our normalization \eqref{quantum weight}.
If $e>1$, then the twisted quantum Bruhat graph of \cite[\S1.1]{FeiginMakedonskyiTwisted17} has vertex set $W_0$, edge labels in $\Phi_0$, and is the quantum Bruhat graph of the dual root system $\Phi_0^\vee$, as remarked in loc.\ cit.; in type $A_{2n}^{(2)}$, one uses the twisted quantum Bruhat graph of type $D_{n+1}^{(2)}$, i.e., of $B_n$ \cite[\S1.2]{FeiginMakedonskyiTwisted17}.
Identifying $\frakh_0\simeq\frakh_0^*$ via the invariant form for which $\langle\alpha,\alpha\rangle=2\operatorname{len}(\alpha)$, the coroot $\alpha^\vee$ corresponds to $\alpha/\operatorname{len}(\alpha)$, so \eqref{Sigma coroot} shows that the root datum $(\Sigma,\Sigma^\vee)$ is obtained from $(\Phi_0^\vee,\Phi_0)$ by rescaling roots by $e/c$ and coroots by $c/e$.
So identifying the label $\alpha$ with $\alpha_\Sigma$, this is the quantum Bruhat graph of $\Sigma$ used in \Cref{sec:qbg-elementary-divisors}.

However, their extended affine Weyl group $W^e=W_0\ltimes P_0$ has translations $t_\lambda$ indexed
 by \emph{weights} $\lambda \in P_0$.
If $e=1$, the corresponding hyperplanes are given by the vanishing of affine \emph{coroots}
  $\alpha^\vee+N\delta$, for which we write $\operatorname{Re}(\alpha^\vee+N\delta)=\alpha^\vee$ and $\deg(\alpha^\vee+N\delta)=N$.
If $e>1$, they are given by the real affine roots $\alpha+N\delta$ of the affine Lie algebra, where $N\in\operatorname{len}(\alpha)\bbZ$ \cite[\S2.1]{FeiginMakedonskyiTwisted17}, and we write $\operatorname{Re}(\alpha+N\delta)=\alpha$ and $\deg(\alpha+N\delta)=N$.
In type $A_{2n}^{(2)}$, these are the real affine roots of type $D_{n+1}^{(2)}$, which are auxiliary: the real affine roots of $A_{2n}^{(2)}$ itself are $\alpha+k\delta$ for $\alpha\in\Phi_0$ and $2\alpha+(2k+1)\delta$ for short $\alpha$, with $k\in\bbZ$, where $\delta$ has current degree $1$ \cite[\S\S2.1, 2.2]{FeiginMakedonskyiTwisted17}.
We write $\widehat{\alpha}_D$ for $-\alpha^\vee+D\delta$, respectively $-\alpha+\operatorname{len}(\alpha)D\delta$, so that $\deg\widehat\alpha_D=\operatorname{len}(\alpha)D$.
Thus $W^e$ is \emph{not} the Iwahori--Weyl group $W=X_*(T)_I\rtimes W_0$ fixed in
  \Cref{sec:qbg-elementary-divisors}.

For $x\in W^e$ with a fixed reduced decomposition, we write
$\beta_1(x),\ldots,\beta_l(x)$ for the associated sequence of affine coroots, respectively affine roots,
\cite[(1.1)]{FeiginMakedonskyi17}, \cite[\S1.2]{FeiginMakedonskyiTwisted17}.  For $u\in W_0$, an element of
$\mathcal{QB}(u,x)$ is a \emph{quantum alcove path}: the datum of a subsequence
$j_1<\cdots<j_s$ of $\{1,\ldots,l\}$, giving rise to
\[
z_0=ux,
\qquad
z_k=z_{k-1}s_{\beta_{j_k}(x)}
\qquad (1\le k\le s),
\]
subject to the requirement that the projection
$\operatorname{dir}(z_0)\to\cdots\to\operatorname{dir}(z_s)$, with labels
$\operatorname{Re}\beta_{j_k}(x)$, be a path in the quantum Bruhat graph
\cite[\S1.2]{FeiginMakedonskyi17}, respectively in the twisted quantum Bruhat graph \cite[\S1.2]{FeiginMakedonskyiTwisted17} (note that when $e>1$ these are called twisted quantum alcove paths and their set is denoted $t\mathcal{QB}(u,x)$ in loc.\ cit.).
In type $A_{2n}^{(2)}$, one requires moreover that no quantum edge comes from a label $\beta_{j_k}(x)$ of odd degree \cite[\S1.2]{FeiginMakedonskyiTwisted17}.  Here
$y=t_{\operatorname{wt}(y)}\operatorname{dir}(y)$ is the canonical decomposition
of $y\in W^e$.  We call $\operatorname{end}(\mathbf p)=z_s$ the endpoint of the
path, $\operatorname{wt}(\operatorname{end}(\mathbf p))\in P_0$ its
weight, and we set
\[\deg(\operatorname{qwt}(\mathbf p))=c\sum\deg(\beta_{j_k}(x)),\] the sum being over
those $k$ for which the corresponding edge of the projection is a quantum edge.
Finally, $W^e$ acts on affine coroots by
$t_\mu(\gamma+N\delta)=\gamma+(N-\langle\mu,\gamma\rangle)\delta$ for $\gamma\in\Phi_0^\vee$
\cite[Lemma~1.14]{FeiginMakedonskyi17}.
If $e>1$, passing to affine coroots $\gamma+\operatorname{len}(\gamma)D\delta\mapsto\gamma^\vee+D\delta$ identifies the real affine roots of the twisted affine algebra with the affine coroots of the untwisted affine algebra of $\frakg_0^\vee$, compatibly with the action of $W^e$ (this is how \cite[Proposition~1.4]{FeiginMakedonskyiTwisted17} is deduced from \cite[Proposition~1.15]{FeiginMakedonskyi17}), so that
$t_\mu(\gamma+\operatorname{len}(\gamma)D\delta)=\gamma+\operatorname{len}(\gamma)(D-\langle\mu,\gamma^\vee\rangle)\delta$.

For an absolute fundamental weight $\chi$ with $\chi_0=c_\chi\varpi_i$ we
set
\begin{equation}\label{rem:two-parabolics}
\lambda=w_0\chi_0=-c_\chi\varpi_{i^\star},
\qquad J=\Delta\setminus\{\alpha_{i^\star}\},
\qquad J^\star=\Delta\setminus\{\alpha_{i}\}.
\end{equation}
Note that we get
\[W_J=\operatorname{Stab}_{W_0}(\lambda)
 \qquad
W_{J^\star}=\operatorname{Stab}_{W_0}(\chi_0),
\]
and correspondingly $w_0W_Jw_0=W_{J^\star}$.

\subsection{The lexicographic chain of a fundamental weight}

We consider the weight $-\lambda=c_\chi\varpi_{i^\star}$ and the set $J=\Delta\setminus\{\alpha_{i^\star}\}$.
Order the Dynkin nodes so that $i_1=i^{\star}$, and let $i_2,\ldots,i_r$ be an arbitrary
ordering of the remaining ones.  For $\alpha\in\Phi_0^{+}$ write
$\alpha^\vee=\sum_k a_{i_k}\alpha^\vee_{i_k}$.
The first coefficient plays a special role, and we write
\[
M_\alpha:=c_\chi a_{i_1}=\langle-\lambda,\alpha^\vee\rangle,
\]
the second equality holding by our choice of node order.

\begin{lemma}
\label{lem:lex-chain-comparison}
\begin{enumerate}
\item For every reduced decomposition of $t_{\lambda}$, the associated sequence of
labels is, as a multiset,
\[
\bigl\{\,\widehat\alpha_D \;:\;
\alpha \in \Phi_0^{+}\setminus \Phi_{0,J}^{+},\;
1\le D\le M_\alpha \,\bigr\}.
\]
\item There is a reduced decomposition of $t_{\lambda}$ for which this sequence is
ordered by the lexicographic order on the vectors
\[
\Bigl(\tfrac{M_\alpha}{D},\ \tfrac{a_{i_2}}{M_\alpha},\ \ldots,\ \tfrac{a_{i_r}}{M_\alpha}\Bigr).
\]
\item If $e=1$, then under the bijection
$\widehat\alpha_D\longleftrightarrow H_{\alpha,-l}$ with
$l=M_\alpha-D$, the order of \emph{(2)} coincides with the
lexicographic $\varpi_{i^\star}$-chain of
\cite[(5.1) and Proposition~5.4]{LenartEtAl17}.
\end{enumerate}
\end{lemma}

\begin{proof}
If $e=1$, parts (1) and (2) are \cite[Proposition~1.18(a) and (b)]{FeiginMakedonskyi17};
the coefficient $a_{i_1}$ appearing there is our $M_\alpha$,
so that (1) reads $1\le D\le M_\alpha$.
If $e>1$ and $c_\chi=1$, they are \cite[Proposition~1.4(a) and (b)]{FeiginMakedonskyiTwisted17}, applied in type $D_{n+1}^{(2)}$ if $\frakg$ is of type $A_{2n}^{(2)}$.
If $c_\chi=2$, then $\varpi_n$ is minuscule for $B_n$, so the labels of $t_{-\varpi_n}$ are the $\widehat\alpha_1$ for $\langle\varpi_n,\alpha^\vee\rangle=1$.
We use the reduced decomposition of $t_{\lambda}=t_{-\varpi_n}t_{-\varpi_n}$ obtained by concatenation, whose labels are, by \cite[Lemma~1.3]{FeiginMakedonskyiTwisted17}, the $\widehat\alpha_2$ followed by the $\widehat\alpha_1$, each in the order of \cite[Proposition~1.4(b)]{FeiginMakedonskyiTwisted17} for $\varpi_n$.
This proves (1), and (2), as $M_\alpha/D\in\{1,2\}$ and the remaining entries of the chain are proportional to those for $\varpi_n$.
Here the first entry of the lexicographic chain is expressed in terms of $D$ as in the untwisted case: by the dictionary of \Cref{FM dictionary}, the label $\widehat\alpha_D$ corresponds to the affine hyperplane $\{\langle x,\alpha^\vee\rangle=D\}$ in $P_0\otimes\mathbb R$, which the segment from $0$ to $-\lambda$ crosses at time $D/M_\alpha$, exactly as for $e=1$.

For (3), the lexicographic $\varpi_{i^\star}$-chain of \cite[(5.1)]{LenartEtAl17} orders the hyperplanes
\[
H_{\alpha,-l},\qquad
\alpha\in\Phi_0^+,\quad
0\le l<\langle\varpi_{i^\star},\alpha^\vee\rangle,
\]
by the lexicographic order on
\[
\frac{1}{M_\alpha}\,(l,\,a_{i_1},\ldots,a_{i_r})
=\Bigl(\tfrac{l}{M_\alpha},\ 1,\ \tfrac{a_{i_2}}{M_\alpha},\ldots,\tfrac{a_{i_r}}{M_\alpha}\Bigr).
\]
The constant second entry may be discarded, and the substitution $l=M_\alpha-D$
turns the first entry into $1-D/M_\alpha$. This is strictly increasing in $M_\alpha/D$,
hence the two lexicographic orders coincide.
\end{proof}

\begin{corollary}
\label{cor:initial-block}
The initial section of the order of \Cref{lem:lex-chain-comparison} consists
exactly of the labels
\[
\widehat\alpha:=\widehat\alpha_{M_\alpha},
\qquad
\alpha\in\Phi_0^+\setminus\Phi^+_{0,J},
\]
one for each such $\alpha$, and the induced order on $\Phi_0^+\setminus\Phi^+_{0,J}$ is
the initial section of the reflection order $<_{\varpi}$ of
\cite[\S6.2]{LenartEtAl17} for the root system $\Sigma$ and its fundamental weight $\varpi$ at the node $i^\star$.
\end{corollary}

\begin{proof}
By \Cref{lem:lex-chain-comparison}(1) one has $D\le M_\alpha$, so the first
lexicographic entry $M_\alpha/D$ is minimal, equal to $1$, precisely when
$D=M_\alpha$ (equivalently $l=0$, by
\Cref{lem:lex-chain-comparison}(3) if $e=1$).

On this initial block the order is the
lexicographic order on $(a_{i_2}/M_\alpha,\ldots,a_{i_r}/M_\alpha)$, which by
\Cref{lem:lex-chain-comparison}(3) is the order of the hyperplanes
$H_{\alpha,0}$ in the lexicographic $\varpi_{i^\star}$-chain if $e=1$.  This is exactly how the
initial section of $<_{\varpi_{i^\star}}$ is defined in \cite[\S6.2]{LenartEtAl17}.
If $e>1$, then by \eqref{Sigma coroot} the coroot $\alpha_\Sigma^\vee$ has the coefficients $a_{i_k}\operatorname{len}(\alpha)/\operatorname{len}(\alpha_{i_k})$ with respect to the simple coroots $(\alpha_{i_k})_\Sigma^\vee$ of $\Sigma$.
Hence passing from $\Phi_0$ to $\Sigma$ multiplies the ratios $a_{i_k}/M_\alpha$ by $\operatorname{len}(\alpha_{i_1})/\operatorname{len}(\alpha_{i_k})$, which does not depend on $\alpha$, and the lexicographic order is unchanged.
\end{proof}

\subsection{The folding lemma}

\begin{lemma}
\label{lem:initial-block-folding}
Let \(\chi\) be an absolute fundamental weight and let \(v,w\in W_0\).  With
$\lambda$ and $J$ as in \eqref{rem:two-parabolics} and $\sigma=ww_0$, there is a path
$\mathbf p\in\mathcal{QB}(\sigma,t_\lambda)$ with
\[
\operatorname{wt}(\operatorname{end}(\mathbf p))=v\chi_0,
\qquad
\deg(\operatorname{qwt}(\mathbf p))
=\left\langle N\chi,\operatorname{wt}(v\Rightarrow w)\right\rangle.
\]
\end{lemma}

\begin{proof}
We start by constructing a finite path in the quantum Bruhat graph of $\Sigma$.
Write $\sigma=\sigma'w_J$ with $\sigma'\in W_0^J$ a
minimal length representative and $w_J\in W_J$, and let $\tau'\in W_0^J$ be the
minimal length representative of the coset $vw_0W_J$.  Applying
\cite[Lemma~6.6(1)]{LenartEtAl17} to the reflection order $<_{\varpi}$ of \Cref{cor:initial-block}
produces a unique path
\begin{equation}
\label{initial-block-finite-path}
\sigma=\sigma_0
\xrightarrow{\alpha_1}\sigma_1
\xrightarrow{\alpha_2}\cdots
\xrightarrow{\alpha_s}\sigma_s=u
\end{equation}
in the quantum Bruhat graph with $u\in\tau' W_J=vw_0W_J$ and with edge labels
$\alpha_k\in\Phi_0^+\setminus\Phi^+_{0,J}$ strictly increasing for
$<_{\varpi}$.  Since $\Phi_0^+\setminus\Phi_{0,J}^+$ is an initial section of
$<_{\varpi}$, this is the unique $<_{\varpi}$-increasing path
from $\sigma$ to $u$, and \cite[Theorem~6.4(2)]{LenartEtAl17} shows that it is
shortest.

Next, we must lift this path to a folding. Fix the reduced decomposition of
$t_\lambda=t_{-c_\chi\varpi_{i^\star}}$ provided by
\Cref{lem:lex-chain-comparison}(2) and write $\widehat\alpha_k$ for the label $\widehat\alpha$ of \Cref{cor:initial-block} attached to $\alpha=\alpha_k$.
By \Cref{cor:initial-block} these labels all belong to the initial block
of the resulting chain, and the order induced on them is precisely
$<_{\varpi}$.  Hence the indices $j_1<\cdots<j_s$ of
$\widehat\alpha_1,\ldots,\widehat\alpha_s$ in the chain are strictly increasing,
so that they define a folding of $\sigma t_\lambda$ in the sense of
\Cref{FM dictionary}.  As $\operatorname{dir}(\sigma t_\lambda)=\sigma$ and
$\operatorname{Re}\widehat\alpha_k=-\alpha_k^\vee$ (respectively $-\alpha_k$ if $e>1$), the projection of this folding
is \eqref{initial-block-finite-path}, which is a path in the (twisted) quantum Bruhat
graph.
In type $A_{2n}^{(2)}$, the labels $\widehat\alpha_k$ have even degree $\operatorname{len}(\alpha_k)M_{\alpha_k}$: this is clear for long roots, and for a short root $\alpha$ of $B_n$, whose coroot is $2\epsilon_l$, the pairing $M_\alpha=\langle-\lambda,\alpha^\vee\rangle$ is even, as $-\lambda\in\{\varpi_1,\ldots,\varpi_{n-1},2\varpi_n\}$.
Therefore the folding defines an element
$\mathbf p\in\mathcal{QB}(\sigma,t_\lambda)$.

Now, we compute the endpoint. By the action recalled in \Cref{FM dictionary}, $t_{-\lambda}$ maps $-\alpha^\vee$ (respectively $-\alpha$) to $\widehat\alpha$,
whence
\[
s_{\widehat\alpha}=t_{-\lambda}s_\alpha t_\lambda,
\qquad\text{equivalently}\qquad
t_\lambda s_{\widehat\alpha}=s_\alpha t_\lambda .
\]
Applying the second identity $s$ times,
\[
\operatorname{end}(\mathbf p)
=\sigma t_\lambda s_{\widehat\alpha_1}\cdots s_{\widehat\alpha_s}
=\sigma s_{\alpha_1}\cdots s_{\alpha_s}\,t_\lambda
=u\,t_\lambda=t_{u\lambda}u
\]
with weight $u\lambda$. Since $u\in vw_0W_J$, and $W_J$ fixes
$\lambda$, this equals $vw_0\lambda=v\chi_0$.

As for the degree, only the quantum edges of
\eqref{initial-block-finite-path} contribute, and
\[
c\deg\widehat\alpha_k
=c\operatorname{len}(\alpha_k)\langle-\lambda,\alpha_k^\vee\rangle
=\langle N(-w_0\chi),(\alpha_k)_\Sigma^\vee\rangle.
\]
Here we apply \eqref{twisted pairing} to the absolute fundamental weight $-w_0\chi$, whose restriction is $-\lambda$.
As
\eqref{initial-block-finite-path} is shortest, the sum of the coroots of $\Sigma$
along its quantum edges is $\operatorname{wt}(\sigma\Rightarrow u)$, so
\begin{equation}
\label{initial-block-degree}
\deg(\operatorname{qwt}(\mathbf p))
=\sum_{k:\,\mathrm{quantum}}\langle N(-w_0\chi),(\alpha_k)_\Sigma^\vee\rangle
=\bigl\langle N(-w_0\chi),\operatorname{wt}(ww_0\Rightarrow u)\bigr\rangle .
\end{equation}

Now, we must identify the right side of \eqref{initial-block-degree} with
$\langle N\chi,\operatorname{wt}(v\Rightarrow w)\rangle$.
Applying \cite[Lemma~3.8(a),(b)]{Schremmer24} on the effect of right multiplication by $w_0$ on the quantum Bruhat graph yields, since $w_0$ commutes with the Galois action,
\[
\bigl\langle N(-w_0\chi),\operatorname{wt}(ww_0\Rightarrow u)\bigr\rangle
=\bigl\langle -w_0N\chi,\,-w_0\operatorname{wt}(uw_0\Rightarrow w)\bigr\rangle
=\bigl\langle N\chi,\operatorname{wt}(uw_0\Rightarrow w)\bigr\rangle .
\]

Finally, $u\in vw_0W_J$ gives $uw_0\in vW_{J^\star}$, where $W_{J^\star}$ fixes $N\chi$.
By
\cite[Lemma~7.2]{LenartEtAl17}, applied to $\Sigma$ and the dominant weight $N\chi$, the
quantity $\langle N\chi,\operatorname{wt}(x\Rightarrow w)\rangle$ depends only on
the coset $xW_{J^\star}$.  Hence
$\langle N\chi,\operatorname{wt}(uw_0\Rightarrow w)\rangle
=\langle N\chi,\operatorname{wt}(v\Rightarrow w)\rangle$, as claimed.
\end{proof}

\subsection{Generalized Weyl modules}

We recall the modules of \cite[\S2.1]{FeiginMakedonskyi17} and \cite[\S2.1]{FeiginMakedonskyiTwisted17}. Let
$\frakg_{0,\alpha}$ be the $\mathfrak h_0$-weight space of $\frakg_0$ for $\alpha \in \Phi_0$, and
$\mathfrak n_+=\oplus_{\alpha \in \Phi_0^+}\frakg_{0,\alpha}$ be the maximal $\mathfrak h_0$-stable nilpotent subalgebra of $\frakg_0$ associated with $\Phi_0^+$.
Write $\mathfrak g[t]^{(e)}=\bigoplus_{k\geq 0}\frakg_{k\bmod e}\otimes t^k\subset \frakg[t]$ for the (twisted) current algebra,
$\mathfrak n^{\mathrm{af}}=(\mathfrak n_+\otimes 1)\oplus\bigoplus_{k\geq1}\frakg_{k\bmod e}\otimes t^k$ for the positive Iwahori current algebra,
both graded by $\deg(x\otimes t^k)=k$, and let $U(\mathfrak n^{\mathrm{af}})_{\le d}$
be the induced increasing filtration.  For $\mu=\sigma\lambda$ with $\lambda\in P_0$
antidominant and $\sigma\in W_0$, the \emph{generalized Weyl module}
$\mathbb W_\mu$ is the cyclic $\mathfrak n^{\mathrm{af}}$-module satisfying the relations of
\cite[Definition~2.1]{FeiginMakedonskyi17}, respectively \cite[Definition~2.1]{FeiginMakedonskyiTwisted17} if $e>1$, respectively \cite[Definition~2.2]{FeiginMakedonskyiTwisted17} in type $A_{2n}^{(2)}$.
It is graded by degree and by $\frakh_0$-weights, with cyclic vector of degree $0$ and weight $\mu$ \cite[Definition~2.5]{FeiginMakedonskyiTwisted17}.
For a
$\mathfrak g[t]$-module $V$ and $z\in\bbC$, we write $V^z$ for $V$
regarded as an $\mathfrak n^{\mathrm{af}}$-module via
$(x\otimes t^i)\cdot m=(x\otimes(t-z)^i)m$ (the shift in \cite[\S2.5]{FeiginMakedonskyiTwisted17} is by $+z$, which amounts to replacing $z$ by $-z$).

Throughout this subsection $\chi$ is an absolute fundamental weight and
$\lambda,J$ are as in \Cref{rem:two-parabolics}. We put $\sigma=ww_0$ and note that
$\sigma\lambda=w\chi_0$, so that the relevant generalized Weyl module
is $\mathbb W_{w\chi_0}$.  Let $W(\chi)$ be the local Weyl module of
the absolute current algebra $\mathfrak g[t]$ with cyclic lowest weight vector $e_{w_0\chi}$, fix
$z\in\bbC^\times$, and put $W_z(\chi)=W(\chi)^z$.  Give $W_z(\chi)$ the
\emph{fusion filtration}
\[
F_d^{\mathrm{FM}}W_z(\chi)
=
U(\mathfrak n^{\mathrm{af}})_{\leq d}\cdot\bigl(\sigma e_{w_0\chi}\bigr)
\]
of \cite[\S2.1]{FeiginMakedonskyi17}, namely the filtration defined in the
paragraph preceding \cite[Definition~2.9]{FeiginMakedonskyi17}, with $M=1$, as $\chi$ is a fundamental weight, respectively the filtration of \cite[\S2.5]{FeiginMakedonskyiTwisted17} if $e>1$.  It
is a filtration by $\mathfrak h_0$-stable subspaces.
If $e>1$, then the $\frakh_0$-weight space $W_z(\chi)_{v\chi_0}$ is the absolute weight space of weight $v\chi$, for $v\in W_0$: every weight of $W(\chi)$ is of the form $v\chi-\eta$ for a sum $\eta$ of absolute roots positive for $B_v$, and these restrict to nonzero sums of relative roots positive for $B_v$.
In particular, $W_z(\chi)_{v\chi_0}$ is one dimensional.

\begin{lemma}
\label{lem:FM-graded-character}
The canonical map of \cite[Lemma~2.11]{FeiginMakedonskyi17}, respectively of the proof of \cite[Lemma~2.6]{FeiginMakedonskyiTwisted17} if $e>1$, is a graded
surjection
\[
\mathbb W_{w\chi_0}
\;\twoheadrightarrow\;
\operatorname{gr}_{F^{\mathrm{FM}}}W_z(\chi).
\]
In addition, $\mathbb W_{w\chi_0}(v\chi_0)$ and $(\operatorname{gr}_{F^{\mathrm{FM}}}W_z(\chi))(v\chi_0)$ are one
dimensional and concentrated in degree
$\langle N\chi,\operatorname{wt}(v\Rightarrow w)\rangle$ for every $v\in W_0$.
\end{lemma}

\begin{proof}
Taking $M=1$ in \cite[Lemma~2.11]{FeiginMakedonskyi17}, respectively by the proof of \cite[Lemma~2.6]{FeiginMakedonskyiTwisted17} if $e>1$, there is a graded surjection
$\mathbb W_{w\chi_0}\twoheadrightarrow\operatorname{gr}_{F^{\mathrm{FM}}}W_z(\chi)$.
If $e>1$, the absolute fundamental weights in loc.\ cit.\ are taken at a fixed representative of each $\tau$-orbit of nodes, and we may choose it to be the node of $\chi$, as the relations used there are symmetric under $\tau$.
In type $A_{2n}^{(2)}$, the weight attached there to the node $n$ is $2\varpi_n$, in accordance with \Cref{FM dictionary}.
The proof of \cite[Theorem~2.21]{FeiginMakedonskyi17}, respectively of \cite[Theorem~2.15]{FeiginMakedonskyiTwisted17}, which follows the same lines, shows that the resulting inequality $\dim\mathbb W_{w\chi_0}\geq\dim W(\chi)$ is an equality.
In type $A_{2n}^{(2)}$, this uses that $\langle\lambda,\alpha_n^\vee\rangle$ is even, which holds as $-\lambda\in\{\varpi_1,\ldots,\varpi_{n-1},2\varpi_n\}$.
Hence the surjection is an isomorphism.
Outside type $A_{2n}^{(2)}$, the character formula of \cite[Corollary~2.22]{FeiginMakedonskyi17}, respectively \cite[Theorem~2.16]{FeiginMakedonskyiTwisted17}, uses our current grading directly.

In type $A_{2n}^{(2)}$, the alcove paths are formed with the auxiliary labels $\widehat\alpha_D$ of type $D_{n+1}^{(2)}$, whose degrees $\operatorname{len}(\alpha)D$ are not current degrees.
By \Cref{lem:A2n-degree-shift} below, at each step of the decomposition procedure of \cite[\S2.6]{FeiginMakedonskyiTwisted17} the kernel is shifted in current degree by $c\deg\beta$ on a quantum step with label $\beta$ and not at all on a Bruhat step, while the weight shifts are those of type $D_{n+1}^{(2)}$, and the steps excluded by the parity condition of \Cref{FM dictionary} contribute no kernel.
Summing over the steps of the procedure, the current degree attached to a path $\mathbf p$ is $c\sum\deg\beta_{j_k}$, the sum over the quantum steps of $\mathbf p$, that is, $\deg(\operatorname{qwt}(\mathbf p))$ as normalized in \Cref{FM dictionary}.

Applying the decomposition and character equality of \cite[Theorems~2.15 and 2.16]{FeiginMakedonskyiTwisted17} with these homogeneous shifts, the graded character of both sides is
\[
C_\sigma^{t_\lambda}
=
\sum_{\mathbf p\in\mathcal{QB}(\sigma,t_\lambda)}
x^{\operatorname{wt}(\operatorname{end}(\mathbf p))}
q^{\deg(\operatorname{qwt}(\mathbf p))}.
\]
The $v\chi_0$-weight space of $W_z(\chi)$ is one dimensional: if $e=1$, as the character of $W(\chi)$ is $W_0$-invariant and its $\chi$-weight space is one dimensional; if $e>1$, as noted before the lemma.
Consequently exactly one
$\mathbf p\in\mathcal{QB}(\sigma,t_\lambda)$ satisfies
$\operatorname{wt}(\operatorname{end}(\mathbf p))=v\chi_0$, and
\Cref{lem:initial-block-folding} exhibits one of degree
$\langle N\chi,\operatorname{wt}(v\Rightarrow w)\rangle$.
\end{proof}

\Cref{lem:A2n-degree-shift} below, used in the proof of \Cref{lem:FM-graded-character}, compares the current grading with the degrees of the auxiliary labels of type $D_{n+1}^{(2)}$ in type $A_{2n}^{(2)}$.
Before we can state it, some preparation is needed. We use the notation of \cite[\S2]{FeiginMakedonskyiTwisted17}, except that the element of $W_0$ and the antidominant weight of loc.\ cit.\ are denoted $\kappa$ and $\nu$, as $\sigma=ww_0$ and $\lambda$ are fixed in this subsection.
For $\alpha\in\Phi_0^+$, let $\widehat\kappa(e_\alpha)$ be the root vector of $\frakn^{\mathrm{af}}$ defined in \cite[\S2.2]{FeiginMakedonskyiTwisted17}, that is,
\begin{equation}\label{hat kappa}
\widehat{\kappa}(e_\alpha)=
\begin{cases}
e_{\kappa\alpha}\in\frakn_+, & \kappa\alpha\in\Phi_0^+,\\
e_{\kappa\alpha+\delta}\in\frakg_1\otimes t, & \kappa\alpha\in\Phi_0^-,\ \alpha\text{ long},\\
e_{2\kappa\alpha+\delta}\in\frakg_1\otimes t, & \kappa\alpha\in\Phi_0^-,\ \alpha\text{ short},
\end{cases}
\end{equation}
where $e_{\gamma+k\delta}$ denotes a generator of the subspace of $\frakg_{k\bmod 2}\otimes t^k$ of $\frakh_0$-weight $\gamma$.
These vectors have $\frakh_0$-weights $\kappa\alpha$, $\kappa\alpha$ and $2\kappa\alpha$, and current degrees $0$, $1$ and $1$, respectively.
The third case is possible as $\frakg_1$ is the irreducible $\frakg_0$-module of highest weight $2\varpi_1$, whose nonzero weights are the roots of $B_n$ and the doubles of its short roots \cite[\S2.1]{FeiginMakedonskyiTwisted17}.

Let $\bar\beta=(\beta_1,\ldots,\beta_r)$ be the sequence of labels of a reduced decomposition of $t_{-\varpi_i}$ in the extended affine Weyl group of type $D_{n+1}^{(2)}$, and, for $0\leq m\leq r$, let $W_{\kappa\nu}(\bar\beta,m)$ be the module of \cite[Definition~2.8]{FeiginMakedonskyiTwisted17}, with cyclic vector $v_m$.
For $\alpha\in\Phi_0^+$, put
\[
l_{\alpha,m}=-\langle\alpha^\vee,\nu\rangle-\bigl|\{j\leq m\mid\operatorname{Re}\beta_j=-\alpha\}\bigr|,
\qquad
p_{\alpha,m}=
\begin{cases}
\lfloor l_{\alpha,m}/2\rfloor, & \alpha\text{ short},\ \kappa\alpha\in\Phi_0^-,\\
l_{\alpha,m}, & \text{otherwise}.
\end{cases}
\]
Then the defining relations of $W_{\kappa\nu}(\bar\beta,m)$ are $\frakh_{k\delta}v_m=0$ for $k\geq1$, $\widehat\kappa(e_{-\alpha+\delta})v_m=0$, and
\begin{equation}\label{FM relation}
\widehat{\kappa}(e_\alpha)^{p_{\alpha,m}+1}v_m=0,
\qquad\alpha\in\Phi_0^+.
\end{equation}
The modules are graded by $\frakh_0$-weight and current degree, and a homogeneous map between them is determined by the image of the cyclic vector.
For a homogeneous map $f$ with $f(v')=y$ for the cyclic vector $v'$ of the source, we call the difference of the $\frakh_0$-weights of $y$ and of $v_m$ the \emph{weight shift} of $f$, and the current degree of $y$ its \emph{degree shift}; these are the monomial $x^{\text{weight shift}}q^{\text{degree shift}}$ by which $f$ multiplies characters, compare \cite[Corollary~2.14]{FeiginMakedonskyiTwisted17}.

\begin{lemma}\label{lem:A2n-degree-shift}
Assume that $\frakg$ is of type $A_{2n}^{(2)}$.
Let $0\leq m<r$, put $\alpha:=-\operatorname{Re}\beta_{m+1}\in\Phi_0^+$ and $l:=l_{\alpha,m}$, and consider the surjection
\[
W_{\kappa\nu}(\bar\beta,m)\twoheadrightarrow W_{\kappa\nu}(\bar\beta,m+1)
\]
of \cite[\S2.6]{FeiginMakedonskyiTwisted17}.
\begin{enumerate}
\item The label $\beta_{m+1}$ has degree $\deg\beta_{m+1}=\operatorname{len}(\alpha)\,l$.
\item If $\kappa\xrightarrow{\alpha}\kappa s_\alpha$ is not an edge of the twisted quantum Bruhat graph of type $D_{n+1}^{(2)}$, or if it is a quantum edge and $\deg\beta_{m+1}$ is odd, then the surjection is an isomorphism.
\item Otherwise, its kernel is nonzero, and it is the image of an injective homogeneous map
\[
W_{\kappa s_\alpha\nu}(\bar\beta,m+1)\hookrightarrow W_{\kappa\nu}(\bar\beta,m)
\]
with weight shift $l\kappa\alpha$, and with degree shift $0$ if $\kappa\alpha\in\Phi_0^+$ and $c\deg\beta_{m+1}=\frac12\operatorname{len}(\alpha)\,l$ if $\kappa\alpha\in\Phi_0^-$.
\end{enumerate}
\end{lemma}

\begin{proof}
Part (1) is the identity $\operatorname{len}(\operatorname{Re}\beta_{m+1})\,l_{-\operatorname{Re}\beta_{m+1},m}=\deg\beta_{m+1}$ in the proof of \cite[Theorem~2.16]{FeiginMakedonskyiTwisted17}, where $\operatorname{len}(\alpha)=2$ for long and $\operatorname{len}(\alpha)=1$ for short $\alpha$.
In particular, $\deg\beta_{m+1}$ is odd if and only if $\alpha$ is short and $l$ is odd.

Passing from $m$ to $m+1$ replaces $l_{\alpha,m}=l$ by $l-1$ and leaves $l_{\alpha',m}$ unchanged for $\alpha'\neq\alpha$.
So the relations \eqref{FM relation} show that the kernel of the surjection is the submodule generated by
\[
y_m:=\widehat{\kappa}(e_\alpha)^{p}\,v_m,
\qquad
p=
\begin{cases}
l/2, & \alpha\text{ short},\ \kappa\alpha\in\Phi_0^-,\ l\text{ even},\\
l, & \text{otherwise}.
\end{cases}
\]
Recall that an edge $\kappa\xrightarrow{\alpha}\kappa s_\alpha$ is a Bruhat edge if $\kappa\alpha\in\Phi_0^+$ and a quantum edge if $\kappa\alpha\in\Phi_0^-$.
If there is no such edge, then $y_m=0$, that is, the surjection is an isomorphism, by \cite[Theorem~2.13(i)]{FeiginMakedonskyiTwisted17}.
If the edge is quantum and $\deg\beta_{m+1}$ is odd, then $\alpha$ is short and $l$ is odd by (1), so $\lfloor l/2\rfloor=\lfloor (l-1)/2\rfloor$, the relation for $\alpha$ does not change, and the surjection is again an isomorphism.
This proves (2).
In the remaining cases, \cite[Theorem~2.13(ii)]{FeiginMakedonskyiTwisted17} and the exactness in \cite[Theorem~2.16]{FeiginMakedonskyiTwisted17} identify the kernel with $W_{\kappa s_\alpha\nu}(\bar\beta,m+1)$, which is nonzero, via the map sending its cyclic vector to $y_m$.
By \eqref{hat kappa}, its weight and degree shifts are $p$ times the $\frakh_0$-weight and current degree of $\widehat\kappa(e_\alpha)$:
\begin{center}
\begin{tabular}{llll}
case & $p$ & weight shift & degree shift\\\hline
$\kappa\alpha\in\Phi_0^+$ & $l$ & $l\kappa\alpha$ & $0$\\
$\kappa\alpha\in\Phi_0^-$, $\alpha$ long & $l$ & $l\kappa\alpha$ & $l=\frac12\cdot2\cdot l$\\
$\kappa\alpha\in\Phi_0^-$, $\alpha$ short, $l$ even & $\frac l2$ & $\frac l2\cdot2\kappa\alpha=l\kappa\alpha$ & $\frac l2=\frac12\cdot1\cdot l$
\end{tabular}
\end{center}
In the last two rows, the degree shift is $\frac12\operatorname{len}(\alpha)\,l=c\deg\beta_{m+1}$ by (1).
\end{proof}

%The edge condition in (2) is not automatic: $\alpha$ is determined by the label $\beta_{m+1}$, but $\kappa$ is the current vertex of the path, which is arbitrary.
%For instance, for $\kappa=1$ there is an edge $1\xrightarrow{\alpha}s_\alpha$ only if $\ell(s_\alpha)=1$, that is, only if $\alpha$ is simple.
%The second case in (2) is the exclusion of quantum edges with labels of odd degree in \cite[Remark~2.12]{FeiginMakedonskyiTwisted17}, see \Cref{FM dictionary}.
%The weight shifts are those of the decomposition procedure in type $D_{n+1}^{(2)}$, so only the degrees are rescaled, by the factor $c$ on quantum steps.

\subsection{The polynomial degree as a filtration degree}
We return to the notation of \Cref{iwahori calc sec}, relating the coefficients $e_\alpha n_b(\alpha,k)$ of $\deg^\chi_{v,w}$ to the twisted current algebra.
We use the embedding $\Lambda\subset\bbC$ and the inverse pinned automorphism $\tau=\tau_\Gamma^{-1}$ fixed in \Cref{FM dictionary}.
By \eqref{tau eigenvectors}, one has $E_{\alpha,k}\in\frakg_{ke_\alpha}$ for all $(\alpha,k)\in v\bar\Psi^+$.
%, including $E_{\alpha,0}\in\frakg_0$, $E_{\alpha,1}\in\frakg_1$ and $E_{2\alpha,1}\in\frakg_1$ for multipliable $\alpha$.
Let $\frakb_0\subset\frakg_0$ be the Borel subalgebra determined by $B_0$, and consider the twisted Iwahori current algebra
\[
\fraki=\frakb_0\oplus\bigoplus_{k\geq1}\frakg_{k\bmod e}\otimes t^k\subset\frakg[t],
\]
graded by the degree in $t$.
For $e=1$, this is the usual Iwahori current algebra.
Then $e_\alpha n_b(\alpha,k)$ is the least degree of $\fraki$ in which $E_{\alpha,k}$ occurs.
Indeed, by \Cref{Iwahori-minor-expansion recall coro}, $e_\alpha n_b(\alpha,k)=ke_\alpha$ if $k>0$, which lies in $\{1,\ldots,e-1\}$ as $k<m_\alpha$,
%(this includes $e_{2\alpha}n_b(2\alpha,1)=1$ for multipliable $\alpha$),
 while $e_\alpha n_b(\alpha,0)$ is $0$ if $\alpha\in\Phi^+$ and $e_\alpha m_\alpha=e$ if $\alpha\in\Phi^-$.
On the other hand, by definition of $\fraki$, a vector of $\frakg_j$ with $0<j<e$ first occurs in degree $j$, and a root vector of $\frakg_0$ occurs in degree $0$ if it lies in $\frakb_0$ and first occurs in degree $e$ otherwise; so the claim follows from $E_{\alpha,k}\in\frakg_{ke_\alpha}$.

\begin{lemma}
\label{lem:minor-degree-one-factor-filtration}
Let $F_dU(\frakg)$ be the image of the part of $U(\fraki)$ of degree at most $d$ under the evaluation $t\mapsto1$.
Define
\[
F_{d,w}V_{v,\chi}
=
F_dU(\frakg)\cdot f_{v,w}^{\chi},
\]
where $V_{v,\chi}=V^H_{\chi,v}\otimes_\Lambda\bbC=\Ind_{B_v^-}(v\chi)$, and $f^\chi_{v,w}$ is the image of the generator $f$ chosen in \Cref{rem:integrality-of-coefficients}.
Its polynomial coefficients in the coordinates $\xi_{\alpha,k}$ are the $u_{\mathbf n}\in\Lambda\subset\bbC$.
Then
\[
\min_{u_{\mathbf n}\neq0}\deg_{v,w}^\chi(\mathbf n)
=
\min\{d\mid f_{v,v}^{\chi}\in F_{d,w}V_{v,\chi}\}.
\]
\end{lemma}

\begin{proof}
Consider the map $\ell:V_{v,\chi}\to \bbC$ given by $\ell(f)=f(1)$.
If \(f\) is a $T$-stable vector of \(T\)-weight \(\eta\), right
equivariance at \(1\) shows that \(\ell(f)=0\) unless
\(\eta=v\chi\).  Thus \(\ell\) detects precisely the
\(v\chi\)-extremal line and kills the remaining $T$-weight spaces.

The filtration $F_\bullet U(\frakg)$ is the multiplicative filtration generated by the increasing filtration of $\frakg$ by the images of the graded pieces of $\fraki$, because $\fraki$ is a graded Lie subalgebra of $\frakg[t]$.
Hence its associated graded algebra is the enveloping algebra of the associated graded Lie algebra, and the filtration is compatible with the Poincar\'e--Birkhoff--Witt decomposition
\begin{equation}\label{big cell dist}
U(\frakn_v^-)\otimes U(\frakh)\otimes U(\frakn_v^+)
\xrightarrow{\sim}U(\frakg),
\end{equation}
where $\frakn_v^\pm$ are the nilradicals of $\Lie(B_v^\pm)$ and $\frakh=\Lie(T)$.
The algebra $U(\frakn_v^-)$ acts under $\ell$ through its augmentation: $\ell(Dm)=0$ for $D$ in its augmentation ideal. Moreover, $U(\frakh)$ acts on weight vectors by scalars.
Hence, by subadditivity of the degree, the minimal $d$ with $f_{v,v}^\chi\in F_dU(\frakg)\cdot f_{v,w}^\chi$
is already attained with $U(\frakn_v^+)$ in place of the full enveloping algebra.

Fix an ordering of $v\bar\Psi^+$, compatible with the ordering of $v\Phi^+$ in \Cref{rem:integrality-of-coefficients}.
This orders the coordinates $\xi_{\alpha,k}$ of \Cref{rem:integrality-of-coefficients}, and the ordered monomials in the divided powers of the $E_{\alpha,k}$ form a basis of $U(\frakn_v^+)$ compatible with the filtration.
It is dual to the monomials in the coordinates $\eta_{\alpha,k}$ of the isomorphism $(\eta_{\alpha,k})\mapsto\prod\exp(\eta_{\alpha,k}E_{\alpha,k})$ onto $U_v$ \cite[I.7.4]{Jan03}.
If $\Phi$ is reduced, then $\eta_{\alpha,k}=\xi_{\alpha,k}$.
If $\alpha$ is multipliable, we order $(\alpha,0),(\alpha,1),(2\alpha,1)$ consecutively.
Then $\eta_{\alpha,k}=\xi_{\alpha,k}$ for $k=0,1$, and the Baker--Campbell--Hausdorff formula in the Heisenberg group, together with $[E_{\tilde\alpha},E_{\tilde\alpha'}]=E_{\tilde\alpha+\tilde\alpha'}$ and the ordering in \eqref{Heisenberg parametrization}, gives $Z_{\tilde\alpha+\tilde\alpha'}=\eta_{2\alpha,1}-\eta_{\alpha,0}\eta_{\alpha,1}+\frac12(\eta_{\alpha,0}^2-\eta_{\alpha,1}^2)$.
Comparing with \eqref{Heisenberg Fourier}, we get $\xi_{2\alpha,1}=\eta_{2\alpha,1}-\eta_{\alpha,0}\eta_{\alpha,1}$.
The monomial $\eta_{\alpha,0}\eta_{\alpha,1}$ has degree $e_\alpha n_b(\alpha,0)+e_\alpha n_b(\alpha,1)\geq1=e_{2\alpha}n_b(2\alpha,1)$, and both this substitution and its inverse do not decrease monomial degrees. They therefore preserve the least degree with nonzero coefficient.
Hence we may compute this minimum in the coordinates $\eta_{\alpha,k}$.
We obtain the coefficient of the monomial
\(\prod \eta_{\alpha,k}^{n_{\alpha,k}}\) in $f^\chi_{v,w}$ as
\[
(-1)^{\sum_{(\alpha,k)}n_{\alpha,k}}\ell\left(
\prod_{(\alpha,k)\in v\bar\Psi^+}
E_{\alpha,k}^{(n_{\alpha,k})}f^\chi_{v,w}\right),
\]
with the product taken in the reverse order,
and its 
filtered degree is exactly $\deg_{v,w}^\chi(\mathbf n)$.
Since the filtration is $\frakh_0$-stable and the $v\chi_0$-weight space is the $v\chi$-extremal line, this identifies the sought minimal $d$ with the minimum of $\deg_{v,w}^\chi(\mathbf n)$ across all $\mathbf n$ such that this coefficient is nonzero.
%, which is the minimum across all $\mathbf n$ with $u_{\mathbf n}\neq 0$ by the above.
\end{proof}

\subsection{Comparison of filtrations}

Let $V(\chi)$ denote the irreducible $\mathfrak g$-module of highest weight
$\chi$, made into a $\mathfrak g[t]$-module by evaluation at $t=0$,
and let $V_z(\chi)=V(\chi)^z$.  Let
\[
\pi_z\colon W_z(\chi)\twoheadrightarrow V_z(\chi)
\]
be the canonical quotient, and endow $V_z(\chi)$ with the quotient filtration
\[
F_dV_z(\chi):=\pi_z\bigl(F_d^{\mathrm{FM}}W_z(\chi)\bigr).
\]
In \Cref{lem:minor-degree-one-factor-filtration}, we defined $F_dU(\frakg)$ as the image of $U(\fraki)_{\leq d}$ under evaluation at $t=1$.
We compare these filtrations below.

\begin{lemma}
\label{lem:quotient-filtration-identification}
Under the identification of $\mathfrak g$-modules $V_z(\chi)\simeq V(\chi)$, we have
\[
F_dV_z(\chi)=F_d\,U(\mathfrak g)\cdot e_{w\chi}.
\]
\end{lemma}

\begin{proof}
Note that $\sigma e_{w_0\chi}=e_{w\chi}$
because $\sigma=ww_0$.
Next,
$\mathfrak n^{\mathrm{af}}$ acts on $V_z(\chi)$ through
\[
x\otimes t^k\longmapsto(-z)^k x \quad(\deg=k,\ x\in\frakg_{k\bmod e},\ k\ge0),
\]
with $x\in\frakn_+$ if $k=0$.
Since $z\neq0$, it follows that $F_dV_z(\chi)$ is spanned by the vectors $x_1\cdots x_l\,e_{w\chi}$ with $x_a\in\frakg_{k_a\bmod e}$, $x_a\in\frakn_+$ if $k_a=0$, and $k_1+\cdots+k_l\leq d$.
In other words, $F_dV_z(\chi)$ is the image of $U(\frakn^{\mathrm{af}})_{\leq d}$ under evaluation at $t=1$, applied to $e_{w\chi}$.
Finally, $\fraki=\frakh_0\oplus\frakn^{\mathrm{af}}$, where $\frakh_0$ sits in degree $0$ and normalizes the graded ideal $\frakn^{\mathrm{af}}$, so that $U(\fraki)_{\leq d}=U(\frakn^{\mathrm{af}})_{\leq d}\cdot U(\frakh_0)$ by the Poincar\'e--Birkhoff--Witt theorem.
As $U(\frakh_0)$ preserves the line spanned by $e_{w\chi}$, we can now identify $F_dV_z(\chi)$ with
the filtration of
\Cref{lem:minor-degree-one-factor-filtration}.
\end{proof}

\begin{proposition}
\label{prop:characteristic-zero-extremal-degree}
For an absolute fundamental weight \(\chi\) and \(v,w\in W_0\), one has
\[
\min_{u_{\mathbf n}\neq0}\deg_{v,w}^\chi(\mathbf n)
=
\left\langle N\chi,
\operatorname{wt}(v\Rightarrow w)\right\rangle.
\]
Here the minimum is over the nonzero universal coefficients $u_{\mathbf n}\in\Lambda\subset\bbC$ of \Cref{rem:integrality-of-coefficients}.
\end{proposition}

\begin{proof}
By \Cref{lem:minor-degree-one-factor-filtration}, the left-hand side is the least
$d$ such that the $v\chi$-extremal line of $V(\chi)$ meets
$F_dU(\mathfrak g)e_{w\chi}$.
By
\Cref{lem:quotient-filtration-identification} the filtration in question is the
quotient filtration $F_\bullet V_z(\chi)$.

The $\frakh_0$-weight spaces $W_z(\chi)_{v\chi_0}$ and $V_z(\chi)_{v\chi_0}$ are both
one dimensional, and equal to the extremal lines of absolute weight $v\chi$, as noted before \Cref{lem:FM-graded-character}.  As $\pi_z$ is surjective and $\mathfrak h$-equivariant, it
restricts to an isomorphism between them.  Since $F^{\mathrm{FM}}$ is a
filtration by $\mathfrak h_0$-stable subspaces,
\[
\bigl(F_dV_z(\chi)\bigr)_{v\chi_0}
=\pi_z\Bigl(\bigl(F_d^{\mathrm{FM}}W_z(\chi)\bigr)_{v\chi_0}\Bigr),
\]
and the two sides vanish for the same values of $d$.  The least $d$ can thus be read off the inclusion $W_z(\chi)_{v\chi_0}\subset F_d^{\mathrm{FM}}W_z(\chi)$.
As $W_z(\chi)_{v\chi_0}$ has dimension $1$, it suffices to compute the degree
of the $v\chi_0$-weight space of $\operatorname{gr}_{F^{\mathrm{FM}}}W_z(\chi)$.  By
\Cref{lem:FM-graded-character} it equals
$\langle N\chi,\operatorname{wt}(v\Rightarrow w)\rangle$.
\end{proof}

\end{appendix}

\newcommand{\etalchar}[1]{$^{#1}$}


\begin{thebibliography}{LRVW26}

\bibitem[AGLR22]{AGLR22}
Johannes Ansch{\"u}tz, Ian Gleason, Jo{\~a}o Louren{\c{c}}o, and Timo Richarz.
\newblock On the {$p$}-adic theory of local models.
\newblock To appear in \textit{Ann.~of Math.}, arXiv:2201.01234, 2022.

\bibitem[ALRR24]{ALRR24}
Pramod Achar, Jo{\~a}o Louren{\c{c}}o, Timo Richarz, and Simon Riche.
\newblock A modular ramified geometric {S}atake equivalence.
\newblock To appear in \textit{J. Inst. Math. Jussieu}, arXiv:2403.10651, 2024.

\bibitem[BHY26]{BessonHongYu:FMCounterexample}
Marc Besson, Jiuzu Hong, and Huanhuan Yu.
\newblock Tangent spaces of spherical schubert varieties and counterexamples to
  the reducedness conjecture.
\newblock Preprint, arXiv:2603.17273, 2026.

\bibitem[BLR90]{bosch_luetkebohmert_raynaud_neron_models}
Siegfried Bosch, Werner L\"utkebohmert, and Michel Raynaud.
\newblock {\em N\'eron models}, volume~21 of {\em Ergebnisse der Mathematik und
  ihrer Grenzgebiete (3) [Results in Mathematics and Related Areas (3)]}.
\newblock Springer-Verlag, Berlin, 1990.

\bibitem[Bou68]{Bou68}
N.~Bourbaki.
\newblock {\em \'{E}l\'{e}ments de math\'{e}matique. {F}asc. {XXXIV}. {G}roupes
  et alg\`ebres de {L}ie. {C}hapitre {IV}: {G}roupes de {C}oxeter et syst\`emes
  de {T}its. {C}hapitre {V}: {G}roupes engendr\'{e}s par des r\'{e}flexions.
  {C}hapitre {VI}: syst\`emes de racines}.
\newblock Actualit\'{e}s Scientifiques et Industrielles, No. 1337. Hermann,
  Paris, 1968.

\bibitem[BS17]{bhatt_scholze_projectivity_of_the_witt_vector_affine_grassmannian}
Bhargav Bhatt and Peter Scholze.
\newblock Projectivity of the {W}itt vector affine {G}rassmannian.
\newblock {\em Invent. Math.}, 209(2):329--423, 2017.

\bibitem[BT72]{BT72}
F.~Bruhat and J.~Tits.
\newblock Groupes r\'{e}ductifs sur un corps local.
\newblock {\em Publ. Math. Inst. Hautes \'{E}tudes Sci.}, (41):5--251, 1972.

\bibitem[BT84]{BT84}
F.~Bruhat and J.~Tits.
\newblock Groupes r\'eductifs sur un corps local. {II}. {S}ch\'emas en groupes.
  {E}xistence d'une donn\'ee radicielle valu\'ee.
\newblock {\em Publ. Math. Inst. Hautes \'Etudes Sci.}, (60):197--376, 1984.

\bibitem[{\v C}es22]{Cesnavicius2022Torsors}
K{\k e}stutis {\v C}esnavi{\v c}ius.
\newblock Problems about torsors over regular rings.
\newblock {\em Acta Mathematica Vietnamica}, 47(1):39--107, 2022.
\newblock With an appendix by Yifei Zhao.

\bibitem[CN90]{CN90}
Ching-Li Chai and Peter Norman.
\newblock Bad reduction of the {S}iegel moduli scheme of genus two with
  {$\Gamma_0(p)$}-level structure.
\newblock {\em Amer. J. Math.}, pages 1003--1071, 1990.

\bibitem[dJ93]{dJ93}
Aise~Johan de~Jong.
\newblock The moduli spaces of principally polarized abelian varieties with
  {$\Gamma_0(p)$}-level structure.
\newblock {\em J. Algebraic Geom.}, (2):667--688, 1993.

\bibitem[DP94]{DP94}
Pierre Deligne and Georgios Pappas.
\newblock Singularit{\'e}s des espaces de modules de {H}ilbert, en les
  caract{\'e}ristiques divisant le discriminant.
\newblock {\em Compos. Math.}, 90(1):59--79, 1994.

\bibitem[DR73]{DR73}
Pierre Deligne and Michael Rapoport.
\newblock Les sch{\'e}mas de modules de courbes elliptiques.
\newblock In {\em Modular functions of one variable II}, pages 143--316.
  Springer, 1973.

\bibitem[Fal97]{Fal97}
Gerd Faltings.
\newblock Explicit resolution of local singularities of moduli-spaces.
\newblock {\em J. Reine Angew. Math.}, 483:183--196, 1997.

\bibitem[Fal01]{Fal01}
Gerd Faltings.
\newblock Toroidal resolutions for some matrix singularities.
\newblock In {\em Moduli of abelian varieties ({T}exel {I}sland, 1999)}, volume
  195 of {\em Progr. Math.}, pages 157--184. Birkh\"{a}user, Basel, 2001.

\bibitem[Fal03]{Fal03}
Gerd Faltings.
\newblock Algebraic loop groups and moduli spaces of bundles.
\newblock {\em J. Eur. Math. Soc. (JEMS)}, 5(1):41--68, 2003.

\bibitem[FM99]{FM:SemiInfinite}
Michael Finkelberg and Ivan Mirkovi\'c.
\newblock Semi-infinite flags. {I}. {C}ase of global curve {$\bold P^1$}.
\newblock In {\em Differential topology, infinite-dimensional {L}ie algebras,
  and applications}, volume 194 of {\em Amer. Math. Soc. Transl. Ser. 2}, pages
  81--112. Amer. Math. Soc., Providence, RI, 1999.

\bibitem[FM17a]{FeiginMakedonskyi17}
Evgeny Feigin and Ievgen Makedonskyi.
\newblock Generalized {W}eyl modules, alcove paths and {M}acdonald polynomials.
\newblock {\em Selecta Math. (N.S.)}, 23(4):2863--2897, 2017.

\bibitem[FM17b]{FeiginMakedonskyiTwisted17}
Evgeny Feigin and Ievgen Makedonskyi.
\newblock Generalized {W}eyl modules for twisted current algebras.
\newblock {\em Theoret. Math. Phys.}, 192(2):1184--1204, 2017.
\newblock Citations to numbered statements refer to arXiv:1606.05219v2.

\bibitem[GD63]{SGA3}
Alexandre Grothendieck and Michel Demazure.
\newblock {\em Sch{\'e}mas en groupes}.
\newblock Publ. Math. Inst. Hautes \'Etudes Sci., 1963.

\bibitem[G{\"o}r01]{Goertz:LinearCase}
Ulrich G{\"o}rtz.
\newblock On the flatness of models of certain {S}himura varieties of
  {PEL}-type.
\newblock {\em Math. Ann.}, 321(3):689--727, 2001.

\bibitem[G{\"o}r03]{Goertz:SymplecticCase}
Ulrich G{\"o}rtz.
\newblock On the flatness of local models for the symplectic group.
\newblock {\em Adv. Math.}, 176(1):89--115, 2003.

\bibitem[G{\"o}r05]{Goertz:RamifiedCase}
Ulrich G{\"o}rtz.
\newblock Topological flatness of local models in the ramified case.
\newblock {\em Math. Z.}, 250(4):775--790, 2005.

\bibitem[Hai]{Haines:Pluecker}
Thomas~J. Haines.
\newblock Drinfeld--pluecker relations for homogeneous spaces, (affine) flag
  varieties, and rapoport-zink models.

\bibitem[Hai18]{Hai18}
Thomas~J. Haines.
\newblock Dualities for root systems with automorphisms and applications to
  non-split groups.
\newblock {\em Represent. Theory}, 22:1--26, 2018.

\bibitem[HH17]{HH17}
Thomas~J. Haines and Xuhua He.
\newblock Vertexwise criteria for admissibility of alcoves.
\newblock {\em Amer. J. Math.}, 139(3):769--784, 2017.

\bibitem[HN02]{HN02}
Thomas~J. Haines and Bao~Ch{\^a}u Ng{\^o}.
\newblock Alcoves associated to special fibers of local models.
\newblock {\em Amer. J. Math.}, 124(6):1125--1152, 2002.

\bibitem[HR08]{HR08}
Thomas Haines and Michael Rapoport.
\newblock On parahoric subgroups.
\newblock {\em Adv. Math.}, 219(1):188--198, 2008.

\bibitem[HV11]{HartlViehmann1}
Urs Hartl and Eva Viehmann.
\newblock The {N}ewton stratification on deformations of local {$G$}-shtukas.
\newblock {\em Journal für die reine und angewandte Mathematik},
  2011(656):87--129, 2011.

\bibitem[Jan03]{Jan03}
Jens~Carsten Jantzen.
\newblock {\em Representations of algebraic groups}, volume 107 of {\em
  Mathematical Surveys and Monographs}.
\newblock American Mathematical Society, Providence, RI, second edition, 2003.

\bibitem[Kac90]{Kac90}
Victor~G Kac.
\newblock {\em Infinite-dimensional {L}ie algebras}.
\newblock Cambridge university press, 1990.

\bibitem[KMW18]{KamnitzerMuthiahWeekes:ReducednessConjecture}
Joel Kamnitzer, Dinakar Muthiah, and Alex Weekes.
\newblock On a reducedness conjecture for spherical {S}chubert varieties and
  slices in the affine {G}rassmannian.
\newblock {\em Transform. Groups}, 23(3):707--722, 2018.

\bibitem[Kot97]{Kot97}
Robert~E. Kottwitz.
\newblock Isocrystals with additional structure. {II}.
\newblock {\em Compos. Math.}, 109(3):255--339, 1997.

\bibitem[KP18]{kisin_pappas_integral_models_of_shimura_varieties_with_parahoric_level_structure}
M.~Kisin and G.~Pappas.
\newblock Integral models of {S}himura varieties with parahoric level
  structure.
\newblock {\em Publ. Math. Inst. Hautes \'{E}tudes Sci.}, 128:121--218, 2018.

\bibitem[KP23]{KP23}
Tasho Kaletha and Gopal Prasad.
\newblock {\em Bruhat-{T}its theory---a new approach}, volume~44 of {\em New
  Mathematical Monographs}.
\newblock Cambridge University Press, Cambridge, 2023.

\bibitem[KPZ26]{KisinPappasZhou_IntegralModels}
Mark Kisin, Georgios Pappas, and Rong Zhou.
\newblock Integral models of {S}himura varieties with parahoric level
  structure, {II}.
\newblock {\em Forum Math. Pi}, 14:Paper No. e14, 89, 2026.

\bibitem[KR00]{KR00}
R~Kottwitz and M~Rapoport.
\newblock Minuscule alcoves for ${G}{L}_n$ and ${G}{S}p_{2n}$.
\newblock {\em Manuscripta mathematica}, 102(4):403--428, 2000.

\bibitem[Lev16]{Lev16}
Brandon Levin.
\newblock Local models for {W}eil-restricted groups.
\newblock {\em Compos. Math.}, 152(12):2563--2601, 2016.

\bibitem[LNS{\etalchar{+}}17]{LenartEtAl17}
Cristian Lenart, Satoshi Naito, Daisuke Sagaki, Anne Schilling, and Mark
  Shimozono.
\newblock A uniform model for {K}irillov--{R}eshetikhin crystals {II}: Alcove
  model, path model, and {$P=X$}.
\newblock {\em Int. Math. Res. Not. IMRN}, 2017(14):4259--4319, 2017.

\bibitem[Lou23]{Lou23}
Jo\~{a}o Louren\c{c}o.
\newblock Grassmanniennes affines tordues sur les entiers.
\newblock {\em Forum Math. Sigma}, 11:Paper No. e12, 65, 2023.

\bibitem[LRVW26]{LRVW:Minuscule}
Jo{\~a}o Louren{\c{c}}o, Timo Richarz, Eva Viehmann, and Torsten Wedhorn.
\newblock Moduli descriptions of v-sheaf local models in the minuscule case,
  2026.
\newblock In preparation.

\bibitem[Luo25]{Luo}
Yu~Luo.
\newblock On the moduli description of ramified unitary local models of
  signature {$(n-1,1)$}.
\newblock {\em Math. Ann.}, 392(4):4661--4738, 2025.

\bibitem[Mil17]{Milne:AlgebraicGroups}
J.~S. Milne.
\newblock {\em Algebraic groups: The theory of group schemes of finite type
  over a field}, volume 170 of {\em Cambridge Studies in Advanced Mathematics}.
\newblock Cambridge University Press, Cambridge, 2017.
\newblock Corrected 2022 printing.

\bibitem[Pap00]{Pap00}
Georgios Pappas.
\newblock On the arithmetic moduli schemes of {PEL} {S}himura varieties.
\newblock {\em J. Algebraic Geom.}, 9(3):577, 2000.

\bibitem[PR03]{PR03}
G~Pappas and M~Rapoport.
\newblock Local models in the ramified case. {I}: {T}he {EL}-case.
\newblock {\em J. Algebraic Geom.}, 12(1):107--145, 2003.

\bibitem[PR05]{PR05}
G.~Pappas and M.~Rapoport.
\newblock Local models in the ramified case. {II}. {S}plitting models.
\newblock {\em Duke Math. J.}, 127(2):193--250, 2005.

\bibitem[PR08]{PR08}
G.~Pappas and M.~Rapoport.
\newblock Twisted loop groups and their affine flag varieties.
\newblock {\em Adv. Math.}, 219(1):118--198, 2008.
\newblock With an appendix by T. Haines and Rapoport.

\bibitem[PR09]{PR09}
Georgios Pappas and Michael Rapoport.
\newblock Local models in the ramified case. {III} {U}nitary groups.
\newblock {\em J. Inst. Math. Jussieu}, 8(3):507--564, 2009.

\bibitem[PR24]{PR21}
Georgios Pappas and Michael Rapoport.
\newblock {$p$}-adic shtukas and the theory of global and local {S}himura
  varieties.
\newblock {\em Camb. J. Math.}, 12(1):1--164, 2024.

\bibitem[PRS13]{PRS13}
Georgios Pappas, Michael Rapoport, and Brian Smithling.
\newblock Local models of {S}himura varieties, {I}. {G}eometry and
  combinatorics.
\newblock In {\em Handbook of moduli. {V}ol. {III}}, volume~26 of {\em Adv.
  Lect. Math. (ALM)}, pages 135--217. Int. Press, Somerville, MA, 2013.

\bibitem[PZ13]{PZ13}
George Pappas and Xinwen Zhu.
\newblock Local models of {S}himura varieties and a conjecture of {K}ottwitz.
\newblock {\em Invent. Math.}, 194(1):147--254, 2013.

\bibitem[PZ22]{Pappas_Zachos:OrthogonalShimuraVarieties}
G.~Pappas and I.~Zachos.
\newblock Regular integral models for {S}himura varieties of orthogonal type.
\newblock {\em Compos. Math.}, 158(4):831--867, 2022.

\bibitem[Rap90]{Rap90}
M.~Rapoport.
\newblock On the bad reduction of {S}himura varieties.
\newblock In {\em Automorphic forms, {S}himura varieties, and {$L$}-functions,
  {V}ol. {II} ({A}nn {A}rbor, {MI}, 1988)}, volume~11 of {\em Perspect. Math.},
  pages 253--321. Academic Press, Boston, MA, 1990.

\bibitem[Rap05]{Rap05}
Michael Rapoport.
\newblock A guide to the reduction modulo $p$ of {S}himura varieties.
\newblock {\em Ast{\'e}risque}, 298:271--318, 2005.

\bibitem[Ric13]{Ri13}
Timo Richarz.
\newblock Schubert varieties in twisted affine flag varieties and local models.
\newblock {\em J. Algebra}, 375:121--147, 2013.

\bibitem[RZ96]{RZ96}
M.~Rapoport and Th. Zink.
\newblock {\em Period spaces for {$p$}-divisible groups}, volume 141 of {\em
  Ann. of Math. Stud.}
\newblock Princeton University Press, Princeton, NJ, 1996.

\bibitem[Sau16]{Sauermann}
Lisa Sauermann.
\newblock On the {$\mu$}-admissible set in the extended affine {W}eyl groups of
  {$E_6$} and {$E_7$}.
\newblock {\em J. Algebra}, 451:526--543, 2016.

\bibitem[Sch86]{Schrijver:LinearIntegerProgramming}
Alexander Schrijver.
\newblock {\em Theory of Linear and Integer Programming}.
\newblock Wiley-Interscience Series in Discrete Mathematics and Optimization.
  John Wiley \& Sons, Chichester, 1986.

\bibitem[Sch17]{Sch17}
Peter Scholze.
\newblock {\'Etale cohomology of diamonds}.
\newblock {\em Preprint, arXiv:1709.07343}, September 2017.

\bibitem[Sch24]{Schremmer24}
Felix Schremmer.
\newblock Affine {Bruhat} order and {Demazure} products.
\newblock {\em Forum of Mathematics, Sigma}, 12:e53, 2024.

\bibitem[Smi11a]{Smithling:PermissibleOn}
Brian~D. Smithling.
\newblock Admissibility and permissibility for minuscule cocharacters in
  orthogonal groups.
\newblock {\em Manuscripta Math.}, 136(3-4):295--314, 2011.

\bibitem[Smi11b]{SmithlingOddUnitary}
Brian~D. Smithling.
\newblock Topological flatness of local models for ramified unitary groups.
  {I}. {T}he odd dimensional case.
\newblock {\em Adv. Math.}, 226(4):3160--3190, 2011.

\bibitem[Smi11c]{Smithling:TopFlatnessEven}
Brian~D. Smithling.
\newblock Topological flatness of orthogonal local models in the split, even
  case. {I}.
\newblock {\em Math. Ann.}, 350(2):381--416, 2011.

\bibitem[Smi14]{SmithlingEvenUnitary}
Brian~D. Smithling.
\newblock Topological flatness of local models for ramified unitary groups.
  {II}. {T}he even dimensional case.
\newblock {\em J. Inst. Math. Jussieu}, 13(2):303--393, 2014.

\bibitem[Smi15]{Smithling15}
Brian Smithling.
\newblock On the moduli description of local models for ramified unitary
  groups.
\newblock {\em Int. Math. Res. Not. IMRN}, (24):13493--13532, 2015.

\bibitem[{Sta}26]{stacks-project}
The {Stacks project authors}.
\newblock The stacks project.
\newblock \url{https://stacks.math.columbia.edu}, 2026.

\bibitem[SW20]{SW20}
Peter Scholze and Jared Weinstein.
\newblock {\em Berkeley lectures on {$p$}-adic geometry}, volume 207 of {\em
  Annals of Mathematics Studies}.
\newblock Princeton University Press, Princeton, NJ, 2020.

\bibitem[Yan25]{Yang:SpinLocalModelsEven}
Jie Yang.
\newblock Topological flatness of orthogonal spin local models.
\newblock Preprint, arXiv:2512.16646, 2025.

\bibitem[YZZ26]{YangZachosZhao:PELTypeD}
Jie Yang, Ioannis Zachos, and Zhihao Zhao.
\newblock On {$p$}-adic integral moduli schemes and local models for {PEL} type
  {D}.
\newblock Preprint, arXiv:2602.23813, 2026.
\newblock Version 2; originally titled ``On quasi-split orthogonal local models
  of PEL type D''.

\bibitem[Zac23]{Zachos}
Ioannis Zachos.
\newblock On orthogonal local models of {H}odge type.
\newblock {\em Int. Math. Res. Not. IMRN}, (13):10799--10836, 2023.

\bibitem[Zhu14]{Zhu14}
Xinwen Zhu.
\newblock On the coherence conjecture of {P}appas and {R}apoport.
\newblock {\em Ann. of Math. (2)}, pages 1--85, 2014.

\bibitem[Zhu17]{Zhu16}
Xinwen Zhu.
\newblock An introduction to affine {G}rassmannians and the geometric {S}atake
  equivalence.
\newblock In {\em Geometry of moduli spaces and representation theory},
  volume~24 of {\em IAS/Park City Math. Ser.}, pages 59--154. Amer. Math. Soc.,
  Providence, RI, 2017.

\end{thebibliography}
\end{document}